%% file: qs-rigidity-arxiv.tex
\documentclass[12pt,oneside]{amsart}
\usepackage{mathtext,amsmath,amsfonts,amssymb,amsthm,bm}
\usepackage{verbatim}
\usepackage{version,etex}
\usepackage{enumitem}
\usepackage[top=1in, bottom=1in,left=1.5in, right=0.5in]{geometry}
\usepackage[pdftex]{graphicx,color}
\usepackage{graphicx}
\usepackage{epsfig}
\usepackage{hyperref}
\usepackage{array}
\usepackage{caption,subcaption}
\DeclareGraphicsExtensions{.pdf_t}
\reversemarginpar
\input prepicte.tex
\input pictex.tex
\input postpict.tex
\def\endpfclaim{$\checkmark$}  

\thanks{
The authors were supported by ERC AdG RGDD No 339523.
Trevor Clark was supported in part by NSF EMSW21-RTG 1045119.
The authors thank Lasse Rempe-Gillen, Daniel Smania and Davoud
Cheraghi for 
several helpful comments, and would like to express their gratitude to
Weixiao Shen for his many insightful remarks.}

\begin{document}
\input{macros.tex}

\begin{abstract}
In the late 1980's Sullivan initiated a programme to prove quasisymmetric rigidity in one-dimensional dynamics: interval or circle maps that are topologically conjugate are  quasisymmetrically conjugate (provided some obvious  necessary assumptions are satisfied). The aim of this paper is to conclude this programme in a natural class of $C^3$ mappings. Examples of such rigidity were established previously, but not, for example, for real polynomials with non-real critical points.

Our results are also new for analytic mappings. The main new ingredients of the proof in the real analytic case are (i) the existence of infinitely many (complex) domains associated to its complex analytic extension so that these domains and their ranges are compatible, (ii) a methodology for showing that combinatorially equivalent complex box mappings are qc conjugate, (iii) a methodology for constructing qc conjugacies in the presence of parabolic periodic points.

For a $C^3$ mapping, the dilatation of a high iterate of any complex extension of the real map will in general be unbounded. To deal with this, we introduce dynamically defined {\em qc\textbackslash bg partitions}, where the appropriate mapping has bounded quasiconformal dilatation, except on  sets with {\lq}bounded geometry{\rq}. To obtain such a partition we prove that we have very good geometric control for infinitely many dynamically defined domains. Some of these results are new even for real  polynomials, and in fact an important sequence of domains turn out to be quasidiscs. This technology also gives a new method for dealing with the infinitely renormalizable case.

We will briefly also discuss why quasisymmetric rigidity is such a useful property in one-dimensional dynamics.



\end{abstract}

\title{Quasisymmetric rigidity in one-dimensional dynamics\\ \today}


\author{Trevor Clark and Sebastian van Strien}
\address{Imperial College London}
\email{t.clark@imperial.ac.uk, s.van-strien@imperial.ac.uk}
\date{\today}
\maketitle
\tableofcontents

\section{Introduction and statements of main results}
\label{sec:introduction}

In this paper, we conclude 
a programme which was initiated by
Sullivan in the late 1980's, see for example 
\cite{Su-ICM}, \cite{Su} and \cite{McM}, 
namely that one has quasisymmetric rigidity in real one-dimensional dynamics.
Around the same time, Herman asked a similar question for critical circle maps
\cite{He2}.
Indeed
we will prove that topologically conjugate smooth mappings satisfying 
certain hypotheses are quasisymmetrically conjugate, see  
Theorem~\ref{thm:main} in Section~\ref{subsection:main theorem}.
For real-analytic
maps this result 
implies the following.

\begin{thm}[QS-rigidity of real-analytic maps]\label{thm:main, analytic}
Let $M$ be either $[0,1]$ or $S^1$ and suppose that
        $f,\tilde f\colon M\to M$ 
are real-analytic, have
at least one critical point or periodic point, 
and are topologically conjugate. 
Moreover, assume that the topological 
conjugacy is a bijection between
		\begin{itemize}
			\item the sets of parabolic periodic points,
			\item the sets of critical points,
		\end{itemize}
		and that the orders of corresponding critical points are the same.
		Then $f$ and $\tilde f$ are quasisymmetrically conjugate.
	\end{thm}

In our approach to extending this result to $C^3$ mappings
$f,\tilde f\colon M\to M$, where $M$ is either $[0,1]$ or $S^1,$
we will make use of complex tools. For these to be applicable, it is essential
to assume that at each critical point $c$ of $f$,
$f(x)=\pm(\phi(x))^\ell+f(c)$,
where $\phi$ is a local $C^3$ diffeomorphism,
$\phi(c)=0$, and $\ell\geq 2$ is an integer, and likewise at 
each critical point of $\tilde f$.
Under these hypotheses, we have the following.
	
\begin{thm}[QS-rigidy of $C^3$ maps with no parabolic cycles]\label{thm:main, C^3 repelling}
Suppose that $f$ and $\tilde f$ are $C^3$ topologically conjugate mappings, as above, with no parabolic cycles.
Assume that $f$ and $\tilde f$ each has
at least one critical point or periodic point, 
that the conjugacy is a bijection between the sets of critical
points of $f$ and $\tilde f$,
and that corresponding critical points have the same orders.
Then $f$ and $\tilde f$ are quasisymmetrically conjugate.
\end{thm}
	
\begin{rem}
	The assumption in 
Theorem~\ref{thm:main, C^3 repelling}
that neither $f$ nor $\tilde f$ has parabolic cycles can be removed
	under some weak additional smoothness and genericity assumptions, see Theorem~\ref{thm:main}.
\end{rem}

\begin{rem}
If $f,\tilde f\colon M\to M$ have periodic points,
then we do not require that they have critical points.
In particular, our theorems hold for interval
maps and for covering maps of the circle of degree $\ne 1$.
They also hold for maps $f,\tilde f\colon [0,1]\to \R,$
which are topologically conjugate on the non-escaping
parts of their dynamics. For example, if $f\colon  \R\to \R$ satisfies $|f(x)|/|x|\to \infty$ as $|x|\to \infty$,
then there exists an interval $[a,b]$ so that $f(\{a,b\})\subset \{a,b\}$ and 
$\lim_{n\to \infty}|f^n(x)|\to \infty$ for each $x\notin [a,b]$.
Our theorems imply qs-rigidity of the non-escaping part of the dynamics.
\end{rem}
	
	\begin{rem}
	Note that the parabolic multiplicity of corresponding parabolic points,
	that is, the order of contact at the parabolic points,
	does \textbf{not} need to be the same for the two maps.
	\end{rem}
	
\begin{rem}
For both analytic and smooth mappings,
we obtain the quasisymmetric conjugacy
as restriction of a quasiconformal mapping of the plane. In
the analytic case, the qc mapping is a 
conjugacy on a necklace neighbourhood of the real line;
however, for smooth mappings this is not the case.
\end{rem}

\begin{rem}
Quasisymmetric mappings are H\"{o}lder, but they are not
necessarily Lipschitz or even absolutely continuous. Indeed one can
give examples of topologically conjugate mappings
as in Theorem~\ref{thm:main, analytic}
for which no conjugacy satisfies these properties,
see for example~\cite{BJ}.
Moreover, it was proved in \cite{MM} that if two mappings with neither
solenoidal nor wild attractors are conjugate by an absolutely
continuous mapping, then they are $C^2$ conjugate, so that any
conjugacy which fails to preserve  multipliers is singular with
respect to Lebesgue measure. The condition that the mappings have
neither solenoidal nor wild attractors is necessary. In 
\cite{MoSm} it was proved that two Feigenbaum mappings with critical
points of the same even order are always conjugate by an absolutely
continuous mapping, and that provided that the degree is sufficiently high,
the same is true for two Fibonacci
mappings.
\end{rem}

	
	The above theorems are essentially optimal: in 
	Subsection~\ref{subsec:necessity of 	assumptions} we will show that the conditions
 	in these theorems are all necessary in general.
In particular, if $f$ is a map of the circle,  
then rigidity in full generality, that is, of all maps, requires that $f$ has  at least one critical point.


\subsection{Previous results}\label{subsec:previous results}
There has been substantial work on quasisymmetric rigidity over the
past three decades.
\subsubsection{A summary of past results}
Let us briefly summarize the previous results on qs-rigidity:
\begin{itemize}
\item For polynomials, rigidity was, in general, only known for real polynomials with all
critical points real.
\item For analytic mappings, rigidity was only known in some special
  cases, and often the results were semi-local.
\item For smooth mappings, rigidity was only known in a few very
  special cases.
\end{itemize}
Let us review these results.
		
\subsubsection{For polynomial maps,} 
		the most general result which is known  about qs-rigidity is the following:

		\begin{thm}[Kozlovski-Shen-van Strien \cite{KSS-rigidity}]\label{thm:KSS-density}
			Let $f$ and $\tilde f$ be real polynomials of
                        degree $d$ with only real critical points,
			which are all of {\em even order}. 
			If $f$ and $\tilde f$ are topologically conjugate
			(as dynamical systems acting on the real line) and corresponding
			critical points have the same order, then 
			they are  quasiconformally conjugate on the
                        complex plane (by a mapping preserving the
                        real line).
		\end{thm}
		
		In \cite{KSS-rigidity}, the assumption that all critical points are real
		is crucial. This assumption also implies that all attracting petals 
		of parabolic periodic points lie along the real line (in fact, 
		in \cite{KSS-rigidity} no details are given how to deal with parabolic points).
		
Prior to that paper, Yoccoz proved rigidity for all non-renormalizable 
quadratic polynomials, see \cite{Hubbard}.
Lyubich \cite{Lyubich-quadratic dynamics}
and Graczyk-{\'S}wi\c{a}tek \cite{GS2} and  \cite{GS3},
proved rigidity for all real quadratic polynomials. For an extension
to the case of real-analytic unimodal maps with a quadratic critical point 
see the  Appendix A of \cite{ALM}.

The results of  \cite{Lyubich-quadratic dynamics}, \cite{GS2} and  \cite{GS3},
rely on ``decay of geometry," a geometric feature of at most finitely
renormalizable quadratic maps that does not
hold even in the case of real polynomials of the form 
$z\mapsto z^d+c$ whenever $d>2$. However, the Quasi-Additivity Law and 
Covering Lemma of \cite{KL} made it possible to bring the theory
of non-renormalizable unicritical mappings in the
 higher
degree case
up to the same level of development as in the quadratic case. Rigidity for
unicritical, finitely renormalizable maps was proved in \cite{AKLS},
for general finitely renormalizable (not necessarily real)
multicritical
polynomials it was proved in \cite{KvS}, and for
real infinitely renormalizable unicritical mappings it was 
proved in \cite{Cheraghi}.

Let $\bm{M}^d$ denote the set of parameters $c$ such that
the Julia set of the mapping
$z\mapsto z^d+c$ is connected.
The set $\bm M^2$ is the Mandelbrot set, and
the set
$\bm M^d$ is often referred to as a {\em Multibrot set}.
In the families $z\mapsto z^d+c$, qc-rigidity is equivalent to the
local connectivity of  $\bm M^d$ at the parameter $c$.
This is due to an open-closed argument going
back to Douady-Hubbard and Sullivan,
and ideas of Yoccoz (see \cite{Hubbard}). 
Local connectivity of $\bm M^d$ remains unknown even in the case $d=2$.		
However, progress on this problem for infinitely renormalizable
quadratic maps $z\mapsto z^2+c$, with $c$ not
necessarily real,  has been made in two complementary
settings: when there is a combinatorial condition on the renormalizations that 
make it possible to prove
\emph{a priori} bounds for the induced polynomial-like maps,
see \cite{Kahn}, \cite{KL2}, \cite{KL3}, \cite{Cheraghi},
and when the (satellite) renormalizations
satisfy certain combinatorial growth conditions,
see \cite{Levin}, \cite{Levin2}, \cite{CSh}.
It is shown in \cite{DH} that the conjecture that the boundary
of the Mandelbrot set is locally connected implies density of 
hyperbolicity in the family $z\mapsto z^2+c$.
See \cite{Schleicher} for the extension to the higher degree case.

The main theorem of our paper shows that the assumptions in 
Theorem~\ref{thm:KSS-density} that all critical points are real,
and that all critical points have even order, 
are both unnecessary to obtain a qs conjugacy on the non-escaping
part of the dynamics on the real line.
When there are complex critical points,
it is not clear how to use the usual Yoccoz puzzle
for a polynomial map, however, our methods make it possible to
construct a complexified Markov partition for the real part of the dynamics. 
We should emphasise that,
nevertheless, our result depends heavily on 
many techniques from \cite{KSS-rigidity}.

\subsubsection{Semi-local results.} 
For real-analytic maps which are not polynomials,
most previous results only assert
the existence of a qs homeomorphism which is a conjugacy restricted
to the post-critical set.   One such result is due to Shen:
\begin{thm}\cite[Theorem 2]{Shen-density}
Let $f,\tilde f$ be real-analytic, topologically conjugate
maps with only hyperbolic repelling periodic points,
non-degenerate (quadratic) critical points and with 
{\em essentially  bounded geometry}.
Then there exists a qs homeomorphism which is a
conjugacy restricted to the post-critical set of
$f$ and $\tilde f$.
\end{thm}

 	 	A map $f$ is said to have \emph{essentially bounded geometry} if
		there exists a constant $C$ such that
		for any recurrent critical point $c$ of $f$
		and any nice interval $I$ containing $c$,
                $|I|/|\mathcal{L}_c(I)|\leq C$, where
                $\mathcal{L}_c(I)$ is the first return domain to $I$
                that contains $c$. See page~\pageref{def:nice} for the
                definitions of nice and $\mathcal{L}_c(I).$
		
		A similar result
		was proved jointly by  Levin and  
		the second author, namely that if $f,\tilde f$ are real-analytic covering maps of the circle
		(of degree $d\ge 2$)
		with at most one critical point, then there exists a
                qs homeomorphism
of $S^1$
		which is a conjugacy restricted to the 
		post-critical sets of $f$ and $\tilde f$. If $f,\tilde f$ have no parabolic
		or attracting periodic points, then in fact $f,\tilde f$ 
		are  qs conjugate. For precise statements
of these results see Theorem~\ref{thm:levstr} in Section~\ref{subsec:circle coverings}.
		We  show that, in fact, all of these qs conjugacies on the post-critical sets
		 extend to globally quasisymmetric conjugacies.

\subsubsection{Critical circle maps.}
For real-analytic circle homeomorphisms with 
precisely one critical point, the first
result about quasisymmetric rigidity was proved in the late
		1980's by  Herman, building on estimates due to {\'S}wi\c{a}tek, 
		see \cite{He2}, \cite{Swi} and also \cite{Petersen_Herman}: 
		 
		\begin{thm}[Herman-{\'S}wi\c{a}tek]  \label{thm:HS}
			A real-analytic unicritical circle mapping
			with an irrational 
			rotation number of bounded type, is  quasisymmetrically conjugate to 
			a rigid rotation of the circle.
		\end{thm}	 
		
		In the late 1990's this result was extended to 
		
\begin{thm}[\cite{dFdM 1} Corollary 4.6]
Suppose that $h\colon  S^1\rightarrow S^1$ is a topological 
conjugacy between two real-analytic
critical circle homeomorphisms $f$ and $\tilde f$,
each with precisely one critical point, $c$, $\tilde c$ respectively,
and $h(c)=\tilde c$.
Then $h$ is quasisymmetric.
\end{thm}

In an early version of this paper, we proposed an approach 
using our methods to treat critical circle maps with several
critical points. Since then, a purely real argument has been found
which solves this problem, so we do not develop this idea further
here.

\begin{thm}\cite{EdF}
\label{thm:EdF}
Let $f$, $\tilde f$ be two $C^r$ ($r\geq 3$) 
multicritical circle homeomorphisms with the same irrational rotation number
and the same number of non-flat critical points.
Let $h\colon S^1\rightarrow S^1$ be a homeomorphism
of $S^1$ conjugating $f$ and $\tilde f$, which maps critical points of
$f$
to critical points of $\tilde f$.
Then $h$ is quasisymmetric.
\end{thm}

		\begin{rem} As mentioned, the presence of at least one critical point is necessary
		in the previous theorem, 
		because for circle diffeomorphisms the analogous statement
		is false. Indeed, 
		one can construct diffeomorphisms for which some sequence of iterates
		has `almost a saddle-node fixed point',
		resulting in larger and larger passing times near these points.
		This phenomenon is also referred to
		as `a sequence of saddle-cascades'.
		It was used by Arnol'd and Herman to construct examples
		of diffeomorphisms of the circle which are conjugate to irrational rotations, but
		where the conjugacy is not absolutely continuous, not qs and for which 
		the map has no $\sigma$-finite measures, see for example Section I.5 in \cite{dMvS}
		and Herman's thesis, \cite{Herman}, Chapters VI and XII.
		In the diffeomorphic case, to get that the conjugacy is quasisymmetric or $C^1$,
		one needs assumptions on the rotation number 
		(to avoid these sequences of longer 			
		and longer saddle-cascades).
		\end{rem}

		\subsubsection{Renormalization, universality and $C^1$-rigidity} Of course our 
		result is related to 
		the famous Feigenbaum-Coullet-Tresser renormalization conjectures,
		which were proved in increasing generality by Sullivan \cite{Su}, 
		McMullen 		\cite{McM},
		 Lyubich \cite{Lyubich-renormalisation}, Avila-Lyubich \cite{AL}. The specific
		 version of this conjecture we are referring to is:

		\begin{thm} Let $f,\tilde f$ be unimodal real-analytic maps which
		are infinitely renormalizable and of bounded type. Assume that $f,\tilde f$ are
		topologically conjugate.   Then there exists a $C^{1+\alpha}$ map which is
		a conjugacy restricted to the post-critical set of
		$f$ and $\tilde f$.
		\end{thm}
Smania has generalised this theorem to real-analytic, infinitely renormalizable
mappings with bounded combinatorics, see \cite{Smania-phase}.

		\medskip
		The previous theorem shows that the conjugacies between maps are even better than 
		quasisymmetric along the post-critical set. We believe that a weaker 
		version of the above $C^1$-rigidity should also hold in the general case: 
		
		\begin{conj} In the setting of the main theorems, whenever  $\omega(c)$ is minimal and 
		has {\lq}bounded geometry{\rq}, the conjugacy is differentiable at each critical point $c$.
		\end{conj}
	
		Without the assumption of bounded geometry, the conjugacy is in general not 
		differentiable: smooth unimodal Fibonacci maps with a critical point of order two
		need not be smoothly conjugate at the critical point,  see \cite{LM}.
			
		 In fact, in the case of critical circle homeomorphisms with 
		 exactly one critical point, 
		 a stronger result is known. 
		It was proved for analytic mappings by Khanin and Teplinsky \cite{khanin_teplinsky} and later extended
		to the smooth category by Guarino and de Melo \cite{GM}.
	
\begin{thm}[Khanin-Teplinsky, Guarino-de Melo]
\label{thm:khanin}
For critical circle mappings with one critical point we have the following:
\begin{itemize}[topsep=1pt]
\item[\cite{khanin_teplinsky}] Let $f$ and $\tilde f$ be two analytic unicritical circle maps with the
same order of critical points and the same irrational rotation number. Then
$f$ and $\tilde f$ are $C^1$-smoothly conjugate to each other.		
\item[\cite{GM}] Any two $C^3$ circle homeomorphisms $f$ and $\tilde f$ with the
same irrational rotation number of bounded type 
and each with precisely one critical point, both
of the same odd order, are $C^{1+\alpha}$-smoothly conjugate,
where $\alpha>0$ is a universal constant.
\end{itemize}
\end{thm}

For $C^4$ mappings the following is known:
\begin{thm}[\cite{GMM}]
 Let $f$ and $\tilde f$ be two $C^4$ unicritical circle maps with the
same order of critical points and the same irrational rotation number.
Let $h$ be the unique topological conjugacy between $f$ and $\tilde f$
that maps the critical point of $f$ to the critical point of $\tilde
f$. Then:
\begin{itemize}
\item $h$ is a $C^1$ diffeomorphism.
\item $h$ is $C^{1+\alpha}$ at the critical point, for a universal
  $\alpha>0.$
\item For a full Lebesgue measure set of rotation numbers, $h$ is $C^{1+\alpha}$.
\end{itemize}
\end{thm}

		By \cite{A} even for analytic maps the theorem of \cite{GM} cannot be
		extended to circle maps with rotation number of unbounded type.
		We should note that these theorems
		build on earlier work of de Faria, de Melo and Yampolsky on renormalization
		of critical circle maps.

		For general interval maps, one cannot expect conjugacies to be 
		$C^1$, because having a $C^1$ conjugacy 
		implies that corresponding periodic orbits have the same multiplier.  
Moreover,
suppose that $f$ and $\tilde f$ are $C^3$ maps of the interval with all periodic points 
hyperbolic repelling and no Cantor attractors.
If the multipliers of corresponding periodic points 
of $f$ and $\tilde f$ are the same, then
we have that $f$ and $\tilde f$ are $C^3$ conjugate, see \cite{Li-Shen Smooth}
and also \cite{MM}. 

		\subsubsection{Covering maps of the
                  circle}\label{subsec:circle coverings}
		
		For the case of covering maps of the circle one has the following theorem: 

\begin{thm}[Levin-van Strien \cite{LevStr:invent}]\label{thm:levstr}
Assume that $f,\tilde f\colon  S^1\to S^1$ are topologically conjugate
real-analytic covering maps of the circle (positively oriented) 
each with exactly one critical point of the same odd order and the conjugacy 
maps the critical orbit to the critical orbit. Then the following hold: 
			\begin{itemize}
				\item if  $\omega(c)$ is minimal,
				there exists a qs homeomorphism which is a conjugacy on the post critical set (and maps
				the critical point to the critical point); 
				\item if $f,\tilde f$ have only repelling periodic orbits, then the conjugacy 
				is quasisymmetric. 
			\end{itemize}
		\end{thm}
	
		The main theorem in our paper shows that the assumption 
		that there exists at most one critical point can be
                dropped, and we prove that there the conjugacy is a global
                quasisymmetric conjugacy even when there are parabolic points.

		\medskip

Note that in our main theorem (see Theorem~\ref{thm:main} or
Theorems \ref{thm:main, analytic},
\ref{thm:main, C^3 repelling} above) we ask for the maps $f$ and $\tilde f$ to have
the same number of critical points. Therefore, the
set-up is  different from the one
in the following ``linearization result'' by Carsten Petersen \cite{Petersen_blaschke}:

\begin{thm}[Petersen]
Assume that $f\colon  S^1\to S^1$ is a degree $d\ge 2$ map which 
is a restriction of a Blaschke product, then $f$ is qs conjugate to $z\mapsto z^d$ if 
\begin{enumerate}
\item  $\omega(c)\cap S^1=\emptyset$ for each recurrent critical point, and 
\item any periodic point in $S^1$ is repelling.
\end{enumerate}
\end{thm}

		Somewhat weaker results were derived in \cite[Exercise 6.2]{dMvS}
		and also in the Cornell thesis of L.K.H. Ma (1994).
		The assumption that $\omega(c)\cap S^1=\emptyset$ for any recurrent critical point $c$
		ensures that $f$ has no `almost saddle-nodes' on the circle. Perhaps one could replace
		this assumption by a kind of `diophantine' condition on the itinerary of the
		critical points of $f$ (so that one does not need to have a globally defined rational map $f$ in this
		theorem).  But in any case, as was remarked in \cite{LevStr:invent},
		see Section \ref{subsec:necessity of assumptions},
		in general a real-analytic covering map of the circle of degree $d$ (with critical points) 
		is {\em not} qs conjugate to  $z\mapsto z^d$ (on the circle $\{z\in \C; |z|=1\}$). 
		To see this, we can use Theorem C
		in \cite{LevStr:invent} to see that within any family, near every map without periodic attractors,
		there exists another map with a parabolic periodic point. Then use the argument 
		given in Section I.5 in \cite{dMvS}.

		If a real-analytic covering map of the circle with degree at least two
		has no critical points, then either
		it has periodic attractors (or parabolic periodic orbits) or some iterate
		of the map is expanding. This follows from a theorem of Ma\~n\'e, for a simplified proof
		see \cite{vS} or  Section III.5 in \cite{dMvS}.
		From this it is easy to see that any such map which
		has no attracting or parabolic periodic orbits is qs conjugate to $z\mapsto z^d$.
		The Main Theorem describes a more general setting in which
		covering maps of the circle can be qs conjugate, even when there
		are parabolic periodic points.
A special case in the analytic setting was considered in \cite{LPS}.

\subsubsection{Quasisymmetric rigidity for smooth maps}
For unimodal Misiurewicz mappings and more generally for geometrically
finite mappings under various conditions, for example negative
Schwarzian derivative or absence of parabolic points,
it was proved in \cite{J}, \cite{Ji} and \cite[Exercise 6.2]{dMvS}.

To tackle the problem of smooth mappings with more general
combinatorics
we use asymptotically holomorphic extensions of smooth mappings.
The specific extension that we use was constructed in
\cite{GSS}. Asymptotically holomorphic extensions were previously used in
\cite{Lyubich-Fibonacci2} to prove rigidity of
smooth mappings with Fibonacci combinatorics and
a quadratic critical point
using exponential decay of geometry.
They were also used in \cite{JS} to show
quasisymmetric rigidity 
for certain unimodal maps ``close to the Chebyshev polynomial
at all levels".

Using entirely real methods, 
Pa{\l}uba proved a Lipschitz rigidity result for
Feigenbaum maps  \cite{Pa}.	

\begin{rem}
As a historical comment,
Jiang \cite{Ji} refers to lecture notes from 1987 where
Sullivan discusses
the problem of quasisymmetric rigidity in one-dimensional dynamics.
Sullivan always had in mind to use qs-rigidity of smooth maps to
analyze the structure of smooth dynamical systems. In this paper
we make a key step in this endeavor, see the next subsection.
\end{rem}

\subsection{Applications of quasisymmetric ridigity}
The main motivation for proving quasisymmetric rigidity stems from the fact 
that this property allows one to use tools from complex analysis, such as the 
Measurable Riemann Mapping 
Theorem to study the dynamics of real-analytic maps.
For example, it turned to be a crucial step towards proving 
the following types of results:
\subsubsection{Quasiconformal surgery}
Theorem~\ref{thm:HS} allows one to use quasiconformal surgery to  prove 
the 	existence of quadratic Siegel 		
polynomials  for which the Julia set is locally connected.
See the work of Ghys, Shishikura, Herman and others,	
see for example~\cite{Petersen_local} and~\cite{BF}.

\subsubsection{Density of hyperbolicity}
Quasisymmetric rigidity is one of the
key tools for proving density of hyperbolic dynamics in dimension-one.
Hyperbolicity is dense within the space of smooth maps. For the unimodal case 
see \cite{GS2}, \cite{GS3}, \cite{Lyubich-quadratic dynamics},
\cite{Koz_ax}, and
for the general (real) multimodal case see 
\cite{KSS-rigidity} and \cite{KSS-density}.
Density of hyperbolicity is also known within certain spaces of real
transcendental maps and in the 
Arnol'd family, see \cite{RvS} and \cite{RvS2}.

\subsubsection{Monotonicity of entropy}
Using quasisymmetric rigidity,
monotonicity of entropy for real polynomials with only real
  critical points was proved in \cite{BvS_mono}.
In the quadratic case this
was shown by quite a few people, see \cite{Douady},
\cite{DH} \cite{Milnor-Thurston}, \cite{Tsuji} and \cite[II. 10]{dMvS}, 
for the cubic case see
\cite{Milnor-Tresser}.
Monotonicity of topological entropy in certain  families of transcendental maps, 
	such as $[0,1]\ni x\mapsto a\sin(\pi x)$ was proved in \cite{RvS2}.

\subsubsection{Hyperbolicity of renormalization}
Quasisymmetric rigidity plays an important role in the proof
of hyperbolicity of renormalization. 
Rigidity is a key ingredient in understanding the
structure of the renormalization horseshoe
and the stable leaves of the renormalization operator, see
\cite{Su}, \cite{McM3}, \cite{Lyubich-renormalisation} and \cite{AL}.
Moreover, in the proof of
hyperbolicity of renormalization,
quasiconformal rigidity is used in conjunction with
the Small Orbits Theorem to prove that renormalization is
exponentially expanding in the direction transverse to the
hybrid class of an infinitely renormalizable unimodal mapping,
\cite{Lyubich-renormalisation} and \cite{Lyubich-regular or
  stochastic}.
This argument requires knowing quasiconformal rigidity of complex
mappings, and is also used in \cite{Smania-Fibonacci} to prove
hyperbolicity of the Fibonacci renormalization.
An infinitesimal version of this argument, which only requires
rigidity for real mappings, is used in
\cite{Smania2} to prove hyperbolicity of renormalization for
multimodal mappings with bounded combinatorics.

\subsubsection{Connectedness of conjugacy classes}
One of the main goals of dynamical systems is to 
classify and understand the conjugacy classes of 
dynamical systems. It is natural to ask whether conjugacy classes of mappings are
connected. This is proved in certain families of polynomials
in \cite{BvS_mono} and \cite{CS}. 

Using the Measurable Riemann Mapping Theorem,
one can prove that the hybrid classes in the space of polynomial-like
germs are
connected, \cite{Lyubich-renormalisation}.
Using
quasiconformal rigidity one can show that the 
topological conjugacy classes coincide with the hybrid conjugacy classes,
and thus obtain that topological conjugacy classes of 
polynomial-like germs are connected.
Using the same argument, this result
also holds for the germs of non-renormalizable persistently recurrent complex
box mappings, \cite{Smania-Fibonacci, KvS}.
For real analytic mappings,
it is proved in \cite{ALM} that the conjugacy class of a
non-renormalizable mapping with a quadratic critical point is
connected. This result exploits the ``decay of geometry''
for non-renormalizable quadratic mappings. 

\medskip

The rigidity result of this paper can be used to prove that
the topological conjugacy class of a real analytic mapping of the
interval with no 
parabolic cycle is connected.
The proof of this and related results will be given in a forthcoming paper.

\subsubsection{Conjugacy classes are manifolds}
Sullivan conjectured
the following:

\medskip
\noindent\textbf{Conjecture.}
Conjugacy classes of analytic (smooth) mappings with all periodic
orbits hyperbolic on the interval are analytic (smooth)
connected manifolds.

\medskip

We expect that our results will be an essential ingredient in the proof.

\medskip

Using transversality results, it is known that 
the set of rational (non-Latt\`es) mappings satisfying critical relations
of the form $f^{k_{i}}(c_i)=f^{k_j}(c_j)$ are
manifolds,
\cite{Epstein},\cite{LSvS1}. Going beyond rational mappings,
this result is proved for certain interval mappings with a flat critical point
in \cite{LSvS2}. These results imply that the topological
conjugacy classes of critically finite rational maps
are analytic manifolds,
namely points.
But in the space of rational maps,
critically finite (non-Latt\`es) mappings are singletons by Thurston
Rigidity
\cite{DH-acta},
so the fact that these are manifolds is rather trivial.

To prove similar results in more general settings rigidity is
vital. For a real quadratic-like mapping
it is shown in \cite{DH}, that the 
conjugacy class of a real quadratic-like mapping is
an infinite dimensional, complex analytic manifold, 
and moreover, it has
co-dimension one in the space of quadratic-like mappings 
\cite{Lyubich-renormalisation}.
The proof that it is an analytic manifold 
uses the identification of the hybrid class
of a polynomial-like
mapping with the space of analytic expanding circle mappings.
Because of
rigidity, the hybrid class is the same as the topological conjugacy class.

Going beyond polynomial-like mappings, in \cite{ALM} it is proved
that the topological conjugacy class of certain
analytic  unimodal
mappings with a non-degenerate critical point
is an analytic manifold; additionally one needs to assume 
either that the mappings are
quasiquadratic
that they have all periodic orbits hyperbolic, see \cite{AM2}.
This result is extended to the higher degree case in \cite{C}.
For smooth mappings, using the fact that topological conjugacy classes 
are stable manifolds of the renormalization operator, 
in \cite{dFdMP} it is
proved that the topological conjugacy classes of unimodal $C^4$ mappings that
are infinitely renormalizable of bounded type
are $C^1$ manifolds.

The previous results rely heavily on complex analytic
machinery,
so
for smooth mappings, the strategy of any proof of Sullivan's
Conjecture
would require
substantial modifications. Indeed, if two smooth mappings
on the interval are topologically conjugate, but one has a repelling
periodic point off the real line where the mapping is conformal
and the other mapping is not conformal at the corresponding point, then
the two mappings cannot be qc-conjugate in a complex neighbourhood
of the interval that contains the repelling point.
(Nevertheless, they are qs-conjugate on the interval.)

Sullivan had another rationale for asking the question about 
quasiconformal or quasisymmetric rigidity, as this fits in his
famous Sullivan dictionary, relating results in complex dynamics, with
results about Kleinian groups and three manifolds, see for example \cite{McM2}. 
Sullivan's point of view is that conjugacy classes of dynamical 
systems can be treated as infinite dimensional Teichm\"uller spaces.
The Teichm\"uller pseudometric on a conjugacy class of an interval
mapping is defined as
$$\dist(f,\tilde f)=\inf \|\mu_{h}\|_{\infty},$$
where $\mu_h$ denotes the dilatation of a qc mapping $h$ and
the infimum is taken over all quasiconformal 
extensions to the plane of qs conjugacies between $f$ and $\tilde f$.

\medskip

\subsubsection{Palis conjecture}
The Palis conjecture is a far reaching conjecture about the
typical dynamics, from the measurable point of view, that occur
in dynamical systems. The Palis Conjecture has been proved 
in dimension-one for analytic 
unimodal mappings. It was first proved in the
celebrated result of Lyubich: almost every real quadratic mapping is
regular or stochastic \cite{Lyubich-regular or stochastic}. 
This 
result was generalized to families of analytic quadratic mappings
in \cite{ALM}.
For higher degree mappings
in \cite{BSvS} it is proved that in a one-parameter family of real
unimodal polynomials, almost every mapping has a unique physical
measure that is supported on the real line.
The regular or stochastic theorems
were generalized
to higher degree unimodal mappings in \cite{ALS} and
\cite{C}.
These results can be improved to obtain that almost every mapping
in a generic analytic family of unimodal mappings,
is regular or Collet-Eckmann, \cite{AM1},\cite{AM2}, and 
it can be shown that
various other
good properties hold for almost every such mapping \cite{AM3}.
Quasisymmetric rigidity plays a crucial role in the proofs of these
results. As we already mentioned,
it is used in the proof of hyperbolicity of 
renormalization, which is the key ingredient in the proof that
infinitely renormalizable mappings have measure zero in 
generic analytic families of unimodal mappings. Moreover,
it plays an important role in the proof that the conjugacy
classes of analytic unimodal mappings are complex analytic
Banach manifolds.
This result is used to transfer metric information
({\em a priori} bounds)
from the dynamical plane to the parameter plane, 
where it is used in parameter-exclusion arguments.

\subsection{Necessity of assumptions}\label{subsec:necessity of assumptions}
In general, all the assumptions in the theorem are necessary to have
a quasisymmetric, or indeed even a H\"older, conjugacy between $f$
and $\tilde f$.

(i) To see that we cannot drop the assumption that the
critical points should have the same order, 
consider two maps $f(z)=z^2+c$ and $\tilde f (z)=z^a+c'$ with $a\geq 4$
and parameters $c$ and $c'$ chosen so that $f$ and $\tilde f$ both have
the same, Fibonacci, combinatorics. (So that the close return times of the critical point
to itself occur at the Fibonacci numbers; these are the fastest possible returns when there are 
no central returns.)
Iterates of the critical point of $f$ and $\tilde f$  accumulate to $0$ at  different rates.
Indeed, let $S_n$ be the Fibonacci sequence $1,2,3,5,8,\dots$, then 
$f^{S_n}(0)$ converges to $0$ at a geometric rate, and $\tilde f^{S_n}(0)$ at a polynomial rate,
see \cite{BKNS}.
So $f$ and $\tilde f$ are topologically conjugate, but not qs conjugate.
Note that it is not known whether two Fibonacci maps 
of the form $f(z)=z^d+c$ and $\tilde f(z)=z^{d'}+c'$ with $d\ne d'$ both large 
are qs-conjugate. 

(ii) Critical points  of {\em odd order} are {\em invisible} from a 'real' point of view; however,
 they can have rather important  consequences for the rate of recurrence. For example,
there exist $f,\tilde f$ degree 2 covering maps of the circle,
$f$ with an odd critical point, and $\tilde f$
everywhere expanding, that are necessarily topologically conjugate 
(because neither of them has wandering intervals), 
but that are not qs conjugate: one can construct the map
$f$, real-analytic, with a critical point
so that it has longer and longer saddle-cascades.

(iii) It is necessary to assume that the topological conjugacy maps
parabolic periodic points to parabolic periodic points, since
 the local escape rate near a hyperbolic periodic point is
completely different from the rate near a parabolic periodic point
 (the former is geometric, the latter is polynomial).

(iv) As mentioned, it is necessary to assume in the case of circle homeomorphisms, that
the maps have at least one critical point. In fact, if $f\colon  S^1\to S^1$ is non-monotone, then it automatically has
periodic orbits, see Proposition~\ref{prop:circlecase}.
		
(v) Although we have no proof that qs-rigidity does not hold for $C^2$
maps, we believe that it 
should be
possible to construct counter examples to our Main Theorem when
one has longer and longer saddle-cascades (in the $C^2$ case, one has 
much weaker distortion and cross-ratio control).		
\medskip

\subsection{Statement of the main theorems}
\label{subsection:main theorem}

\subsubsection{The class of maps}\label{subsec:class C}
We let $\mathcal{C}=\mathcal{C}(M)$\label{def:class of maps} 
denote the set of maps $f\colon  M\rightarrow M$ with the following properties:
\begin{enumerate}[label=(\arabic*)]
\item $f$ is $C^3$ on (a neighbourhood of) $M$;
\item $f$ has finitely many critical points all with {\bf integer orders}; that is, 
if $c_i$ is a critical point of $f$, then we can express
$$f(x)=\pm(\phi_i(x))^{d_i}+f(c_i),$$
where $\phi_i$ is a local
diffeomorphism near $c_i$ of class $C^3$, 
$\phi_i(c_i)=0,$ and $d_i\in\{2,3,\dots\};$
\item if $p$ is a periodic point of period $n$ and multiplier $\lambda=\pm 1$, then the
Taylor  expansion of $f^n$ at $p$ is given by
$$f^n(x)=p+\lambda(x-p)+a(x-p)^{d+1}+R(|x-p|),$$
in a real neighbourhood of $p$ and $a\in\mathbb{R},$
where $R=o(|x-p|^{d+1})$,
and as usual $o(t)/t\rightarrow 0$ as $t\rightarrow 0$.
We call $d$ the \emph{parabolic multiplicity}
of the parabolic periodic point.
\end{enumerate}		
The set $\mathcal{C}$ includes all $C^3$ maps without
parabolic periodic points, all $C^3$ maps with parabolic points which all have parabolic multiplicity 2, 
and all analytic maps.
Notice that condition (3) implies that $f$ can have at most finitely
many
parabolic cycles.

\subsubsection{Statement of the main theorem}
\begin{thmx}\label{thm:main}
Suppose that $f,\tilde f\colon [0,1]\rightarrow [0,1]$,
or alternatively that $f,\tilde f\colon  S^1 \to S^1$, 
each has at least one periodic point, are in $\mathcal{C}$.
Suppose that $f$ and $\tilde f$ are topologically conjugate and 
moreover, assume that the topological conjugacy is a bijection between
\begin{itemize}
\item the  sets of critical points,
\item the set of parabolic periodic points,
\end{itemize}
and the orders of corresponding critical points are the same.
Then $f$ and $\tilde f$ are quasisymmetrically conjugate. 
\end{thmx}
As we observed in Theorem~\ref{thm:EdF}
this theorem holds for critical circle maps without periodic
points. Because of \cite{EdF}, we do not give an independent proof of
this fact.

\medskip

This result implies Theorem 3.8 of \cite{RvS2}, where it was stated
for analytic mappings of the interval $[a,b]$ such that 
$f(\{a,b\})\subset\{a,b\}$.

\medskip

A straightforward application of the techniques 
in this paper implies the following result
that is useful in applications:
Suppose that $f,\tilde f\colon M\rightarrow M$ are 
mappings of class $\mathcal C$ as in the 
statement of Theorem~\ref{thm:main}.
Suppose that $N$ is a union of intervals in $M$
and that $N'$ is a union of intervals
$J'\subset N$ such that for each interval $J'$,
there exists $k_{J'}\geq 1$ such that $f^{k_{J'}}(J')\subset N$
and that $f^{k_{J'}}$ does not have a critical point on 
$\partial J'$. Suppose that we have the same objects for 
the mapping marked with a tilde. 
Define $F\colon N'\rightarrow N$ by $F|_{J'}
=f^{k_{J'}}|_{J'}$ and 
$\tilde F\colon\tilde N'\rightarrow \tilde N$ by $\tilde F|_{\tilde J'}
=\tilde f^{k_{\tilde J'}}|_{\tilde J'}$.

\medskip
\noindent\textbf{Theorem A$\bm{'}$.}
\emph{If $F\colon N'\rightarrow N$ and $\tilde F\colon \tilde N'\rightarrow \tilde N$
are topologically conjugate by a mapping $h\colon M\rightarrow M$ that is order preserving
on $M$ and the conjugacy is
a bijection between 
the sets of parabolic periodic points, and the sets of critical points
and corresponding critical points have the same order, then
$F$ and $\tilde F$ are qs-conjugate.}

\medskip

We believe that the techniques that we develop
will be useful in other contexts, so even though some of these
theorems are quite technical, we state them here.

\subsubsection{Complex bounds theorems: the general case}

To prove Theorem~\ref{thm:main}, we make use of various complex extensions
of $f$. Close to a critical point of $f$, we use usual
\emph{complex box mappings}, see page~\pageref{def:box mapping}
for the definition. These were constructed for asymptotically holomorphic
extensions of smooth mappings at persistently recurrent critical
points in~\cite{CvST}, and 
at reluctantly
recurrent critical points
for analytic mappings with only even
critical points in \cite{KSS-density}. See
page~\pageref{def:gapextension} for the definitions of the gap and
extension properties.

\begin{figure}
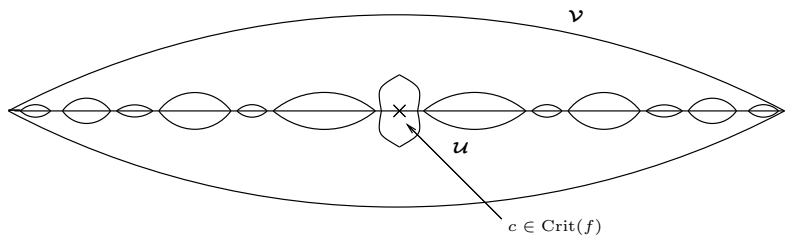
\caption{A complex extension of the real return mapping where the
 range is a Poincar\'e lens domain.}
\label{fig:lensmapping}
\end{figure}

In the \textit{persistently  recurrent} case we build on \cite{CvST},
we prove the following,
see Theorem~\ref{thm:box mapping persistent infinite branches}:
\begin{thmx}
Suppose that $c_0$ is a persistently recurrent critical point of $f$,
at which $f$ is at most finitely renormalizable,
and let $\Omega_1$ denote the set of critical points in $\omega(c_0)$. 
Then there exist arbitrarily small, combinatorially defined real
neighbourhoods $\mathcal I$ of $\Omega$ and $m\in\mathbb N$, such that
the return mapping
$R_{\mathcal I}^m\colon\Dom^*(\mathcal I)\rightarrow \mathcal I$
extends to a complex box mapping
$$F_{it}\colon\bm{\mathcal{U}}_{it}\rightarrow \bm{\mathcal{V}}_{it}$$ such
that $\bm{\mathcal{V}}_{it}$ is a union of  Poincar\'e lense domains,
$(\bm{\mathcal{V}}\setminus\bm{\mathcal{U}})\cap\mathbb H^+$ is a
quasidisk
and this complex box mapping satisfies the gap and extension properties.
\end{thmx}
Let us point out that in this theorem,
the box mapping has infinitely many branches, since we do
not
restrict to components of the domain that intersect the post-critical
set, and the range is a Poincar\'e lens domain, see
Figure~\ref{fig:lensmapping}.

In the \textit{reluctantly recurrent} case, building on \cite{KSS-density} we
construct complex box mappings about reluctantly recurrent critical
points. This result is new even in the analytic case, since we allow
for odd critical points. Let us give a rough statement of the
theorem here,
see Theorem~\ref{thm:box mapping reluctant} for a precise statement.
\begin{thmx}
Let $\Omega_r$ denote the set of reluctantly recurrent critical points
of $f$. Then there exist arbitrarily small combinatorially defined real
neighbourhoods
$\mathcal I$ of $\Omega_r$ such that the return mapping to $\mathcal I$
extends to a complex box mapping
$F\colon\bm{\mathcal{U}}\rightarrow\bm{\mathcal{V}}$ where
$\bm{\mathcal{V}}$ is a union of Poincar\'e lens domains, 
$(\bm{\mathcal{V}}\setminus\bm{\mathcal{U}})\cap\mathbb H^+$ is a
quasidisk and this complex box mapping satisfies the gap and extension
properties.
\end{thmx}
The key results used in the proof of this theorem
are the real bounds at reluctantly recurrent
critical points, Proposition~\ref{prop:real bounds reluctant}.

\medskip

This result, together with the complex bounds of \cite{CvST},
implies Theorem 3.9 of \cite{RvS2}.

\subsubsection{Touching box mappings}
The previous two theorems provide us with control near the critical points of
$f$, but we need a tool that will make it possible to construct a qs
conjugacy on all of $M=[0,1]$ or $S^1$. To do this, we make use of 
\emph{touching box mappings}, see page~\pageref{def:tbm} for the
definition. Figure~\ref{fig:touchingglobal} depicts a touching box mapping.

	\begin{figure}[htp] \hfil
		\beginpicture
		\dimen0=0.3cm
		\setcoordinatesystem units <\dimen0,\dimen0>
		\setplotarea x from -9 to 30, y from -5 to 4
		\setlinear
		\plot -8 0 26 0 /
		\circulararc 120 degrees from 5 0  center at  0 -3 
		\circulararc -120 degrees from 5 0  center at  0 3 
		\put {$\partial \VV_T$} at 5 2 
		\circulararc 110 degrees from 1 0  center at  -2 -2 
		\circulararc -110 degrees from  1 0 center at  -2 2 
		\circulararc 90 degrees from 5 0  center at  3 -2 
		\circulararc -90 degrees from  5 0 center at  3 2 
		\circulararc 115 degrees from 14 0  center at  9.5 -3 
		\circulararc -115 degrees from 14 0  center at  9.5 3 
		\circulararc 115 degrees from 23 0  center at  18.5 -3 
		\circulararc -115 degrees from 23 0  center at  18.5 3 
		\circulararc 115 degrees from 23 0  center at  20 -2 
		\circulararc -115 degrees from  23 0 center at 20 2 
		\circulararc -90 degrees from 14 0  center at  15.5 -1.5 
		\circulararc 90 degrees from  14 0 center at 15.5 1.5 
		\put {$\partial \UU_T$} at 1.5 1.3 
		\put {$I$} at 9 -1 
		\put {$*$} at 9 0 
		\endpicture
		\caption{A touching box mapping $F_T\colon  \UU_T\to \VV_T$: the domain does not contain critical points of $f$
		(marked with the symbol $*$). \label{fig:touchingglobal}}
		\end{figure}

Informally a touching box mapping is a complex box mapping
$F_T\colon\bm{\mathcal{U}}_T\rightarrow\bm{\mathcal{V}}_T$ such that
$F_T$ is a diffeomorphism on each component of its domain,
$\bm{\mathcal{V}}_T$ is a union of finitely many 
Poincar\'e lens domains that form a
necklace neighbourhood of $M$, $\bm{\mathcal{U}}_T$
has finitely many components,
$\bm{\mathcal{U}}_T\subset\bm{\mathcal{V}}_T$,
but not compactly, and the intersection points of 
$\partial \bm{\mathcal{U}}_T$ and $\partial \bm{\mathcal{V}}_T$
are only along the real line at points where components of 
$\bm{\mathcal{U}}_T$ and $\bm{\mathcal{V}}_T$ meet tangentially.
If the mapping has a parabolic periodic point, then we choose
$\bm{\mathcal{V}}_T$ so that it contains the parabolic cycle in
its boundary. This makes is possible to use a local analysis near the
parabolic points to construct a conjugacy there.
We construct touching box mappings in
Theorem~\ref{thm:touching box map}. 

One consequence of this Theorem is rigidity away from critical points,
see Corollary~\ref{cor: pullback argument for touching box mappings}.
\begin{thmx}
Suppose that $f$ and $\tilde f$ are of class $\mathcal C$ and
that $h$ is a topological conjugacy between $f$ and $\tilde f$ as in 
Theorem~\ref{thm:main}. Let $\mathcal I,$ a union of intervals,
be a combinatorially defined
neighbourhood of $\crit f$ and let $\tilde{\mathcal{I}}$ be the corresponding
neighbourhood for $\tilde f$. Let $E(\mathcal I)$ denote the set of points
whose orbits avoid $\mathcal I$. Then $h|_{E(\mathcal I)}$ is quasisymmetric.
\end{thmx}
This Theorem
immediately implies qs-rigidity for mappings with periodic points
and no critical points, for example, topologically expanding maps of the circle 
with degree at least 2; here $E(\mathcal I)$ may contain parabolic cycles. 

\begin{rem}
After proving that we can construct an initial external conjugacy between $f$
and $\tilde f$ for the complex box mappings near the critical points,
see Theorem~\ref{thm:external conjugacy}, we can 
almost argue as in \cite{KSS-rigidity} to prove 
Theorem~\ref{thm:main, analytic}, qs-rigidity, in the analytic case.
\end{rem}

\subsubsection{Quasidisks}
Long cascades of central
returns are one of the main reasons why
the argument from the analytic case
does not go through for smooth mappings.
One step towards dealing with this issue
is proving that certain puzzle pieces are
quasidisks; these puzzle pieces are obtained by a 
modification of the
construction of the enhanced nest, \cite{KSS-rigidity}.
See Section~\ref{subsec:improving} and
Proposition~\ref{prop:quasidisks} for the details.
We also prove that puzzle pieces from the enhanced nest are
quasidisks,
something that was neither proved nor needed in \cite{KSS-rigidity}.

\begin{thmx}\label{thmx:good central}
There exists $\delta>0$ such that the following holds.
\begin{itemize}
\item Suppose that $$\V^m\supset \V^{m+1}\supset\dots\supset \V^{m+n+1}$$ is 
long, maximal, cascade of central returns 
for a complex box mapping
$F\colon\bm{\mathcal{U}}\rightarrow\bm{\mathcal{V}}.$
Then there exists a puzzle piece $\bm{E}$ such that $\V^{m+2}\subset\bm
E\subset \V^m$ and $\bm E$ is a $\bm E$ is $\delta$-nice and 
a $1+1/\delta$-quasidisk.
\item If $\bm E$ is a puzzle piece from the enhanced nest, then 
 $\bm E$ is a $1+1/\delta$-quasidisk.
\end{itemize}
\end{thmx}
See page~\pageref{def:rhonice} for the definition of $\delta$-nice.
\subsubsection{QC control along real pullbacks}

Since we work with asymptotically holomorphic extensions of
smooth mappings, we need a tool to the dilatation. 
Let us give
informal statements of our results here.
We control the dilatation of pullbacks along
the real line that do not pass through long cascades of central
returns, see Propositions \ref{prop:sum of squares, puzzle pieces} 
and \ref{prop:squares}. 
Here we build on the work of \cite{Li-Shen},
and we use the real bounds of \cite{vSV}
as we pullback from one non-central return to
the next.
Roughly Proposition~\ref{prop:sum of squares,
  puzzle pieces}
can be stated as:
\begin{thmx}\label{thmx:sq}
Suppose that $F\colon\bm{\mathcal{U}}\rightarrow\bm{\mathcal{V}}$
is a complex box mapping.
Suppose that $\{\bm G_j\}_{j=0}^s$ is a chain of puzzle pieces that
does not contain a long cascade of central returns, then
$$\sum_{j=0}^s\diam(\bm G_j)^2<C,$$
and  $C\rightarrow 0$ as $\max_{\V\subset\bm{\mathcal{V}}}\diam
(\V)\rightarrow 0,$ where the maximum is taken over connected components 
$\V$ of $\bm{\mathcal{V}}.$
\end{thmx}
A central return in a chain of puzzle pieces occurs when 
there exists a critical point $c$ of $f$ such that
$c\in \bm G_{i_2}\subset\bm G_{i_1}\subset \bm G_{i_0}$,
$\bm G_{i_2}$ is a first return domain to $\bm G_{i_1}$,
$\bm G_{i_1}$ is a first return domain to $\bm G_{i_0}$ and
$R_{\bm G_{i_0}}|_{\bm G_{i_1}}(c)\in \bm G_{i_1}$.
In other words, the return time of $\bm G_{i_2}$ to
$\bm G_{i_1}$ is the same as the return time of 
$\bm G_{i_1}$ to $\bm G_{i_0}.$

Assume that $V^m\supset V^{m+1}\supset \dots\supset V^{m+n+1}$ is
a central cascade of intervals.
Suppose that $(a_0,b_0)$ is a component of $\mathrm{int}(V^m\setminus V^{m+1})$, 
and let $(a_k,b_k), k\leq n$ be a component of
$\mathrm{int}(V^{m+k}\setminus V^{m+k+1})$,
such that $R_{V^m}^k\colon(a_k,b_k)\rightarrow (a_0,b_0)$ is a
diffeomorphism. 
An important step in the proof of
Theorem~\ref{thmx:sq} is controlling the sum 
$\sum_{j=0}^k\diam(\bm G_j)^2$ when $\bm G_0\cap (a_k,b_k)\neq
\emptyset$
and $k$ is arbitrarily large.  
To do this we use Theorem~\ref{thmx:good central}, which implies that
$\bm G_k$ is contained in a Poincar\'e disk $D_{\theta}(a_0,b_0)$ 
with $\theta$ bounded away from zero, and then the bound on 
$\sum_{j=0}^k\diam(\bm G_j)^2$ follows.

\subsubsection{The good nest.} We adapt the construction of
good nest of puzzle pieces \cite{AKLS} to the multicritical setting
and prove complex bounds for this nest of puzzle pieces, see the
definition on page~\pageref{def:good nest}, Proposition~\ref{prop:good
bounds}, Corollary~\ref{cor:Wnice} and Proposition~\ref{prop:goodbg}. 
The next theorem, which summarizes these results,
gives complex bounds on a much finer nest of puzzle
pieces than the enhanced nest.

\begin{thmx}
There exists $\delta>0$ such that the following holds.
Suppose that $F:\bm{\mathcal{U}}\rightarrow\bm{\mathcal{V}}$ is a
complex box mapping. Let
$$\V^0\supset \W^0\supset\V^1\supset\W^1\supset \V^2\supset\dots$$
be the good nest of puzzle pieces about $\crit{F}$. Then for each
$n\geq 0$ and each critical point $c$ of $F$ the following hold:
\begin{itemize}
\item $\mod(\V^n_c\setminus \W^n_c )\geq \delta$,
\item $\V^n_c$ has bounded geometry at the critical point,
\item $\W^n$ is $\delta$-nice.
\end{itemize}
\end{thmx}
See page~\pageref{def:bddgeo} for the definition of $\delta$-bounded
geometry and page~\pageref{def:rhonice} for the definition of
$\delta$-nice.

\subsubsection{QC\textbackslash BG partitions}

Motivated by the QC Criterion of \cite{KSS-rigidity}, we 
say that a mapping is $(K,\delta)$-qc\textbackslash bg
if it is $K$-\emph{qc} except on a set where one has \emph{bg}
(bounded geometry, depending on $\delta$).
See page~\pageref{page:def qcbg} for the
precise definition.

\medskip

\begin{figure}
\centering
\begin{subfigure}{.5\textwidth}
  \centering
  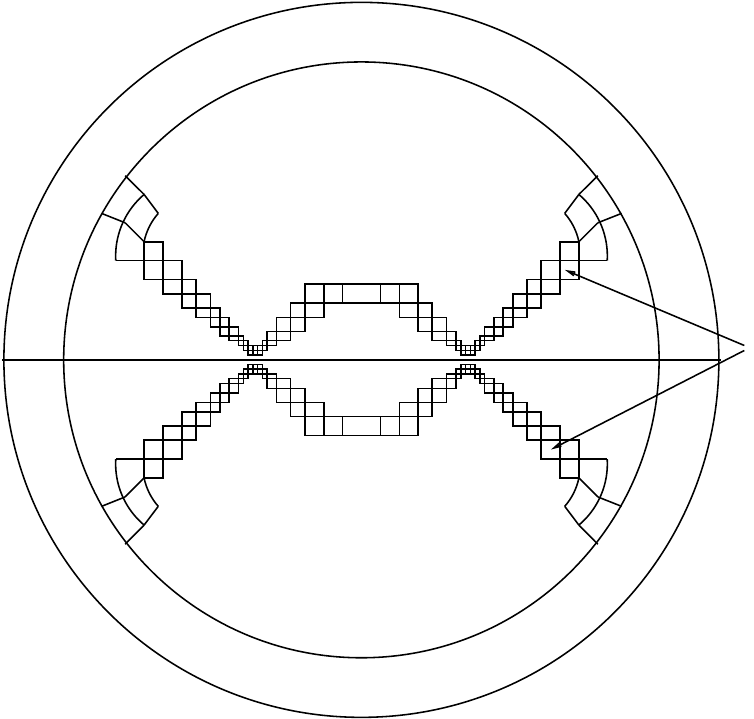
  \caption{A qc\textbackslash bg partition of a saddle-node
  cascade. }
\end{subfigure}%
\begin{subfigure}{.5\textwidth}
  \centering
  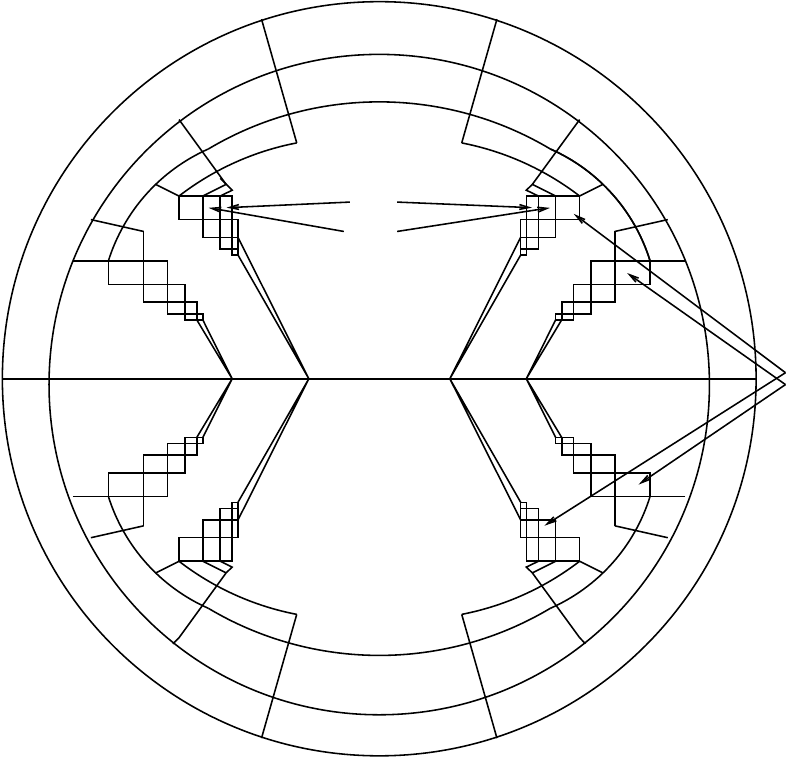
  \caption{A qc\textbackslash bg partition after a renormalization,\\
    $F\colon\V^1\rightarrow\V^0$.}
 \label{fig:renormalization}
\end{subfigure}
\caption{$X_1$ are pieces with bounded geometry, $X_2$ are pieces
  where the qc mapping is defined dynamically and $X_3$ are pieces
  that are $K_1$-quasidisks. The thin triangles in (B) are built from small
    topological squares, but we do not show that level of detail.
In (B), we can avoid using the set $X_1$
    since we can show that the boundary maps are uniformly quasisymmetric}
\label{fig:qcbg}
\end{figure}

Theorem~\ref{thmx:good central} is an important step 
in the proof of a QC\textbackslash BG Partition of a central cascade.
Let us give an informal statement here.
See Theorem~\ref{thm:central cascades} for a precise
statement and Figure~\ref{fig:qcbg} for an illustration.
A similar partition is also used to glue together
conjugacies that are made between renormalization
levels of an infinitely renormalizable mapping,
see Lemmas \ref{lem:pl decomp} and \ref{lem:renintersection}.




\begin{thmx}[QC\textbackslash BG Partition of a Central Cascade]
\label{thmx:qcbgcen}
There exist $K_1,K\geq 1$ and $\delta>0$ such that the
following holds:
Let $\bm E\owns c_0$ be a puzzle piece given by Theorem~\ref{thmx:good
  central}, and let $\widetilde{\bm{E}}$ be the corresponding piece for
$\tilde f$. There exists a $(K,\delta)$-qc\textbackslash bg mapping 
$H\colon\bm E\rightarrow\widetilde{\bm{E}}$
and a partition of $\bm E$ into subsets $X_0\cup X_1\cup X_2\cup X_3\cup
\bm{\mathcal{O}},$ 
and similarly for the objects marked with a tilde,
such that the following holds:
\begin{itemize}
\item $H(X_i)=\tilde X_i$ for $i=0,1,2,3$, and
  $H(\bm{\mathcal{O}})=\widetilde{\bm{\mathcal{O}}}.$
\item $\bm{\mathcal{O}}$ is a neighbourhood of the set of critical
  points of $R_{\bm E}|_{\mathcal{L}_{c_0}(\bm E)},$ and $H|_{\bm{\mathcal{O}}}$ is $K$-qc.
\item $X_0$ has measure 0 and (\ref{eqn:qc bound 0}) holds at each
  $x\in X_0$.
\item If $\bm P$ is a connected component of $X_1$ and
  $\widetilde{\bm{P}}=H(\bm P)$ then both $\bm P$ and $\widetilde{\bm{P}}$ have
  $\delta$-bounded geometry and 
$$\mod(\bm E\cap\mathbb H^+\setminus\bm P)\geq \delta\quad\mbox{and}\quad
\mod(\widetilde{\bm{E}}\cap\mathbb H^+\setminus\widetilde{\bm{P}})\geq\delta.$$
\item $H|_{X_2}$ is $K$-quasiconformal, $H\circ F(z)=\tilde F\circ
  H(z)$ for $z\in X_2,$ and $X_2$ contains $E=\bm E\cap\mathbb R\setminus\bm{\mathcal{O}}.$
\item Each component $\bm W$ of $X_3$ is a $K_1$-quasidisk and there
  exists a $K$-qc mapping from $\bm W$ onto $\widetilde{\bm{W}}.$
\end{itemize}
\end{thmx}

\subsubsection{A real, qc/bg, speading principle}
Since we cannot control the dilatation of high iterates of the
complex extension of $f$, we are unable to use the usual
Spreading Principle, see \cite{KSS-rigidity}, to glue together 
qc mappings defined near the post-critical set into
global qc mappings. Instead, we start with a
$K_0$-qc external
conjugacy, see page~\pageref{def:extconj},
$H_0\colon\mathbb C\to\mathbb C$ that conjugates
$F\colon\UU\rightarrow\VV$ and 
$\tilde F\colon\widetilde{\UU}\rightarrow\widetilde{\VV}$
on $\partial\UU$ and refine it by pulling back.
In order to pull $H_0$ back, we have to ``move''
the points $H_0(F(c)), c\in\Crit(F)$ so they match
$\tilde F(\tilde c)$. When we do this, the dilatation of 
resulting map could be arbitrarily large. 
However, by the complex bounds for the good nest
proved in Section~\ref{sec:goodnest},
we can show that the region where the dilatation is uncontrolled
has bounded geometry and moduli bounds. This is done in the proof
of Proposition~\ref{prop:realspreading}, which we summarize now:

\begin{thmx}[A Real Spreading Principle]\label{thmx:rsp}
There exist arbitrarily small nice
open complex neighbourhoods $\bm{\mathcal{O}}$
of $\Crit(F)$ and a mapping $H_{\bm{\mathcal{O}}}\colon\mathbb C\rightarrow\mathbb C$
such that 
\begin{itemize}
\item $H_{\bm{\mathcal{O}}}(\bm{\mathcal{O}})=\widetilde{\bm{\mathcal{O}}};$
\item $H_{\bm{\mathcal{O}}}$ agrees with the
boundary marking on $\partial\bm{\mathcal{O}};$
\item $H_{\bm{\mathcal{O}}}(\Dom'(\bm{\mathcal{O}}))=\Dom'(\widetilde{\bm{\mathcal{O}}});$
\item for each component $\U$ of $\Dom'(\bm{\mathcal{O}}),$
$H_{\bm{\mathcal{O}}}\circ F(z)=\tilde F\circ H_{\bm{\mathcal{O}}}(z)$ for $z\in\partial\U$.
\end{itemize}
Moreover, 
$H_{\bm{\mathcal{O}}}\colon\mathbb{H}^+
\rightarrow \mathbb{H}^+$
is a $(K,\delta)$-qc\textbackslash bg mapping.
\end{thmx}

By the QC Criterion, see Theorem~\ref{thm:qc-criterion},
Theorem~\ref{thmx:rsp}, implies that there exists
a $K'=K'(K,\delta)$-qc mapping of the plane whose restriction to
the real line
agrees with $H_{\bm{\mathcal{O}}}|_{\mathbb R}.$

\section{Ingredients of the paper and a sketch of the proof}
Before we get into the details, let us point out 
that in Section~\ref{sec:guide} we give a short guide to the paper
for the reader interested in specific results.

\subsection{Old and new ingredients of this paper}
	Our paper not only 
	builds on earlier work on rigidity of polynomials, \cite{KSS-rigidity},
	\cite{KSS-density} and \cite{Koz_ax}, but also relies heavily on ideas from
	\cite{LevStr:invent} and  uses that one has complex and real bounds in the generality required 
	in this paper, see the recent paper \cite{CvST}.
	Let us list the issues that were not present in earlier works on this topic
	for which we will need to develop new tools to obtain our
        results.

\begin{table}[h]
\begin{tabular}{|p{1.5in}|p{1.75in}|p{2.75in}|}
\hline
\multicolumn{1}{|p{1.5in}}{\center\textbf{Polynomials}} & 
\multicolumn{1}{|p{1.75in}}{\center\textbf{Analytic mappings}}  & 
\multicolumn{1}{|p{2.75in}|}{\center\textbf{Smooth mappings}}\\
\hline
\multicolumn{1}{|p{1.5in}}{\center B\"ottcher Coordinate} &            
\multicolumn{2}{|p{4.5in}|}{\center Touching box mappings}\\
\hline

\multicolumn{1}{|p{1.5in}}{\center Yoccoz puzzle pieces} &            
  \multicolumn{2}{|p{4.5in}|}{\center Complex box mappings}
    \\
\hline

\multicolumn{2}{|p{3.25in}}{\center The enhanced nest} &            
  \multicolumn{1}{|p{2.75in}|}{\center
    The good nest to prove a spreading principle,
 and the enhanced nest over a long central
                                                      cascade}
\\
\hline
\multicolumn{2}{|p{3.25in}}{\center Puzzle pieces in the enhanced nest have complex bounds} &            
  \multicolumn{1}{|p{2.75in}|}{\center    
Puzzle pieces in the good nest have complex bounds
}

\\

\hline
\multicolumn{2}{|p{3.25in}}{\hspace{1in}} &            
  \multicolumn{1}{|p{2.75in}|}{\center   
Puzzle pieces in the enhanced nest over a long central cascade have complex bounds and they are quasidisks
}
\\
\hline

\multicolumn{2}{|p{3.25in}}{\center The Spreading Principle} &            
  \multicolumn{1}{|p{2.75in}|}      
   {\center A Real Spreading Principle, and a QC\textbackslash BG Partition of a central cascade}
     \\
\hline

\multicolumn{3}{|p{6in}|}{\center A QC\textbackslash BG Partition for a
(quasiregular) polynomial-like mapping} \\ 
\hline
\end{tabular}%
\vspace{0.25in}
\caption{This table indicates the changes to the proof
of rigidity for polynomials in \cite{KSS-rigidity} to prove rigidity
for analytic and smooth mappings.}
\end{table}

\subsubsection{New challenges in the polynomial case}
\begin{enumerate}[leftmargin=*] 
\item \textit{Odd critical points.} 
Lack of dynamical symmetry around an odd critical point
	complicates the construction of complex box mappings, even
	near a reluctantly recurrent critical point; we treat this situation in this paper.
	In the persistently recurrent case, we rely on the complex bounds proved in
	\cite{CvST}, which also deals with odd critical points.
		
\item \textit{Polynomials with non-real critical points.} 
Even if the maps are real polynomials with non-real critical points,
	qs-rigidity on the real line was not known for infinitely
        renormalizable maps. 
        In this
	case, it is not known (without using complex bounds)
if it is possible to find critical puzzle pieces that contain
	only real critical points, see \cite{KvS}.
Even in the non-renormalizable case the maps
	$f$ and $\tilde f$ will, in general, not be conjugate on the complex plane,
	and therefore the proof of \cite{KvS} does not immediately apply.
	However,
	we prove qs-rigidity of the real part of dynamics.
\end{enumerate}

\subsubsection{New challenges in the analytic case}

\begin{enumerate}[leftmargin=*] 
\item \textit{No B\"ottcher Coordinate.} If the maps are real-analytic, but not polynomial, then there
are no B\"ottcher coordinates, and therefore we cannot use the usual
construction of the Yoccoz puzzle. 
In the analytic case, the strategy of the proof is to combine the proof of 
qc-rigidity of \cite{KSS-rigidity} and the complex bounds of
\cite{CvST} with the proof of
rigidity for maps with a
non-minimal post-critical set of
\cite{LevStr:invent} to provide topologically
conjugate maps $f$ and $\tilde{f}$ with compatible 
``external structures."
We use
touching box mappings,
see Section~\ref{subsec:global touching box mappings},
to play the role of the B\"ottcher coordinates
and complex box mappings to play the role of the Yoccoz puzzles 
in our setting. Gluing local information together is more difficult,
since the puzzle pieces given by
complex box mappings do not globally have the same Markov properties
as Yoccoz puzzle pieces. In particular, we need to construct complex box mappings
whose ranges are Poincar\'e lens domains to make it possible to transfer
the global information given by a touching box mapping close to a critical point.

\item\textit{Yoccoz puzzle pieces are automatically nested or
    disjoint; this is not the case for puzzle pieces constructed by
    hand.} 
For polynomials, the case when one critical point $c$ accumulates on 
another critical point $c'$, for example on a persistently recurrent critical point,
does not cause difficulties. While here, we need to consider complex box mappings
containing all, possibly infinitely many, components of the domain of the first entry map to a puzzle piece
containing $c'$.  To deal with this,
we need to control the geometry of infinitely many puzzle pieces
near $c'$ at a given level.

	\item\textit{Parabolic points.} To construct qc-conjugacies near parabolic points, we need to consider
	touching box mappings where the domain and range touch
        tangentially,
see Section~\ref{subsec:global touching box mappings}. 
	If we did not consider maps with parabolic points, we could
	use regions that are quasi-circles to play the role of
	B\"ottcher coordinates 
	as in \cite{LevStr:invent}.
	Here, we will have to use a precise local analysis near parabolic points
	to ensure that we obtain a qc-conjugacy, respecting the puzzle pieces, near a
	parabolic point.
\end{enumerate}

\subsubsection{New challenges in the $C^3$ case}
\begin{enumerate}[leftmargin=*] 
\item\textit{The usual spreading principle fails.} For maps that are $C^3$, but not analytic, we make use of asymptotically
holomorphic extensions on the real-line. Under high iterations, these mappings
are far from holomorphic. 
One consequence of this is that we cannot use the usual
Spreading Principle, 
Proposition~\ref{prop:spreading analytic}.
We use Propositions \ref{prop:sum of squares, puzzle pieces} and
\ref{prop:squares} to control the dilatation of $f$
along certain chains,
and we use Theorem~\ref{thm:central cascades},
to treat long cascades of central returns. Even when there are no long
cascades of central returns, we cannot proceed as in
\cite{KSS-rigidity}, since we can only the dilatation of $f$ 
through real pullbacks.

\item \textit{The good nest.} Without the usual Spreading Principle,
it is difficult to ``glue together'' conjugacies defined on the landing
domain to a neighbourhood of the critical points into a global qc
mapping. We develop the good nest for multicritical
mappings. It was introduced for unicritical mappings in \cite{AKLS}.
We are able to use this nest 
to prove a Real Spreading Principle for smooth mappings,
see Theorem~\ref{thmx:rsp},
since between any
two levels of the good nest, there can be at most one central cascade.
We prove complex bounds for the good nest in Section~\ref{sec:goodnest}.

Even when there are no long cascades of central returns, the mapping
(a pseudo-conjugacy)
that we construct will have unbounded dilatation. However, the set where 
the dilatation cannot be controlled will have good geometric
properties - it is contained in a union of puzzle pieces with bounded
geometry and appropriate moduli bounds. In particular, such puzzle 
pieces which do not intersect the real line are well-inside the upper
half plane, but this is not the case for every puzzle piece,
see Figure~\ref{fig:wellinside}.
\begin{figure}
\resizebox{0.8\textwidth}{!}{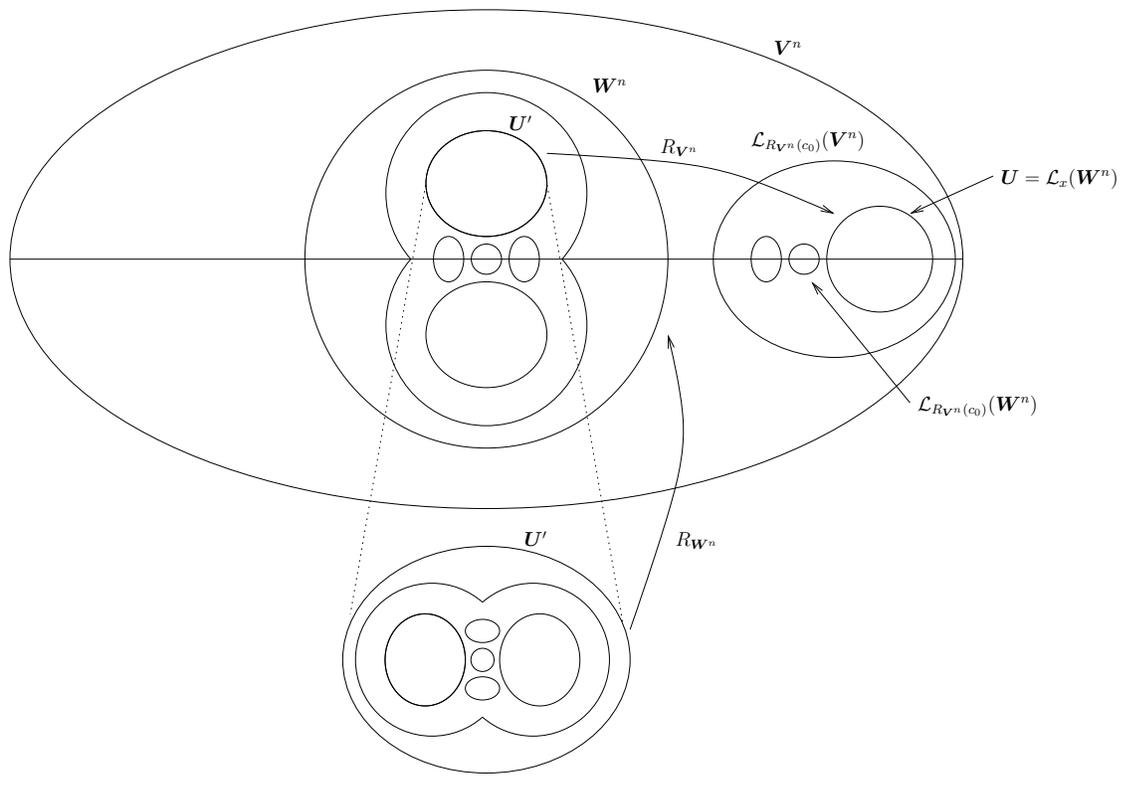}
 \caption{The domain $\U'=(R_{\V^n}|_{\W^n})^{-1}(\U)$
need not be well-inside $\mathbb H^+$;
however, the components of $(R_{\W^n}|_{\U'})^{-1}(\Dom(\W^n))$
are well-inside $\U'$, since $\W^n$ is $\delta$-nice.}
 \label{fig:wellinside}
\end{figure}

\item \textit{Long cascades of central returns.}
Let us describe the problem caused by long cascades of
central returns.
Suppose that $\V^0\supset \V^1\owns c_0$
and
$R_{\V^0}\colon\V^1\rightarrow\V^0$ has a long 
saddle-node cascade. 
Let $\V^{i+1}=\comp_{c_0}R_{\V^0}^{-1}(\V^i)$ be the central
puzzle pieces in the principal nest about $c_0$.
Then there are no 
fixed points of $R_{\V^0}$ on the real line, and
$R_{\V^0}(c_0)\in\V^N$ for some $N$ large.
We have that $R_{\V^0}$ has periodic orbits that are far from the
real line.
\begin{figure}[htbp]
\begin{center}
\includegraphics[width=3in]{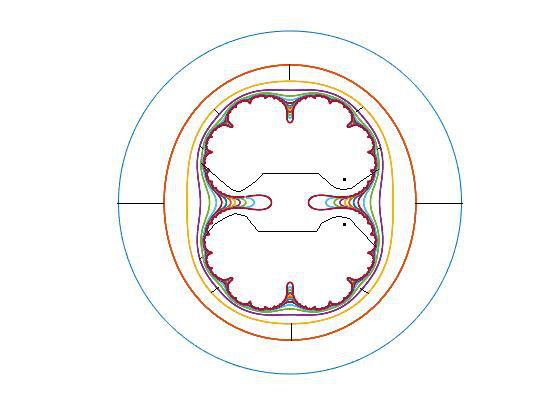}
\end{center}
\caption{Pullbacks of a disk by a qr mapping that is close to
  $z\mapsto z^2+1/4+\varepsilon$. The return mapping is 
non-conformal at the complex fixed points, and this non-conformality 
will blow up under high iterates. The fixed points in the complex
plane and the boundaries of the components of $X_3$ are indicated,
compare
Figure~\ref{fig:qcbg}.}
\label{fig:nearpara}
\end{figure}
It is useful to think of 
$z\mapsto z^2+1/4+\varepsilon$, see Figure~\ref{fig:nearpara}.
Because of these periodic points,
the diameters of the $\V^i$ are bounded from below, and so
the dilatation of the mapping
$R_{\V^0}^i\colon(\V^i\setminus\V^{i+1})\rightarrow \V^{0}\setminus \V^1$,
$0\leq i<N$ can diverge for large $i$ and $N$.
If we are in the situation (which may arise infinitely
many times in the proof of rigidity) 
that $H\colon \mathbb C\rightarrow\mathbb
C,$ 
with $H(\V^0)=\widetilde{\V}^0,$ and
$H(\Dom^*(\V^0))=\Dom^*(\widetilde{\V}^0)$ is a $K$-qc mapping.
If we define 
$\hat H\colon \V^i\setminus
\V^{i+1}\rightarrow\widetilde{\V}^i\setminus\widetilde{\V}^{i+1}$
for $0\leq i\leq N$ by the formula  $H\circ R^i(z)=\tilde R^i\circ\hat
H(z), z\in \V^i\setminus \V^{i+1}$,
then the dilatation of $\hat H$ could diverge as $N\rightarrow
\infty$.
This means that we cannot construct a qc conjugacies between $f$ and
$\tilde f$
by ``pulling back'' (even if we knew that the external conjugacy was
a conjugacy on the post-critical set).
Because of this, a new strategy for proving rigidity in the smooth
case is required. We use a geometric construction of a
qc\textbackslash bg mapping,
see Theorem~\ref{thm:central cascades}.
\end{enumerate}

\subsubsection{Different nests of puzzle pieces}
For analytic mappings our proof builds on the
proof of qc rigidity of real polynomials, \cite{KSS-rigidity}.
In that paper, the authors introduced the enhanced nest.
The enhanced nest is defined for return mappings 
$F:\UU\to\VV$ defined
in 
a neighbourhood of 
$\crit(f)\cap\omega(c_0)$ where $c_0$ is
a persistently recurrent critical point.
The enhanced nest goes to deep levels very efficiently,
and a puzzle piece $\II_n$ in the enhanced nest 
has the property of of being $\delta$-{\em free}: 
there exist $\delta>0$, an annulus $\bm A^+$ 
with inner boundary
$\partial\II_n$ and an annulus $\bm A^-$ with 
outer boundary $\partial\II_n$ such that 
$\mod(\bm A^+),\mod(\bm A^-)\geq \delta,$ and
$\overline{\bm A^+\cup\bm A^-}\cap\PC(F)=\emptyset.$
This control on the post-critical set makes it possible to prove
complex bounds for the nest of puzzle pieces, \cite{KSS-rigidity, CvST}. 
Moreover,
it makes it possible to prove that the puzzle pieces in the
enhanced nest are quasidisks, Proposition~\ref{prop:quasidisks}.
It is easier to deal with critical points that are reluctantly 
recurrent, that is, not persistently recurrent.
If $f$ is reluctantly recurrent at $c_0$, then 
there are arbitrarily small neighbourhoods of $c_0,$
which are mapped up to some definite scale with bounded degree.
As in \cite{KSS-rigidity} using these two nests of 
puzzle pieces about the critical points of $f$, we
prove rigidity for analytic mappings.

For smooth mappings it seems to be difficult to use the 
enhanced nest to prove rigidity. Between any two levels of
the enhanced nest there can be many distinct
cascades (different return mappings), and
the geometry of these intermediate puzzle pieces 
can not be controlled. The control we have for
the enhanced nest is eventually lost, since the degree
of the mapping from a pullback of $\II_n$ nested between
$\II_n$ and $\II_{n+1}$ could be arbitrarily large.
To deal with this
we introduce the good nest for multicritical mappings,
see \cite{AKLS} for the unimodal case, and prove
complex bounds for this nest. This nest goes to 
deeper levels much more slowly than the enhanced
nest. Indeed there can be at most one central cascade
between any two levels of the good nest.

It does not seem to be possible for us to use the good nest 
alone. In order to prove {\em a priori} bounds for the good nest,
we need to be able to control the dilatation of the mapping from
an arbitrarily deep puzzle piece in the good nest up to the top level.
We are able to prove such a statement for the enhanced nest since they
are controlled by Poincar\'e disks:
there exist $\theta\in(0,\pi/2)$ and $\delta>0$ so that
$\II_n\subset D_{\theta}((1+2\delta)I_n)$.
Before proving moduli bounds,
this is the only tool that we have to control the dilatation
of high iterates of $f.$ Since critical values of $f$ 
can be arbitrarily close to the boundaries of puzzle pieces in the 
good nest construction, puzzle 
pieces in the good nest may not be contained in such Poincar\'e disks.
In order to overcome this, we prove Proposition~\ref{prop:sum of squares, puzzle pieces} which gives us the control on the
dilatation that we need. Proposition~\ref{prop:sum of squares, puzzle pieces} relies on 
knowing that we can control the dilatation of arbitrarily 
long pullbacks through a real branch of a central cascade, 
and to do this 
we need to be able to control puzzle pieces at the top of a 
central cascade with Poincar\'e disks. This is given by
Proposition~\ref{prop:good geometry for central cascades}.
The puzzle piece given by 
Proposition~\ref{prop:good geometry for central cascades}
comes from a nest that is closely related to the enhanced nest,
which
we refer to as the {\em enhanced nest over a central cascade}.
These puzzle pieces have similar geometric properties to the puzzle 
pieces in the enhanced nest; for example,
they are also $K$-quasidisks.
We also
use Proposition~\ref{prop:good geometry for central cascades}
to carry out a geometric construction of a qc\textbackslash bg
mapping between a puzzle piece at the top
of central cascade of $f$ and 
the corresponding puzzle piece for $\tilde f,$ see
Theorem~\ref{thm:central cascades}.
In this construction it is important that the puzzle piece
at the top of the central cascade is a
$K$-quasidisk.

\begin{table}[h]
\begin{tabular}{|p{1.5in}|p{0.8in}|p{0.4in}|p{0.9in}|p{2in}|}
\hline
\multicolumn{1}{|p{1.5in}}{\center{}}&
\multicolumn{1}{|p{0.8in}}{\center $\delta$-nice and
                                                $\delta$-bounded
                                                geometry}&
\multicolumn{1}{|p{0.4in}}{\center $\delta$-free} &
\multicolumn{1}{|p{0.9in}}{\center $K$-quasidisk}&
\multicolumn{1}{|p{2in}|}{\center Notes}\\
\hline
\multicolumn{1}{|p{1.5in}}{\center Enhanced Nest}&
 \multicolumn{1}{|p{0.8in}}{\center \checkmark}&
\multicolumn{1}{|p{0.4in}}{\center \checkmark}&
\multicolumn{1}{|p{0.9in}}{\center \checkmark} &
\multicolumn{1}{|p{2in}|}{Goes to deep levels efficiently.
Used to prove qs rigidity for
analytic complex box mappings at a persistently
                                                 recurrent critical point.}
\\
\hline

\multicolumn{1}{|p{1.5in}}{\center Nest about reluctantly recurrent critical points}&
 \multicolumn{1}{|p{0.8in}}{\center \checkmark}&
\multicolumn{1}{|p{0.4in}}{\center \text{\sffamily X}}&
\multicolumn{1}{|p{0.9in}}{\center \checkmark} &
\multicolumn{1}{|p{2in}|}{Used to prove rigidity for
analytic
                                                 complex box mappings
                                                 at a reluctantly
                                                 recurrent critical point.}\\
\hline

\multicolumn{1}{|p{1.5in}}{\center Good Nest}&
 \multicolumn{1}{|p{0.8in}}{\center \checkmark}&
\multicolumn{1}{|p{0.4in}}{\center \text{\sffamily X}}&
\multicolumn{1}{|p{0.9in}}{\center \text{\sffamily X}} &
\multicolumn{1}{|p{2in}|}{Goes to deep levels slowly. Used to
                                                         prove a
                                                         spreading
                                                         principle for
                                                         $C^3$ mappings.}
\\
\hline
\multicolumn{1}{|p{1.5in}}{\center Enhanced nest over a central cascade}&
 \multicolumn{1}{|p{0.8in}}{\center \checkmark}&
\multicolumn{1}{|p{0.4in}}{\center \text{\sffamily X}}&

\multicolumn{1}{|p{0.9in}}{\center \checkmark} &
\multicolumn{1}{|p{2in}|}{Used for $C^3$ mappings to control dilatation of
                                                 pullbacks along a
                                                 real branch of a long
                                                 central cascade
                                             and in the construction of a
                                                 qc\textbackslash bg
                                                 partition of a
                                                 central cascade.}
\\
\hline

\end{tabular}
\vspace{0.25in}
\caption{Properties of the different nests of puzzle pieces.}
\end{table}

\subsubsection{Gluing in the infinitely renormalizable case.}
The gluing argument in \cite{KSS-rigidity} Section 7 takes a significant amount
of effort beyond what is needed in the finitely renormalizable case.
Our method of dealing with the infinitely renormalizable case is an 
extension of our qc\textbackslash bg technique for 
dealing with long central cascades. Figure~\ref{fig:renormalization}
indicates how this is done.

\subsection{The QC Criterion and qc\textbackslash bg mappings.}
We recall the following result from \cite{KSS-rigidity}:
\begin{thm}[QC Criterion \cite{KSS-rigidity} Lemma 12.1]
\label{thm:qc-criterion}
For any constants $K\geq 1$ and $\varepsilon>0$,
there exists a constant $K'\geq 1$ with the following property.
Let $\phi\colon \bm \Omega\rightarrow\widetilde{\bm{\Omega}}$
be a homeomorphism between two Jordan domains, which extends
continuously to the boundary.
Assume that $X_0\cup X_1\cup X_2$ and 
$\tilde X_0\cup\tilde X_1\cup \tilde X_2$ 
are partitions of $\bm \Omega$ and $\widetilde{\bm{\Omega}}$, respectively, 
such that the following hold:
\begin{enumerate}
\item $\phi(X_i)=\tilde X_i$ for $i=0,1,2$.
\item If $\bm P$ is a component of $X_1$, then both
$\bm P$ and $\phi(\bm P)$ have $\varepsilon$-bounded geometry,
and moreover
$$\mathrm{mod}(\bm \Omega\setminus \bm P)\geq \varepsilon,\;
\mathrm{mod}(\widetilde{\bm{\Omega}}\setminus\phi(\bm
P))\geq\varepsilon.$$
\item $X_0$ has zero measure and for each $x\in X_0$,
\begin{equation}\label{eqn:qc bound 0}
\liminf_{r\rightarrow 0}\frac{\sup_{|y-x|=r}
    |\phi(y)-\phi(x)|}
{\inf_{|y-x|=r} |\phi(y)-\phi(x)|}<\infty.
\end{equation}
\item For every $x\in X_2$
\begin{equation}\label{eqn:qc bound 2}
\liminf_{r\rightarrow 0}\frac{\sup_{|y-x|=r}
    |\phi(y)-\phi(x)|}
{\inf_{|y-x|=r} |\phi(y)-\phi(x)|}\leq K
\end{equation}
\end{enumerate}
Then there exists a $K'$-qc map 
$\psi\colon \bm \Omega\rightarrow\widetilde{\bm{\Omega}}$
such that $\psi=\phi$ on $\partial\bm\Omega$.
	\end{thm}
\begin{rem}
We refer to \cite{KSS-rigidity} Lemma 12.1 as the QC Criterion.
\end{rem}
With the QC Criterion as motivation, we say that 
a mapping $\phi\colon\bm\Omega\rightarrow\widetilde{\bm{\Omega}}$ is a
$(K,\varepsilon)$-\emph{qc\textbackslash bg mapping}
\label{page:def qcbg}
if there is a partition, as in Figure~\ref{fig:qcbg}, 
of $\bm \Omega$ into sets $X_0,X_1,X_2,$
and a partition of  $\widetilde{\bm{\Omega}}$ into sets $\tilde X_0,\tilde
X_1,\tilde X_2,$
 such that
\begin{itemize}
\item $\phi(X_i)=\tilde X_i$ for $i=0,1,2$,
\item $X_0$ has measure zero, and satisfies (\ref{eqn:qc bound 0}) on $X_0$.
\item $X_1$ can be partitioned into open sets $\bm P$ such that 
 both
$\bm P$ and $\phi(\bm P)$ have $\varepsilon$-bounded geometry,
and moreover
$$\mathrm{mod}(\bm \Omega\setminus \bm P)\geq \varepsilon,\;
\mathrm{mod}(\widetilde{\bm{\Omega}}\setminus \phi(\bm P))\geq\varepsilon.$$
\item $\phi$ satisfies (\ref{eqn:qc bound 2}) on $X_2$.
\end{itemize}

We refer to such a partition as a \emph{qc\textbackslash bg partition.}
We believe that qc\textbackslash bg mappings will be a useful idea in 
other settings. 
We will say that 
Suppose that $\bm \Omega'\subset\bm \Omega$. We say that
$\phi\colon\bm \Omega\rightarrow\widetilde{\bm{\Omega}}$ 
is a \emph{qc\textbackslash bg mapping with moduli bounds in $\bm\Omega'$}
if it is a qc\textbackslash bg mapping and 
$\mathrm{mod}(\bm \Omega'\setminus \bm P)\geq \varepsilon$ and
$\mathrm{mod}(\phi(\bm{\Omega}')\setminus \phi(\bm P))\geq\varepsilon$
for each component $\bm P$ of $X_1$ which is contained in $\bm \Omega'$.

One may read the name qc\textbackslash bg as \emph{quasiconformal
  except on a set with bounded geometry.}
The QC Criterion implies that if 
$\phi\colon\bm\Omega\rightarrow\widetilde{\bm{\Omega}}$ is a 
$(K,\varepsilon)$-qc\textbackslash bg mapping,
then there exists a $K'$-qc mapping
$\phi'\colon\bm\Omega\rightarrow\widetilde{\bm{\Omega}}$ that has the same
boundary values as $\phi$. 
To prove Theorem~\ref{thm:main} for smooth mappings
we are going to construct a
sequence of
$(K,\varepsilon)$-qc\textbackslash bg mappings
on the upper half-plane, whose extensions to the real line
converge to a conjugacy between $f$ and $\tilde f$.
The qc\textbackslash bg mappings we construct will not be defined
dynamically everywhere in their domains. Up to sets of measure zero,
their domains 
will consist of a set where the mappings are dynamically defined, 
a set that can be decomposed into quasidisks, and a set separating
these two sets whose
components have bounded geometry and moduli bounds.
This level of flexibility makes it possible for us to deal with 
central cascades. We will describe this
in the next subsection.

\subsection{Sketch of the proof and organization of the paper}
\label{subsec:sketch}
Let us now outline the structure of the paper, and give
more detail about some of the aspects of the proof.
 		
\subsubsection{Developing the required tools, Sections
  \ref{sec:definitions} and \ref{sec:real dynamics}}
In Section \ref{sec:definitions} we give some background definitions.
Specifically we discuss quasiconformal and quasisymmetric mappings,
asymptotically holomorphic extensions of $C^3$ maps, and
Poincar\'e disks and their pullbacks by asymptotically holomorphic
mappings.
In Section \ref{sec:real dynamics}
we deal with the necessary real dynamics. 
In Subsection \ref{subsec:decomp} we give a hierarchical decomposition of the 
set of critical points of a map $f$ and describe a nice real neighbourhood of 
the set of critical points. We will use these often. Finally, we present the
real bounds that we will need to construct complex box mappings.

\subsubsection{Compatible complex box mappings,
Section~\ref{sec:complex bounds}}
A few types of complex box mappings play vital roles in the proof.
When combined, they
perform the same function as the Yoccoz puzzle for polynomial mappings.
Using the hierarchical decomposition of the critical points
and
a chosen nice real neighbourhood $\mathcal I$ of the critical set of $f$, see 
Section \ref{subsec:decomp}, we build complex box mappings about
certain critical points. 
The complex box mapping at
each critical point needs to be compatible
with the box mappings around other critical points.
This compatibility is obtained using touching box mappings,
see below.
\medskip

We will make use of the usual complex box mappings which are complex
extensions of the real
return maps to neighbourhoods of the set of critical points of $f$.
These were constructed in \cite{CvST} about critical points $c_0$ 
with the property that $f$ is persistently recurrent on $\omega(c_0)$.
In the reluctantly recurrent case, we construct complex box mappings
in Subsection \ref{subsubsec:box mappings reluctant}. This was
done for analytic maps with all critical points of even order in
\cite{KSS-density},
and here we allow for odd critical points.
In the reluctantly recurrent case, the  complex bounds are obtained from 
the real bounds at reluctantly recurrent critical points,
and
the fact that we can go from any small scale up to a fixed
scale by a mapping of bounded degree,
see Proposition~\ref{prop:real bounds reluctant}
and Theorem~\ref{thm:box mapping reluctant}.
In the persistently recurrent case, the necessary bounds were proved in
\cite{CvST}, see Theorem~\ref{thm:box mapping persistent}
and Theorem~\ref{thm:box mapping persistent infinite branches}.

Complex box mappings are very useful for pulling back qc-conjugacies;
however, on their own, they do not show that the conjugacies glue together
to give a global conjugacy.
For this we need:
		
\subsubsection{Touching box mappings, Section~\ref{sec:touching box mappings}}
To obtain a global conjugacy, we construct in
Subsection \ref{subsec:global touching box mappings}
``global touching box mappings." These mappings are the analogue of B\"ottcher
 coordinates in our setting. 
Here the presence of parabolic points causes some additional
difficulties.
		
Some care is needed to make sure that the ranges of the complex box mappings
are compatible with with ranges of the touching box mappings.
In the reluctantly recurrent case the real bounds
(see Proposition~\ref{prop:real bounds reluctant})
imply, among other things, that the return domains to certain small intervals $I$
 about a critical point are very deep inside of $I$. Such bounds make it possible 
to extend the return map to $I$ to a complex box mapping whose range is a 
Poincar\'e lens domain. In general, we do not have this sort of geometric control in the
persistently recurrent case. Nevertheless, in 
Theorem~\ref{thm:box mapping persistent infinite branches}, for maps that
are at most finitely renormalizable, 
we construct a complex box mapping associated to an iterate of a
return map whose range is a Poincar\'e lens domain.
		
\begin{figure}
\resizebox{1.2\textwidth}{!}{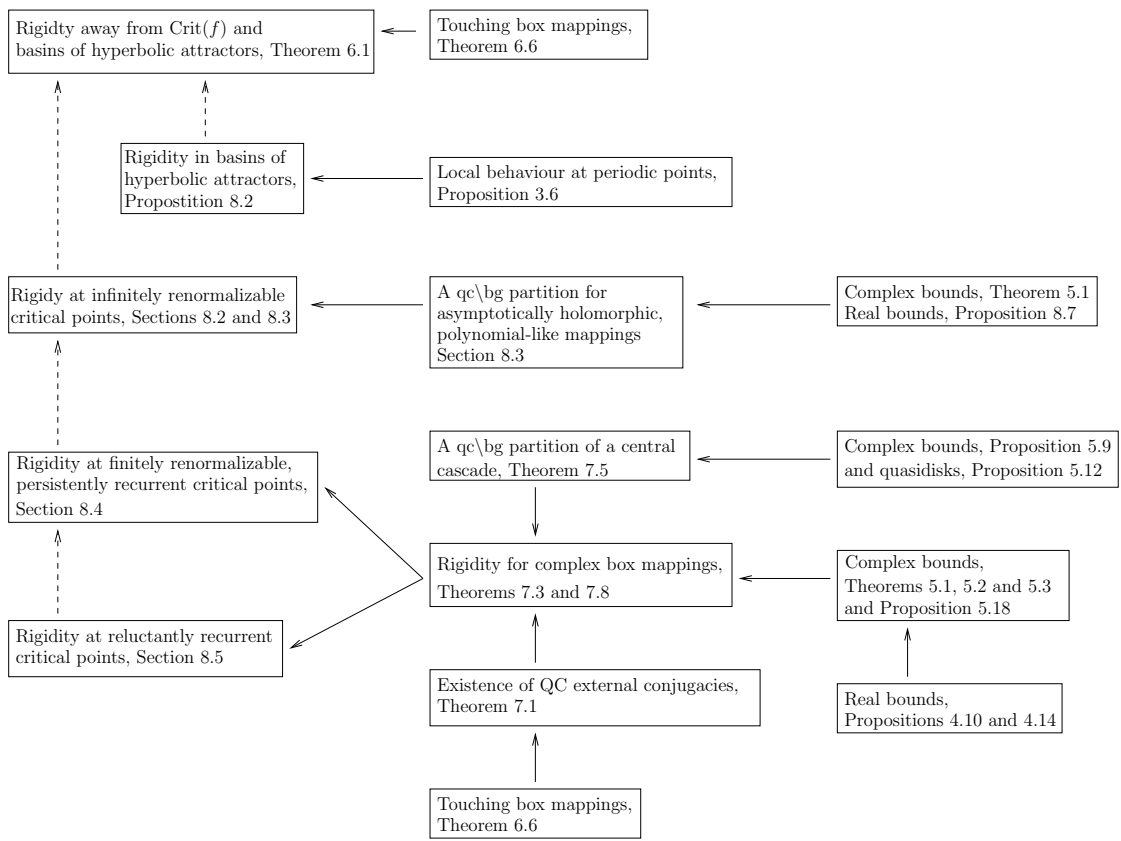}
 \caption{Scheme for proving rigidity: the vertical dashed arrows indicate the order in which we prove rigidity
 (starting at the top and ``pulling back"),
 and the solid arrows indicate key steps in the proof.}
 \label{fig:proofscheme}
\end{figure}
		
Consider two global touching box mappings
$F_T\colon \bm{\mathcal{U}}_T\rightarrow
\bm{\mathcal{\V}}_T,$ 
$\tilde{F}_T\colon \widetilde{\UU}_T\rightarrow\widetilde{\bm{\mathcal{V}}}_T$,
associated to $f$ and $\tilde f$.
These are mappings whose range is a finite ``necklace neighbourhood"
of the interval, each component of the range is a Poincar\'e
lens domain, and each component of 
the domain is mapped diffeomorphically onto a component 
of the range (we omit a neighbourhood of the critical points from the
domain of  a touching box mapping).
In Subsection \ref{subsec:global touching box mappings} we show that it is 
possible to
construct a qc-conjugacy between $F_T$ and $\tilde{F}_T$, and pull it back.
This yields a qs-conjugacy between $f$ and $\tilde{f}$ on the set of
points that avoid a real neighbourhood of
$\crit(f)$ and the attracting basins.
Constructing the conjugacy in attracting basins is not hard
since the local dynamics near attracting cycles is well understood
for maps in class $\mathcal{C}$. We explain how this is done
in Section \ref{sec:attracting}.
In the remaining part of the paper, we
focus on how to extend the conjugacy to the real 
neighbourhood $\mathcal I$ of $\crit(f)$.

\subsubsection{Rigidity of non-renormalizable complex box mappings,
  Section~\ref{thm:external conjugacy}-
Section~\ref{sec:cascades}}
Let us start with critical points at which $f$ is
at most finitely renormalizable .
To prove rigidity of complex analytic box mappings (on their real traces),
we use complex bounds 
together with  the Spreading Principle, Proposition
\ref{prop:spreading analytic}, 
and the QC Criterion of \cite{KSS-rigidity},
see Theorem~\ref{thm:qc-criterion}.
The nice geometric properties of the complex box mappings
allow us to construct qc external conjugacies,
see Theorem~\ref{thm:external conjugacy}, between
two topologically conjugate complex box mappings. 
The external conjugacy replaces the conjugacy 
given by the B\"ottcher coordinate for
polynomials. With these tools in hand, 
it is possible to repeat the proof of 
qc-rigidity for non-renormalizable 
mappings in \cite{KSS-rigidity} to prove
rigidity for analytic complex box mappings.
		
In the smooth case, since it is impossible to control the
dilatation of all pullbacks of puzzle pieces in the Yoccoz puzzle,
we cannot proceed exactly as in the analytic case. In particular the
usual Spreading Principle cannot be applied, and we need to use a 
different strategy to show that we obtain a qc mapping of the plane
with bounded dilatation
at each stage of our construction.
If we ignore controlling the dilatation of the conjugacy,
then our construction is reminiscent of the one used in
\cite{Lyubich-Fibonacci2, Lyubich-quadratic dynamics}.
We start with the
external conjugacy, see page~\pageref{def:extconj},
conjugating the dynamics of
$F\colon\bm{\mathcal{U}}\rightarrow\bm{\mathcal{V}}$
and $\tilde F\colon \tilde{\bm{\mathcal{U}}}\rightarrow
\tilde{\bm{\mathcal{V}}}$ on $\partial \bm{\mathcal{U}}$,
and pullback, modifying the mapping before we pullback through
critical branches so that it maps critical values to corresponding critical values.
However, because we do not have growth of moduli, the dilatation of
the resulting mapping can be arbitrarily large.
To deal with this, we prove complex bounds for the good nest,
see Section~\ref{sec:goodnest}.
These results make it possible to confine the region where the
dilatation of
the mapping that we obtain is uncontrolled to a union of 
puzzle pieces with bounded geometry and \emph{a priori} bounds.
When there are no long central cascades, 
after applying the QC Criterion to the puzzle pieces
in $\Dom'(\W^n)$, this
gives us the following:
There exists $K_1\geq 1 $ and $\delta_1>0$, such that
for any $n\geq 0$, there exists a $(K_1,\delta_1)$-qc\textbackslash bg
mapping $H_n:\mathbb H^+\rightarrow\mathbb H^+$
such that for any component $\U$ of $\Dom'(\W^n),$ 
$z\in\partial \U$, we have that $H_n\circ F(z)=\tilde F\circ H_n(z)$.
Quasisymmetric rigidity of complex box mappings on their real traces
follows from this together with the QC Criterion, applied on the upper
half-plane, and the compactness
of $K$-qc mappings.

To prove the complex bounds, we prove, roughly, that if a chain
$\{\bm G_j\}_{j=0}^s$ of
puzzle pieces such that $\bm G_0$ intersects the real line
does not contain a long central cascade, then
$\sum_{j=0}^s\diam(\bm G_j)^2$ is bounded,
see Proposition~\ref{prop:squares} for more general statement.
This estimate gives us that there are exactly two cases where we cannot
control the dilatation along a chain: along pullbacks of puzzle pieces
that are disjoint from the real line and nested inside of puzzle
pieces that are disjoint from the real line, and along long central
cascades. 
By not refining
the conjugacy inside of complex puzzle pieces that do not intersect
the real line, our complex bounds for the good nest give us a solution
to the first problem.
Our solution to the second problem requires 
us to find a puzzle piece near
the ``top'' of a long central cascade that is a quasidisk,
see Proposition~\ref{prop:good geometry for central cascades}.
Using this together with
a careful analysis of central cascades, 
we construct a QC\textbackslash BG Partition of 
a central cascade, see
Theorem~\ref{thm:central cascades},
which allows us to avoid the problem of the dilatation 
blowing up through long central cascades.

\subsubsection{Application of qc\textbackslash bg mappings to central cascades,
Section~\ref{sec:cascades}}
Suppose that we are in the situation described above where
$R_{\bm E^0}\colon \bm E^1\rightarrow \bm E^0$ has a long central cascade and
$\bm E^0$ is a $(1+1/\delta)$-quasidisk and $\delta$-nice,
see Figure~\ref{fig:qcbg}(A).
Using that $\bm E^0$ is a quasidisk and our complex bounds,
we can decompose $\bm E^0$ into combinatorially defined,
artificial puzzle pieces each with either 
good dynamical or geometric properties. The set 
with good dynamical properties is $X_2$ and the sets
with good geometric properties are $X_1$ and $X_3$.

Very roughly,we use a ``real'' argument to show that
if we start with a $K$-quasiconformal mapping
$H_0\colon \bm E^{0}\rightarrow\widetilde{\bm{E}}^{0}$ that maps
the domain of the return mapping to $\bm E^{0}$
to the domain of the return mapping to $\widetilde{\bm{E}}^{0}$
and conjugates $F$ and $\tilde F$ on the boundaries of these puzzle
pieces,
then we can pull $H_0$ back along the real branches of $R\colon E^1\rightarrow
E^0$
and $\tilde R\colon\tilde E^1\rightarrow \tilde E^0$ with 
the loss of dilatation controlled by $\diam(\bm E^0).$
This gives us, at least, a necklace neighbourhood of $E^0$ contained in $X_2$ 
where 
$H$ is defined. 
However, this neighbourhood could be very pinched
near the ``middle'' of the cascade; this is the case when there is
 a long saddle-node cascade, and in the high return case, it is
a ``necklace neighbourhood'' (as opposed to an open) 
neighbourhood of the interval.

The complement of $X_2\cup X_0,$
where $X_0$ is the measure zero set in the 
definition of a qc\textbackslash bg mapping, 
is decomposed into two types of regions: 
\begin{itemize}
\item The first type, corresponding to $X_3$, 
has geometry that we understand
very well because of the Yoccoz Lemma or the linear behaviour near a
repelling periodic point. We will show that if $\bm W$ is a connected
component 
of this type for $R$ and $\widetilde{\bm{W}}$ is a the corresponding component for
$\tilde R$, 
then $\bm W$ and $\widetilde{\bm{W}}$ are $K_1$-quasidisks, and
we
construct a
$\hat K$-qc mapping from $\bm W$ to $\widetilde{\bm{W}}$ that respects
some minimal dynamical information of their boundaries.
\item The remainder of $\bm E^0$ is $X_1$.
We will show that
that this subset, and its image under $H$, have the property that
their connected
  components have
bounded geometry and moduli bounds in both
  $\bm E^0$ and $\bm E^0\cap\mathbb H^+$.
\end{itemize}
This allows us to apply the QC Criterion. Let us point out that
we use the QC Criterion in two ways: in the way that is familiar from
\cite{KSS-rigidity} to obtain a qc mapping on a puzzle piece that
agrees with the boundary marking on the puzzle piece, and on the upper half-plane 
to obtain a qc mapping on the upper half-plane that agrees with 
the mapping that conjugates the dynamics on the real line outside
of a neighbourhood of the critical points.

\subsubsection{Proving quasisymmetric rigidity, Section~\ref{sec:qc rigidity}}
Figure \ref{fig:proofscheme} indicates how we proceed.

We start by constructing a conjugacy on immediate basins of attraction
at infinitely renormalizable critical points, 
and then we spread the conjugacy around 
to critical points $c$ at which $f$ is at most finitely renormalizable and 
so that $f$ is persistently recurrent on $\omega(c)$,
and finally to reluctantly recurrent critical points.

For real analytic infinitely renormalizable mappings
we refer to the argument of \cite{KSS-rigidity}
to construct a qs conjugacy in a real neighbourhood of
each critical point at which $f$ is infinitely renormalizable,
and use the touching box mappings to show that this conjugacy is
compatible
with the conjugacies constructed around other critical points.
The argument of \cite{KSS-rigidity} makes use of a 
lengthy real ``gluing'' argument, which we can avoid rather neatly by
using qc\textbackslash bg partitions, see Section~\ref{subsec:C3 inf renorm rigidity}.

For a smooth mapping $f\in\mathcal{C}$, we have
that a sufficiently high renormalization of $f$ extends to a
quasiregular polynomial-like mapping $F\colon\U\rightarrow\V$
(we define these on page~\pageref{def:box mapping}).
From the original mapping $F\colon\U\rightarrow\V$ it is 
impossible to construct arbitrarily small puzzle pieces.
Using a method of \cite{Levin-Przytycki} we construct rays for
$F\colon\U\rightarrow\V$. By doing this for each subsequent 
renormalization of $F$, we obtain arbitrarily small
puzzle pieces. (For analytic polynomial-like 
mappings this can be done using the Douady-Hubbard
Straightening Theorem, but
in general, a smooth mapping may not be
topologically conjugate to a polynomial.)
Using these puzzle pieces and complex bounds
we construct qc\textbackslash bg partitions
in the infinitely renormalizable setting, see Lemmas
\ref{lem:pl decomp} and  \ref{lem:renintersection}.
Thus obtaining a qc\textbackslash bg mapping, with moduli bounds in
the upper half-plane,
on a neighbourhood of the real line, which is a conjugacy
on the real line.

\subsection{A guide to the paper}
\label{sec:guide}
Sections~\ref{sec:definitions} and \ref{sec:real dynamics} contain definitions,
background material, and results in real dynamics. The reader can refer to these sections as needed.

\begin{itemize}
\item Whether a reader is interested in the analytic or the class 
$\mathcal C$ case,
the following sections are relevant:
\begin{itemize}
\item Sections \ref{subsubsec:real bounds - persistent}, 
\ref{sec:rbinfren} and \ref{subsubsec:real bounds - reluctant} for the
real bounds for persistently recurrent mappings,  infinitely renormalizable mappings,
 and reluctantly recurrent mappings, respectively.
\item Sections~\ref{subsec:complex bounds} and \ref{subsubsec:box
    mappings reluctant}
 for the complex bounds.
\item Section~\ref{sec:touching box mappings} where we construct touching
box mappings.
\item Section~\ref{sec:external conjugacy} where we construct an external conjugacy
between complex box mappings.
\end{itemize}

\item A reader who is interested in the case of analytic non-renormalizable
mappings should refer to
Section~\ref{sec:rigidity analytic} where we prove rigidity for
  analytic box mappings.
\end{itemize}
In Section~\ref{sec:qc rigidity}, these results are combined to show that
$f$ and $\tilde f$ are qs-conjugate. A reader who is only interested in
the non-renormalizable analytic case could start reading
that section at Section~\ref{subsec:qc rigidity-persistent}.

\begin{itemize}
\item A reader who is additionally interested in the analytic, infinitely
renormalizable case, should see
\begin{itemize}
\item Section~\ref{subsec:complex bounds} for the complex bounds.
\item Section~\ref{subsec:qc rigidity - infinitely renormalizable} where these
critical points are included in the proof of rigidity.
\end{itemize}

\item A reader who is focussing on the case of smooth 
non-renormalizable mappings
without long cascades of central returns ({\it i.e.} Fibonacci
mappings) should go on to read the following:
\begin{itemize}
\item Section~\ref{sec:diam sq} where we prove that we can control
dilatation as long as we do not pull back through long central
cascades, note that 
this
is easier when there are no long cascades of central returns whatsoever.
\item Section~\ref{sec:goodnest} where we prove complex bounds for the
  good nest of puzzle pieces.
\item Section~\ref{subsec:restricted spreading} for a real spreading
  principle. When there are no long cascades of central returns, the
  reader needs only read steps one
through five of the proof.
\end{itemize}
\item A reader who is interested in general smooth mappings should
  addtionally read:
\begin{itemize}
\item Section~\ref{subsec:improving} where we prove that we can find
  puzzle pieces which are quasidisks at the start of long central 
cascades.

\item Section~\ref{sec:diam sq}, with the case of long central
  cascades in mind.

\item Section~\ref{subsec:restricted spreading} for a real spreading
  principle. We treat long central cascades in Step 6 of the proof.

\item Section~\ref{sec:cascades} for the proof of
  Theorem~\ref{thmx:qcbgcen}, a qc\textbackslash bg partition of a
  central cascade.

\item Section~\ref{subsec:C3 inf renorm rigidity} where we treat the
  case of infinitely renormalizable smooth mappings.
\end{itemize}
\end{itemize}

\section{Preliminary definitions, remarks and some tools}\label{sec:definitions}

\subsection{Basic definitions and
  terminology}\label{subsec:terminology}
When it will not cause confusion, we let $M$ denote $[0,1]$ or
$S^1$. We let $\N,$ $\R,$ $\R^2$ and $\C$ denote the natural numbers ($\{1,2,3,\dots\}$),
the real line, the real plane and the complex plane, respectively.
	As usual,  we will consider $S^1$ as $\R\; \mathrm{mod}\, 1$.
Let $\H^+$ and $\H^-$ denote the upper and lower half-planes in $\mathbb{C}$.
	Let $B(a,r)$ denote the ball of radius $r$ centered at $a$ in $\mathbb{C}.$
	Let $\mathbb{D}_r=B(0,r)$ and $\mathbb{D}=\mathbb{D}_1$.
	An annulus is a doubly connected domain in $\mathbb{C}$.
	Any annulus $A$ is uniquely uniformized by a round annulus
	$\mathbb{D}_R\setminus\overline{\mathbb{D}}$, $R>1$,
	where $R=\infty$ is allowed (this case corresponds to the punctured disk).
	We define $\mathrm{mod}(A)=\log R$, when $R$ is finite,
	and in the case of the punctured disk, set $\mathrm{mod}(A)=\infty$.
	
If $X,Y$ are sets in $\mathbb{R}$ or $\mathbb{C}$,
we will let $\mathrm{diam}(X)$ denote the diameter of $X$
measured in the Euclidean metric, and 
$$\mathrm{dist}(X,Y)=\inf_{x\in X, y\in Y}|x-y|,$$
denote the distance between $X$ and $Y$ measured in the Euclidean metric.

For $k=0,1,2,\dots,$ $C^k(I)$ denotes continuous functions that are
	$k$-times differentiable on $I$ with $f^{(k)}\in C^0(I)$.
\medskip

We will adopt the convention that if a complex domain, say
$\bm P$, is named with a bold symbol, then 
$P=\bm P\cap\mathbb R.$

\medskip

	We shall denote the connected component that contains $x$ 
	of a set $P$ by $P(x)$ or by $\comp_x P$ when it improves clarity.

\medskip

We shall say that a set $\U\subset \C$ is {\em real-symmetric} if it is symmetric
with respect to the real line. A mapping  $F\colon  \U\to \V$
is {\em real-symmetric} if $\U$ and $\V$ are both real-symmetric
and if $F$ commutes with $z\mapsto \bar z$.
Unless otherwise stated, all constructions in the complex plane
are carried out real-symmetrically.

\medskip

If $I$ is a union of disjoint intervals, $|I|$ denotes the sum of their lengths.
If $I=(a-x,a+x)$ is any interval and $c>0$ we define $cI=(a-cx,a+cx).$
If we say that $I$ is \emph{well-inside}\label{well-inside}
of $I'$, we mean that $I'\supset (1+2\delta)I$
for a universal $\delta>0$;
that is, $\delta$ does not depend of $I,I'$. 
	
If $\U\subset \V$
are domains in $\mathbb{C}$ we say that $\U$ is
\emph{well-inside} $V$ if $\mod(\V\setminus \U)>\delta$, again
with $\delta>0$ universal.
If $\U$ is a topological disk in $\mathbb{C}$ and $\delta>0$,
we say that $\U$ has $\delta$-\emph{bounded geometry at} $x$
if \mbox{$B(x,\delta\cdot\mathrm{diam}(\U))\subset \U$.}
We say that $\U$ has $\delta$-\emph{bounded geometry}
\label{def:bddgeo}
if for some $x\in \U$ it has $\delta$-bounded geometry at $x$.

\medskip

If $f$ is differentiable,
we let $\mathrm{Crit}(f)$ denote the set of critical points of $f$
and $\mathrm{PC}(f)$ denotes the set of
strict foward images of $\mathrm{Crit}(f)$.
For a critical point $c$, we let 
$\mathrm{Forward}(c)$\label{page:forward}
denote all critical points in the closure
of $\{f^n(c)\}_{n=0}^\infty$, and we let
$\mathrm{Back}(c)$ denote the set of all critical points $c'$
such that $c\in\mathrm{Forward}(c')$. Critical points of
even order will be called \emph{folding points.}

\medskip

By the \emph{forward orbit of a point} $x$, we mean the set of points
$f^{n}(x)$ for $n\geq 0$. 
As usual, we let $\omega(x)$ denote the \emph{omega-limit set of x}.
We say that $\omega(x)$ is {\em minimal}, if the orbit of
every $x'\in \omega(x)$
is dense in $\omega(x)$.

\medskip

	We say that a real mapping $f$ is {\em renormalizable} at a non-periodic point $x$
	if there exists an 
	interval $K\owns x$ and $s\ge 2$ so that $K,\dots,f^{s-1}(K)$ have
	disjoint interiors and $f^s(K)\subset K$.
	In this case, we call $K$ a \emph{properly periodic interval} of period $s$.
	We say that $f$ is {\em infinitely renormalizable} at some non-periodic 
	point $x$ if there exist properly periodic intervals containing $x$
	of arbitrarily high period.
	This implies that there exists a sequence of {\em renormalization 
	intervals} around $x$ with periods $s_n\to \infty$, i.e. intervals 
	$K_n\ni x$ and integers $s_n\to \infty$ so that 
	$K_n,\dots,f^{s_n-1}(K_n)$ have disjoint interiors, 
	$f^{s_n}(K_n)\subset K_n$ and $f^{s_n}(\partial K_n)\subset \partial K_n$. 
If $f$ is infinitely renormalizable at $c$,
then $\omega(c)$ is a minimal set.
\medskip

	We say that a set  $P$ is {\em nice}\label{def:nice} if 
	for each $x\in \partial P$, $f^n(x)\notin \mathrm{int} P$ for all $n>0$. 
	It is easy to see that pullbacks of nice open intervals are nested or disjoint.
	
\medskip

Given a set $P$, let  
$$\Dom(P)=\{x\colon  \exists k>0 \mbox{ so that } f^k(x)\in P\}.
\label{notation:Dom}$$

Observe that $\Dom(P)$ is collection of \emph{all} landing domains to
$P$ together with first return domains to $P$;
it is not restricted to those that intersect $\omega(c)$ for some
critical point $c$.
We will refer to $\Dom(P)$ as the \emph{domain of the 
first entry mapping to} $P$.
	If $P$ is nice, then the boundary of a component of $\Dom(P)$ is mapped
	into the boundary of $P$. 
	If $\Omega\subset\crit(f)$, we let 
	 $\Dom_{\Omega}(P)\label{notation:Dom_Omega}$ denote the components of $\Dom(P)$ that
	 intersect $\cup_{c\in\Omega}\omega(c)$.
	\label{notation:D(I)}
	
If $x\in\Dom(P)$ the smallest $k>0$ such that
$f^{k}(x)\in P$ is called the \emph{entry time of x to P}. 
In this case we define $\mathcal{L}_x(P)=
\mathrm{Comp}_x f^{-k}(P)
=\comp_x \Dom(P).\label{notation:entry domain}$
We define 
	\begin{equation*}\label{landing domain}
		\hat{\mathcal{L}}_{x}(P)= \left\{ 
		\begin{array}{rl}
		P, & \mathrm{if\ } x\in P,\\ 
		\mathcal{L}_{x}(P)&\mathrm{if}\  x\notin P.
		\end{array}\right.
	\end{equation*}
The first entry mapping restricted to 
$\Dom^*(P):= P\cap \Dom(P)$ is the \emph{first return mapping to P};
\label{page:return}
the \emph{first landing mapping to} $P$ is defined on 
$\Dom'(P):=\mathrm{Dom}(P)\cup P$
and is given by the first entry mapping to $P$ on
$\mathrm{Dom}(P)\setminus P$ and the identity on $P$.

	\medskip

A set $X$ is \emph{invariant} if $f^{-1}(X)= X$
and \emph{forward invariant} if $f(X)\subset X.$
If $I$ is a union of intervals we let $E(I)$
denote the set of all points whose forward orbits never enter $I$.

\medskip

If $I$ is an interval and $J$ is a pullback of $I$,
we define the \emph{chain}, $\{G_j\}_{j=0}^s$ associated to the pullback by
setting $G_s=I$ and $G_j=\mathrm{Comp}_{f^{j}(J)}f^{-(s-j)}G_s$
for $j=0,\dots, s-1$.
The \emph{order} of chain $\{G_j\}_{j=0}^s$ is the number of intervals
$G_j,$ $j=0,\dots,s-1$ that contain critical points. We define
chains and the order of a chain
for complex puzzle pieces (defined on page~\pageref{page:defpuzzle}) 
in the same way. 

\medskip

Suppose that $\Omega\subset\mathrm{Crit}(f)$.
We say that a real neighbourhood $\mathcal I$ of $\Omega$ is
\emph{admissible} if 
\begin{itemize}
\item $\mathcal I$ has exactly $\#\Omega$ components,
each containing an element of $\Omega$; and 
\item for each connected component $J$ 
of $\mathrm{Dom}(\mathcal I)$, either $J$ is a component of $\mathcal I$ or $J$
is compactly contained in a component of $I$.
\end{itemize}
 
\medskip


Suppose that $\mathcal I$ is a union of nice open intervals
and that $\mathcal J$ is a union of intervals such that for each
connected component $J$ of $\mathcal J$, there exists $k_J$ such that
$f^{k_J}(J)\in \mathcal I$ and $f^{k_J}(\partial J)$ is contained in
the boundary of a component of $\mathcal I$.
Define $g\colon\mathcal J\rightarrow \mathcal I$ by 
$g(x)=f^{k_J}(x)$ whenever $x$ is in the connected component $J$ of
$\mathcal J$. We say that a complex box mapping (see
page~\pageref{def:box mapping})
$G\colon\bm{\mathcal U}\rightarrow\bm{\mathcal{V}}$ \emph{extends}
$g\colon\mathcal J\rightarrow \mathcal I$ if
\begin{itemize}
\item the critical points of $g$ and $G$ are the same,
\item for each component $J$ of $\mathcal J$, there exists a component
$U$ of $\bm{\mathcal{U}}$ that contains $J$ and $g|_J=G|_J,$
\item each component of $\bm{\mathcal{U}}$ contains a unique
  component of $\mathcal J$ and each component of $\bm{\mathcal{V}}$
contains a unique component of $\mathcal I$.
\end{itemize}

\subsection{Quasiconformal and quasisymmetric maps}\label{subsec:qc maps}
We say that a homeomorphism $h\colon [0,1]\rightarrow[0,1]$ is
$\kappa$-\emph{quasisymmetric} (abbreviated $\kappa$-$qs$),
\label{def:qs}
if
$$\frac{1}{\kappa}\leq\frac{h(x+t)-h(x)}{h(x)-h(x-t)}\leq \kappa$$
whenever we have both $x-t,x+t\in[0,1].$ 
The analogous definition holds for a homeomorphism $h\colon  S^1\rightarrow S^1.$
A homeomorphism is called \emph{quasisymmetric} if it is 
$\kappa$-quasisymmetric for some $\kappa\geq 1$.
If $X\subset M$, we say that $h\colon  X\rightarrow M$ is $\kappa$-qs 
if it has an extension to a $\kappa$-qs map
$h\colon  M\rightarrow M$.
	
Let $\U\subset\C$ be a domain.
A mapping $h\colon  \U\rightarrow \C$ is $K$-\emph{quasiconformal} ($K$-\emph{qc}) if it is a
homeomorphism onto its image, and if for any annulus $\bm A\subset \U$,
$$\frac{1}{K}\mod(\bm A)\leq\mod(h(\bm A))\leq K\mod(\bm A).$$
		The minimum such $K$ is called the \emph{dilatation} of $h$.
		Suppose that $X\subset [0,1].$ Any $\kappa$-quasisymmetric
                mapping $h\colon  X\rightarrow [0,1]$ has an extension to
                a
$K(\kappa)$-quasiconformal mapping \mbox{$h\colon\C\rightarrow\C$.}

We say that a Jordan disk is a \emph{(K-)quasidisk},
if it is the image of 
$\mathbb{D}$ under a ($K$-)qc mapping $h\colon\mathbb C\rightarrow\mathbb C$.
We say that a Jordan curve is a ($K$-)\emph{quasicircle},
if it is the boundary of a ($K$-)quasidisk. We will
use the Ahlfors-Beurling Criterion
to prove that certain topological disks
are quasidisks:
\begin{lem}[Ahlfors-Beurling Criterion]\label{lem:Ahlfors-Beurling}
For any $C>0$ there exists a constant $K=K(C)>1$,
such that the following holds.
Suppose that $T$ is a Jordan curve such that for any two points
$z_1, z_2\in T$, there exists an arc $\gamma\subset T$
containing $z_1$ and $z_2$ such that
$\diam(\gamma)\leq C\dist(z_1,z_2)$.
Then $T$ is a $K$-quasicircle.   
\end{lem}

Let us collect some useful results about 
quasisymmetric mappings on the real line.
Suppose that $X\subset\mathbb R$.
We say that a mapping $h\colon X\rightarrow h(X)$ is $(C,p)$-quasisymmetric, 
if there exists a constant $C>0$ and $p\in\mathbb N$, if
whenever $x,a,b\in X$  we have
$$\frac{|h(x)-h(a)|}{|h(x)-h(b)|}<
C\max\{\Big(\frac{|x-a|}{|x-b|}\Big)^p,\Big(\frac{|x-a|}{|x-b|}\Big)^{1/p}\}.$$
In this definition, we allow $X$ to be any,
even a totally disconnected, subset of $\mathbb R$.

\begin{lem}\cite[Lemma 2.3]{Vellis}\label{lem:qs gluing}
For any $C\geq 0, p\geq 1,$ $\kappa\geq 1$ and $\delta>0$
there exists $\hat\kappa\geq 1$ such that the following holds.
Let $I$ be a bounded interval.
Suppose that $h\colon\mathbb R\setminus I\rightarrow\mathbb R$ is
$(C,p)$-quasisymmetric, $h_I\colon(1+2\delta)I$ is
$\kappa$-quasisymmetric and that for each $h|_{(1+2\delta)I_i\setminus
I_i}=h_I|_{(1+2\delta)I_i\setminus I_i}$. Then 
$$\hat H(z)=\left\{
\begin{array}{ll}
H(z),& z\in\mathbb R\setminus I\\
h(z), &z\in I\end{array}\right.
$$
is $\hat\kappa$-quasisymmetric.
\end{lem}

For any $M\geq 1$,
we say that a set $S$ is $M$-relatively
connected if for any $x\in S, r>0$ such that 
$\overline{B(x,r)}\neq S$ we have that either 
$\overline{B(x,r)}=\{x\}$ or that
$ \overline{B(x,r)}\setminus B(x,r/M)\neq\emptyset.$
It follows from the definitions that for 
any $M,\kappa \geq 1$ there exists $M'\geq 1$ such that
if 
$S$ is $M$-relatively connected and
$h$ is $\kappa$-qs, then $h(S)$ is $M'$-relatively connected.

\begin{prop}\cite[Theorem 1.2]{Vellis}\label{prop:Vellis}
For any $M,\kappa\geq 1$ there exists $\kappa'\geq 1$ such that
the following holds.
If $S$ is an $M$ relatively connected subset of 
$\partial\mathbb D$ and 
$h\colon S\rightarrow\partial\mathbb D$ is $\kappa$-quasisymmetric mapping,
then $h$ extends to a $\kappa'$-qs mapping
$\hat h\colon\partial\mathbb D\rightarrow\partial\mathbb D$.
\end{prop}

\subsection{Asymptotically holomorphic mappings}
	\label{sec:dynamics of asymptotically conformal maps}

The asymptotically holomorphic extensions of smooth maps
that we use were constructed in \cite{GSS}, and we
refer the reader to that paper for the background on these extensions.

Recall that 
$$\frac{\partial}{\partial\bar{z}}=
\frac{1}{2}(\frac{\partial}{\partial x}+i\frac{\partial}{\partial y}).$$

Let $K\neq\emptyset$ be a compact subset of $\mathbb{R}^{2}$, 
$\bm U$ an open neighbourhood of $K$ in the plane and 
$H\colon  \bm U\rightarrow\mathbb{C}$ a $C^{1}$ map. 
We say that $H$ is $asymptotically\ holomorphic$ of order $t$, 
$t\geq 1,$ on $K\subset\mathbb{R}^{2}$ if for every $(x,y)\in K$ 
$$\frac{\partial}{\partial\bar{z}}H(x,y)=0,$$
and 
$$\frac{\frac{\partial}{\partial\bar{z}}H(x,y)}{d((x,y),K)^{t-1}}\rightarrow 0$$
uniformly as $(x,y)\rightarrow K$ for $(x,y)\in \bm U\setminus K$.

In our applications of asymptotically conformal extension, 
$K$ will be an interval contained in the real line, and 
$$\mu(z)\equiv\frac{\bar{\partial}H}{\partial H}=o(|y|^{t-1}) $$

A $C^{n}$ function $H$ is \emph{asymptotically holomorphic}
of order $n$ at a point $z_{0}=x_{0}+iy_{0}$ is equivalent to $H$ having a
complex Taylor expansion at $z_{0}$ of order $n$:
$$H(z)=H(z_{0})+\sum_{k=1}^{n}\frac{\partial^{k}H}{\partial
  x^{k}}(z_{0})
\frac{(z-z_{0})^{k}}{k!}+R(x,y),$$
where $z=(x,y)$, the remainder $R(x,y)$ is $C^{n}$
and $D^{j}R(x_{0},y_{0})=0$ for $j=0,\dots,n$,
and $R(x,y)=o(\|(x,y)-(x_0,y_0)\|^n)$.
	

\begin{lem}\cite[Lemma 2.1]{GSS}
Let $I\subset\mathbb{R}$ be a compact interval with non-empty interior
and $h\colon  I\rightarrow\mathbb{C}$ a $C^{n}$ function
($n\in\mathbb{Z}^{+})$.
Then there exists a $C^{n}$ extension
$H$ of $h$ to the complex plane that is
asymptotically holomorphic on $I$ of order $n$.
The mapping from $h$ to $H$ is linear and continuous in the $C^{n}$
topology (for $H$ we use the $C^{n}$ topology
on compact subsets of the complex plane)
and commutes with affine rescaling.
\end{lem}

\begin{lem}\cite[Lemma 2.2]{GSS}
Let $h\colon  I\rightarrow\mathbb{R}$ be a $C^{1}$ diffeomorphism
onto its image,
where $I\subset\mathbb{R}$ is compact and has non-empty interior.
If $H$ is any $C^{1}$ extension of $h$ to a neighbourhood of $I$
in the complex plane which is asymptotically holomorphic on
$I$ of order 1,
then $H$ is an orientation preserving diffeomorphism
on some (possibly smaller) complex neighbourhood of $I$.
\end{lem}

It will be useful to observe the following.
Let $K\neq\emptyset$ be a compact subset of $\mathbb{C}$ and
$H$ a $C^{1}$ complex-valued function defined in a open neighbourhood
of $K$
in the complex plane.
If $H$ is asymptotically holomorphic of order $t$, $(t\geq 1)$ on $K$, then
\begin{enumerate}
\item if $\phi$ is a holomorphic function defined in a neighbourhood of $H(K)$,
	then $\phi\circ H$ is asymptotically holomorphic of order $t$ on $K$;
\item if $H$ is a diffeomorphism onto its image,
	then $H^{-1}$ is asymptotically conformal of order $t$ on $H(K)$.
\end{enumerate}

The local dynamics near an attracting fixed point
for an asymptotically holomorphic mapping
can be studied in the same way as 
for a holomorphic mapping:	
\begin{prop}\label{prop:local dynamics}
Suppose that $f\in C^3$ and  
that $0$ is a fixed point of $f$. 
Let $\lambda=f'(0)$ and suppose $|\lambda|<1$.
Let $F$ denote the asymptotically holomorphic extension of $f$ of order 3.
\begin{enumerate}
\item If $0<|\lambda|$, then there exists $r_0>0$, so that
for any $r\in(0,r_0),$
$F$ is conjugate to
$z\mapsto\lambda z$ in the
neighbourhood $\mathbb D_r=\{z:|z|<r\}$
of zero by a $K(r)$-quasiconformal mapping and
$K(r)\rightarrow 1$ as $r\rightarrow 0$.
\item If $|\lambda|=0$,  then there exists $r_0>0$, so that
for any $r\in(0,r_0),$
$F$ is conjugate to
$z\mapsto z^d$ at 0, where $d$ is the order of the critical point,
in the neighbourhood $\mathbb D_r$
by a $K(r)$-quasiconformal mapping and
$K(r)\rightarrow 1$ as $r\rightarrow 0$.
\end{enumerate}
\end{prop}
\begin{pf}
The proof is the same as for holomorphic maps
(see \cite{Milnor}).
First suppose that $\lambda\neq 0$.
Let $c$ be so that $c^2<|\lambda|<c$.
Then by Taylor's Theorem, there exists a neighbourhood 
$\mathbb D_{r_0}$
of the origin so that $|F(z)|<c|z|$ for all $z\in \mathbb D_{r_0}$, and 
fixing any $0<r<r_0$, there exists
$k>0$ so that $|F(z)-\lambda z|<k|z|^2$ for all $z\in \mathbb D_r$.
Then given any starting point $z_0\in \bm D_r$, the orbit of $z_0$,
$z_0,z_1,\dots$, satisfies $|z_n|<rc^n$,
so that $|z_{n+1}-z_n|< kr^2c^{2n}.$
Then the sequence of points $w_n=z_n/\lambda^n$ satisfies
$$|w_{n+1}-w_n|\leq\frac{kr^2}{|\lambda|}\Big(\frac{c^2}{|\lambda|}\Big)^n.$$
Thus $|w_{n+1}-w_n|$ converges geometrically to $0$.
Hence the sequence of functions $z_0\mapsto w_n(z_0)$ converges
uniformly in $\mathbb D_r$.
Let $\phi\colon [-r,r]\rightarrow[-r,r]$ denote the limit of this sequence.
It fixes the end points of the interval and 0. Arguing as we just did,
there exists a constant $k'>0$ such that $|\phi(z)-z|<k'|z|^2$,
so $\phi$ is a local diffeomorphism.
Moreover, there exists $K'\geq 1$ (independent of $n$) so that
the mapping $z_0\mapsto F^n(z_0)/\lambda^n$ is $K'$-qc.
Hence, by compactness of $K'$-qc mappings, the mapping $\phi$
is $K'$-qc. That $\phi$ satisfies $\phi\circ F(z)=\lambda\phi(z)$
is automatic.
			
		\end{pf}

\subsection{Poincar\'e disks}\label{sec:poincare disks}

Let $I$ be an open interval and fix $\theta\in(0,\pi)$. 
Let $\bm D$ denote a round disk intersecting the real line with
$\bm D\cap\mathbb R=I$ with
the property that $\bm D^{+}=\bm D\cap\mathbb{H}^{+}$ intersects the
real line with external angle $\theta$.
Let $\bm D^-$ denote the
reflection of $\bm D^+$ about the real line.
Let $D_{\theta}(I)=\bm D^{+}\cup\bm D^{-}\cup I$; 
$D_{\theta}(I)$ is the \emph{Poincar\'e disk with angle}
\label{page:poincare} 
$\theta$
based on $I$.
We will
refer to Poincar\'e disks with angle $\theta$
close to $\pi$ as \emph{lens domains}.

	 \begin{prop}[\cite{GSS} Proposition 2]\label{prop:Almost Schwarz Inclusion}
		Let $h\colon  I\rightarrow\mathbb{R}$ be a $C^{3}$ diffeomorphism from a compact interval $I$ 
		with non-empty interior into the real line.
		Let $H$ be a $C^{3}$ extension of $h$ to a complex neighbourhood of $I$,
		with $H$ asymptotically holomorphic of order $3$ on $I$. 
		Then there exist $K>0$ and $\delta>0$ such that 
		if $a,c\in I$ are distinct, $0<\alpha<\pi$ and $\diam(D_{\alpha}(a,c))<\delta,$
		then
		$$H(D_{\alpha}(a,c))\subset D_{\tilde{\alpha}}(h(a),h(c)),$$
		where
                $\tilde{\alpha}=\alpha-K|c-a|\diam(D_{\alpha}(a,c))$. 
Moreover, $\tilde{\alpha}<\pi$.
	\end{prop}
		
	This proposition implies a version of the Almost Schwarz Inclusion Lemma from \cite{dFdM} when
	$F$ is analytic:
	\begin{lem}[Almost Schwarz Inclusion \cite{dFdM}]
		There exist $K<\infty$, $a_{0}>0$  and a function 
		$\theta\colon (0,a_{0})\rightarrow(0,\infty)$ satisfying
		$\theta(a)\rightarrow 0$ and $a/\theta(a)\rightarrow 0$ as 
		$a\rightarrow 0$ such that the following holds.
		Let $F\colon \mathbb{D}\rightarrow\mathbb{C}$ be univalent and real-symmetric,
		and assume that $I\subset\mathbb{R}$ is an interval containing 0 and
		$|I|<a\in(0,a_{0})$. Let $I'=F(I)$. Then
		\begin{enumerate}[label=(\alph*)]
			\item for all $\theta\geq\theta(|I|)$, we have 
				$$F(D_{\theta}(I))\subset D_{(1-K|I|^{1+\delta})\theta}(I'),$$ where
				$0< \delta<1$ is a universal constant;
			\item for all $\theta\in(\pi/2,\pi)$ we have $$F(D_{\pi-\theta}(I))\subset D_{\pi-K|I|\theta}(I').$$
		\end{enumerate}
	\end{lem}

	\begin{lem}[\cite{CvST}, Lemma 4.2]\label{lem:zd}
		Let $\ell\geq 2$ be a natural number and let $\theta\in(0,\pi).$
		Let $P(z)=z^{\ell}$.
		\begin{itemize}
		\item Suppose $\ell$ is even and let $K\geq 1$.
		Then there exists $\lambda=\lambda(K,\ell)\in(0,1)$ such that
				$$P^{-1}(D_{\theta}(-K,1))\subset D_{\lambda\theta}(-1,1)).$$
			
			\item Suppose that $\ell$ is an odd integer. Let $K>0$. Then
			there exists $\lambda=\lambda(\ell)\in(0,1)$ such that 
				$$P^{-1}(D_{\theta}(-K^{\ell},1))\subset D_{\lambda\theta}((-K,1)).$$
		\end{itemize}
	\end{lem}

	\begin{lem}\label{lem:zd lower bound}
		Let $\ell\in\mathbb{N}$ and let $P(z)=z^\ell$. 
		For any $A>0$ and any $\theta\in(0,\pi)$,
		there exists $\theta'\in(0,\pi)$ such that
		\begin{itemize}
			\item If $\ell\geq 2$ is even, then 
				$$P^{-1}(D_{\theta}((-A,1)))\supset D_{\theta'}((-1,1)).$$
			\item If $\ell\geq 3$ is odd,
				$$P^{-1}(D_{\theta}((-A^\ell,1)))\supset D_{\theta'}((-A,1)).$$
		\end{itemize}
	\end{lem}

\begin{lem}\label{lem:z2'}
Let $\ell\in\mathbb N$ and let $P(z)=z^\ell$.
\begin{itemize}
\item Suppose that $\ell$ is even. For any $K\geq 1$, there exists
$\theta_0\in(0,\pi/2)$ and $C=C(K,\ell)$ such that for any
$\theta\in(0,\theta_0)$
we have that 
$$P^{-1}(D_{\theta}(-K,1))\subset D_{\theta}(-1,1)\cup D_{\pi/2}(-C,C).$$
\item Suppose that $\ell$ is odd. There exists
$\theta_0\in(0,\pi/2)$ and $C=C(\ell)$ such that for any
$\theta\in(0,\theta_0)$
we have that 
$$P^{-1}(D_{\theta}(-K^\ell,1))\subset D_{\theta}(-K,1)\cup
D_{\pi/2}(C\cdot(-K,1)).$$
\end{itemize}
 \end{lem}
In the even case, this lemma follows from:
\begin{lem}\cite[Appendix]{LS-local}
Let $P_\ell(z)=z^\ell$. Let $K>1$. There exists $\theta_0=\theta_0(K)\in(0,\pi)$ such that,
for all $\theta\in(0,\theta_0)$, the boundaries of $P_{\ell}D_{\theta}((-1,1))$ and
$D_{\theta}((-K,1))$ intersect each other at a point $Z(K,\theta)$ and its complex
conjugate. Furthermore, 
$$Z(K,\theta)\rightarrow K^2\in\mathbb R,\mbox{ as }\theta\rightarrow 0.$$
Hence, the difference
 $$\Delta_\ell(K,\theta):=D_{\theta}((-K,1))\setminus P_{\ell}D_{\theta}((-1,1))$$
tends to the interval $[1,K^2]$, as $\theta\rightarrow 0$.
\end{lem}
\begin{pf}
Compare the proof of
Lemma~\ref{lem:zd}, see \cite{CvST}.
Suppose that $z\in \partial P^{-1}(D_{\theta}(-K^\ell,1))\cap \partial
D_{\theta}(-K,1)$.
Then since $z\in \partial D_{\theta}(-K,1)$,
we have $$\arg\frac{z-1}{z+K}=\theta,$$
and since $z^\ell\in\partial D_{\theta}(-K^\ell,1),$
we have that $$\arg\frac{z^\ell-1}{z^\ell+K^\ell}=\theta.$$
But now expressing $z=re^{i\theta},$
we have
$$\frac{r(K+1)\sin t}{r^2+r(K-1)\cos t
  -K}=\tan\theta\quad\mathrm{and}\quad\frac{r^\ell(K^\ell+1)\sin\ell
  t}{r^{2\ell}+r^\ell(K^\ell-1)\cos\ell t - K^\ell}=\tan\theta.$$
Dividing these expressions we have that
$$\frac{r^\ell\sin\ell t(K^\ell+1)(r^2+r(K-1)\cos
  t-K)}{(r^{2\ell}+r^{\ell}(K^\ell-1)\cos \ell t-K^\ell)r\sin t(K+1)}=1,$$
Expanding this and setting $RK=r,$
we have
\begin{eqnarray*}
1&=&\frac{\sin \ell t}{\sin t}\cdot\frac{r^{\ell-1}(K^\ell+1)(r^2+
rK\cos t -r\cos t -K)}{(K+1)(r^{2\ell} +r^\ell K^\ell\cos\ell
t-r^\ell\cos \ell t -K^\ell)}\\
&=&\frac{\sin\ell t}{\sin t}\cdot\frac{r^{\ell+1}K^\ell + r^\ell K^{\ell
     +1} \cos t - r^\ell K^\ell \cos t - r^{\ell-1}K^{\ell +1} +
     r^{\ell+1} + r^\ell K\cos t-r^\ell \cos t - K r^{\ell-1}}
{r^{2\ell}K+r^\ell K^{\ell+1}\cos\ell t-r^\ell K\cos\ell
    t-K^{\ell+1}+r^{2\ell}+r^\ell K^\ell\cos\ell t -r^\ell\cos\ell t
    -K^\ell}\\
&=& \frac{\sin\ell t}{\sin t}\cdot\frac{
\frac{r^{\ell+1}}{K^{\ell+1}}+\frac{r^\ell}{K^\ell}\cos
  t-\frac{r^\ell}{K^{\ell+1}}\cos
  t-\frac{r^{\ell-1}}{K^\ell}+\frac{r^{\ell+1}}{K^{2\ell+1}}+\frac{r^\ell}{K^{2\ell}}\cos
  t - \frac{r^\ell}{K^{2\ell+1}}\cos
    t-\frac{r^{\ell-1}}{K^{2\ell}}
}
{
\frac{r^{2\ell}}{K^{2\ell}} 
+\frac{r^\ell}{K^\ell}\cos\ell t
-\frac{r^\ell}{K^{2\ell}}\cos\ell t
-\frac{1}{K^{\ell}}
+\frac{r^{2\ell}}{K^{2\ell+1}}
+\frac{r^\ell}{K^{\ell+1}}\cos\ell t
-\frac{r^{\ell}}{K^{2\ell+1}}\cos\ell t
-\frac{1}{K^{\ell+1}}
}\\
&=&\frac{\sin\ell t}{\sin t}\cdot
\frac{
R^{\ell+1}+R^{\ell}\cos t
-R^{\ell}\frac{\cos t}{K}-R^{\ell-1}\frac{1}{K}
+R^{\ell+1}\frac{1}{K^\ell}+
R^\ell\frac{\cos t}{K^\ell}
-R^\ell\frac{\cos t}{K^{\ell+1}}
-R^{\ell-1}\frac{1}{K^{\ell+1}}
}
{
R^{2\ell}
+R^\ell\cos\ell t
-R^\ell\frac{\cos \ell t}{K^\ell}
-\frac{1}{K^\ell}
+R^{2\ell}\frac{1}{K}
+R^{\ell}\frac{\cos\ell t}{K}
-R^{\ell}\frac{\cos \ell t}{K^{\ell+1}}
-\frac{1}{K^{\ell+1}}
}
\end{eqnarray*}
For large $K$, the last term is dominated by 
$$\frac{\sin\ell t}{\sin t}\cdot\frac{R^{\ell+1}+R^\ell\cos
  t}{R^{2\ell}+R^{\ell}\cos\ell t}.$$
So either $R$ is bounded, in which case we are done, 
or $R$ is large, which implies that $\frac{R^{\ell+1}+R^\ell\cos
  t}{R^{2\ell}+R^{\ell}\cos\ell t}$ is very small, 
 and $\frac{\sin \ell t}{\sin t}$ 
is very large, but $\frac{\sin \ell t}{\sin t}$  is bounded from above.
So we can assume that $K$ is bounded,
but then we have that
$$1\asymp\frac{\sin\ell t}{\sin t}\cdot
\frac{R_n^{\ell+1}+R_n^\ell\cos t_n -\frac{R_n^\ell}{K}\cos t_n
  -\frac{R_n^{\ell-1}}{K}}{R_{n}^{2\ell}+R_n^\ell\cos\ell
  t_n-\frac{R_n^\ell}{K^\ell}\cos\ell t_n-\frac{1}{K^\ell}},$$
and again we have that either $R$ is bounded, in which case we are
done or we have that $R$ is very large, in which case
$$\frac{R_n^{\ell+1}+R_n^\ell\cos t_n -\frac{R_n^\ell}{K}\cos t_n
  -\frac{R_n^{\ell-1}}{K}}{R_{n}^{2\ell}+R_n^\ell\cos\ell
  t_n-\frac{R_n^\ell}{K^\ell}\cos\ell t_n-\frac{1}{K^\ell}}$$
is very small, but then we would need to have $$\frac{\sin\ell t}{\sin
  t}$$
large, which is impossible.
\end{pf}

\begin{lem}[\cite{CvST}, Lemma 5.10]\label{lem:angle control}
Suppose that 
$f\colon M\rightarrow M$ is $C^3$
and let $f$ also denote an asymptotically holomorphic extension of order three of $f$.
Then
for any $\delta>0$, $N\in\mathbb{N}\cup\{0\}$ and $\theta\in(0,\pi/2),$
there exist 
$\varepsilon>0,$ $\eta>0$, and $\tilde\theta\in(0,\pi/2)$
 such that the following holds.
Suppose that $J_s$
is a nice interval with  $|J_s|<\varepsilon,$
which does not contain a parabolic periodic point in its boundary,
and that
$\{J_j\}_{j=0}^s$ is a chain with order bounded by $N.$
Assume that $J_j\setminus(1+2\delta)^{-1}J_j$ 
does not contain an even critical point for $j=0,1,\dots,s-1$.
Let $\U_s=D_{\theta}(J_s)$, and set
$$\U_{j}=\comp_{J_j}(f^{-(s-j)}(\U_s)).$$
Then $\U_0\subset D_{\tilde \theta}(J_0)$
and $f^s\colon\U_0\rightarrow\U_s$ is $(1+\eta|\mu(J_s)|^{1/2})$-quasiregular.
\end{lem}

When we apply this lemma, we will use Lemmas \ref{lem:sum of lengths} and \ref{lem:sum of nice lengths}
to obtain the estimate on the sum of the lengths of intervals in the chain.

The following corollary follows easily:
\begin{cor}\label{cor:angle control}
For each $\delta>0$, there exist $\delta'>0$ and $\lambda\in(0,1)$
such that the following holds for all sufficiently small
$\varepsilon>0$.
Let $I$ be a nice interval with $|I|<\varepsilon$,
$J$ a domain of the first entry mapping to $I$, and $s\geq 1$ the minimal
number so that $f^{s}(J)\subset I$.
Let $\{H_j\}_{j=0}^s$ be a chain with $(1+2\delta)H_s\subset I$ and
$H_0\subset J$, then there exists an interval $H'$ with $H_0\subset
H'\subset(1+2\delta')H'\subset J$ so that
$$\comp_{H_0}f^{-s}D_{\theta}(H_s)\subset D_{\lambda\theta}(H').$$
\end{cor}

To treat diffeomorphic pullbacks of lens domains, we have
\begin{lem}\label{lem:lens angle control}
For any $\theta\in(0,\pi/2)$, there exists $\varepsilon_0>0$ and $\theta'\in(0,\pi/2)$ such that
the following holds.
Suppose that $H_s$ is a nice interval with $|H_s|<\varepsilon_0$ 
that does not have a parabolic periodic point in its boundary
and that $H_0$ is diffeomorphic pullback of $H_s$ by $f^s$.
Then $$\comp_{H_0}f^{-s}(D_{\pi-\theta}(H_s))\subset D_{\pi-\theta'}(H_0).$$
Moreover as $|H_s|\rightarrow 0$, $|\theta'-\theta|\rightarrow 0$.
\end{lem}
\begin{pf}
By Lemma~\ref{lem:sum of nice lengths}, there exists $\eta>0$, such that
\begin{equation}\label{eqn:length est}
\sum_{j=0}^{s}|H_j|^{3/2}<\eta,\mathrm{\ so\ that\ }
\sum_{j=0}^{s}|H_j|^{2}\leq\frac{|H_s|^{1/2}}{C}\eta.
\end{equation}
Let $K$ be the constant from Proposition~\ref{prop:Almost Schwarz Inclusion}.
For $\theta\in(0,\pi/2),$ $\diam(D_{\pi-\theta}(H_{j}))=|H_j|$,
so that by Proposition~\ref{prop:Almost Schwarz Inclusion}
$$\comp_{H_{s-1}}f^{-1}(D_{\pi-\theta}(H_{s})\subset D_{\pi-\theta-K|H_s|^2}(H_{s-1}).$$
Repeating this $s$ times, we have that
$$\comp_{H_{0}}f^{-s}(D_{\pi-\theta}(H_{s}))\subset D_{\pi-\theta-K\sum_{j=0}^s |H_j|^2}(H_{0}).$$
Estimate~(\ref{eqn:length est}) shows that the quantity
$K\theta\sum_{j=0}^s |H_j|^2$ can be made as small as we like by taking $|H_s|$ small.
\end{pf}

\section{Real dynamics}\label{sec:real dynamics}

We refer to \cite{dMvS} for background in real dynamics.
Let $M$ be either the interval or the circle, and suppose that
$f\colon M\rightarrow M$ is a mapping in $\mathcal C$.
Let us remark that since $f$
has no flat critical points, by
Theorem A of Chapter IV of [dMvS], $f$ does not 
possess wandering intervals;
and by Theorem B of that chapter, $f$ has
at most finitely many non-repelling cycles.
We also have by Ma\~n\'e's Theorem
that the set of points that avoids a (real) neighbourhood of the 
critical points of $f$ has measure zero.

Suppose that $p$ is a periodic point $p$ of period $s$.
Let $\lambda=Df^s(p)$ be the \emph{multiplier} of the periodic orbit.
A periodic orbit is called 
\emph{attracting}, if $|\lambda|<1$;
\emph{parabolic}, if $|\lambda|=1$,
and \emph{repelling}, if $|\lambda|>1$.
An attracting periodic point is called \emph{super-attracting},
if $\lambda=0$,
and in this case the periodic cycle contains a critical point of $f$.

	\begin{lem}\cite[Theorem C]{vSV}\label{real Koebe}
		Suppose that $f$ is a $C^3$ mapping of either the interval or the circle with 
		at most finitely many critical points $c_i$, $0\leq i\leq b-1$, and that 
		at each $c_i$ we can express
		$$f(x)=\pm (\phi(x))^{d_i}+f(c_i),$$
		where $\phi$ is $C^3$, $\phi(c_i)=0$, and $d_i>1$.	
		Then one has the following properties:
		\begin{enumerate}[label=(\arabic*)]
			\item \textit{Improved Koebe Principle.} For each $S>0$, $\delta>0$ and $\xi>0$,
				there exists $K>0$ such that the following holds. Suppose $J\subset T$ are intervals
				with $f^{n}|_T$ a diffeomorphism, $f^{n}(J)$ $\xi$-well-inside $f^{n}(T)$ and either
				\begin{enumerate}[label=(\alph*)]
					\item $\sum_{i=0}^{n-1}|f^{i}(J)|\leq S$ or
					\item $f^{n}(T)\cap B_{0}(f)=\emptyset$ and
						$\mathrm{dist}(f^{i}(T), \mathrm{Par}(f))\geq \delta, i=0,\dots, n-1.$
				\end{enumerate} 
				Then $f^{n}|_J$ has bounded distortion, i.e. for any $x,y\in J,$
				$$\frac{|Df^{n}(x)|}{|Df^{n}(y)|}\leq K.$$
				Here $B_{0}(f)$ is the union of immediate basins of periodic (possibly parabolic) attractors
				and $\mathrm{Par}(f)$ is the set of parabolic periodic points of $f$.
			\item Negative Schwarzian derivative. For each critical point $c$ that is not in the basin of a periodic attractor,
				there exists a (real) neighbourhood $U$ of $c$ such that whenever
				$f^{n}(x)\in U$ for some $x\in I$ and $n\geq 0$, the Schwarzian derivative of $f^{n+1}$ at $x$ is negative.
		\end{enumerate}
	\end{lem}

		\begin{lem}[Improved Macroscopic Koebe Principle]\label{lem:koebe}
		Let $f$ be as in Lemma~\ref{real Koebe}.
	\begin{enumerate}
		\item For each $\xi>0$ there exists $\xi'>0$ so that if $(1+2\xi)J\subset I$ 
			where $J$ and $I$ are both nice intervals, then for any $x$ and any $k$ with $f^k(x)\in J$,
			$$(1+2\xi')\comp_x f^{-k}(J) \subset \LL_x(I).$$

		\item For each $\xi>0$ there exists $\xi'>0$ so that if $(1+2\xi)J_s\subset G_s$ 
		and $J_s\supset \LL_x(G_s)$ for some $x\in J_s$, the following holds.
			Let $J_i\subset G_i$ be pullbacks of $J_s\subset G_s$. Then $(1+2\xi')J_0\subset G_0$.
			Here $\xi'(\xi)\to \infty$ when $\xi\to \infty$.
			\end{enumerate}
	\end{lem}
	
	\begin{pf}[Proof of Lemma~\ref{lem:koebe}]
		Statement (1) of this lemma is Theorem C(1) of \cite{vSV}, see also
		\cite{vSV-erratum}. The second part of this lemma is a special case of
		 \cite{CvST} Lemma 3.12
		\end{pf}
		
		As in \cite{dMvS}, for intervals $J\subset T$ so that $L,R$ are the components
of $T\setminus J$ define
$$C(T,J)=\dfrac{|T||J|}{|L||R|} \mbox{ and }
C(f,T,J)=\dfrac{C(f(T),f(J))}{C(T,J)}.$$

\begin{lem}[Cross-ratio inequalities]\label{lem:cross-distortion}
 Assume that $f$ is $C^3$ (weaker assumptions are sufficient, see
\cite{dMvS}). Then for each $K<\infty$ there exists $C<\infty$  and $\kappa>1$
such that the following hold:
\begin{enumerate}
\item If $J\subset T $ are sufficiently small intervals as above with 
$f^s|_T$ a diffeomorphism and $T,\dots, f^{s-1}(T)$ have intersection multiplicity $\le K$
(this means that each point is contained it at most $K$ of these intervals), 
then 
$$C(f^s,T,J)\ge 1- C\max_{i=0,\dots,s-1}|f^i(T)|.$$
\item 
If $f$ has no parabolic periodic orbits, then above statement holds even without the assumption 
that $T,\dots, f^{s-1}(T)$ have intersection multiplicity $\le K$.
\item 
 If $J\subset T $ are sufficiently small intervals as above satisfying 
 $|T|\ge (1/K) \dist(T,c)$ where $c$ is a critical point of $f$, and $|L|,|R|\geq |T|/K,$
then 
$$C(f,T,J)\ge \kappa>1.$$
\end{enumerate}
\end{lem}
\begin{pf}
See Sections IV.1 and IV.2 of \cite{dMvS}. 
The proof of (3) is the same as the proof of the Second Expansion Principle,
\cite[Chapter IV, Theorem 2.2]{dMvS}. Note that in the 
proof of Theorem 2.2 the assumption that $|L|,|R|\geq |T|/K$
is used (and should have been assumed) rather than $T\supset (1+2\tau)J$.
\end{pf}

We will use the following lemmas to control the sums of squares of the lengths of intervals.
	\begin{lem}[Theorem B \cite{Li-Shen}]\label{lem:sum of lengths}
		Let $f$ be $C^3$ with all periodic orbits hyperbolic repelling. Then for any $\alpha$ there exists $C=C(\alpha)$ such that
		for any interval $T$ and any $s\in\mathbb{N}$, if $f^s\colon  T\rightarrow f^s(T)$
		is a diffeomorphism, then
		$$\sum_{i=0}^{s}|f^i(T)|^{1+\alpha}<C.$$
	\end{lem}
	Under more restrictive conditions, we can remove the condition that $f$ have no parabolic points.
	\begin{lem}\label{lem:sum of nice lengths}
	Let be $f\in\mathcal{C}$ and let $N\in\mathbb{N}\cup\{0\}$.
	Let $\mathcal{P}'$ be a partition of $M$ given by Lemma \ref{lem:startingpartition}.
	Let $T'$ be a component of $f^{-N}(\mathcal{P}')$ that does not contain a 
	parabolic periodic point in its boundary.
	Then for any $\alpha$ there exists $C=C(\alpha)$ such that
		for any interval $T$ and any $s\in\mathbb{N}$, if $f^s\colon  T\rightarrow f^s(T)=T'$
		is a diffeomorphism, and $f^i(T)$ avoids the immediate basins of periodic attractors of $f$,
		then
		$$\sum_{i=0}^{s}|f^i(T)|^{1+\alpha}<C.$$
	\end{lem}
	\begin{pf}
We will explain how to adjust the proof of Lemma~\ref{lem:sum of lengths} in this situation.
In the proof of 
Proposition 2 of \cite{Li-Shen} in Subcase 1.1 on page 1490,
using the notation of that paper,
we replace their estimate
$$|J_n|\leq C\lambda ^{-n}|I|$$
by
$$\sum_{i=0}^{q}|f^{ir}(T)|^{1+\alpha}\leq \max_{0\leq i<q}|f^{ir}(T)|^{\alpha}\sum_{n=0}^{q}|f^{ir}(T)|
\leq \frac{|J_0|^{\alpha}}{C}|I|\leq\frac{|I|^{1+\alpha}}{C},$$
where the first inequality is trival, the second follows from the fact that
each interval $f^{ir}(T)$ is contained in $J_0=J$
and the fact that $f^{ir}(T)\cap f^{jr}(T)=\emptyset$ for $i\neq j$,
and the third from the fact that $J_0\subset I$.
It is worth remarking that in Subcase 1.1 of  \cite{Li-Shen}
we assume that the monotone mapping $f^r|_K$
has no fixed point.
The remaining parts of the proof are not affected by the
presence of parabolic cycles.
	
Let us briefly comment on the remaining ingredients of the proof.
Proposition 1 is Theorem C(1) from \cite{vSV}, and does not require that $f$ not have parabolic points.
Proposition 3 gives an estimate in the case when a return domain $J$ to a nice interval $I$
has length comparable to $I$. It is Proposition 5.2 of \cite{Li-Shen Hausdorff}, and
follows from \cite{vSV} Lemmas 2 and 3, and \cite{Shen-Todd} Theorems 1 and 2,
all of which hold for $C^3$ maps with no flat critical points.
Proposition 4 is Theorem A from \cite{vSV}, which does not depend on $f$ not having 
parabolic points.	
The final ingredient is the 
Yoccoz Lemma, which does not require that $f$ have no parabolic cycles. 
	\end{pf}
	
	For maps $f\in\mathcal{C}$, we will need to allow for the
	boundaries of the domains to contain 
	parabolic points. This will be important in 
	Section \ref{sec:touching box mappings}.
	To handle this we will use the following estimate.

	\begin{lem}\label{lem:parabolic rate}
		Let $f\in\mathcal{C}$. 
		Suppose that 0 is a parabolic periodic point with multiplier 1 and repelling side $[0,x).$
		Suppose that $x_0$ is so small that $f$ is a diffeomorphism on 
		$[0,x_0]=T$ and we can approximate $f$ by
		$$f(x)=x+ax^{d+1}+o(x^{d+1})$$ in $T$.
		For any $\alpha>0$, there exists a constant $C=C(d,\alpha)$ such that		
		$$\sum_{i=0}^{\infty}|f^i(T)|^{d+\alpha}<C.$$
	\end{lem}
	The exponent on the left is among the key reasons that we
	require higher regularity near parabolic periodic points.
	
	\begin{pf}
	We may assume $a=1$.
	Approximate $f$ by $x+x^{d+1} +o(x^{d+1})$ near 0.
	Let $x_0>0$ be a boundary point of $T$, define
	a sequence by
	 $x_{n}=x_{n+1}+x_{n+1}^{d+1}+o(x^{d+1})$.
	Then $$x_{n+1}-x_{n}=-x_{n}^{d+1}-o(x_{n}^{d+1}).$$
	So we can approximate $x_{n}$ by the solution of the differential equation:
	$$y'(t)=-y(t)^{d+1}-o(y(t)^{d+1}), y(0)=x_0.$$
	Since
        $y''(t)=-(d+1)y'(t)^d=-(d+1)(-y(t)^{d+1}-o(y(t)^{d+1})^d>0$,
for $y(t)>0$ sufficiently small,
	we have that the time one mapping eventually satisfies
	$y(t+1)\leq y(t)-ay(t)^{d+1}$, for some $\eta\in(0,\infty)$.
	Solving the differential equation 
	$$\phi'(t)=-\eta\phi(t)^{d+1}, \phi(0)=x_0$$
	gives
	$$\phi(t)=\Big(\frac{1}{\frac{1}{x_0^d}+\eta dt}\Big)^{1/d},$$
	so that
	$$y(t)\leq\Big(\frac{1}{\frac{1}{x_0^d}+\eta dt}\Big)^{1/d},$$
	and we have
	$$\sum_{i=0}^{\infty}|f^i(T)|^{d+\alpha}\leq\sum_{i=0}^{\infty}
        \Big(\frac{1}{\frac{1}{|x_0|^d}+\eta di}\Big)^\frac{d+\alpha}{d} <C.$$	
	\end{pf}

	\subsection{A real partition}\label{subsec:real partition}
              Suppose that either
		$f\colon  [0,1]\to [0,1]$ or the $f\colon  S^1\to S^1$ is not a 
		homeomorphism. 
		In the latter case $f$ has periodic orbits, in view of the following proposition
		(which, surprisingly,  we could not find in the literature).
		

\begin{prop}[Non-monotone circle maps have periodic orbits]
\label{prop:circlecase}
If $f\colon  S^1\to S^1$ is of class $\mathcal{C}$ and 
has no periodic points,
then $f$ is an orientation preserving homeomorphism.
\end{prop}
		
\begin{pf}
If the degree of $f$ is not equal to $1$, then $f$ has a fixed point.
So assume that $\deg(f)=1$.
Let $S^1=\R\, \;\mod\ 1$ and define $\pi\colon  \R\to S^1$
by $\pi(x)=x\, \mod\ 1$. Let $F\colon  \R\to \R$ be the the lift of $f$
(i.e. $\pi \circ  F=f\circ \pi$). 
Let $\rho(f,z)$ equal to the set of limit points of the sequence 
$\dfrac{F^n(x)-x}{n}$ where $\pi(x)=z$. It is well-known that
$\rho(f)=\{\rho\colon  \rho=\rho(f,z)\mbox{ for some }z\in S^1\}$
is a connected set,
called the rotation interval of $f$. In the classical case when $f$ is a
homeomorphism $\rho(f)$ is a singleton. If the
interval $\rho(f)$ contains a rational number $p/q$, then $f$
has a a periodic point of period $q$. So we are done unless $\rho(f)$
is equal to an irrational number. 
It is also well known that the rotation interval is equal to
the interval $[\rho(f_-),\rho(f_+)]$ where $F_-(x)=\inf_{y\ge x} F(y)$
and $F_+(x)=\sup_{y\le x} F(y)$ are the lower and upper maps 
associated to $F$, and $f_-$ and $f_+$ are the corresponding
circle maps.
Moreover,  let $F_t\colon  \R\to \R$, $t\in [0,1]$
be a family of maps depending
continuously on $t$ with $F_0=F_-$, $F_1=F_+$,
$F_-\le F_t\le F_+$  and so that $F_t=F$ except
on some intervals on which $F_t$ is constant.
Choose $F_t$ so that it is increasing in $t$, and
let $f_t$ be the corresponding circle map, see Figure~\ref{fig:circlemap}.

\input circlemap.tex
	
If $f$ is not a homeomorphism, then it is
has local extrema, and then $F_+\ge F_-$ with strict inequality in some places.

\noindent
\medskip
{\bf Claim:} If $\rho(f_-)$ is irrational, then  $\rho(f_+)>\rho(f_-)$
and there exists $t>0$ arbitrarily close to $0$ so that $f_t$ has a periodic point.

To prove this claim, we first observe that the maps $f_t$ do not have
wandering intervals.  Indeed, in \cite{MMS} 
monotone degree one circle maps $g\colon  S^1\to S^1$ 
were considered  with the following additional property:
(i) $g$ is smooth outside a finite number of points,
(ii) $g$ has a finite number of plateaus and critical points, 
and (iii) near each of these points $x_i$ the mapping $g$ is 
(up to smooth coordinates) of one of the following forms
$$x\mapsto \left\{\begin{array}{ll} -|x-x_i|^{\alpha_i}, &  x\le x_i \\
0 & x \ge x_i \end{array} \right. \,\,\, \mbox{ or }\,\,\, x\mapsto
\left\{\begin{array}{ll} 0 & x \le x_i  \\
|x-x_i|^{\alpha_i}, &  x\ge x_i 
\end{array} \right.$$
where $\alpha_i\ge 1$. It is shown in \cite{MMS}
that such a mapping has no wandering intervals,
i.e. for each interval $J$ so that $g^n(J)$, $n=0,1,2,\dots$
are disjoint there exists 
$n\ge 0$  so that $g^n(J)\subset \mbox{Plat}(g)$. 
Here $\mbox{Plat}(g)$ is the union of the plateaus of $g$.
In particular, $f_t$ has no wandering intervals for each $t\in [0,1]$.

Take a plateau $K$ of $F_-$ which is not eventually mapped
into any other plateau and let $\hat x$ be the left endpoint of $K$. 
Since $\rho(f_-)$ is irrational and $f_-$ has no wandering intervals,  there
exist iterates $n_i\to \infty$ so that $(f_-)^{n_i}(K)$
accumulates from the left to $\hat x$.
Therefore there exists a sequence $p_i$
of positive integers such that 
$(F_-)^{n_i}(\hat x) <\hat x + p_i$ for all $i$ and 
$(F_-)^{n_i}(\hat x) -\hat x - p_i\to 0$ as $i\to \infty$.
Notice that for each $t>0$ there exist $\alpha(t)>0$  and $N(t)$ so that
$(F_t)^n(x)\geq (F_-)^n(x) + \alpha(t)$
for all $n\ge N(t)$ and all $x\in \R$.
Indeed, there exists $\epsilon(t)>0$ and an open set
$U_t\subset \mbox{Plat}(f_t)\setminus \mbox{Plat}(f_-)$
so that  $F_t(x)>F_-(x)+\epsilon(t)$ for all $x\in U_t$. 
Since $U_t$ is not contained
in a plateau of $F_-$ and $\rho(f_-)$
is irrational, there exists $N(t)$ so that for each $x\in \R$
there exists $0\le k<N(t)$ so that $(F_-)^k(x)\in U_t \mod\ 1$.
It follows that $(F_t)^{k+1}(x)\ge (F_-)^{k+1}(x)+\epsilon(t)$.
Moreover, there exists $\tau(t)>0$ so that 
$(F_t)^m(y+\epsilon(t))\ge (F_-)^m(y)+\tau(t)$
for each $0\le m\le N(t)$ and each $y\in \R$.
It follows that for $\alpha(t)=\min(\epsilon(t),\tau(t))>0$ so that
$(F_t)^n(x)\geq (F_-)^n(x) + \alpha(t)$ for $n \ge N(t)$ and all $x\in \R$.
This, $(F_-)^{n_i}(\hat x) -\hat x - p_i\to 0$ and  the
Intermediate Value Theorem implies that there exists $t'\in (0,t)$, so that
$(F_{t'})^{n_i}(\hat x) -\hat x-p_i=0$.
Hence $\hat x$ is a periodic point of $f_{t'}$.
\end{pf}
With this lemma in hand, we will assume that $f$ has
periodic orbits.

Suppose that $f$ is a $C^1$ mapping of $M$ with
$b<\infty$ critical points.
We say that an $f$-forward invariant set
$Z\subset I$ is \emph{admissible with respect to} $f$
if it is a finite set that is disjoint from
$\mathrm{PC}(f)$.
From now on $Z$ will always denote an admissible set.
Observe that if $Z$ is an admissible set,
and $\Omega\subset\crit(f),$
then the union of the components 
of $M\setminus Z$ that intersect $\Omega$ is
an admissible neighbourhood of $\Omega$.
We call an interval $I$ a 
\emph{(real) puzzle piece with respect to} $Z$
if it is a component of 
$f^{-n}(Y)$ where $Y$ is a component of $[-1,1]\setminus Z.$
In this case we define the \emph{depth} of the puzzle piece to be $n$. 
We will use the notation $Y^Z_n(x)$ to denote the puzzle piece
with respect to $Z$ of depth $n$ that contains $x$.
We will let $\mathcal{Y}^Z_n$ denote the union of puzzle pieces of level $n$.
When it will not be confusing, we will omit the $Z$ in the notation. 
We call a puzzle piece \emph{critical} if it contains a critical point.

		If $P$ and $Q$ are critical puzzle pieces, we say that $Q$ is a \emph{child}  of $P$,
		if it is a unicritical pullback of $P$; that is, for some $n$, $Q$ is a component of 
		$f^{-n}(P)$ and the mapping $f^{n-1}\colon  f(Q)\rightarrow P$ maps a 
		neighbourhood of $f(Q)$ diffeomorphically onto its image.
		Notice that we do not require that $Q$ be contained in $P$.
		Let $c_0\in\mathrm{Crit}(f)$. 		 
		A mapping $f$ is called \emph{persistently recurrent on $\omega(c)$}\label{persistently recurrent}
		if each critical point in $\omega(c)$
		is non-periodic and each $\omega(c)$-critical puzzle piece has only finitely many children.
		Under these circumstances we will also say that $f$ is persistently recurrent at $c$.
		If $f$ is persistently recurrent on $\omega(c)$, then $\omega(c)$ is a minimal set.
		We say that $f$ is \emph{reluctantly recurrent} at $c$ if $c$ is recurrent and
		and there exists an $\omega(c)$-critical puzzle piece with infinitely many children.

		Let $c$ be a critical point of $f$ and let $Z$ be an admissible set for $f$.
		We say that $f$ is $Z$-recurrent at $c$ 
		if for any $n\geq 0$, there exists some $k\geq 1$ such that $f^{k}(c)\in Y^{Z}_{n}(c)$. 
		We say that $f$ is $Z$-\emph{renormalizable at} $c$, or that $c$ is $Z$-\emph{renormalizable} 
		if there exists a positive integer $s$, such that $f^{s}(c)\in Y^{Z}_{n}(c)$ for any $n\geq 0,$ 
		and the minimal positive integer $s$ with this property is 
		called the $Z$-\emph{renormalization period} of $c$.
		
		We will specify a choice of starting partition in Lemma~\ref{lem:startingpartition}.

		\begin{lem}[\cite{KSS-rigidity}, Fact 6.1]\label{lem:shrinking puzzle pieces}
			Let $Z$ be an admissible set. 
			Then for each $c\in\Crit(f)$, if $f$ is $Z$-recurrent, 
			but not $Z$-renormalizable at $c$, then 
			$|Y^{Z}_{n}(c)|\rightarrow 0$ as $n\rightarrow\infty.$
		\end{lem}

		Let $\mathrm{Crit}_{trival}(f)$ denote the subset of $\mathrm{Crit}(f)$ 
		which are either  contained in an attracting basin 
		or that have a finite forward orbit. 
		We will call a mapping $f$ trivial if
		$\mathrm{Crit}_{trival}(f)=\mathrm{Crit}(f)$.

\subsection{Renormalizable maps}\label{subsec:renormalizable maps}
If $c$ is $Z$-renormalizable,we define 
$$A^{Z}(c)=\cap_{n=0}^{\infty}\cup_{i=0}^{s-1}Y^{Z}_{n}(f^{i}(c))\cap\Crit(f),$$
where $s$ stands for the $Z$-renormalization period of $c$. 
The set $\cap_{n=0}^{\infty}\cup_{i=0}^{s-1}Y^{Z}_{n}(f^{i}(c))$
is a union of intervals $I$, which contains $\omega(c)$, and such that 
for each such $I$,
$f^s(I)\subset I$ and $f^{s}(\partial I)\subset \partial I;$
that is, it is the union of periodic intervals of period $s$, containing $\omega(c)$.
Note that any critical point $c'\in A^{Z}(c)$ is also renormalizable with period $s$ 
and that $A^{Z}(c)=A^{Z}(c').$

If $I$ is a nice interval containing a critical point $c$,
we define the \emph{principle nest} with $I^0=I$ at $c$ by
$$I=I^0\supset I^1\supset I^2\supset\dots,$$
where $I^{i+1}=\mathcal{L}_c(I^i)$.
We say that an interval $I$ is \emph{terminating} if 
for each critical point $c$ of $R|_{I^1}\colon  I^1\rightarrow I,$ 
$R^i(c)\in I^i$ for all $i\in\mathbb N$.
Otherwise we say that $I$ is \emph{non-terminating}.
If $I$ is terminating, then $R_I\colon  I^1\rightarrow I$ is renormalizable and
$I^{\infty}=\cap_i I^i$ is a periodic interval.

If $f$ is renormalizable (that is $Z$-renormalizable)
at a critical point $c$, then there is an even
 critical point $c_0\in\omega(c).$ Fix such a $c_0$.
There exists an interval $T$ and an involution
$\tau\colon  T\rightarrow T,$ such that $f\circ\tau=f$ holds on $T$.  
Let $B_1\supset B_2\supset B_3\supset \dots$ be all the symmetric
		(with respect to $\tau$) properly periodic intervals that contain $c_0$,
		and let $1\leq s_1<s_2<s_3<\dots$ be the corresponding periods.
		Define $\beta_n$ by $B_n=(\beta_n,\tau(\beta_n))$.
		Let $\alpha_n$ be the orientation reversing fixed point of
		$f^{s_n}$ closest to $c_0$,
		and set $A_n=(\alpha_n,\tau(\alpha_n))$.

	\subsection{A starting partition}
		\begin{lem}[Choice of a starting partition on the real line]\label{lem:startingpartition}
			Assume that  $f$ is of class $\mathcal{C}$ and $f\colon  [0,1] \to \R$ or $f\colon  S^1\to S^1$,
			and that $f$ has periodic orbits.  Then there exists a finite, non-empty, 
			forward invariant set $Z$ so that the partition $\P$ induced by $Z$ has the following properties:
			\begin{enumerate}[leftmargin=*]
				\item The boundary points of  $B_0$ are
					contained in the set $Z$. Here $B_0$ is the immediate
					basin of periodic attractors (possibly one-sided).
				\item $Z$ contains all parabolic periodic orbits.
				\item If $f$ is not infinitely renormalizable at a critical point $c$, 
					then  $\mathcal{P}(c)$ is the smallest periodic interval containing $c$.
				\item If $f$ is infinitely renormalizable at a critical point $c$, then 
					$\P(c)$ is a renormalization interval and each critical point $c'$
					which is contained in the orbit of $\P(c)$ has the same $\omega$-limit set. 
				\item Suppose that $x$ is not a boundary point
                                  of $[0,1]$. If $x$ is not eventually mapped
					into the interior of a renormalization interval $I$ of period $s$
					for which $f^s\colon  I\to I$ is infinitely renormalizable
					and also is not eventually mapped into $\clos(B_0)$, then
					there exists a sequence of points $y_n,z_n\in f^{-n}(Z)$
					with $y_n\downarrow x$ and
                                        $z_n\uparrow x$ as $n\to
                                        \infty$.
\item Suppose that $f\colon[0,1]\rightarrow \R$. Then there exists an
  extension of $f$ to a neighbourhood $W$ of $[0,1]$. In this case
(5) holds for $x\in\{0,1\}$ for the extension of $f$ to $W$,
provided that the orbit of $x$ is not contained in $\{0,1\}.$
If the orbit of $x\in\{0,1\}$ is contained in $\{0,1\},$ 
$x\in Z$.
			\end{enumerate}
		\end{lem}
		
		\begin{pf}
			This lemma is the analogue of Lemma 6.1 in \cite{KSS-rigidity}.
			Take $Z_a$ to be the set of boundary points of immediate basins
			of periodic attractors, and $Z_{par}$ be the set of parabolic periodic orbits of $f$.
			Note that $f$ has only a finite number of periodic attractors, 
			see \cite{MMS} and also \cite[Theorem B, Chapter IV]{dMvS}.
			These finite sets consist of periodic and pre-periodic points.
			Next for each critical point $c$ at which $f$ is not infinitely renormalizable, 
			let $I_c$ be the smallest periodic interval (i.e. $f^s(I)=I$ and 
			$f^s(\partial I)\subset \partial I$) containing $c$. Let $Z_{nr}$ be the set
			of boundary points of such intervals $I_c$. 
			For each critical point at which $f$ is infinitely renormalizable, take a periodic interval $I\ni c$
			which is so small that for each critical point $c'$ in the orbit of $I$
			one has $\omega(c)=\omega(c')$. This can be done since $f$ has only a finite number
			of critical points.
			Finally, in order to ensure that $Z$ is non-empty,
			let $Z_f$ be the set of all fixed points of
                        $f$ together with any $x\in\{0,1\}$ such that
                        the orbit of $x$ is contained in $\{0,1\}$ in the interval case
			and the set of all periodic points of period $N$ in the circle case
			(where $N$ is the minimal period of a periodic point of $f$).
			The set $Z=Z_a\cup Z_{par} \cup Z_{nr}\cup Z_{ir}\cup Z_f$
			is non-empty and has the required properties.
			Finally, whenever $I$ is a smallest renormalization interval,
			include in $Z$ a periodic point
			which lies in the interior of $I$.
			Properties 5 and 6 claimed in the lemma
			follow from the following fact:
			
			\noindent
			{\em Fact:} Assume that $f$ is of class $\mathcal{C}$ and has a fixed point.
					Moreover, assume that $x$ is a point which is not contained
					in a renormalization interval
					and that $x$ is not eventually mapped into $\clos(B_0)$.
					Then there are backward iterates of some fixed point 
					$a$ which accumulate to $x$ from the left
					(and also backward iterates which accumulate to $x$ from the right).

			\noindent
			{\em Proof of fact:}  Take $x$ 
					and a closed one-sided neighbourhood $J$ of $x$. 
					Assume by contradiction
					that no point in $J$ eventually is mapped into any fixed point. 
					The mapping $f$ is of class $\mathcal{C}$, and therefore has no wandering intervals.
					Since $J$ is not contained in the basin of a periodic attractor,
					it follows that there exists $n<n'$
					so that $g^n(J)\cap g^{n'}(J)\ne \emptyset$. 
					Hence $\hat J=\cup_{k\ge 0}f^{k(n'-n)}(J)$ is connected.
					Since $\hat J$ does not contain any fixed points of $f$, it cannot be equal
					to $S^1$ and so $\hat J$ is an interval. The interval
					$\hat J$ has the property that $f^{n'-n}(\hat J)\subset \hat J$.
					If $n'-n>1$, then $\hat J$ is a renormalization interval for $f$, so
					since $x$ is not eventually mapped into a renormalization interval, 
					it follows that $f(\hat J)\subset \hat J$ 
					and so $\hat J$ contains a fixed point of $f$,  a contradiction.
			\end{pf}

\subsection{A decomposition of the set of critical points}
\label{subsec:decomp}	
Let us define a partial ordering on the set of critical points as follows.
We say that $c_1\geq c_0$ if $c_0\in\omega(c_1)$ 
or if $c_0\in\mathrm{Forward}(c_1).$ 
We define $[c]=\{c\}$ if $c$ is non-recurrent 
and otherwise define $[c]$ to be the set of critical points 
$c'$ with $\omega(c)=\omega(c')$. 
If $c_1\geq c_0$, then $\omega(c_1)\supset \omega(c_0).$ 
This implies that if $c_1\geq c_0,$ 
$c_1$ is recurrent and $c_0$ is reluctantly recurrent, 
then $c_1$ is reluctantly recurrent too.
We associate to the partial ordering a graph 
whose vertices are the critical points
and where there is a directed edge from $c$ to $c'$ if $c'\leq c$.
Let $\Omega$ be any connected component of this graph. 
We decompose $\Omega$ as follows.			
Let $\Omega_0$ be the set of critical points at which 
$f$ is infinitely renormalizable
or which tend to a periodic, possibly parabolic,
attracting cycle. 
For $c'\in\Omega_0$ let $I(c')$ be the component of the basin
of the periodic attractor that contains $c'$ 
or a periodic interval that contains $c'$
of period high enough that for any critical point
$c\in\crit(f)$ such that either $c'\notin\omega(c)$ or
$c'\neq f^j(c)$ for any $j\in\mathbb{N}$,
we have that 
$f^k(c)\notin I(c')$ for all $k=0,1,2,\dots$.
In certain circumstances, when $f$ is infinitely renormalizable at $c'$,
we may ask that $I(c')$ be selected smaller, but this is always possible
by real \emph{a priori} bounds
(see for instance Theorem 3.1 of \cite{CvST}).

		Let $\Omega_1\subset\Omega$ be the set of critical points $c$
		such that $\omega(c)$ is persistently recurrent, and $f$ is not infinitely renormalizable at $c$.
		Choose a smallest subset $\Omega_1'\subset\Omega_1\label{notation:Omega}$ such that 
		$$\cup_{c\in\Omega_1'}\omega(c)=\cup_{c\in\Omega_1}\omega(c).$$
		Then each $\omega(c)$ contains exactly one critical point in $\Omega_1'$.
		
		Let $\Omega'\neq\Omega$ be any other connected component of the graph.
		Observe that since $\Omega'$ is a connected component of the graph, 
		for each critical point $c\in\Omega$ we can, and will,
                choose a nice real neighbourhood
		$I(c)$ so that all iterates of critical points in 
		$\Omega'$ remain outside $I(c)$.
		First, choose a nice neighbourhood $I'_1$ of $\Omega_1'$ and let
		$$I_1=\cup_{c\in\Omega_1\setminus\Omega_1'}\mathcal{L}_{c}(I'_1)\cup I'_1.
		\label{notation:critical neighbourhoods}$$

		Let $\Omega_2=\Omega\setminus(\Omega_0\cup\Omega_1)$; it is a
		set of critical points $c\in\Omega$ such that
		$f$ is reluctantly recurrent at $c$ or such that $c$ is nonrecurrent.
		Let $\Omega_2'\subset\Omega_2$
		be the set of critical points $c\in\Omega_2$ such that $\mathcal{L}_c(I(c_0))$
		is empty for each $c_0\in\Omega_1\cup\Omega_0.$
		If $c\in\Omega_2\setminus\Omega_2',$
		then we take $I(c)=\mathcal{L}_{c}(I(c_0))$
		where $c_0\in\Omega_0\cup \Omega_1$ and $I(c_0)$ is the first neighbourhood of 
		such a critical point that the orbit of $c$ enters.
		For $c\in\Omega_2'$. we take $I(c)$ to be some arbitrary nice interval about $c$
		that is so small that no critical point in $\Omega_1$ enters it.
		We define
		$$I_2=\bigcup_{c\in\Omega_2}I(c).$$

		Note that we can choose $I'_1$ so small that each critical point
		$c'\notin I'_1$ that eventually enters $I'_1$ is contained in $\Omega_2$.
		Since the real puzzle pieces shrink, we can choose $I$ so that it is a nice set.

	\subsection{Real bounds, old and new}\label{subsec:real bounds}

Suppose that $f$ is a $C^{3}$ mapping of the interval
with at most $b$ critical points of finite order. 
That is, for each critical point $c_i$ of $f$ there exists $\ell_i<\infty$,
so that near $c_i$ we can write
$f(x)=\pm(\phi(x))^{\ell_i}+f(c_i)$,
where $\phi$ is a local $C^{3}$ diffeomorphism
and $\phi(0)=0$.
We will denote this set of maps by
$\mathcal{A}^{3}_{\underline{b}}$,
where $\underline{b}=(\ell_{1},\dots,\ell_{b})$.  
We say that a mapping of
$f\colon  I\rightarrow I$ in $\mathcal{A}^{3}_{\underline b}$ is
$\tau$-extendible, if it extends to mapping 
$\tilde{f}\colon (1+\tau)I\rightarrow \mathbb{R}$
with $\tilde{f}\in\mathcal{A}^{3}_{\underline b}$,
i.e.,  so that $\tilde f$ has no additional critical points. 

A nice interval $I$ is called 
\begin{equation*}\begin{array}{ll}\label{interval geometry}
\mbox{$\rho$-\emph{nice}}&\mbox{ if for each $x\in I\cap\omega(c_{0})$ },
(1+2\rho)\mathcal{L}_{x}(I)\subset I,\\
\mbox{$\rho$-\emph{externally free}} &\mbox{ if there exists a nice interval
$J\supset (1+2\rho)I$
so that $J\supset I$ is a nice pair}\\&\mbox{ and }J\cap \omega(c_{0})\subset I,\\
\mbox{$\rho$-\emph{internally-free}} &
\mbox{ if there exists a nice interval 
$J'\subset (1+2\rho)^{-1}I$ so that $I\supset J'$ is a nice pair}\\ 
&\mbox{ and }I\cap \omega(c_{0})\subset J',\\
\mbox{$\rho$-\emph{free}}&
\mbox{ if $((1+2\rho)I\setminus(1+2\rho)^{-1}I)\cap\omega(c_{0})=\emptyset$.}
\end{array}\end{equation*}


\subsubsection{Real bounds in the persistently recurrent, finitely renormalizable case}
\label{subsubsec:real bounds - persistent}
The following is Part (a) of Theorem 3.1 of \cite{CvST}.
It asserts that we have certain real bounds for arbitrarily small
nice intervals about critical points at which $f$ is at most finitely
renormalizable and persistently recurrent. These intervals are given by
the enhanced nest, see Section 2 of \cite{CvST} or Section 8 of \cite{KSS-rigidity}
for the definition of this sequence of real puzzle pieces.

\begin{prop}[Real geometry of the enhanced nest, \cite{CvST} Theorem
  3.1]
\label{prop:real bounds persistent}
Suppose that $f\in \mathcal{A}^3_{\underline b}$. 
There exists $\varepsilon_f>0$ such that the following holds.
Assume that $c_0$ is a critical point at which $f$ is persistently
recurrent, and that either $c_0$ has even order or
that every critical point in $\omega(c_0)$ has odd order.
Suppose that $I_0\owns c_0$ is a nice interval with $|I_0|<\varepsilon_f.$ 
Let $I_{0}\supset I_{1}\supset\dots$ be the enhanced nest
for $f$ at  $c_0$.
Then
there exists $\rho>0$ such that if $I_n$ is non-terminating, 
then $I_n$ is $\rho$-nice. In addition, if $I_{n-1}$ is
non-terminating then
$I_n$ is $\rho$-externally and $\rho$-internally free,
where the externally free space is given by an interval $J\supset
(1+2\rho)I_n$
and the internal free space is given by an interval
$J'\subset (1+2\rho)J'\subset I_n.$  
Moreover, if $c_0$ is even $|J'|\geq\rho|I_n|$; 
if $c_0$ is odd then for each $\nu>0$ there exists 
$\rho'>0$ so that if $|I_{n-1}|/|I_n|<\nu$ then $|J'|\geq\rho'|I_n|$.
\end{prop}

\subsubsection{Real bounds in the infinitely renormalizable case}
\label{sec:rbinfren}		

\begin{prop}[Real bounds in the infinitely renormalizable case] 
\label{prop:realboundsinf}
For each $D$,  there exist $\kappa>1$ and $\delta>0$ with the following properties.
Let $f$ be infinitely renormalizable with $D$ critical points.
Assume that $K$ is a renormalization interval for $f$ of period $s$ with $s$  sufficiently large. 
Then there exists $\hat K',\hat K$ so that for $\hat K'\supset \hat K\supset K$
 the following holds. Let $\hat K_i'\supset \hat K_i\supset f^i(K)$ are the pullbacks of $\hat K_s'=\hat K'$
and $\hat K_s=\hat K$. Also let  $J^{\pm}_s$ be the components of $\hat K_s'\setminus \hat K_s$
and let $J_i\subset \hat K_i'$ be its pullback. Then
\begin{enumerate}
\item the intersection multiplicity of $\hat K_s',\dots,\hat K_0'$ is at most $6$.
\item  $f^s$ is a diffeomorphism on  $J_0^+$ and on $J_0^-$. 
\item $|J^+_s|=|J^-_s|\ge  \kappa |J^+_0|$ and 
$|J^+_s|=|J^-_s|\ge  \kappa |J^-_0|$.
\item $|\hat K_s'|\ge \kappa |\hat K_0'|$.
\item each component of $\hat K_s'-\hat K_0'$ has length $\delta|K|$ where $\delta\in (0,1)$.
\end{enumerate}
\end{prop}
\begin{pf}
This follows from the so-called `smallest interval argument', 
see for example \cite[Lemma VI.2.1]{dMvS} and Lemmas 2 and 3 of \cite{vSV}.
Indeed, choose $m\in \{0,\dots,s-1\}$ so that $f^m(K)$ is the smallest interval among
the intervals $K,\dots,f^{s-1}(K)$. We can take $s$ so large that there exists an interval
$M_m\subset [0,1]$ which is a $0.9$-scaled neighbourhood of $f^m(K)$ (and which
therefore does not (strictly) contain any other interval from the collection
$K,\dots,f^{s-1}(K)$ other than $f^m(K)$). Moreover, if $s$ is large enough,
then each $f^i(K)$, $i=0,\dots,s-1$ contains at most one critical point.
Now pull back $M_m$ to a neighbourhood $M_0$ of $K$. By the real Koebe Lemma,
there exists $\delta$ (which does not depend on $K$ and $s$) so that
$M_0$ contains a $\delta$-neighbourhood of $K$ (the pull back has intersection multiplicity $\le 6$).
Let $T_s:=M_0\supset K \supset f^s(K)$.  Then $T_s$ does not 
contain any other interval from the collection
$K,\dots,f^{s-1}(K)$ other than $K$. Let $T_i\supset f^i(K)$ be the pullback of $T_s$. Again, by construction, 
$T_0,\dots,T_s$ has intersection multiplicity $\le 6$, and so 
$T_0$ contains a $\delta'$-neighbourhood of $K$ (where again $\delta$ does not depend on $K$ and $s$).

In fact, $f^s$ is monotone on $T_0^+$ and on $T_0^-$ where
$T^{\pm}_0$ are the components of $T_0\setminus K$.
Indeed,  otherwise there exists $j<s$ so that $f^j(T_0^\pm)$ contains a turning point,
and then $f^{j+1}(T_0^{\pm})$ contains $f^{j+1}(K)$ which implies
that $T_s\setminus K$ contains one of the intervals $f^i(K)$, $i=1,\dots,s-1$, 
a contradiction.  If $f$ has no critical points of odd order then we are done.
Simply take $\hat K'=T_s$, $\hat K=K$ and $J_s=T_s\setminus K$.
A simple argument then shows all the required properties are satisfied.
Since we need to deal with the general case anyway, we will not give the details here.

So let us allow for the possibility that $f$ has critical points of odd order (which cannot be `seen'
when $f$ is considered as a mapping on the real line). Then
we cannot expect the previous argument to prove that $f^s$ is a diffeomorphism on 
each component of $T_0^+$ or $T_0^-$. However, $f^s|_{T_0^+}$ and $f^s|_{T_0^-}$
has  at most $D$ critical points (all of odd order), 
where $D$ is the number of critical points of $f$ of odd order.
(In fact, one can easily show that if the orbit of $K$ contains all critical points 
of $f$, that  $f^s|_{T_0^+}$ and $f^s|_{T_0^-}$ each contain at most
one inflection point, but this observation does not significantly simplify 
the arguments below, so we shall not use this observation.)
To complete the proof of the proposition and 
choose  $\hat K_s\supset \hat K_s^*\supset K$ with the required properties,
we will need the following two lemmas.

\begin{lem}\label{lem:slopes}
For each $c>1$  there exist $\lambda\in (1,\infty)$ and $\kappa>1$ so that the following holds.
Assume that $g\colon [0,1]\to \R$ is monotone increasing
with $g(0)=0$, $g(1)>1$ and with no fixed point in $(0,1)$.
Assume that $t=(x,y)$ with $0\leq x<y<1$ is so that both of the following assumptions are satisfied
\begin{enumerate}[label=(\alph*)]
\item either $x'=0$ or 
$|y'-x'|/|x'|\ge \lambda$ where $x'=g(x)$ and $y'=g(y)$;
\item  whenever $l,j,r\subset t$ are adjacent intervals (ordered from left to right)
so that $l'=g(l)$, $j'=g(j)$ and $r'=g(r)$ all have length $|g(t)|/3$,
then 
$(|g(t)|/|t|)(|g(j)|/|j|) \ge c \cdot (|g(l)|/|l|)(|g(r)|/|r|)$.
\end{enumerate}
Then there exists an interval 
$m=(a,b)\subset t$ with $0<a<b$ and 
$$|g(m)|/|m|\ge \kappa\mbox{ and }|g(m)|\ge |g(t)|/3\mbox{ and }g(b)/b\ge \kappa.$$
\end{lem}
\begin{pf}
Take $m\subset t$ with $|g(m)|\ge |g(t)|/3$, and let $\alpha$ be the component of $[0,1]\setminus m$
which contains $0$. Write $m'=g(m)$ and $\alpha'=g(\alpha)$. Note that $|\alpha|\le |\alpha'|$
because $g|\alpha$ has no fixed points other than $0$. Also $|\alpha'|\le 5|m'|$ (since
$|y'-x'|/|x'|\ge \lambda>1$ and since $|m'|\ge |g(t)|/3$).
Hence $|m'|/|m|\ge \kappa$ implies
$|\alpha'\cup m'|/|\alpha\cup m| \ge |\alpha'\cup m'|/|\alpha'\cup m| \ge (5|m'|+|m'|)/(5|m'|+|m|)\ge 
 (6|m'|)/(5|m'|+(1/\kappa)|m'|)=
  6/(5+1/\kappa)> 1 $.
  So the lemma follows if there exists an interval $m\subset t$ with $|g(m)|\ge |g(t)|/3$ and $|m'|/|m|\ge \kappa$.

In this paragraph we show that such an interval $m$ always exists.
Let $l,j,r\subset t$ be as in assumption (b) of the lemma.
By the previous paragraph we may assume that  $|j'|\le \kappa |j|$ and $|t'|\le \kappa |t|$.
Since $(|g(t)|/|t|)(|g(j)|/|j|) \ge c \cdot (|g(l)|/|l|)(|g(r)|/|r|)$, 
it follows that either $|l'|/|l|\le \kappa/\sqrt{c}$ or $|r'|/|r|\le \kappa/\sqrt{c}$.
Let us show that neither is possible when $\lambda$ is large and $\kappa>1$ is close to $1$.
Let $\alpha$ be the component of $[0,1]\setminus t$ containing $0$ and let $\alpha'=g(\alpha)$.
If $|l'|/|l|\le \kappa/\sqrt{c}$ then (by assumption (a)),
$|\alpha'\cup l'|/|\alpha\cup l| \le |\alpha'\cup l'|/|l| \le (|\alpha'\cup l'|/|l'|)(\kappa/\sqrt{c})\le
((3/\lambda)+1)(\kappa/\sqrt{c})<1$, where, using assumption (b) when 
$\alpha'\neq 0$, the last inequality holds if $\lambda$ is large enough
and $\kappa>1$ sufficiently close to $1$.
Hence $g$ has a fixed point in $\alpha\cup l$, a contradiction.
If  $|r'|/|r|\le \kappa/\sqrt{c}$ then (by assumption (a) and 
because by the previous paragraph we can assume $|l'|\le \kappa |l|$, $|j'|\le \kappa |j|$), 
$|\alpha'\cup l'\cup j'\cup r'|/|\alpha\cup l \cup j \cup r| \le 
 |\alpha'\cup l'\cup j'\cup r'|/| l \cup j \cup r| \le 
(1+1/\lambda)( | l'\cup j'\cup r'|/ l \cup j \cup r|) \le (1+1/\lambda) (3/(1/\kappa + 1/\kappa + \sqrt{c}/\kappa))<1
 $, provided $\lambda$ is large enough
and $\kappa>1$ sufficiently close to $1$. Hence $g$ has a fixed point in $\alpha\cup l\cup j \cup r$, a contradiction.
\end{pf}

\begin{lem}\label{lem:crossexpans}
There exists $c>1$ and for each $\delta>0$ there exists $c(\delta)>1$ with the following properties:
Assume that $j\subset t$ are sufficiently small intervals
and denote the components of $t\setminus j$ by $l$ and $r$.
Write $l'=f^s(l)$, $j'=f^s(j)$ and $r'=f^s(r)$ and assume that $|l'|=|j'|=|r'|$.
Moreover assume that $f^s|_t$ is a diffeomorphism 
and that there exists an interval $\hat t\supset t$ so that $\hat t,\dots,f^s(\hat t)$ has 
intersection multiplicity $\le 6$ and so that $f^s(t')$ contains a $1$-scaled neighbourhood of $f^s(t)$.
Then the following properties hold.
\begin{enumerate}[label=(\alph*)]
\item If $f^s|_t$ has a critical point in its boundary, then one has the following cross-ratio inequality
$(|f^s(t)|/|t|)(|f^s(j)|/|j|) \ge c \cdot (|f^s(l)|/|l|)(|f^s(r)|/|r|)$.
\item if the interval $t$ is so that $|t|\ge \delta \dist(t,c)$ for some critical point $c$ of $f$, 
then one has
$(|f^s(t)|/|t|)(|f^s(j)|/|j|) \ge c(\delta) \cdot (|f^s(l)|/|l|)(|f^s(r)|/|r|)$.
\end{enumerate}
\end{lem}
\begin{pf}
Let $l_i.j_i,r_i$ be the pullbacks of $l',j,r'$ contained in $f^i(t)$.
By assumption of statement (b) and the real Koebe lemma, 
it follows that there exists a universal constant $K<\infty$ (which does not depend on the interval $t$)
so that $1/K\le |l_i|/|j_i|,|r_i|/|j_i|\le K$.
If $t$ is so that $|t|\ge \delta \dist(t,c)$, then the result follows immediately
from Lemma~\ref{lem:cross-distortion}.
Alternatively there exists $i$ so that $t_i$ has a critical point in its boundary, and then
again $(|f(t_i)|/|t_i|)(|f(j_i)|/|j_i|) \ge c \cdot (|f(l_i)|/|l_i|)(|f(r_i)|/|r_i|)$ and
$(|f^s(t)|/|t|)(|f^s(j)|/|j|) \ge  \cdot (|f^s(l)|/|l|)(|f^s(r)|/|r|)$ follow from  Lemma~\ref{lem:cross-distortion}.
\end{pf}

Let us now continue with the proof of Proposition~\ref{prop:realboundsinf}.
Let $T_0^\pm$ denote the components of $T_0\setminus K$.
Let $u$ be the fixed point of $f^s$ in $\partial K$ and take
$T_0^+$ to be the component of $T_0\setminus K$ containing $u$
in its boundary. There is no loss in generality in assuming that $u$ is the right endpoint
of $\partial K$ and that $f^s|T_0^+$ is increasing.
Take $t$ be the interval in $T_0^+$
containing $u$ so that $|f^s(t)|=(1/2) |f^s(T_0^+)|$. Note that 
since $T_s$ contains a $\delta$-neighbourhood of $K$, 
we have  $|f^s(t)|\ge (\delta/2)|K|$. By Koebe, there exists $\delta'>0$ so that
$|t|\ge \delta'|K|$ and therefore
$|t|\ge \delta'\dist(t,c)$ for some critical point $c$ of $f$.

Divide $t$ into intervals $t_0,\dots,t_m$
so that $u$ is in the closure of $t_0$, so that $t_i$ is closer to $u$ than $t_{i+1}$
and so that the common boundary points of $t_i$ and $t_{i+1}$
are critical points of $f^s$ (so $m$ is at most the number of critical points of odd order
of $f$).  Let $t_i'=f^s(t_i)$ be the image of these intervals in $T_s\setminus K$. So $|t_0'\cup \dots \cup t_s'|\ge (\delta/2)|K|$.
Moreover, let $t_{-i}$ be the corresponding intervals in the other component
of $T_0\setminus K$ with $f(t_{-i})=f(t_i)$. 
Since  $K$ contains a critical point of even order, $|t_i|/|t_{-i}|$
is approximately equal to one, and $f^s(t_{-i})=f^s(t_i)=t_i'$.

Now take the universal constant $c>1$  as in Lemma~\ref{lem:crossexpans}.
Next take the universal constants $\lambda>1$  and $\kappa>1$ associated to $c>1$ as in 
Lemma~\ref{lem:slopes}.
Next take $i\in \{0,1,\dots,m\}$ minimal so that $|t_i'|\ge \lambda |t_0'\cup \dots \cup t_{i-1}'|+\epsilon$
where $\epsilon:=(1+\lambda)^{-(m-1)}(\delta/4)|K|$.
If no such $i$ exists then we get by induction
that  $|t_0'\cup \dots \cup t_i'|\le (1+\lambda)^i|t_0'|+(1+\lambda)^{i-1}\epsilon$ for each $i=1,\dots,m$
and then  $|t_0'|\ge (1+\lambda)^{-m} |t_0'\cup \dots \cup t_m'|-(1+\lambda)^{-1}\epsilon\ge(1+\lambda)^{-m} (\delta/4)|K|$.
Hence, according to Lemma~\ref{lem:crossexpans}, there exists a universal constant 
$c(\hat \delta)>1$ (corresponding to $\hat \delta=(1+\lambda)^{-m}(\delta/4)|K|$) so that
the corresponding cross-ratio inequality holds. This means that we can apply Lemma~\ref{lem:slopes}
taking for $x,y$ in   the left resp. the right endpoint of $t_0$ (so $x=u$). Hence 
there exist a universal constant $\kappa>1$ 
and an interval $j=(a,b)\subset t_0$ with $|g(j)|/|j|\ge \kappa$, $|g(j)|\ge \hat \delta/3$
and $g(b)|/b\ge \kappa$. This means that in this case we can take $\hat K=(a_{-},a)$ and $\hat K'=(b_{-},b)$,
where $a_-,b_-$ are in  $T_0^-$ with
$f(a_-)=f(a)$ and $f(b_-)=f(b)$. Note that  if $c$ is the turning point in $K$
then $|a_- - c|/|a-c|$ and $|b_--c|/|c-b|$ are close to one, when $s$ is large (and therefore
$T_0$ is small). By redefining $\kappa$
the required properties in Proposition~\ref{prop:realboundsinf} hold.

On the other hand, if there exists $i$ with $|t_i'|\ge \lambda |t_0'\cup \dots \cup t_{i-1}'|+\epsilon$
then, because $|t_i'|\ge \lambda |t_0'\cup \dots \cup t_{i-1}'|$ 
we can apply the first part of Lemma~\ref{lem:crossexpans} and 
Lemma~\ref{lem:slopes} to obtain again an interval $j_i=(a,b)\subset t_i$ with 
$|j_i'|/|j_i|\ge \kappa$, $|j_i'|\ge |t_i'|/3$
and $|f^s(b)-u|/|b-u|\ge \kappa$ where $b$ is the right endpoint of $j_i$. Note that
$|t_i'|\ge \epsilon$ and so $|j_i'|\ge \epsilon/3$.
So we can take in this case  $\hat K=(a_{-},a)$ and $\hat K'=(b_{-},b)$,
where as before $a_-,b_-$ are in $T_0^-$ with $f(a_-)=f(a)$ and $f(b_-)=f(b)$.
Again by redefining $\kappa$
the required properties in Proposition~\ref{prop:realboundsinf} hold.
\end{pf}

		\subsubsection{Real bounds in the reluctantly recurrent case}
			\label{subsubsec:real bounds - reluctant}
			A subset $\Omega\subset\crit(f)$ is called a \emph{block} if it is contained in a connected component of 
			the graph of critical points, as defined in Subsection \ref{subsec:decomp} and if
			$\mathrm{Back}(c)\subset\Omega$ for all $c\in \Omega$.
			A block of critical points is called \emph{non-trivial} if it is disjoint from the basins
			of periodic attractors and there exists $c\in\Omega$ with infinite forward orbit.
			Let $\Omega_{r}\label{notation:Omega_reluct}$ be a non-trivial block of critical points 
			such that each recurrent critical point in $\Omega_{r}$
			is reluctantly recurrent, and if $\Omega'$ is the component of the graph of critical points that contains
			$\Omega_{r}$, then $f$ is not infinitely renormalizable at any $c'\in\Omega'$.
			
			Let $\Dom_{\Omega_r}(I)$ denote the union of the components of the domain of 
			the return mapping to $I$ which intersect the orbit of $c$ for some $c\in\Omega_r$.
			Let $\Dom'_{\Omega_r}(I)=\Dom_{\Omega_r}(I)\cup I$
			and $\Dom_{\Omega_r}^*(I)=\Dom_{\Omega_r}(I)\cap I$.
			We let $R_I\colon \Dom_{\Omega_r}^*(I)\rightarrow I$ denote the first return mapping to $I$
			restricted to $\Dom_{\Omega_r}^*(I)\rightarrow I$.
			Let $I$ be a real neighbourhood of $\Omega_r$ such that
			each component of $I$ contains exactly one critical point in $\Omega_r$
			and for each component $J$ of $\Dom_{\Omega_r}(I)$ either $J$ is a component of $I$
			or $J$ is compactly contained in a component of $I$. 
			
			Let 
			$$C_1(\Omega_r)=\Omega_r\setminus \mathrm{Dom}_{\Omega_r}(I) \mbox{ and }
			C_2(\Omega_r)=\{c'\in\Omega_r: I(c')\subset\mathrm{Dom}_{\Omega_r}(I)\}.$$

			For intervals $J\subset T$, we let $L$ and $R$ denote the components of $T\setminus J$ and define 
			$$\gap(L,R)=\frac{1}{\Space(T,J)}:=\frac{|T||J|}{|L||R|}.$$
			We define $\gap(I)$ by
			$$\gap(I)=\inf_{(J_{1},J_{2})}\gap(J_{1},J_{2}),$$
			where $(J_{1},J_{2})$ runs over all distinct pairs of intervals in $\Dom_{\Omega_r}'(I)$.
			
			Let $$I^{*}=\cup_{c'\in C_{2}(I)}I(c'),\quad I^{\#}=I\setminus I^{*}.$$
			We define $\Space(I)$ by
			$$\Space(I)=\inf_{J}\,\Space(\Comp_{J}I,J),$$
			where the infimum is taken over all components $J$ of $\Dom_{\Omega_r}^*(I)\cap I^{\#}$.

\begin{figure}
\resizebox{0.9\textwidth}{!}{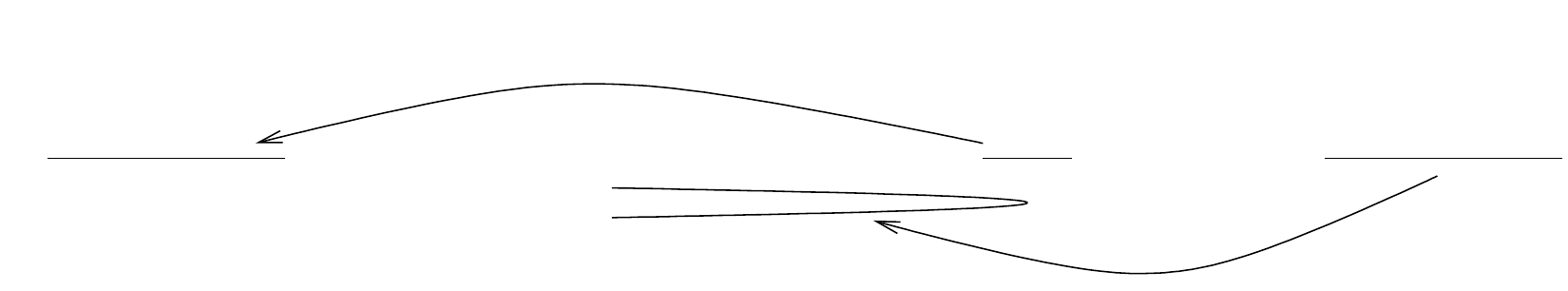}
 \caption{When $\mathrm{Cen}_1(I)$ is small.}
 \label{fig:cen1}
\end{figure}

\begin{figure}
\resizebox{0.9\textwidth}{!}{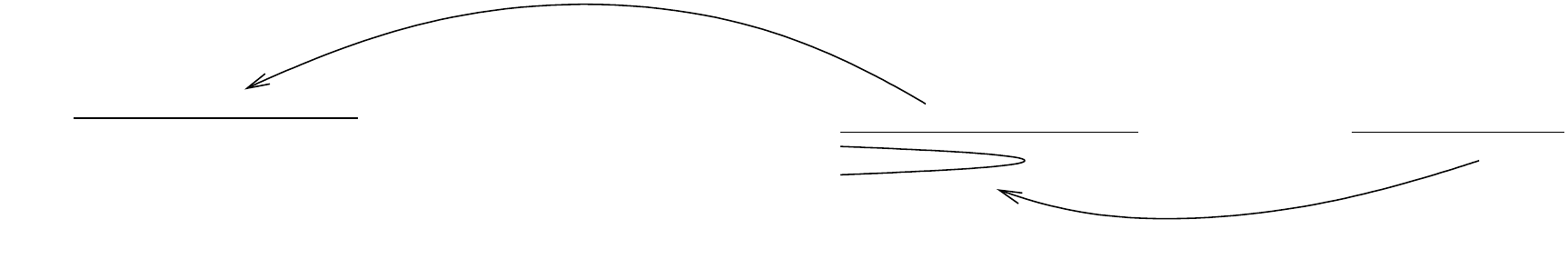}
 \caption{When $\mathrm{Cen}_2(I)$ is small.}
 \label{fig:cen2}
\end{figure}

			Let $\Omega_{r, even}$ denote
			the even critical poins in $\Omega_r$.
			For any $c\in\Omega_r$, let $\hat{J}(c)$ be the component of
			$\Dom_{\Omega_r}'(I)$ which contains $f(c)$ and define
			$$\cen_{1}(I)=\max_{c\in\Omega_{r,even}\setminus C_{2}(I)}\frac{|\hat{J}(c)|}{|f(I(c))|},$$
			$$\cen_{2}(I)=\max_{c\in C_{2}(I)\cap \Omega_{r,even}}\Bigg(\Bigg|\frac{|\hat{J}(c)|}{|f(I(c))|}-2\Bigg|\Bigg),$$
			and $\cen(I)=\max(\cen_{1}(I),\cen_{2}(I))$. Figures \ref{fig:cen1} and \ref{fig:cen2} indicate the
			arrangement of intervals when these quantities are small. It is important to notice that
			in the definition of $\cen_1(I)$ and $\cen_2(I)$ we only take the 
			$\max$  over the even critical points. For odd critical points, having $\Space(I)$
			big implies that $\cen_1(I)$ is small, and we do not require the bound on
			$\cen_2(I)$ to control the loss of angle when
                        we pull back a Poincar\'e disk
			through an odd critical point, see Lemma \ref{lem:zd}.
			
			We decompose $\Omega_r$ it into two parts:
			\begin{itemize}\label{notation:reluct decomp}
				\item  $\Omega_{r,e}$ is a set of critical points in $\Omega_{r}$
					so that if $c\in\Omega_{r,e}$ has odd order,
					then there exists $c'\in\Omega_{r,e}$
					with even order so that $c\geq c'$,
					where we use the ordering from Subsection \ref{subsec:decomp}.
				\item $\Omega_{r,o}=\Omega_{r}\setminus\Omega_{r, e}$.
			\end{itemize}
			Observe that $\Omega_{r,e}$ is not the same as $\Omega_{r,even}$.
			It follows from the definition that for each 
			$c\in\Omega_{r,o}$ and each critical point $c'\in\omega(c)$ we have that either
			\begin{itemize}
				\item $c'$ is a critical point of odd order,
				\item $\omega(c')$ is persistently recurrent,
				\item$c'$ is in the basin of a periodic attractor or
				\item $c'$ is eventually mapped into a 
periodic interval where the mapping is infinitely renormalizable.
			\end{itemize}
			
\begin{prop}[Prop. 3 of \cite{KSS-density} ]\label{prop:real bounds reluctant}
				Let $\Omega$ be either $\Omega_{r,e}$ or $\Omega_{r,o}$.
				For any $\varepsilon>0$ and $C>0$, there exists an arbitrarily small admissible neighbourhood
				$I$ of $\Omega$ such that $\gap(I)>C, \Space(I)>C$ and $\cen(I)<\varepsilon$.
			\end{prop}
			\begin{pf}
				We will only give an outline, and we will refer heavily to \cite{KSS-density}.
				In the notation of \cite{KSS-density}, one should take as
				$\Omega$ either the set $\Omega_{r,e}$
				or the set $\Omega_{r,o}$.
				
				As we shall see below, the proof for the case that we take 
				$\Omega_{r,o}$ is unchanged. Indeed critical points  $c\in \Omega_{r,o}$ 
				do not accumulate on other reluctantly recurrent critical points $c'$, such that 
				$c'$ has even order. Of course $c$ can accumulate
				on a critical point $c'$ for which $\omega(c')$ is persistently recurrent
				and which has even order, but these critical points $c'$ will not cause troubles.

				So let us concentrate the discussion on the case that $ \Omega_{r,e}$
				(and implicitly explain why the original proof goes through for $\Omega_{r,o}$).
				So define an admissible open set $I$ of $\Omega_{r,e}$ as in \cite{KSS-density}.
				By the definition of $\Omega_{r,e}$, each critical point of odd
				order in $\Omega_{r,e}$ is eventually mapped into $I$.
				So for each $c\in \Omega_{r,e}$ of odd degree we can define 
				$\hat I_c=\LL_c(I)$ (here $\LL$ is defined as in Section~\ref{sec:definitions}) and
				\begin{equation}
				\hat I= \bigcup_{c} \hat I_c \,\,  \bigcup \,\,  I\,  ,\label{eq:evenodd}
				\end{equation}
 				where the union runs over all critical point of
				odd degree in $\Omega_{r,e}$. Note that
				$\hat I$ is an admissible open neighbourhood of $\Omega_{r,e}$.
				
				Note that in \cite{KSS-density}, Proposition 3, $I$ is of the form $T_N=\A^N(T)$, see
				page 159 of \cite{KSS-density}, and that for each critical point $c$,
				\begin{equation}
					\mbox{$T_N$ is well-inside $T_{N-1}$ and $T_N\supset \comp_c (\Dom(T_N))$},
					\label{eq:extension}
				\end{equation}
				see  line 3, page 158 on  \cite{KSS-density}).
				
				The proof of Proposition 4, given in \cite{KSS-density} is for
				maps with a non-recurrent critical point. This goes through as before
				(since the results in \cite{vSV} work for critical points of all orders).
				So we suppose that $c$ is reluctantly recurrent at a critical point $c$.
				We need to show that the following claim holds: for any $C>0$ and any 
				$c\in\Omega$, there exists an arbitrarily small $C$-nice interval.
				By Lemma 6.5 of \cite{KSS-rigidity}, for every critical puzzle piece 
				$I(c)$ there exist arbitrarily large positive integers $m$ such that
				$f^m(c)\in I(c)$, and letting $I^m(c)=\comp_c f^{-m}(I(c))$, we have that
				the mapping $f^m\colon  I^m(c)\rightarrow I(c)$ has bounded degree.
				Since $f$ is at most finitely renormalizable at $c$, 
				by Theorem A of \cite{vSV},
				we can assume that
				$I(c)$ is $\delta$-nice for some $\delta>0$, universal.
				Then by Lemma 3.11 of \cite{CvST}, there exists $\eta>0$ so that if
				$(1+2\eta)I^m(c)\subset I(c)$, $I^m(c)$ is $C$-nice.
				Taking $m$ sufficiently big finishes the proof.		
								
				In the definition of $\A(T)$  on 156-157 in \cite{KSS-density} it is important
				that $\Omega$ only contains critical points of even order.
				Let us first explain why a change is needed. If the critical points are of even order, then one can ensure as in
				\cite{KSS-density} that a nice interval $I$ around a critical point is symmetric around $c$ (in the sense
				that $f(\partial I)$ consists of one point). This means that if $\A(T)(c)$ is much smaller than $T(c)$,
				then $\A(T)(c)$ is roughly in the middle of $T(c)$ (and this will be why one
				can get $\cen_2(I)$ small). This is no longer true for a nice interval
				that is not symmetric around a critical point (which certainly can happen 
				for critical points
				of odd order). Let us make some comments about
				the definition of $\mathcal{A}(T)(c)$. When  $c\in C_1(T)$,
				we should ensure that $\spa(T(c),\A(T)(c))$ is large. Since we want
				$\A(T)$ to be a nice  neighbourhood of $\Omega$ we do this by taking
				$\A(T)(c)$ so that its boundary points  are also 
				boundary points of domains of $R_T$.  (This is possible since $R_T$ has infinitely many domains 
				which accumulate from both sides to $c$.) Note that this part of the construction is not combinatorially defined,
				but if two maps are conjugate, then one can make a choice a combinatorial choice
				which works for both maps. When $c\in C_2(T)$ then take $\mathcal{A}(T)(c)$ to be the maximal
				interval containing $c$ such that
$\frm_T^{k(c)}(\mathcal{A}(T)(c))$ is contained in a component of $\Dom(T)$.
				(What was written here in \cite{KSS-density} was a typo.)

				Then Lemma 3 holds as before except that in statements (2) and (3)  `diffeomorphically'
				should be replaced by `homeomorphically'. Note that $k(c)$ is bounded by
				$\#\Omega$ because otherwise one of the components
				of $I$ would be periodic, and also note that  the pullback of $T$ to $J(c)$ meets
				at most
				$\#\Omega$ of critical points (all of odd order). 
				Then Lemma 4 goes through as before.
				As stated in \cite{KSS-density},
				Lemma 5 is only valid for maps without parabolic periodic orbits
				(this is the context in which Theorem 3 in \cite{KSS-density} was used in that paper).
				However, all we need is Lemma~\ref{lem:koebe} of this paper which holds even for maps
				with parabolic periodic points. The proof of Lemma 6 goes through as before
				(note that not only $|\A(T)(c)|/|T(c)|$ but also $\spa(T(c),\A(T))$ is large). 
				Lemma 7 also goes through: in the final
				paragraph it is shown that $\cen_2(\A(T))$ is small, which holds
				because $c'$ is a critical point of even order and so
				$\A^2(T(c')$ and $\A(T)(c')$ are symmetric around $c'$.
				This allows us to argue as before: although the presence of critical points of odd order
				implies that $f^s\colon  J\to f^s(J)$ extends merely to a homeomorphism onto $T(c')$,
				but as the order of the corresponding chain is still bounded by the number of critical points,
				we can use part (2) of
                                Lemma~\ref{lem:koebe} to pull back the space.	
			\end{pf}

\section{Compatible complex box mappings and complex bounds}\label{sec:complex bounds}
Let $\bm U$ and $\bm V$ be open Jordan disks in $\mathbb{C}$. 
We say that a mapping $F\colon \bm U\rightarrow \bm V$ is a 
\emph{$K$-qr branched covering}
if it can be decomposed as $F=G\circ H$,
where $H\colon \bm U\rightarrow \bm U$ is 
$K$-qc homeomorphism and 
$G\colon \bm U\rightarrow \bm V$
is a holomorphic branched covering.

A mapping $F\colon \bm{\mathcal{U}}\rightarrow \bm{\mathcal{V}}$ is a 
\emph{holomorphic box mapping}\label{def:box mapping} 
if $F$ is holomorphic and 
$\bm{\mathcal{U}}$ and $\bm{\mathcal{V}}$ are 
open subsets of the complex plane where
\begin{itemize}
\item $\bm{\mathcal{V}}$ is a union of finitely 
many pairwise disjoint Jordan disks;
\item every connected component 
$\bm V$ of $\bm{\mathcal{V}}$ is either a connected component of
$\bm{\mathcal{U}}$ or $\bm V\cap\bm{\mathcal{U}}$ is a union of
 Jordan disks with pairwise disjoint closures
that are compactly contained in $\bm V$; 
\item for each component $\bm U$ of $\bm{\mathcal{U}}$, $F(\bm U)$ 
is a component of
$\bm{\mathcal{V}}$ and $F|_{\bm U}$ is a proper mapping with at most one critical point;
\item each connected component of $\bm{\mathcal{V}}$ contains 
at most one critical point of $F$;
\item each component of $\bm{\mathcal{U}}$ is either a component of 
$\bm{\mathcal{V}}$ or it is compactly contained in 
$\bm{\mathcal{V}}$.
\end{itemize}
If $F$ is as in the definition of holomorphic box mapping,
except that $F$ is not holomorphic, but instead there exists $K\geq 1$
such that
for each component $\bm U$ of $\bm{\mathcal{U}}$,
$F|_{\U}$ is a $K$-qr branched covering onto its image, 
a connected component of $\bm{\mathcal{V}}$, we call
$F$ a $K$-\emph{qr box mapping}.
We will use the terminology \emph{complex box mapping} to
refer to both holomorphic box mappings and
to $K$-qr box mappings, and unless otherwise stated,
we assume that covering maps are quasiregular.
A $K$-qr box-mapping whose domain has exactly one component is called a
$K$-\emph{qr polynomial-like mapping}.
	We define the \emph{filled Julia set} of a $K$-qr box mapping (which includes the case of a $K$-qc
	polynomial-like mapping)
	$F\colon\bm{\mathcal{U}}\rightarrow\bm{\mathcal{V}}$ to be
	$$K(F)=\{z\in\bm{\mathcal{U}}:f^n(z)\in \bm{\mathcal{U}}\mbox{ for all } n=0,1,2,\dots\},$$
	and we define the \emph{Julia set} of
        $F\colon\bm{\mathcal{U}}\rightarrow
\bm{\mathcal{V}}$
	to be $J(F)=\partial K(F)$.
The filled Julia set of a complex box mapping is
compactly contained in its domain.

We shall need to use the \textit{gap condition} and the \textit{extension condition}.
The {\em gap condition} \label{def:gapextension} 
states that there exists $\delta>0$ so that
for each two components $\U_i$ and $\U_j$ of $\bm{\mathcal{U}}$ there exist
$\U'_i\supset \U_i$ and $\U'_j\supset\U_j$ so that 
\begin{enumerate}
\item $\mathrm{mod}(\U'_i\setminus \U_i)\ge \delta$ and 
$\mathrm{mod}(\U'_j\setminus \U_j)\ge \delta$ and
\item either $\U'_{i}\cap \U_{j}=\emptyset$ or $\U'_{j}\cap \U_{i}=\emptyset$;
\item in addition, if $\U_i$ is contained in the component $\V_i$ of
$\bm{\mathcal{V}}$,
then either $\V_i=\U_i$ or 
$\mathrm{mod}(\V_i\setminus \U_i)\ge \delta$.
\end{enumerate}
The {\em extension condition} requires that there exist $\delta>0$
and $\bm{\mathcal{V}}'\supset \bm{\mathcal{V}}$ so that for each component
$\U_i$ of $\bm{\mathcal{U}}$ there exists $\U_i'\supset \U_i$ 
so that
\begin{enumerate}
\item $\mathrm{mod}(\V'_j\setminus \V_j)\ge \delta$ for each component 
$\V_j'$  of $\bm{\mathcal{V}}'$
and with $\V_j$ the component of 
$\bm{\mathcal{V}}$ contained in $\V_j'$, and
\item
if $F(\U_i)=\V_j$ then there exists a $K'$-quasiregular extension of
$F$ to $\U_j'$ 
so that
$F\colon \U_i'\to \V_j'$ is a branched covering with the same degree
as
$F\colon \U_i\to \V_j$ (usually, $K'\leq 1+\mu(\V_j')$).
\end{enumerate}
These conditions will enable us to control the shape of the domains
$\U_i$,
which helps us in constructing
quasiconformal homeomorphisms on 
$\bm{\mathcal{V}}\setminus\bm{\mathcal{U}}$,
see Theorem~\ref{thm:external conjugacy}.
		
\medskip
		
A component of $F^{-n}(\V)$ of a component $\V$ of the
                range 
$\bm{\mathcal{V}}$ will be called a  
{\em puzzle piece of level $n$}\label{page:defpuzzle}.
It is important to observe that
puzzle pieces are either nested or disjoint.

\begin{definition}[Renormalizable box mapping]\label{def2}
A box mapping $F\colon\bm{\mathcal{U}}\to\bm{\mathcal{V}}$
is called {\em renormalizable} at $c\in\crit(F)$
if there exists $s>1$ (called the {\em period}) and 
a puzzle piece $\bm W$ 
containing a critical point $c$ of $f$ such that
$f^{ks}(c)\in \bm W$, $\forall k\ge 0.$
\end{definition}

\subsection{Compatible complex box mappings and 
bounds at persistently recurrent points}\label{subsec:complex bounds}
	
Suppose that $F\colon\bm{\mathcal{U}}\rightarrow\bm{\mathcal{V}}$ is a 
$K$-qr box mapping.
Let $\rho>0$. 
Let $\V$ be a component of $\bm{\mathcal{V}}$.
We say that $\V$ is
$\rho$-\emph{nice}\label{def:rhonice}, if for each $x\in \V\cap\omega(c_0)$, we
have that 
$\mod(\V\setminus \mathcal{L}_x(\V))>\rho$. We say that
$\V$ is $\rho$-free if there exist puzzle pieces $\bm P^+$ and $\bm P^-$
such that $\bm P^-\subset \V\subset \bm P^+$, and both
$\mod(\V\setminus\bm P^-), \mod(\bm P^+\setminus \V)>\rho.$
See page~\pageref{def:bddgeo} for the definition of
bounded geometry.

The following is the main theorem of \cite{CvST}.

\begin{thm}[Complex bounds at persistently recurrent critical points]
\label{thm:box mapping persistent}
There exist $\rho>0$ and $C>0$, universal,
such that for each $n$ sufficiently large the following holds. 
Suppose that $f\in\mathcal{A}^3_{\underline b}(M)$.
Suppose that $c_{0}$ is a critical point such that
$f$ is persistently recurrent on $\omega(c_{0})$. 			
\begin{enumerate}
\item If $f$ is finitely renormalizable, 
there exists an arbitrarily small,
combinatorially defined, nice real neighbourhood $I$ of
$\crit(f)\cap\omega(c_0)$ such that
the first return mapping to $I$ extends to a $K_n$-qr box mapping
$F\colon\bm{\mathcal{U}}\rightarrow\bm{\mathcal{V}}$ with
$\crit(F)=\crit(f)\cap\omega(c_0)$ that is
$\rho$-extendible and such that the components of $\bm{\mathcal{V}}$ are
$\rho$-nice, $\rho$-free, the components of $\bm{\mathcal{U}}$ and 
$\bm{\mathcal{V}}$
have $\rho$-bounded geometry, and there exists $\theta'\in(0,\pi)$
such that each component $\U$ of $\bm{\mathcal{U}}$ or
$\bm{\mathcal{V}}$ there exists an interval $\tilde U$
such that $\U$
is contained in $D_{\theta'}(\tilde U)$, $U$ is well-inside $\tilde U$
and $(\tilde U\setminus U)\cap\omega(c_0)=\emptyset$.
\item If $f$ is infinitely renormalizable, then for all $s$ sufficiently big,
if $J$ is a periodic interval for $f$ of period $s$, then $f^s\colon J\rightarrow J$
extends to a $K(\V)$-qr polynomial-like mapping
$F\colon \U\rightarrow \V$ with
$\crit(F)=\crit(f)$,
$\mod(\V\setminus\U)>\rho,$ $\U$ has $\rho$-bounded geometry
and $\mathrm{diam}(\U)<C|J|$.
\end{enumerate}
\end{thm}

Here $K_n$ and $K(\V)$ are defined (up to a constant multiple) as
follows.  
Let $\mu_n=\max\diam(\J)$, where the maximum is taken over all
pullbacks $\J$ of $\bm{\mathcal{V}}$ that intersect
$\omega(c_0)$. Set $K_n=(1+\mu_n)/(1-\mu_n).$
Observe that $\mu_n$  is small when $|I_n|$ is small.
Similarly, if $F\colon \U\rightarrow \V$
is a polynomial-like mapping that extends the return map
$f^s\colon J\rightarrow J$, where $J$ is a periodic interval, we define
$$\mu(\V)=\max_{1\leq i\leq s}\diam(\comp_{f^i(c_0)} f^{-(s-i)}(\V))
\mbox{ and } K(\V)=(1+\mu(\V))/(1-\mu(\V)).$$
Again $\mu(\V)$ is close to zero when $\V$ is small.

In the finitely renormalizable (and persistently recurrent) case,
the intervals $I(c_0)$ in the theorem are given by the enhanced nest about $c_0$
(see \cite{KSS-rigidity} and \cite{CvST} where the infinitely
renormalizable case was considered).

			\subsubsection{Extending complex box mappings when $\omega(c)$ is persistently recurrent.}
				\label{subsubsec:persistent infinitely many branches}
				
				By Theorem~\ref{thm:box mapping persistent},
				if $c$ is a persistently recurrent critical point
				at which $f$ is at most finitely renormalizable,
				we have that
				there exist arbitrarily small nice real neighbourhoods $I_1$ 
				of $\crit(f)\cap\omega(c)$ for which there exists a 
				complex box mapping that extends
				$R_{I_1}\colon  \Dom_{\omega(c)}(I_1)\to I_1$.
				In this subsection we will still assume that $c$ is persistently recurrent,
				but  want to consider also  (infinitely many) branches of 
				$R_{I_1}$ which do {\em not} intersect $\omega(c)$.
				In fact, we need to treat different omega limit sets, $\omega(c),$ 
				on which $f$ is persistently recurrent simultaneously. 
				
So let $\Omega$ be any connected component
of the graph associated to the set of critical points. We proceed exactly as in
Subsection~\ref{subsec:decomp}:
First, decompose $\Omega=\Omega_0\cup\Omega_1\cup\Omega_2$.
Recall that $\Omega_0$ are critical points that are either
contained in a basin of attraction of $f$ or at which $f$
is infinitely renormalizable, $\Omega_1$ are the critical points 
at which $f$ is persistently recurrent and at most finitely renormalizable
and $\Omega_2$ are the critical points at which $f$ is 
reluctantly recurrent or nonrecurrent.
Now, identify a smallest subset $\Omega'_1$ with
$$\cup_{c\in\Omega'_1}\omega(c)=\cup_{c\in\Omega_1}\omega(c).$$
Choose a neighbourhood $I'^{\Omega}_1$ of $\Omega'_1$ and 
 choose neighbourhoods $I_0^{\Omega}$, $I_1^{\Omega}$, $I_2^{\Omega}$ of $\Omega_0$,
$\Omega_1$, $\Omega_2$, respectively, as in Subsection~\ref{subsec:decomp}.
Let 
$$\I=\bigcup_{\mathrm{connected\ components}\ \Omega} (I_0^{\Omega}\cup I_1^{\Omega}\cup I_2^{\Omega}).
\label{notation:all intervals}$$
Recall that we choose the intervals around the critical points so that $\I$ is a nice set. 
From now on, we will drop the superscript $\Omega$ from the notation. 
				
Let 
$$\Dom_{\mathrm{diff}}(I_1)\subset I'_1 \label{notation:dom_diff}$$ 
be the union of the components $J$ of $\Dom(I_1)$ which are {\em inside} $I'_1$
and which 
\begin{enumerate}
\item do {\em not} intersect $\bigcup_{c\in \Omega_1}\omega(c)$;
\item if  $s$ is the return time of $J$ to $I_1$,
then $f^i(J)\cap  \I=\emptyset$ for all $i=0,1,\dots,s-1$.
\end{enumerate}
So $f^s|_J$ is a diffeomorphism.
It will be important to consider the first entry map
to $I_1\cup I_0$ restricted to the above sets
because we will need to control orbits of 
critical points $c'\notin \Omega_1$ which linger for long periods near
$\cup_{c\in \Omega_1} \omega(c)$.
Let
$$R_{\hat I_1}\colon  \Dom_{\Omega_1}(I_1)\cup
\Dom_{\mathrm{diff}}(I_1)\cup (\Dom(I_0)\cap I'_1)\to I_1\cup I_0.
				\label{pageref:F_it}$$
			 (We write $R_{\hat I_1}$ for this first entry map
				because this mapping is not defined on $I_0$.)
				The mapping $R_{\hat I_1}$ will be used to  control the pieces of the orbit of $c'$ 
				during which $c'$ remains close to $\omega(c)$.
				Note that this first return mapping can have infinitely many components

\begin{thm}[Complex box mappings in the persistently recurrent
case with infinitely many branches included]
\label{thm:box mapping persistent infinite branches}
There exists $\theta'\in(0,\pi/2),$
such that
for each $\theta\in(0,\theta')$,
there exists an integer $m>0$, $\delta>0$ and a complex box mapping
$F_{it}\colon \bm{\mathcal{U}}_{it}\rightarrow \bm{\mathcal{V}}_{it},$
which is an extension of the $m$-th iterate of
$R_{\hat{I}_{i}}$. 
For each $c\in\Omega'_{1}\cup \Omega_0$, the component
$\V(c)$ of $\bm{\mathcal{V}}_{it}$
is a Poincar\'e lens domain based on $I(c)$, that is,
$$\V(c)=D_{\pi-\theta}(I(c)).$$
Moreover, $F_{it}\colon \bm{\mathcal{U}}_{it}\rightarrow \bm{\mathcal{V}}_{it}$
satisfies the extension and gap properties,
and if $\bm D$ is any component of either $\bm{\mathcal{U}}_{it}$ or
$\bm{\mathcal{V}}_{it}$, then $\bm D$ has $\delta$-bounded geometry.	
\end{thm}
				
\begin{pf}
Let $F_{\Omega_1}\colon \bm{\mathcal{U}}_{\Omega_1}
\rightarrow \bm{\mathcal{V}}_{\Omega_1}$
be the box mapping given by Theorem~\ref{thm:box mapping persistent}
which extends the real return mapping to $I_1$. 
Note that we are free to choose
$I_1$ so that this extension exists.
First, observe that by Propositions~\ref{prop:real bounds persistent}
and \ref{prop:realboundsinf},
the following real analogues of the
extension and gap properties hold for
$R_{\hat{I}_1}$:
\begin{enumerate}[label=(\arabic*)]
\item there exists a constant $\rho>0$ such that
if $J$ is a component of $\Dom(R_{\hat{I}_1}),$
then either $J$ is a component of $I_1$, 
$J$ is well-inside $I_1$, or $J$ is disjoint from $I_0\cup I_1$ and
$J\cap\Big(\bigcup_{c\in \Omega_1}\omega(c)\Big)\neq\emptyset$;
\item for any components $J_1$ and $J_2$
of $\Dom(R_{\hat{I}_1})$, there exist
$J'_i\supset J_i$ 
such that $J_i$ is well-inside $J_i'$ is for $i=1,2;$ 
\item either $J_1\cap J_2'=\emptyset$ or $J_2\cap J_1'=\emptyset;$
\item there exists a set $I'$ with $I'\supset I_0\cup I_1$
such that for each component $I(c)$ of $I_0\cup I_1$,
there exists a component $I'(c)$ of $I'$ such that
$I(c)$ is well-inside $I'(c)$;
\item for each component $J$ of $\Dom(R_{\hat{I}_1})$,
there exists an interval $J'\supset J$ such that if
$R_{\hat{I}_1}\colon J\rightarrow I(c)$,
then $R_{\hat{I}_1}|_J$ extends to a mapping
from $J'$ to $I'(c)$ with the same degree as 
$R_{\hat{I}_1}|_J$, where $I'(c)$ is the interval from (4).
\end{enumerate}
Let $c\in \Omega_0\cup\Omega_1$.
Notice that we cannot rule out that there are 
critical points $c'\in\Omega_2$ which 
enter $(1+2\delta)(I(c))$, where $I(c)$ is the component of 
$I_0\cup I_1$ that contains $c$.
However, any first entry domain 
$J=\mathcal{L}_x((1+2\delta)I(c))$ with $x\in\omega(c)$,
does not contain any critical point $c'\notin\omega(c)$
provided that $I(c)$ is sufficiently small:
if $c'\notin\omega(c)$ was a critical point 
with $c'\in J$ then $J$ is an interval
that contains both $c'$ and $x\in\omega(c)$, 
which is a first entry domain to $(1+2\delta)I(c)$.
This is impossible by the absence of wandering intervals.

Let us now show that the required complex extension exists.
Let $\theta_0,\theta_1\in(0,\pi/2)$.
Let $\theta=\max\{\theta_0,\theta_1\}$.
For each $c\in\Omega_0$
let $\V(c)=D_{\pi-\theta_0}(I(c))$,
and for each $c\in\Omega_1'$ let
let $\V(c)=D_{\pi-\theta_1}(I(c)).$
For each $c\in\Omega_1\setminus\Omega_1'$
let $s_c>0$ be smallest so that $f^{s_c}(c)\in I(c')$
for some $c'\in\Omega_1'\cup\Omega_0$, and set 
$\V(c)=\comp_{c}f^{-s_c}(\V(c')).$
Provided that $I_0\cup I_1$ is sufficiently small, the
domains $\V(c)$, $\V(c'), c\neq c',$ are disjoint.
Let $\V=\cup_{c\in\Omega_0\cup\Omega_1}\V(c)$.
Taking $\theta_0,\theta_1$, smaller if necessary (see below),
this will be the range of the complex box mapping
we construct. 
					
For each component $J$ of 
$\Dom_{\Omega_1}(I_1)\cup
\Dom_{\mathrm{diff}}(I_1)\cup(\Dom(I_0)\cup I'_1)$
Let $s_J>0$ be minimal so that $f^{s_J}(J)\subset I_1\cup I_0$,
let $c$ denote the critical point contained in the component
of $I_0\cup I_1$ containing $f^{s_J}(J),$ and set
Let $\U'_J=\comp_Jf^{-s_J}(\V(c)).$
Define a complex extension of $R_{\hat{I}_1}$,
$\hat {F}\colon \UU'\rightarrow \VV$
by $\hat{F}(\U'_J)=f^{s_J}(\U'_J)$.
Note that $\hat {F}\colon \UU'\rightarrow \VV$
is not a box mapping, since we cannot ensure that 
the components of the domain are compactly contained in the range.

By \cite{CvST}, Corollary 10.2,
there exist $\theta'\in(0,\pi)$ and an interval $\tilde J\supset J$
such that $\U_{\Omega_1}(J)\subset D_{\theta'}(\tilde J)$.
Moreover, $J$ is well-inside $\tilde J$ and
$(\tilde J\setminus J )\cap\omega(c_0)=\emptyset.$					
This implies that the components of the 
complex pullbacks lie close to their real traces.
Moreover,
$R_{\hat{I}_1}$ is a diffeomorphism on each component of 
$\Dom_{\mathrm{diff}}(I_1)\cup(\Dom(I_0)\cap I'_1),$
which implies
by Lemma~\ref{lem:lens angle control}
that pulling back along these branches results in little loss of angle.

\textit{Claim 1:} 
There exists $m$ such that each component of
$\hat{F}^{-m}(\V)$ which intersects the real line is well-inside of $\V$.

\textit{Proof of Claim 1:} First, consider the case when 
$\hat{F}^m(J)\cap I_0=\emptyset.$
Then $\hat{F}^m$ corresponds to $m$ iterations of $\hat F$.
Assume that $m$ is large and
consider a component $J$ of $R^{-m}_{\hat{I}_1}(I_1)$
and the corresponding domain $\U^{m}_{J}=\comp_J f^{-m}(\V).$ 
Let $I_J$ be the component of $I_1$ that contains $J$.
Let $0\leq k\leq m$ be maximal so that
$(R_{\hat{I}_1})^k(J)$ is contained in $\Dom_{\mathrm{diff}}(I_1).$
Then $(R_{\hat{I}_*})^{k+1}$ is a diffeomorphism on $J$.
This holds because when $c\in\Omega_1$,
then $\omega(c)$ does not intersect $\Dom_{\mathrm{diff}}(I_1)$, and,
moreover, by the definition of the domain of $R_{\hat{I}_1}$,
the iterates of $J$ under $R_{\hat{I}_1}$ up to the $k+1$-st,
avoid critical points $c$ in $\Omega_2$.
Indeed, one can find an interval $J_k\supset J$ such that
$(R_{\hat{I}_1})^{k+1}(J_k)$ is a component of $I_1$ and 
$(R_{\hat{I}_1})^{k+1}$ is a diffeomorphism on $J_k$. 
Since $f$ has no wandering intervals, $|J_k|/|I_J|$
is small if $k$ is large. 
Since each component of the domain of $R_{\hat{I}_1}$
is well-inside the component 
of $I_1$ that contains it, it follows that there exists a large constant $C$ so that 
$(1+C)J_k\subset I_J$.
Since $\hat{F}^{k+1}(\U^m_{J})$ is contained in some 
$\U'_{J'}$ a component of the domain of the first 
return mapping to $\V$,
we can conclude that $\U^m_{J}$ is well-inside $\V$ when $k$ is large.

Now we consider the case when $k$ is not large. Then $m-k$ is large. 
Then since the diameters of puzzle pieces tend to zero, each component of
$\hat{F}^{-(m-k-1)}(\V)$ (that is, each puzzle piece) has very small diameter
compared to the component $I$ of $I_1$ that intersects 
$\hat{F}^{-(m-k-1)}(\V).$
So we can take an interval $K$ with
$|K|/|I|$ small and a $\theta'\in(0,\pi)$
that is independent of the choice of $J$,
so that 
$\mathrm{Comp}_{\V}\hat{F}^{-(m-k-1)}(\V)\subset D_{\theta'}(K)$,
but then since $J_k$ is well-inside $I$, we also get that
$\U^m_{J}=\mathrm{Comp}_J\hat{F}^{-m}(\V)$
is well-inside $\V$.

We now consider the remaining case that 
$\hat{F}^m(J)\subset I_0$. 
Notice that $I_0\cap \Dom(R_{\hat I_1})=\emptyset$.
Let $k\leq m$ be minimal so that
$\hat{F}^k(J)\subset I_0$. 
Then since $\hat{F}^k|_J$ is a diffeomorphism 
(as critical points in $\Omega_1$ are 
minimal), and components of 
$\V$ with real trace, $I_{0,i}$
are of the form $D_{\pi-\theta_0}(I_{0,i}),$
we can take $\theta_0$ so small that
$\mathrm{Comp}_J\hat{F}^{-k}(\V)$ is well-inside $\V.$
This completes the proof of the claim.
\endpfclaim

					\begin{figure}[htp] \hfil
					\beginpicture
					\dimen0=0.2cm
					\setcoordinatesystem units <\dimen0,\dimen0>
					\setplotarea x from -9 to 30, y from -8 to 8
					\setlinear
					\plot -8 0 26 0 /
					\circulararc 120 degrees from 5 0  center at  0 -3 
					\circulararc -120 degrees from 5 0  center at  0 3 
					\circulararc 300 degrees from 6 0 center at 7 -2
					\circulararc -300 degrees from 6 0 center at 7 2
					\endpicture
					\caption{Controlling the geometry of $F_{it}$, see the proof
					of Theorem~\ref{thm:box mapping persistent infinite branches}. \label{fig:gapconditionl}}
					\end{figure}
					
Taking $\bm{\mathcal{V}}_{it}=\V$, chosen as above.
For each component $J$ of $R^{-m}_{\hat{I}_1}(I_0\cup I_1),$
let $\V_J$ be the component of $\bm{\mathcal{V}}_{it}$ that contains
$R_{\hat{I}_1}^m(J)$, and let 
$\U_J=\comp_J R_{\hat{I}_1}^{-m}(\V_J)$.
Let $\bm{\mathcal{U}}_{it}=\cup \U_J$ where the 
union is taken over all components $J$ of the domain of 
$R_{\hat{I}_1}^m$.
Define the mapping $$F_{it}\colon \bm{\mathcal{U}}_{it}
\rightarrow \bm{\mathcal{V}}_{it}$$
by $F_{it}(\U_J)=R^m_{\hat{I}_1}(\U_J)$.
Claim 1 implies that if $m$ is chosen sufficiently large
that $F_{it}$ is a complex box mapping.

\medskip
\noindent\textit{Claim 2:} 
The gap and extension conditions hold for the map
$$F_{it}\colon  \bm{\mathcal{U}}_{it}\rightarrow \bm{\mathcal{V}}_{it}.$$
Moreover, any component $\bm D$
of either $\bm{\mathcal{U}}_{it}$ or $\bm{\mathcal{V}}_{it}$ has $\delta$-bounded
geometry.
				
\medskip
\noindent\textit{Proof of Claim 2:}
That $F_{it}$ satisfies the extension condition
follows from Claim 1 and the
real bounds at persistently recurrent critical points
and at infinitely renormalizable critical points.
To see that gap condition holds. Suppose that
$\U_1$ and $\U_2$ are two connected components of
$\bm{\mathcal{U}}_{it}$. If $F_{it}|_{\U_1}$ and 
$F_{it}|_{\U_2}$ are both diffeomorphisms, then 
$\U_{1}$ and $\U_{2}$ are both contained 
in Poincar\'e lens domains and the gap condition 
follows from the real bounds. Otherwise
let $\{\U^1_{j}\}_{j=0}^{s_1}$ and 
$\{\U^2_{j}\}_{j=0}^{s_2}$
be the chains so that for $i=1,2$ we have that 
$\U^i_0=\U_i$ and $\U^i_{s_i}$ is a component of
$\bm{\mathcal{V}}_{it}$. We can assume that
$s_1\leq s_2$. Then $F_{it}^{s_1}(\U^1_0)$ is a component $\V$ of
$\bm{\mathcal{V}}_{it}$ and $F_{it}^{s_1}(\U_2)$ is
disjoint from $\V$. Evenmore, by the real bounds,
there exists $\theta\in (0,\pi/2)$ so that
$F_{it}^{s_1}(\U_2)\subset D_{\theta}(\tilde U_2),$
where $U_2$ is well-inside $\tilde U_2$.
The gap condition follows.

Propositions~\ref{prop:real bounds persistent}
and
Lemma~\ref{lem:angle control} imply that 
$\bm D$ has $\delta-$bounded geometry. 
\endpfclaim
\end{pf}

\subsection{Compatible complex box mappings and bounds at reluctantly recurrent points}
\label{subsubsec:box mappings reluctant}
In this subsection, we will construct a complex box mapping associated
to the set of reluctantly recurrent critical points of $f$.
Note that this includes the case that $c$ is recurrent and $\omega(c)$ is non-minimal. 

Let us define 
$\Omega_{r}$ as on page~\pageref{notation:Omega_reluct}. As before,
this set can be decomposed into two parts:
\begin{enumerate}
\item
$\Omega_{r,e}$ is a set of critical points so
that if $c\in \Omega_{r,e}$ has odd order, 
then there exists $c'\in \Omega_{r,e}$
with even order so that $c\ge c'$ (here we use the ordering defined in
Subsection~\ref{subsec:decomp}), and
\item $\Omega_{r,o}$.
\end{enumerate}
			For each  $c\in \Omega_{r,o}$,  
			and each critical point $c'\in \omega(c)$
			we have that  $c'$  is a critical point of odd order in 
			$\Omega_{r,o}$,  that  $f$ is persistently recurrent on $\omega(c')$,
			that $c'$ is in the basin of a periodic attractor
			or that $c'$ is eventually mapped into a renormalization interval
			(where the mapping is infinitely renormalizable).
			
\begin{thm}[Existence of complex box mappings in the reluctantly
  recurrent case] 
\label{thm:box mapping reluctant}
 Let $f\colon  [0,1]\to \R$ or $f\colon  S^1\to S^1$ with periodic orbits
be in $\mathcal{C}$.
  Let $\Omega$ be either $\Omega_{r,e}$  or  $\Omega_{r,o}$,
  as defined above.
Then, for any $\theta\in(0,\pi/2)$
there exists  arbitrarily small real  
  neighbourhoods $I$ of $\Omega$, such that the first return map
  $R_I\colon \Dom_{\Omega}^*(I)\to I$ extends to a complex box mapping
$F\colon  \UU\to \VV$ with the following properties:
  \begin{itemize}
  \item For each $c\in\Omega$ the set $\VV(c)$ is contained in $D_{\pi-\theta}(I(c))$.
  	\item Let $\VV(c)^\pm=\VV(c)\cap \mathbb{H}^{\pm}$.
 	 Then $\VV(c)^+\setminus \UU$ and $\VV(c)^-\setminus \UU$ are quasidiscs.
  \item This complex box mapping satisfies the extension and gap conditions
	defined on page~\pageref{def:gapextension}.
  \end{itemize}
  \end{thm}

The proof of this theorem closely follows the proof of part II of Theorem 3 in
\cite[pages 150-159]{KSS-density},
and therefore we only discuss what changes
we have to make in that proof
(for the terminology used in the proof we also refer to \cite{KSS-density}).
 The main difference is that we allow for the existence of
critical points of odd order, but we also
have to check that the presence of parabolic points does not cause difficulties. 

Note that the last two assertions in
Theorem~\ref{thm:box mapping reluctant} are new. They are
needed because we will want to argue as in the proof of Theorem A' and
B'
of \cite{LevStr:invent}.

This theorem follows immediately from the following proposition
from \cite{KSS-density}, which holds for critical points of all orders 
without any essential changes to the proof, and
the real bounds for reluctantly recurrent maps,
Proposition~\ref{prop:real bounds reluctant}.
The construction given \cite{KSS-density} was done
in the non-minimal case; however,
it depends only on certain real geometric conditions,
which hold in the reluctantly recurrent case as well.
\begin{prop}[\cite{KSS-density} Propostion 1]
\label{prop: big space box mapping}
There exists $\varepsilon_{0}>0$, $C_{0}>0$ and $\theta\in(0,\pi/2,)$
with the following properties. Let $I$ be an admissible neighbourhood of
$\Omega$ such that $\cen(I)<\varepsilon_{0},\Space(I)>C_{0}$ and $\gap(I)>C_{0}$.
Additionally, assume that $\max_{c\in\Omega}|I(c)|$ is sufficiently small.
Then there exists a real-symmetric complex box mapping
$F\colon \UU\rightarrow \VV$ whose real trace is a real box mapping $R_{I}$
and for each critical point $c$, $\VV(c)$ is contained in $D_{\pi-\theta}(I(c))$. 
\end{prop}
			
\medskip
\noindent\textit{Proof of Theorem~\ref{thm:box mapping reluctant}.}
Let us make a few remarks.
Throughout, we should take for $\Omega$ in the proof in \cite{KSS-density},
the set  $ \Omega_{r,e}$ or the set $\Omega_{r,o}$.
It is enough to show that one can construct a box mapping $F\colon \UU\to \VV$
extending  $R_I\colon \Dom_{\Omega}^*(I)\to I$, because this can be easily extended
to a box mapping extending
$R_{\hat I}\colon \Dom_{\Omega}^*(\hat I)\to \hat I$: just define
for each critical point $c\in \Omega_{r,e}$ of odd order 
the sets $\UU(c)=\VV(c)=\LL_c(\VV)$.
In the statement of  Theorem~\ref{thm:box mapping reluctant}, 
$I$ should then be taken as $\hat I$.

%
%
				
The only place where an assumption $\cen_2(I)$  is needed is in Lemma 2 (through 
the first part of the proof of Proposition 1 in \cite{KSS-density});
there the argument goes through (even more easily,
and without using any assumption on $\cen_2(I)$)
each time the pullback passes through a critical point 
of odd order. We need the assumption that $\cen_2(I)$ is small
for critical points of even order to ensure that when we pull back to through an
even critical point, the size of the domain on the
real line is not very small compared to the 
diameter of the part that pulls back to the complex plane.
Note that the pullback is disjoint,
so each critical point is only passed through at most once.
The assumption that $\cen_1(I)$ is small
(which means  that the critical value $f(c)$ is roughly 
in the middle of $\hat J(c)$) is used in the proof of
Proposition 1 of \cite{KSS-density}
is at the bottom of page 153
(only in the case of critical points of even order).
So that we can choose for the range $\VV(c)$ 
a lens-shaped set of the form $\VV(c)=D_{\pi-\theta}(I)$,
where $\theta\in(0,\pi/2)$ can be chosen small, but bounded away from zero.
As in the statement of Theorem 3 in \cite{KSS-density}
we get control of the shape of each domain. In fact, because of
(\ref{eq:extension}) one has the following extension property.
There exists $\delta>0$ and, for each component $\V_j$ of $\VV,$
there exists a neighbourhood  $\hat{\V}_j\supset \V_j,$ 
with $\mod(\hat{\V}_j\setminus \V_j)\ge
\delta$ 
such that the following hold:
\begin{enumerate}
\item If $\U_i$  is a component of $\UU$ and
$F|_{\U_i}=f^{s_i}$
so that $\U_i$ maps onto $\V_j$, 
then there exists a neighbourhood $\hat{\U}_i$ of $\U_i$
so that $f^{s_i}\colon  \hat{\U}_i\to \hat{\V}_j$ has the same degree as 
$f^{s_i}\colon \U_i\to \V_j$.
\item If $\U_i$ is compactly contained in $\VV$, then $\hat{\U}_i$ is contained in $\VV$.
\end{enumerate}
Here one can take for $\hat{\V}_j$ a Poincar\'e disc based on the set $\hat T_{N-1}$
defined in (\ref{eq:extension}) (where $\hat T_{N-1}$ is the 
set associated to $T_N$ in (\ref{eq:evenodd}).
Because $T_N$ is a return domain to $T_{N-1}$,
and all critical points are contained in $T_N$ the above property follows.
From this we get that $\comp_{\U_i}f^{-s_i}(\hat{\V}_j)$ is a domain
with real trace
$\hat U$, and we can arrange it so that the above properties hold. 
All this implies in particular that each $\U_i$ has a smooth boundary without cusps
and with controlled shape (they have bounded geometry and are 'roughly round').

To show that $\bm A^+:=\VV(c)^+\setminus \UU$ is a quasidisc we
use the Ahlfors-Buerling Criterion in a similar way as was done
on page 440 of \cite{LevStr:invent}  
This criterion says that $\partial \bm A^+$ is a quasicircle provided there exists
$L<\infty$ with the following property.
For any two points $z_1,z_2$ in  $\partial(\bm A^+)$,
take one of the two arc $C$ in $\partial  \bm A^+$ connecting $z_1$ and $z_2$ 
is contained in a disc of diameter $L\cdot |z_1-z_2|$. If $z_1,z_2$ are in 
one component of $\UU$ this follows from the control of $\U_i$ shown above.
If $z_i$ and $z_2$ are in two distinct domains $\U_i$ and $\U_j$ then
this follows 
from the fact that $z_1,z_2$ are 'roughly round' and from the gap property
of the admissible neighbourhood $I$:   the  gap between two domains of $R_I$ 
in one component of $I$ is large 
(i.e. the gap between two domains is much large  than one of the two domains).

The last item in Theorem~\ref{thm:box mapping reluctant}
follows then from Lemma 3(3) in \cite{KSS-density}.
\qed

			\begin{figure}[htp] \hfil
			\beginpicture
			\dimen0=0.3cm
			\setcoordinatesystem units <\dimen0,\dimen0>
			\setplotarea x from -9 to 30, y from -5 to 4
			\setlinear
			\plot -12 0 26 0 /
			\circulararc 120 degrees from 0 0  center at  -5 -3 
			\circulararc -120 degrees from 0 0  center at  -5 3 
			\put {$V$} at 5 2 
			\put {$*$} at 7 0
			\put {$*$} at 20 0
			
			\circulararc 50 degrees from -3 0  center at  -4 -2 
			\circulararc -50 degrees from  -3 0 center at  -4 2 
			\circulararc -90 degrees from -9 0  center at  -7.5 -1.5 
			\circulararc 90 degrees from  -9 0 center at  -7.5 1.5 
			
			\circulararc 75 degrees from 10 0  center at  7 -4 
			\circulararc -75 degrees from 10 0  center at  7 4 
			
			\arrow <3mm> [0.2,0.67] from 9 1 to 13 3

			\circulararc 130 degrees from 25 0  center at  18.5 -3 
			\circulararc -130 degrees from 25 0  center at  18.5 3 
			
			\circulararc -90 degrees from 14 0  center at  15.5 -1.5 
			\circulararc 90 degrees from  14 0 center at 15.5 1.5
			 \arrow <3mm> [0.2,0.67] from 15.5 1 to 17 3
			\circulararc 125 degrees from 22 0  center at  20 -1 
			\circulararc -125 degrees from  22 0 center at 20 1 
			\arrow <3mm> [0.2,0.67] from 18.5 -1 to 8.7 -1
			\circulararc 50 degrees from 24.5 0  center at  23.5 -2 
			\circulararc -50 degrees from  24.5 0 center at  23.5 2 
			 \arrow <3mm> [0.2,0.67] from 23.5 0.5 to 22 2.5
			\endpicture
			\caption{The box mapping associated to the non-minimal case has infinitely many components
			of its domain. \label{fig:imal}}
			\end{figure}

\subsection{Complex bounds for central cascades}
\label{subsec:improving}
In this section, we will prove complex bounds for certain
puzzle pieces near the top of central cascades.
These results are only used in the proof
quasisymmetric rigidity
for smooth mappings; however, they are new even for polynomials,
and are of independent interest.
Suppose that $F\colon\bm{\mathcal{U}}\rightarrow\bm{\mathcal{V}}$ is a real-symmetric
complex box mapping given by Theorem~\ref{thm:box mapping persistent},
Theorem~\ref{thm:box mapping persistent infinite branches} or Theorem~\ref{thm:box mapping reluctant}
with $\diam(\bm{\mathcal{V}})$ sufficiently small.

We will use a
nest of puzzle pieces whose construction
is closely related to the construction of the
enhanced nest.
Recall that if $F$ is renormalizable, and the 
first renormalization $F^s\colon J\rightarrow J$ of $F$ is not of intersection type, 
then the enhanced nest construction gives a terminating puzzle piece
$\V'$ for $F$ such that the first return mapping
$F\colon\U'\rightarrow\V'$ is a polynomial-like mapping, 
$\mod(\V'\setminus\U')$ is bounded away from zero and 
$\V'$ has bounded geometry. We are going to modify the
construction of the enhanced nest
to produce a puzzle piece that has the same properties as 
$\V'$ near the top of a sufficiently long
central cascade.
We can obtain complex bounds for this nest using the 
initial bounds for the box mapping
$F\colon\bm{\mathcal{U}}\rightarrow\bm{\mathcal{V}}$
and the method of \cite{CvST} or using the argument of
\cite{KvS}. We will explain how to use the argument of
\cite{KvS}.

We will say that a puzzle piece is 
$\delta$-\emph{excellent}
\label{page:defgood} if it is $\delta$-nice,
has $\delta$-bounded geometry
and is a $(1+1/\delta)$-quasidisk for some $\delta>0$.
Note that we exclude $\delta$-free from this definition.
From the complex bounds of \cite{CvST}, we have that
every puzzle piece in the enhanced nest is $\delta$-nice and 
has $\delta$-bounded geometry. We will prove in this section
that there exists $\delta>0$, 
depending on the 
geometry of $\bm{\mathcal{V}}$, such that
each puzzle piece in the enhanced nest is a
$(1+1/\delta)$-quasidisk.
\begin{rem}\label{rem:good}
Observe that for any $\delta>0,$ $d\in\mathbb N$ and $K\geq 1$,
there exists $\delta'>0$ such that the following holds.
Suppose that 
$\bm I$ is a $\delta$-excellent puzzle piece, and that
$\bm I'$ is a pullback of $\bm I$ with
degree bounded by $d$ and dilatation bounded by $K$.
Assume there exist sets $\bm I^+,\bm I^-$ so that if we set
$\bm A^+=\II^+\setminus \II$ and
$\bm A^-=\II\setminus \II^-$ we have that
$\mod(\bm A^{+}),\mod(\bm A^-)>\delta$ 
and that the mapping from $\bm I'$ to
$\bm I$ extends to a mapping onto $\II^+$
with no critical values in $\bm A^+\cup\bm A^{-}.$
Then $\bm I'$ is $\delta'$-excellent.
\end{rem}

Let us suppose now that the 
principal nest starting with $\Z^0\owns c_0:$
$$\bm Z^0\supset\bm Z^1\supset\bm Z^2\supset\dots
\supset\bm Z^N\supset\bm  Z^{N+1}$$
is a central cascade; 
that is, $\bm Z^1=\mathcal{L}_{c_0}(\bm Z^0)$ and
we have that every critical value of $R_{\bm Z_0}|_{\bm Z^1}$
is contained in $\bm Z^N$, and there exists a critical point, $c$
of $R_{\bm Z_0}|_{\bm Z^1}$ such that 
$$R_{\bm Z_0}|_{\bm Z^1}(c)\in \bm Z^N\setminus \bm Z^{N+1}.$$
Consider the principal nest starting with a puzzle piece $\Z^{-1}\owns
c_0:$
$$\Z^{-1}\supset \bm Z^0\supset\bm Z^1\supset\bm Z^2\supset\dots
\supset\bm Z^N\supset\bm  Z^{N+1},$$
where $\bm Z^0\supset\bm Z^1\supset\bm Z^2\supset\dots
\supset\bm Z^N\supset\bm  Z^{N+1}$ is a central cascade.
If
the return time of $\Z^0$ to $\Z^{-1}$
is strictly less than the return time of 
$\Z^1$ to $\Z^{0},$
then we say that the central cascade 
$\bm Z^0\supset\bm Z^1\supset\bm Z^2\supset\dots
\supset\bm Z^N\supset\bm  Z^{N+1}$
is \emph{maximal}.
Assume that $\bm Z^0\supset\bm Z^1\supset\bm Z^2\supset\dots
\supset\bm Z^N\supset\bm  Z^{N+1}$ is a maximal central cascade.
Let $s$ be the return time of $\bm Z^1$ to $\bm Z^0$ under $F$ and
suppose that $N\geq 3$. 
In these circumstances, for $N\in\mathbb N$,
we say that a puzzle piece $\bm J$ is
\emph{combinatorially close to }$\bm Z^0$
if $\bm J\supset \Z^2 \owns c_0$ and
the return mapping from $\Z^3$ to 
$\Z^{2}$ extends to a mapping of the same degree onto
$\bm J$.

\subsubsection{Modifying the construction of the enhanced nest
to produce puzzle pieces containing a central cascade}
For each critical point $c\in\Crit(F),$
let $s_c$ be the first entry time of $c$ to $\bm Z^0.$ 
We define 
$$PC(\bm Z^0)=\bigcup_{c\in\Crit(F)}\bigcup_{j=0}^{s_c}F^j(c).$$ 
Notice that we include the critical point in this set.

Suppose that $\bm I=\II^0\owns c_0$ is a critical puzzle piece.
Let 
$k_{c_0}$ be such that
$R_{\bm I}^{k_{c_0}}(c_0)\in\bm I\setminus \bm I^1$
we define
$$\hat{\mathcal A}(\bm I)=\comp_{c_0}(R_{\bm I}^{-k_{c_0}}\mathcal{L}_{R^{k_{c_0}}_{\bm I}(c_0)}(\bm I)).$$

\begin{lem}\label{lem:short step}
Assume that
$\bm Z^0\supset \bm Z^1\supset \bm Z^2\supset\dots\supset \bm
Z^N\supset \bm Z^{N+1}\owns c_0$
is a maximal central cascade with $N\geq 2$.
Suppose that the return time of $\Z^1$ to $\Z^0$ is
strictly greater than the return time of $\II^1$ to $\II$ and
$$\bm I\supset \bm Z^0\supset \hat{\mathcal A}(\bm I).$$ 
Then there exists a puzzle piece $\bm J$ such that 
$\bm J$ is a pullback of $\bm I$ with order bounded by $2b$ and
$\bm J$ is combinatorially close to $\bm Z^0$.
\end{lem}
\begin{pf}
Suppose first that $\bm Z^0\supset \bm I^1$. Then 
$\bm I^1\supset \bm Z^1\supset \bm I^2\supset \bm Z^2.$
It is immediate that $\bm I^1$
is combinatorially close to $\Z^0$.

Now, assume that $\bm I^1\supset \bm Z^0$.
Suppose first that the return of
$\bm I^1$ to $\bm I$ has the property that
$R_{\bm I}(c_0)\notin\bm I^1.$
Then $R_{\bm I}(\bm Z^1)$ is contained in a
non-central return domain to $\bm I$.
By assumption, the pullback, $\bm J=\hat{\mathcal A}(\bm I)$,
of this domain under $R_{\bm I}$
is contained in $\bm Z^0$.
So we have that
$\bm Z^1\subset \bm J\subset \bm Z^0,$
and hence is combinatorially close to $\bm Z^0$.

Now suppose that the return of
$\bm I^1$ to $\bm I$ is central for $c_0$.
Let 
$$\bm I^1\supset\bm I^2\supset \bm I^3\supset\dots$$
be the principal nest. 
Let $k>0$ be minimal so that
$R_{\bm I}(\bm Z^1)\subset \bm I^k\setminus \bm I^{k+1}.$
Then $R^{k+1}_{\bm I}(\bm Z^1)$ is contained in a non-central component of the 
domain of the return mapping to $\bm I.$ 
Setting
$$\bm J=\hat{\mathcal A}(\bm I)\subset \bm Z^0$$
we have that $\bm Z^1\subset \bm J\subset \bm Z^0,$
which again finishes the proof.
\end{pf}

\begin{prop}[Compare \cite{KSS-rigidity}, Lemma 8.2.]
\label{lem:nu}
Assume that $\Z^0\supset \Z^1\supset\dots\supset \Z^N\supset \Z^{N+1}$ is
a maximal central cascade with $N\geq 8.$
Suppose that $\II\supset \Z^0\owns c_0$ is a pair of puzzle pieces such that
$\hat{\mathcal{A}}(\II)\supset \Z^0$.
Then there exists a sequence of combinatorially defined puzzle pieces
$$\II=\II_0^{\Z^0}\supset \II_1^{\Z^0}\supset\dots\supset\II_n^{\Z^0}\supset\Z^2$$
such that the following holds: 
\begin{itemize}
\item Each $\II_{i+1}^{\Z^0}$ is a pullback of $\II_i^{\Z^0}$ by $F^{p_i}$ of
  bounded degree.
\item There exist sets $\II^+_i,\II_i^-$ such that 
$\II_i^-\subset \II_i^{\Z^0}\subset \II_i^+$,
and the critical values of $F^{p_i}|_{\II_{i+1}^{\Z^0}}$ are contained in
$\II_{i}^-$.
\item $\Z^2\subset\II_n^{\Z^0}\subset \Z^0$.
\end{itemize}
\end{prop}

\begin{rem}
In the course of the proof of this
proposition, we will define the
\emph{enhanced nest above a (maximal) central cascade}.
\end{rem}

\begin{pf}
\medskip
\noindent\textit{Step 1.}
We are going to inductively construct a sequence of pairs of 
sets of puzzle pieces: 
$\bm T_n,\bm J_n.$ 
This construction terminates either 
\begin{enumerate}
\item[(1)] when we obtain a 
puzzle piece $\bm U_0$ that is combinatorially close to $\Z^0$
or
\item[(2)]  when we find a collection of critical puzzle pieces
$\bm P_c\Supset\bm P'_c$ that contains $\Crit(F)$
such that each $\bm P_c$ is a pullback of $\bm I$
with order bounded by $b^2-b$, and for each $c\in \Crit(F)$
and any
$z\in PC(\Z^0)\cap(\bm P_c\setminus \bm P_c')$ there exist
a critical point $\hat c$ and a puzzle piece $\bm Z\owns z$
such that the landing mapping from $\bm Z$ to $\bm P_{\hat c}$
is a diffeomorphism. 
\end{enumerate}

Let $\bm T_0=\bm I$ and $\bm J_0=\mathcal L_{c_0}(\bm T_0).$
If the return time of $\Z^1$ to $\bm I$ is 
the same as its return time to $\bm Z^0,$ then we
set $\bm U_0=\mathcal{L}_{c_0}(\bm I)$, and we are done.
So suppose that
$\bm I$ is not combinatorially close to $\bm Z^0.$
If for any critical point 
$c \in\Crit(F)$, $c\neq c_0$
we have that $R_{\bm I}(c)\in\bm J_0,$
set $\bm P_c=\hat{\mathcal L}_c(\bm T_0)$
and $\bm P_c'=\hat{\mathcal L}_c(\bm J_0).$
Otherwise, there exists a critical point $c_1\in\Crit(F)$
such that $R_{\bm I}(c_1)\in\bm T_0\setminus \bm J_0.$
Let $\bm T_{1}=\bm J_0\cup
\comp_{c_1}R_{\bm T_0}^{-1}(\mathcal{L}_{R_{\bm T_0}(c_1)}(\bm T_0)),$
and let 
$\bm J_{1}=
\mathcal L_{c_0}(\bm T_{1})\cup\mathcal{L}_{c_1}(\bm T_{1}).$

Suppose first that
$\bm Z^0\supset \bm J_1(c_0).$
Observe the either 
\begin{itemize}
\item $\bm J_1(c_0)=\mathcal L_{c_0}(\bm J_0(c_0))$
or 
\item $\bm J_1(c_0)=\mathcal L_{c_0}(\comp_{c_1}R_{\bm
  T_0}^{-1}(\mathcal{L}_{R_{\bm T_0}(c_1)}(\bm T_0))),$
\end{itemize}
In the first case we have that 
$$\bm J_0\supset \bm Z^0\supset \bm
J_1(c_0)\supset\mathcal L_{c_0}(\bm Z^0),$$
and $\bm J_1(c_0)$ is combinatorially close to $\Z^0$.
In the second case, since $\bm J_1(c_0)$ is a first entry
domain to $\bm T_1(c_0)$ the orbit of $c_0$ enters
$\comp_{c_1}R_{\bm T_0}^{-1}(\mathcal{L}_{R_{\bm T_0}(c_1)}(\bm T_0))$
before returning to $\bm T_1(c_0)=\bm J_0(c_0).$
So the orbit of $\mathcal L_{c_0}(\bm Z^0)$ enters
$\comp_{c_1}R_{\bm T_0}^{-1}(\mathcal{L}_{R_{\bm T_0}(c_1)}(\bm T_0))$
before returning to $\bm Z^0$, pulling this back to $c_0$,
we have that
$$\mathcal L_{c_0}(\bm Z^0)\subset  \mathcal L_{c_0}(\comp_{c_1}R_{\bm
  T_0}^{-1}(\mathcal{L}_{R_{\bm T_0}(c)}(\bm T_0)))\subset
\bm Z^0$$
and we have that
$\mathcal L_{c_0}(\comp_{c_1}R_{\bm
  T_0}^{-1}(\mathcal{L}_{R_{\bm T_0}(c)}(\bm T_0)))$ is
combinatorially close to $\Z^0$.

So let us suppose that
$\bm J_1\supset\bm Z^0.$
If for any critical point 
$c\in\Crit(F)$, $c\neq c_0, c_1$
we have that $R_{\bm T_1}(c)\in\bm J_1,$
set $\bm P_c=\hat{\mathcal L}_c(\bm T_1)$
and $\bm P_c'=\hat{\mathcal L}_c(\bm J_1).$
Otherwise, there exists a critical point 
$c_2\in\Crit(F)$ such that
$R_{\bm T_1}(c_2)\in\bm T_1\setminus \bm J_1.$
Let $\bm T_{2}=\bm J_1\cup
\comp_{c_2}R_{\bm T_1}^{-1}(\mathcal{L}_{R_{\bm
    T_1}(c_2)}(\bm T_1)),$ 
and let 
$\bm J_{2}=
\mathcal L_{c_0}(\bm T_{2})\cup\mathcal L_{c_1}(\bm T_{2})
\cup \mathcal{L}_{c_2}(\bm T_{2}).$

Suppose first that
$\bm Z^0\supset \bm J_2(c_0).$
There are three possibilities:
\begin{itemize}
\item[(a)] $\bm J_2(c_0)=\mathcal L_{c_0}(\bm T_2(c_0))$ and
  the orbit of $c_0$ enters $\bm T_2(c_0)$ before entering $\bm
  T_2(c_1)\cup\bm T_2(c_2)$,
\item[(b)] $\bm J_2(c_0)=\mathcal L_{c_0}(\bm T_2(c_1))$  and
  the orbit of $c_0$ enters $\bm T_2(c_1)$ before entering $\bm
  T_2(c_0)\cup\bm T_2(c_2)$  or
\item[(c)] $\bm J_2(c_0)=\mathcal L_{c_0}(\bm T_2(c_2))$ and
  the orbit of $c_0$ enters $\bm T_2(c_2)$ before entering $\bm
  T_2(c_0)\cup\bm T_2(c_1)$,
\end{itemize}
In the Case (a), we immediately have that
$$\bm T_2(c_0)\supset \bm Z^0\supset \bm J_2(c_0)=\mathcal
L_{c_0}(\bm T_2(c_0))\supset
\mathcal L_{c_0}(\bm Z^0).$$

In the Case (b), we have that
$$\bm T_2(c_1)=\bm J_1(c_1)=\mathcal L_{c_1}(\bm
T_1)=\mathcal L_{c_1}(\bm J_0\cup\comp_{c_1}R^{-1}_{\bm T_0}(\mathcal
L_{R_{\bm T}(c_1)}(\bm T_0))).$$
So there are two possibilities, either
$\bm J_2(c_0)=\mathcal L_{c_0}\mathcal L_{c_1}(\bm J_0)$,
or
$$\bm J_2(c_0)=\mathcal L_{c_0}\mathcal L_{c_1}(\comp_{c_1}R^{-1}_{\bm T_0}(\mathcal
L_{R_{\bm T_0}(c_1)}(\bm T_0))).$$
In the second case, it is important to observe that 
after entering $\mathcal L_{c_1}(\comp_{c_1}R^{-1}_{\bm T_0}(\mathcal
L_{R_{\bm T_0}(c_1)}(\bm T_0))),$
the orbit of $\mathcal L_{c_0}(\bm Z^0)$
returns to $\comp_{c_1}R^{-1}_{\bm T_0}(\mathcal
L_{R_{\bm T}(c_1)}(\bm T_0))$
before returning to $\bm J_0$,
so we have that 
$$\bm T_2(c_0)\supset \bm Z^0\supset \bm J_2(c_0)\supset
\mathcal L_{c_0}(\bm Z^0).$$

In Case (c), 
$$\bm J_{2}(c_0)=\mathcal L_{c_0}(\comp_{c_2}R_{\bm
  T_1}^{-1}(\mathcal{L}_{R_{\bm T_1}(c_2)}(\bm T_1))),$$
and the orbit of $\mathcal L_{c_0}(\bm Z^0)$ enters
$\comp_{c_2}R_{\bm
  T_1}^{-1}(\mathcal{L}_{R_{\bm T_1}(c_2)}(\bm T_1))$
before returning to $\bm T_2(c_0)$.
So again
$$\bm T_2(c_0)\supset \bm Z^0\supset \bm J_2(c_0)\supset
\mathcal L_{c_0}(\bm Z^0),$$
and we are done.

So we may assume that $\bm Z^0\subset \bm J_2(c_0)$.
Continuing in this fashion, 
we find that either
\begin{enumerate}
\item[A.] There exists $\bm U_0=\bm
  T_m$ that is combinatorially close to $\Z^0$ and is a pullback of
  $\II$ of order bounded by $b^2-b$ or
\item[B.]
\begin{enumerate}
\item there is a collection $\bm T_m, m<b$
of puzzle pieces around some critical points, and this collection is strictly nice;
\item for any other critical point $c'\in\Crit(F),$ $R_{\bm T_m}(c')\in\bm J_m;$
\item any puzzle piece of $\bm T_m$ is a pullback of $\bm I$
  with order bounded by $b^2-b$.
\end{enumerate}
\end{enumerate}
For each $c\in\crit_{\bm Z^0}(F),$ let $\bm P_c=\hat{\mathcal{L}}_c(\bm T_m)$ and
$\bm P_c'=\hat{\mathcal{L}}_c(\bm J_m).$
By the construction, we immediately see that either
(1) or (2) holds.

It will be useful later to make the following definition.
\begin{defn}
Given a pair of critical puzzle pieces $\bm I\supset\bm
Z^0\owns c_0$ we define
$$\mathcal T(\bm I)=\bm T_m(c_0).$$
\end{defn}

\medskip
\noindent\textit{Step 2.}
Let $k$ be the maximal integer such that for some $c\in\crit_{\Z^0}(F)$ the set
$F^{-k}(\bm P_c)$ 
has a component $\bm Q$
such that
\begin{itemize}
\item $\bm Q$ is a child of $\bm P_c$; that is, $\bm Q$ is
  a unicritical pullback of $\bm P_c$ that contains a critical point, and
\item $\hat{\mathcal{L}}_{c_0}(\bm Q)\supset\bm Z^0$.
\end{itemize}
Let $\bm U_0=\hat{\mathcal{L}}_{c_0}(\bm Q),$
and let $\nu$ be the positive integer such that $F^{\nu}(\bm
U_0)=\bm I$.
Let $k_1\geq 0$ be so that $F^{k_1}(\bm U_0)=\bm Q$,
and
let $l\geq 0$ be the integer such that $\bm P_c$ is a pullback of 
$\bm I$ by $F^l$.
Observe that since $\mathcal{L}_{F^{\nu}(c_0)}(\bm I)$ is a first return
domain to $\bm I$ and $F^l(\bm P'_c)\cap
\mathcal{L}_{F^{\nu}(c_0)}(\bm I)\neq\emptyset$, we have that
$F^l(\bm P_{c}')\subset \mathcal{L}_{F^{\nu}(c_0)}(\bm I),$
which immediately implies
$\bm P_c'\subset\comp_{F^{k+k_1}(c_0)}F^{-l}(\mathcal{L}_{F^{\nu}(c_0)}(\bm I)).$
Thus to show that a point 
$x\in \comp_{c_0} F^{-\nu}(\mathcal L_{F^\nu(c_0)}(\bm I))$
it is sufficient to show that
$F^{k_1+k}(x)\in \bm P'_c.$

Let $\{\bm G_j\}_{j=0}^s$  be the chain with
$\Z^0=\bm G_s$ to $\Z^1=\bm G_0$. 
Suppose first that $x\in\mathrm{PC}(\Z^0)$,
and that $x\in \bm G_j$, for some $0\leq j<s.$
Suppose that $F^k(x)\in\bm P_c\setminus \bm P'_c$.
Then by construction, there exists a puzzle
piece $\bm W\owns F^k(x)$ and a critical point $\hat c$,
so that $\bm W$ is mapped diffeomorphically
onto $\bm P_{\hat c}$. Since $F^k\colon\bm Q\rightarrow\bm P_{c}$
is unicritical, 
we have that $\bm G_{k+j}\subset \bm W\subset \bm P_c\setminus \bm
P_c'$. This contradicts the maximality of $k$.
If $\hat{\mathcal{A}}(\U_0)\subset\Z^0$, then 
by Lemma~\ref{lem:short step}, we have that there is a
pullback of $\U_0$ with order bounded by $b$ that is 
combinatorially close to $\Z^0$. So we can suppose that
$\hat{\mathcal{A}}(\U_0)\supset\Z^0.$
Let $\{\bm{H}_{j}\}_{j=0}^r$ denote the chain given by
$\bm{H}_r=\bm{U}_0$ and
$\bm{H}_0=\mathcal{L}_{c_0}(\bm{U_0})$, and let
$\{\hat{\bm{H}}_{j}\}_{j=0}^{\hat{r}}$
be the chain with
$\hat{\bm{H}}_{\hat r}=\U_0$ and 
$\hat{\bm{H}}_0=\hat{\mathcal{A}}_{c_0}(\U_0).$
Let us show that neither of these two chains intersects $x.$

Suppose that $x\in\mathrm{PC}(\Z^0)\cap \bm Q$ and
$F^k(x)\in\bm P_c\setminus \bm P'_c.$
Then there exists a puzzle piece $\bm W\owns F^k(x)$
and a critical point $\hat c$ so that
$\bm W$ is mapped diffeomorphically onto $\bm P_{\hat c}$.
If for some $j,0\leq j\leq r$, 
$\bm{H}_j\owns x$,
then $\bm{H}_{j+k}$ contains $f^k(x)\in\bm{P}_c\setminus\bm
P_{c}'$ and $\bm G_{j'}$ for some $0\leq j'\leq s$.
Since $\bm G_{j'}\subset \bm P_c'$ and puzzle pieces are nested or disjoint,
we have that $\bm{H}_{j+k}\supset \bm P_c'$.
But this is impossible since $\hat{\mathcal{A}}_{c_0}(\U_0)\supset\Z^0.$
Arguing identically, we have that
the chain $\{\hat{\bm{H}}_j\}_{j=0}^{\hat r}$
does not intersect $x$.

\begin{rem} Each critical puzzle piece $\II$ 
such that $\mathcal{L}_{c_0}(\II)\supset \Z^0$
has a smallest successor
that contains $\bm Z^0$, whether $c_0$ is reluctantly or persistently recurrent.
\end{rem}

Whenever $\U_0$
satisfies 
$\hat{\mathcal{A}}(\U_0)\supset\Z^0$,
we
let $\Gamma_{\Z^0}(\U_0)$ 
denote the successor 
of $\U_0$ with the following properties.
Let
$\{\bm H_j'\}_{j=0}^{r'}$ be the chain with
$\bm H'_{r'}=\U_0$, $\bm H'_0=\Gamma_{\Z^0}(\U_0)$.
Assume that 
\begin{itemize}
\label{page:gammaprops}
\item[(1)]for all $j,$ $\bm H'_j$ is disjoint from the set of all $x\in\bm U_0\cap\PC(\Z^0)$ 
with $f^{k+k_1}(x)\in\bm P_c\setminus\bm P'_c$,
\item[(2)]$\bm H'_0\supset\bm Z^1$,
\item[(3)] $r'$ is maximal over all chains satisfying (1) and (2).
\end{itemize}
Repeating this,
we have that either $\hat{\mathcal{A}}(\Gamma_{\Z^0}(\U_0))\subset\Z^0$. or
or $\Gamma_{\Z^0}(\U_0)$ has two successors with properties (1) and (2), and in
this case, we set $\Gamma^2_{\Z^0}(\U_0)$ to be the last such successor.

Let $\bm I\supset \bm Z^0$ be a puzzle piece.
Suppose that $\hat{\mathcal{A}}\mathcal{T}(\bm I)\supset\bm Z^1$.
	Define 
	$$\mathcal A_{\bm Z^0}(\bm I)=\comp_{c_0} (F^{-\nu}(\mathcal L_{F^{\nu}(c)}(\bm I)))$$
	$$\mathcal B_{\bm Z^0}(\bm I)=\comp_c(F^{-\nu}(\bm I)).$$

\medskip

\begin{defn}
We define the \emph{enhanced above }$\bm Z^0$ as follows:
\label{page:nest over}
\begin{itemize}
\item If $\bm E^{\bm Z^0},$ is defined, then we stop the construction.
\item Let $\bm I^{\bm Z^0}_0=\bm I.$ Suppose that $\bm I^{\bm Z^0}_{n-1}$ has been constructed.
\item 
If $\mathcal T(\bm I^{\bm Z^0}_{n-1})=\mathcal B(\bm I^{\bm Z^0}_{n-1})$,
set $\bm E^{\bm Z^0}=\mathcal T(\bm I^{\bm
  Z^0}_{n-1})$ (this occurs when we are in case (1) of Lemma~\ref{lem:nu}).
\item Otherwise, we can assume that (2) holds. If
$\mathcal T\mathcal A_{\bm Z^0}(\bm I^{\bm Z^0}_{n-1})
= \mathcal{B}_{\bm Z^0}\mathcal{A}_{\bm Z^0}(\bm I^{\bm Z^0}_{n-1})$,
set $\bm E^{\bm Z^0}=\bm I^{\bm Z^0}_n=\mathcal T \mathcal A(\bm I^{\bm Z^0}_{n-1})$. 
From now on, we can assume that
$$\bm Z^0\subset\mathcal B_{\bm Z^0}\mathcal A_{\bm Z^0}(\bm I^{\bm Z^0}_{n-1})
\subset \mathcal T\mathcal A_{\bm Z^0}(\bm I^{\bm Z^0}_{n-1}).$$
\item If $\hat{\mathcal A}\mathcal{B}_{\bm Z^0}\mathcal{A}_{\bm Z^0}(\bm I^{\bm Z^0}_{n-1})
\subset \bm Z^0,$ set
$$\bm E^{\bm Z^0}=\bm I^{\bm Z^0}_n$$
to be the pullback of $\mathcal{B}_{\bm Z^0}\mathcal{A}_{\bm Z^0}(\bm
I^{\bm Z^0}_{n-1})$
that is combinatorially close to $\Z^0$ that is given by
Lemma~\ref{lem:short step}.
\item Otherwise, $\bm Z^0\subset \hat{\mathcal A}\mathcal{B}_{\bm Z^0}\mathcal{A}_{\bm Z^0}(\bm I^{\bm Z^0}_{n-1})$.
If for some $T, 0\leq T<5b$,
$$\hat{\mathcal A}\Gamma_{\bm Z^0}^T\mathcal{B}_{\bm Z^0}\mathcal A_{\bm Z^0}(\bm I^{\bm Z^0}_{n-1})\subset \bm Z^0,$$
then let $T<5b$ be minimal with this property, and let 
$$\bm I^{\bm Z^0}_{n}=\bm E^{\bm Z^0}$$
be the pullback of $\Gamma_{\bm Z^0}^{T}\mathcal{B}_{\bm Z^0}
\mathcal{A}_{\bm Z^0}(\bm I^{\bm Z^0}_{n-1})$
that is combinatorially close to $\bm Z^0$ that is given by
Lemma~\ref{lem:short step}.
\item If no such $T$ exists, set $T=5b,$ and define
$$\bm I^{\bm Z^0}_{n}=\Gamma_{\bm Z^0}^T\mathcal B_{\bm Z^0}\mathcal A_{\bm Z^0}(\bm I^{\bm Z^0}_{n-1}).$$
\end{itemize}
Notice that this last case, is the only one where $\bm E^{\bm Z^0}$ is
not defined, and in this case, we continue on to define $\bm I^{\bm
  Z^0}_{n+1}$.
\end{defn}

Concluding the proof, we set
$\II^+_n=\Gamma_{\Z^0}^{5b}\mathcal{B}_{\Z^0}\mathcal{B}_{\Z^0}(\II^{\Z^0}_{n-1})$
and $\II^-_n=\Gamma_{\Z^0}^{5b}\mathcal{A}_{\Z^0}\mathcal{A}_{\Z^0}(\II^{\Z^0}_{n-1}).$
\end{pf}

Suppose that 
$\bm I_0^{\bm Z^0}\supset\bm  I_1^{\bm Z^0}\supset\bm  I_2^{\bm Z^0}\supset\dots\supset \bm I_{n+1}^{\bm Z^0}$
are defined and that $\hat{\mathcal A}(\bm I_{n+1})$
contains $\bm Z^0$.
Let $p_n$ be so that
$F^{p_n}(\bm I_{n+1}^{\bm Z^0})=\bm I_n^{\bm Z^0}$.
If $\II$ is a puzzle piece, we
let $r(\bm I)$ be the minimal return time of a point in $\bm I$ to $\bm I$.
Let $\bm M_{n,i}=\Gamma^i_{\bm Z^0}\mathcal{B}_{\bm Z^0}\mathcal{A}_{\bm Z^0}(\bm I_n)$.
Let $q_{n,i}$ be so that $F^{q_{n,i}}(\bm M_{n,i})=\bm M_{n,i-1}.$

\begin{lem}[See \cite{KSS-rigidity} Lemma 8.3, Transition and return time relations]\label{lem:transition times}
Suppose that 
$\bm I_0^{\bm Z^0}\supset\bm  I_1^{\bm
  Z^0}\supset\bm  I_2^{\bm Z^0}
\supset\dots\supset \bm I_{n+1}^{\bm Z^0}$
are defined and that $\hat{\mathcal A}(\bm I_{n+1}^{\bm Z^0})$ 
contains $\bm Z^0$. Then
\begin{enumerate}
\item $3r(\bm I_{n+1}^{\bm Z^0})\geq p_n$,
\item $p_{n+1}\geq 2 p_n.$
\end{enumerate}
\end{lem}
\begin{pf}
The proof follows the proof of \cite{KSS-rigidity} Lemma 8.3.
The proofs of these inequalities follow from the 
first three inequalities of \cite{KSS-rigidity} Lemma 8.3.
These are:
\begin{itemize}
\item $2b^2 r(\mathcal{A}_{\Z^0}(\bm I_n))\geq s_n\geq r(\II_n^{\Z^0}),$
\item $b^2 r(\mathcal{B}_{\Z^0}\mathcal{A}_{\Z^0}(\II_n^{\Z^0}))\geq t_n
\geq r(\mathcal{A}_{\Z^0}(\bm I_n^{\Z^0})),$
\item $r(\bm M_{n,i})\geq q_{n,i}\geq 2r(\bm M_{n,i}),$
\end{itemize}
where $s_n$ and $t_n$ are such that
$$\mathcal{A}_{\Z^0}(\II_n)=\comp_{c_0}F^{-s_n}(\II_n^{\Z^0})$$
and 
$$\mathcal{B}_{\Z^0}\mathcal{A}_{\Z^0}(\II_n^{\Z^0})=\comp_{c_0}F^{-t_n}(\mathcal{A}_{\Z^0}(\II_n^{\Z^0})).$$
The proofs of the first and second inequalities are unchanged from \cite{KSS-rigidity}.
We need only comment on the proof of the third:
Since $F$ is non-renormalizable, each puzzle piece $\bm M_{n,i-1}$ 
has at least two children satisfying conditions (1) and (2)
on page~\pageref{page:gammaprops}.
Moreover, by construction, for each $i$, the second child of 
$\bm M_{n,i-1}$ contains $\bm Z^0,$
and hence also contains $\bm M_{n,i}.$
Thus $R_{\bm M_{n,i-1}}(\bm M_{n,i})\cap \bm M_{n,i}=\emptyset$,
so $q_{n,i}\geq 2r(\bm M_{n,i-1}).$
Finally observe that for
$0<j<q_{n,i}$, $F(\bm M_{n,i})\not\owns c_0$,
so $F(\bm M_{n,i})\cap \bm M_{n,i}=\emptyset$.
This implies that $r(\bm M_{n,i})\geq q_{n,i}.$
\end{pf}

Lemma~\ref{lem:transition times} implies that
there are at most six returns of $c_0$ to $\Z^0$ under 
$F^{p_{n}+p_{n-1}+\dots +p_{i}}$.
Since we chose the central cascade to have length at least eight,
this implies that for each $i$, the critical values of
$F^{p_{n}+p_{n-1}+\dots +p_{i}}|_{\II_{n+1}^{\Z^0}}$
are contained in $\Gamma^{5b}_{\Z^0}\mathcal{A}_{\Z^0}\mathcal{A}_{\Z^0}(\II^{\Z^0}_{i-1})$.

\subsubsection{Complex bounds for the enhanced nest above a central cascade}
We will use the following lemma which give us
``small distortion of thin annuli.''
\begin{lem}[Compare \cite{KvS} Lemma 9.1]\label{lem:small distortion of thin annuli}
For every $K\in(0,1)$, $\delta>0$,
there exists $\kappa>0$ such that the following holds.
Suppose that $F\colon\bm U\rightarrow \bm Z$ is $(1+\delta)$-quasiregular,
real-symmetric, branched covering map,
of degree $D$ with all critical values real. 
Assume that $\bm A\subset \bm U,$ $\bm B\subset \bm Z$ are simply connected domains,
symmetric with respect to the real axis,
and that the domain $\bm A$ is the
connected component of $F^{-1}(\bm B)$
with respect to the real-line
and the degree of $F|_{\bm A}$ is $d$.
Moreover, suppose that $F$ can be decomposed as a composition of mappings
$f_1\circ\dots\circ f_n$ with all maps $f_i$ real
and either a $(1+\delta_i)$-quasiconformal diffeomorphism,
or a $(1+\delta_i)$-quasiregular branched covering map
with a unique real critical point, where
$$\prod_{i=1}^n(1+\delta_i)\leq (1+\delta),$$ 
then
$$\mod(\bm U\setminus \bm A)\geq \frac{K^D}{2d(1+\delta)}\min\{\kappa,\mod(\bm Z\setminus \bm B)\}.$$
\end{lem}

The previous lemma is an adaptation of Lemma 9.1 in \cite{KvS},
where the mappings are assumed to be analytic,
to the asymptotically holomorphic setting, and the proof can be 
repeated verbatim.
In the real-analytic setting, it implies and sharpens the
\emph{Covering Lemma} of Kahn and Lyubich \cite{KL}.

With Lemma~\ref{lem:small distortion of thin annuli} in hand,
one can repeat the arguments of
\cite{KvS} to prove complex bounds for complex bounds for 
the enhanced nest above $\bm Z^0$:

\begin{prop}[see \cite{KvS} Proposition 10.1]
\label{prop:modified delta nice}
For any $\varepsilon>0,$ there exists $\delta>0$ such that the following holds. Suppose that $F$ be a complex box mapping.
If $\bm I^{\bm Z^0}_0$ is $\varepsilon$-nice, and $\bm I^{\bm Z^0}_0\supset\bm I^{\bm Z^0}_1\supset \bm I^{\bm Z^0}_2\supset\dots\supset \bm I^{\bm Z^0}_n$ is the
enhanced nest above $\bm Z^0$, and that $n \geq 8$.
then all $\bm I_i^{\bm Z^0}$, $n=0,1,2,\dots, n$ are $\delta$-nice.
\end{prop}

\begin{pf}
The proof is a slight modification of the proof of Proposition 10.1 of \cite{KvS}.
Let us omit the $\bm Z^0$ in the notation.
Suppose that $0\leq j\leq n$.
Let $\delta_j\geq 1$ be such that $\II_j$ is $\delta_j$-nice.
Fix $M>4$, to be chosen later, and suppose that $j>M+1$.
Let $\bm A$ be any component of the domain of the first return mapping to 
$\II_j$ and let $r$ be its return time.
\begin{rem}
Let us remark that we do not require that $\bm A$ intersect the
post-critical set of $F$.
\end{rem}
Let  $P_{j,M}=p_{j-1}+p_{j-2}+\dots+p_{j-M}$ and let 
$\bm B=F^{P_{j,M}}(\bm A).$

First, we show that $\bm B\subset\mathcal{L}_x(\II_{j-4}).$ 
Observe that $$r\geq r(\II_j)\geq
\frac{1}{3}p_{j-1}\geq\frac{16}{3}p_{j-5}\geq p_{j-5}+p_{j-6}+\dots
\geq P_{j-4,M-4},$$
so $s=r-P_{j-4,M-4}\geq 0$. Since $F^r(\bm A)=\II_j=F^s\circ
F^{P_{j-4,M-4}}(\bm A),$
we have that 
$$F^s(F^{P_{j,M}}(\bm A))=F^{p_{j-1}}\circ F^{p_{j-2}}\circ
F^{p_{j-3}}\circ F^{p_{j-4}}\circ F^s\circ F^{p_{j-4,M-4}}(\bm
A)=\II_{j-4}.$$
Thus $F^{P_{j-M}}(\bm A)$ is contained in some
$\mathcal{L}_x(\II_{j-4})$.

Let $\II_j^+=\Gamma^T\mathcal{BB}(\II_{j-1}).$ Then
$\II_j\subset\II_j^+\subset\II_{j-1}$ and
$\mod(\II_j^+\setminus\II_j)\geq K_1\delta_{j-1}.$

For each $\II_j^+$
let $\mu_j=\mu(\II_j^+)=\max\diam(\U),$ where the maximum is taken over
components of $\Dom'(\II_j^+)$.
Fix any $x\in F^{p_{n,M}}(\bm A)$ and $k\in\{n-M,n-M+1,\dots, n-4\}.$
Let $\nu$ be so that $F^{\nu}(x)=R_{\II_{k}}(x)$.
Let $\bm A_k$ denote the pullback of $\II_{k}^+\setminus\II_k$
by $F^\nu$ that contains $x$ in the bounded component of its complement. 
The degree of the mapping $F^\nu\colon\mathcal{L}_x(\II_k)\rightarrow\II_k$
is bounded from above by a universal constant depend only the critical
points of $F$ and their orders, and by Lemma~\ref{lem:angle control},
there exists a constant $C>0$ so that
the dilatation of $F^\nu|_{\bm A_k}$ is at most
$C\mu_k$.
By property (1) in the construction of 
$\Gamma_{\Z^0}$, see page~\pageref{page:gammaprops},
we have the same for 
$F^\nu\colon\comp_{x}F^{-\nu}(\II_k^+)\rightarrow\II_k^+.$
Hence there exists a constant $K_2$ such that 
$$\mod(\bm A_k)\geq \frac{K_1}{K_2}\delta_{k-1}.$$

The annuli $\bm A_k$ are nested
and surround $F^{p_{n,M}}(\bm A)$, which implies that
$$\mod(\II_{n-M},F^{p_{n,M}}(\bm A))\geq
\frac{K_1}{K_2}(\delta_{n-M-1}+ \delta_{n-M}+\dots+\delta_{n-5}).$$

The degree of the mapping $F^{P_{n,M}}|_{\bm A}$
is bounded from above by a constant independent of $M$, and the 
degree of $F^{p_{n,M}}|_{\II_n}$ is bounded from above by a constant
$K_4(M)$.

Let $\V=\II_{n-M},$ and let 
$\bm B$ and $\bm B'$ be the domains bounded by the 
inner and outer boundary of $\bm A_{n-4}$, respectively.
Then
$$\mod(\V,\bm B)\geq 
\frac{K_1}{K_2}(\delta_{n-M-1}+\delta_{n-M}+\dots\delta_{n-5}).$$
Let $\delta_{n-M-1,n-5}=\min_{n-M-1\leq i\leq n-5}\{\delta_{i}\}.$
We may assume that $K_1/K_2<1.$ 
Take $M$ to be the smallest integer 
greater than $4d(1+\delta)K_2/K_1$ and set
Let $D=K_4(M)$ and take $K=(\frac{4d(1+\delta)K_2}{K_1M})^{1/D}.$
Let $\kappa>0$ be the number associated to $K$
given by Lemma~\ref{lem:small distortion of thin annuli}.
By  Lemma~\ref{lem:small distortion of thin annuli}
we have that 
$$\mod(\II_{n}\setminus\bm A)
\geq \frac{K^{D}}{2d(1+\delta)}\min\{\kappa,\mod(\V\setminus\bm B)\}
$$
$$\geq \min\{ \frac{K^{D}}{2d(1+\delta)}\kappa ,
\frac{K^{D}}{2d(1+\delta)} \frac{K_1}{K_2}M\delta_{n-M-1,n-5}\}$$
$$\geq\min\{\kappa', 
\frac{4d(1+\delta)K_2}{2d(1+\delta)K_1M}\frac{K_1}{K_2}M\delta_{n-M-1,n-5}\}
=\min\{\kappa',2\delta_{n-M-1,n-5}\},$$
where $\kappa'=\frac{K^{D}}{2d(1+\delta)}\kappa.$
We also have that there exists $K$ so that $\delta_{j}\geq K\delta_{j-1}.$
This completes the proof.
\end{pf}
Using the construction of the enhanced nest above $Z^0$ we immediate deduce:
\begin{prop}\label{prop:space}
 There exists $\delta>0$ such that for all $j,0\leq j\leq n+1$,
$$\mod(\II_j^+\setminus
\II_j^{\Z^0})\geq\delta,\quad\mbox{and}\quad
\mod(\II_j^{\Z^0}\setminus\II_j^-)\geq\delta.$$
\end{prop}
Since for all $j\in\{0,\dots,n-1\}$ we have that
the critical values of $F^{P_{j}}$ are contained 
in $\II_{j}^-,$ and $\mod(\II_j^{\Z^0}\setminus \II_j^-)$
is bounded away from 0, independently of $j$, the same argument
given in \cite{KvS} implies:
\begin{prop}[See \cite{KvS} Proposition 11.1]
\label{prop:modified bounded geometry}
For any $\varepsilon,\rho>0$, there exists $\delta>0$ such that the following holds.
Let $F\colon\bm{\mathcal{U}}\rightarrow \bm{\mathcal{V}}$ be a real-symmetric complex box mapping and suppose that
$\bm I^{\bm Z^0}_0\supset \bm I^{\bm Z^0}_1\supset\dots\supset \bm I^{\bm Z^0}_n$
is the enhanced nest above $\bm Z^0$.
Assume that $\bm I^{\bm Z^0}_0$ be $\varepsilon$-nice and that it has $\rho$-bounded geometry 
with respect to $c_0$. Then all $\bm I^{\bm Z^0}_j,0\leq j\leq n$ have $\delta$-bounded geometry.
\end{prop}

Now we are in a position to prove our new result, that the puzzle
pieces in the enhanced nest or an enhanced nest over a 
central cascade are all quasidisks with bounded dilatation.

\begin{prop}\label{prop:quasidisks}
For any $K\geq 1$ there exists $K'\geq 1$ such that the following
holds.
Let $\bm I=:\II_0^{\bm Z^0}\supset\bm Z^0$ be puzzle pieces.
Suppose that $\bm I$ is a $K$-quasidisk. Then for all $j,0\leq
j\leq n$, $\bm I^{\bm Z^0}_j$ is a $K'$-quasidisk.
\end{prop}
\begin{rem}
In our applications, 
$F\colon\bm{\mathcal{U}}\rightarrow \bm{\mathcal{V}}$
is the complex box 
box mapping given by either 
Theorem~\ref{thm:box mapping persistent infinite branches}
or Theorem~\ref{thm:box mapping reluctant}, so that 
$\II_0^{\bm Z^0}$ is a Poincar\'e lense domain $D_{\theta}(I_0)$.
In particular, it is a $K(\theta)$-quasidisk.
\end{rem}
\begin{pf}
We will give the proof for the enhanced nest. The proof for the
enhanced nest above a long cascade of central returns is identical.
	
We will verify that the Ahlfors-Beurling Criterion,
see Lemma~\ref{lem:Ahlfors-Beurling}),
is satisfied for constant $C$ that depends only on the geometry of
$\bm{\mathcal{V}}\setminus\bm{\mathcal{U}}$
and the constant $\rho$ from
Theorem~\ref{thm:box mapping persistent}.
By Theorem~\ref{thm:box mapping persistent}, there exists $\rho>0$, such that
for any $k> 0$, there exist puzzle pieces
$\bm P^+_k\supset \bm I_k \supset \bm P^-_k$
such that 
$\mod(\bm P^+_k\setminus  \bm I_k)$ and 
$\mod(\bm I_k\setminus \bm P^-_k)>\delta,$
and such that each $\bm I_n$ has $\rho-$bounded geometry at the critical point.

\medskip
\noindent\textit{Claim:}
There exists a constant $\eta_0=\eta_0(\delta)>0$ such
$$\dist(\partial \bm I_k,\bm P^-_k) \geq \eta_0\cdot\diam(\bm I_k)\mbox{ and }
\dist(\partial \bm P^+_k, \bm I_k)\geq \eta_0\cdot\diam(\bm P_k^+).$$
	
\medskip
\noindent\textit{Proof of claim:}
We will need the following:
	
\medskip
\noindent\textit{Fact 1.} For any $\delta>0$,
there exists a constant $c>0$ such that if
$\mod(\V\setminus \U)>\delta$, then $\dist(\partial \V,\U)>c\cdot \diam(\U)$.
	
\medskip
\noindent\textit{Fact 2.}
There exists a constant $M_0>0$ such that if $\U\Subset \V$ are topological disks with
$\mod(\V\setminus \U)> M_0$, then for any $z_0\in \U$, there exists a
round annulus $\bm A=\{z:r<|z-z_0|<R\}\subset \V\setminus \U$ 
with $\mod(\bm A)\geq \mod(\V\setminus \U)-M_0$.

\medskip
	
We will prove the first estimate
in the claim, the second is similar.
First suppose that $\mod(\bm I_k \setminus \bm P^-_k)> 3M_0$.
Then by the above fact, there exists a round annulus
$\bm A=\{z:r<|z-c_0|<R\}$
contained in $\bm I_k\setminus \bm P^-_k$
with $\mod \bm A\geq \mod(\bm I_k\setminus \bm P^-_k) -M_0$.
Since $\bm I_k$ has $\delta$-bounded geometry with respect to $c_0$,
we can assume that $R$ is comparable to $\diam(\bm I_k)$. But now since 
$r$ is much smaller than $R$, $\dist(\partial\bm I_k, \bm P^-_k)$ is comparable to
$\diam(\bm I_k).$

Now suppose that $\mod(\bm I_k \setminus \bm P^-_k)\leq 3M_0$,
Then since $\bm I_k$ has $\delta$-bounded geometry at $c_0$, 
there exists a constant $c'$, depending only on $M_0$, such that
$\mathrm{diam}(\bm P^-_k)\geq c'\cdot\diam(\bm I_k).$
Then since there exists a constant $c$ such that 
$\mathrm{dist}(\bm P^-_k,\partial\bm I_k)>c\cdot\diam(\bm P_k^-)$
we have that $\mathrm{dist}(\bm P^-_k,\partial\bm I_k)>cc'\cdot\diam(\bm I_k).$
\endpfclaim

With the claim in hand, let $\bm I_n$ be a critical puzzle piece
in the enhanced nest, and choose
points $x,y\in\partial \bm I_n$. For each $j\geq 0$, let $p_j$ be defined so that
$F^{p_j}|_{\bm I_{j+1}}$ maps $\bm I_{j+1}$ onto $\bm I_j$ as a 
quasi-regular
branched
covering mapping with bounded degree and bounded dilatation,
and define $s_k=p_{n-1}+p_{n-2}+\dots p_{k}$. 
Let $x_j=F^{s_j}(x)$ and $y_j=F^{s_j}(y)$.
For each $j\leq n$, let $\gamma_j$ be the shortest path in
$\partial \bm I_j$ connecting $x_j$ and $y_j$.
Suppose that there exists $k\geq 0$, such that
$|x_k-y_k|>(\eta_0/2)\diam(\bm I_k)$,
and let $0\leq k_0<n$ be maximal with this property.
then $|x_{k_0}-y_{k_0}|>(\eta_0/2)\diam(\gamma_{k_0})$.
Pulling back by $F^{p_{k_0}}\colon I_{k_0+1}\rightarrow I_{k_0}$,
we have that there exists a constant $C$ such that
$C|x_{k_0+1}-y_{k_0+1}|>\diam(\gamma_{k_0+1}),$
and $|x_{k_0+1}-y_{k_0+1}|<(\eta_0/2)\diam(\bm I_{k_0+1})$.
Since $|x_{k_0+1}-y_{k_0+1}|<(\eta_0/2)\diam(\bm I_{k_0+1})$,
by the claim, the line segment $l_{k_0+1}$ connecting $x_{k_0+1}$ and
$y_{k_0+1}$ is well-inside $\bm P^+_{k_0+1}\setminus \bm P_{k_0+1}^-$,
and $\dist(\gamma_{k_0+1},\partial \bm P^{+}_{k_0+1}\cup\partial \bm P^{-}_{k_0+1})
\geq\eta_0\cdot\diam(\bm I_{k_0+1})$.
By real-symmetry and since $\partial\bm I_{k_0+1}$
intersects the real line exactly twice, we can assume that
$\gamma_{k_0+1}$ intersects the real line at most once. 
So we can remove a real ray that does not intersect either 
$l_{k_0+1}$ or $\gamma_{k_0+1}$ and contains the post-critical set of $F$.
Choosing the branch of
$F^{s_{k_0+1}}\colon\bm I_n\rightarrow \bm I_{k_0+1}$
that maps $x$ to $x_{k_0+1}$ and $y$ to $y_{k_0+1}$,
and pulling back by this map,
which has bounded dilatation, we have that there exists a 
constant $C>0$, independent of $n$ such that 
$C|x_n-y_n|\geq \diam(\gamma_n)$.
Suppose now that no such $k$ exists. Then $|x_0-y_0|<(\eta/2)\diam(\bm I_0),$
and there exists a constant $C>0$ such that $C|x_0-y_0|>\diam(\gamma_0)$,
so we can argue as in the previous case.
	\endpfclaim
	\end{pf}

From the definition of the enhanced nest above $\bm Z^0$ we immediately have
\begin{prop}[Good geometry for central cascades]
\label{prop:good geometry for central cascades}
There exists $\delta>0$ such that
if $\bm{Z}^0\supset \bm{Z}^1\supset\dots\supset\bm{Z}^{N+1}$ is a
maximal central cascade,
then there exists a $\delta$-excellent puzzle piece 
$\bm G^0$ which is combinatorially close to $\bm P^0$. See
page~\pageref{page:defgood}
for the definition of $\delta$-excellent.
\end{prop}

\subsection{Chains without long central cascades}
\label{sec:diam sq}
In this subsection, we will bound
the sum of the squares of the diameters of the puzzle pieces in any
chain of puzzle pieces that intersects the real-line and does not 
contain a long central cascade. These estimates are used to control
the dilatation of certain restrictions of $f^n$.
We will use the notation $\bm Y^i$ to denote
puzzle pieces in a principal nest.
\subsubsection{Combinatorial depth}
Let $c$ be a critical point and suppose that
$\bm Y\owns c$ is a critical puzzle piece.
Let $\bm Y=\bm Y^0\supset \bm Y^1\supset \bm Y^2\supset \dots,$
where $\bm Y^{i+1}=\mathcal{L}_c(\bm Y^i)$, be the principal nest
about $c$. Let $m>0$, where we allow $m=\infty,$ be minimal so that
$R_{\bm Y^0}|_{\bm Y^1}(c)\in \bm Y^{m-1}\setminus \bm
Y^m$.
If $m$ is finite,
define $\mathcal{C}_c(\bm Y)=\bm Y^{m}$.
If $m$ is infinite, let $\mathcal{C}_c(\bm Y)=\II^\infty.$

If $\J$ is a return domain to an arbitrary complex puzzle piece 
$\II$, and $\{\bm G_{i}\}_{i=0}^{r}$ is the chain with $\bm G_{r}=\II$
and 
$\bm G_{0}=\J$ where $r$ is the return time of  $\J$ to $\II$, we define
$$\mathrm{Crit}(\II;\J)=\Big(\bigcup_{i=0}^{r-1}\bm G_{i}\Big)\cap\mathrm{Crit}(F).$$
Similarly, if $\mathbb{G}=\{\bm G_{j}\}_{j=0}^{s}$ is an arbitrary
chain
such that the pullbacks of $\bm G_s$ and $\II$ are either nested or disjoint, 
$\bm G_{0}\subset \II$ and $0=n_{0}<n_{1}<\dots<n_{p}=s$
are the integers with $\bm G_{n_{i}}\subset \II$, we define
$$\mathrm{Crit}(\bm I;\mathbb{G})=\bigcup_{i=0}^{p-1}
\mathrm{Crit}(\II;\mathcal{L}_{\bm G_{n_{i}}}(\II)).$$\label{crit symbol}

For any complex puzzle piece $\II$ and any critical point $c$ we define 
$$k_{c}(\II,\mathbb{G})=\inf\{k_{c}\ge 0:\mbox{ there exists no }
j=1,2,\dots,s-1\mbox{ with }\bm G_{j} \subset  
\mathcal{C}^{k_{c}}_c(\hat{\mathcal{L}}_{c}(\II))
\},$$
where we take  $\mathcal{C}^0_c(\mathcal{\hat L}_{c}(\II))= \hat{\mathcal{L}}_{c}(\II)$
and
$$k(\II,\mathbb{G})=\sum_{c\in\mathrm{Crit}(\II;\mathbb{G})}k_{c}(\II,\mathbb{G}).$$
The \emph{combinatorial depth of the chain $\mathbb{G}$ with respect
  to $\II$} 
is defined to be $k(\II,\mathbb{G})$. 
Note that $k(\II,\mathbb{G})$ is well-defined even if 
$\II$ does not contain a critical point. 
These definitions are the same if we replace $\II$ and $\J$ by nice intervals.

Let us recall the following result:
\begin{prop}\cite[Proposition 7.1]{CvST}\label{prop:combinatorial depth angle control}
For each $\delta>0, k\geq 0$ and $N\geq 0$ there exist $\mu(k,N,\delta)\in(0,1)$ and $\delta'>0,$ and for each $\theta\in(0,\pi)$ 
there exists $\varepsilon>0$ so that 
the following holds. Let $I$ be a $\delta$-nice interval with $|I|<\varepsilon$.
Suppose that $\mathbb{G}:=\{G_i\}_{i=0}^s$ is a chain
such that $G_0, G_s$ are nice intervals contained in $I$,
the pullbacks of $G_s$ and $I$ are nested or disjoint, the intervals 
$G_0,\dots,G_{s-1}$ are pairwise disjoint
and $G_0\cap\omega(c_0)\neq\emptyset.$
Assume that 
$$k(I,\mathbb{G})\leq k\;\; and\; \;\#\mathrm{Crit}(I;\mathbb{G})\leq N.$$
Let $\hat{G}_{s}$ be an interval with $G_{s}\subset \hat{G}_{s}\subset (1+2\delta)\hat{G}_{s}\subset I$. Let 
$\V=D_{\theta}(\hat{G}_{s})\cap\mathbb{C}_{G_{s}}$ and
$\U_{i}=\mathrm{Comp}_{G_{i}}F^{-(s-i)}(\V)$ for $i=0,\dots,s$.
Then, there exists an interval $\hat{I}\supset G_{0}$ with $(1+2\delta')\hat{I}\subset I$ such that 
$$\U_{0}\subset D_{\mu(k,N,\delta)\theta}(\hat{I}).$$ 
\end{prop}

For any point $x$, let $\bm Y^0_x\supset \bm Y_x^1\supset
\bm Y^2_x\supset \dots$ be the principle nest
constructed about $x$. 
Let $\hat m>0$, where we allow $\hat m=\infty,$ be minimal
$R_{\bm Y^0_x}|_{\bm Y^1_x}(x)\in \bm Y^{\hat m-1}\setminus \bm
Y^{\hat m}$.
If $\hat m$ is finite,
define 
$\tilde{\mathcal{C}}_x(\bm Y)=\bm Y^{\hat m}_x$ and
$\hat{\mathcal{C}}_x(\bm Y_x)=\bm Y^{\hat m+1}_x$,
By \cite[Theorem A]{vSV}
we have that $\hat{\mathcal{C}}_x(\bm Y_x)$ is well-inside
$\tilde{\mathcal{C}}_x(\bm Y_x)$, and that 
there exists a constant $\delta>0$ so that 
$\hat{\mathcal{C}}_x(\bm Y)$ is $\delta$-nice.
If $\hat m$ is infinite, let
$\tilde{\mathcal{C}}_x(\bm Y)=\hat{\mathcal{C}}_x(\bm Y)=\II^\infty.$

Given a puzzle piece $\U\owns x$,
let $\hat k(\bm Y^0_x,\U)=\hat k$ be minimal so that
$\tilde{\mathcal{C}}_x^{\hat k}(\bm Y^0_{x})\subset\U$, and we
define the \emph{depth} of a puzzle piece $\U\owns x$ to be
$\hat k(\bm Y^0_{\U},\U)$.

Let $\{\U_j\}_{j=0}^s$ be a chain of complex puzzle pieces.
Suppose that $\U_{j_1}$ is a first return domain to $\U_{j_0}$
and $\U_{j_2}$ is a first return domain to $\U_{j_1}$, 
$0\leq j_2<j_1<j_0\leq s$.
We say that $\{\U_j\}_{j=0}^s$
is \emph{free from central returns} if
we have that $R_ {\U_{j_0}}|_{\U_{j_1}}(\U_{j_2})\subset
\U_{j_0}\setminus \U_{j_1}$.
A central cascade in a chain $\{\U_j\}_{j=0}^s$
is a sequence of critical puzzle pieces
$$\U_{j_0}\supset \U_{j_1}\supset
\U_{j_2}\supset\dots\supset\U_{j_N}\supset \U_{j_{N+1}}$$
from the chain so that for each $i=0,\dots, N,$
$\U_{j_{i+1}}$ is a first return domain to $\U_{j_{i}}$
and return time of $\U_{j_{i+1}}$ to $\U_{j_{i}}$
is the same as the return time of $\U_{j_1}$
to $\U_{j_0}$. The maximum number $N$ so that this occurs
for a sequence of first returns in the chain
is the
\emph{length of the central cascade}.
We say that  $\{\U_j\}_{j=0}^s$ \emph{does not have central
cascades of length greater than} $N$
if any central cascade in the chain $\{\U_j\}_{j=0}^s$ has length
less than $N$.

\subsubsection{Controlling the diameters of
complex puzzle pieces}
For convenience, let us restrict to the case that 
$F\colon\bm{\mathcal{U}}\rightarrow\bm{\mathcal{V}}$
is a complex box mapping given by either
Theorem~\ref{thm:box mapping persistent infinite branches}
or Theorem~\ref{thm:box mapping reluctant}.

Let $M\in\mathbb N_{\geq 2},$ and
let $\{\U_j\}_{j=0}^s$ be a chain of complex puzzle pieces.
Let us construct the \emph{M-shielding intervals}:
For each
$U_j:=\U_j\cap\mathbb R$,
let $k$ 
be the depth of $\U_j$. Assume that $k$ is at least $M+1$.
Let $x\in\U_j$.
Then
$U_{j}\subset\tilde{\mathcal{C}}_x^{k-1}(Y_x^0)$.
For $i,$ $0\leq i\leq M-2$,
let $\mathcal{S}_i(U_j)=\hat{\mathcal{C}}_x^{k-(M-i)}(Y_x^0).$

\begin{lem}\label{lem:flower lemma}
There exist constants $M\geq 2, \lambda\in(0,1)$ and $\theta_0\in(0,\pi/2)$, such that
for any
$N\geq M$
the following holds:
Suppose that $\U$ is a puzzle piece
with $\hat k(\bm Y^0_{\U},\U)=N$
and $|\mathcal{S}_0(U)|$ sufficiently small.
Let
$\mathcal{S}_0(U)\supset \mathcal{S}_1(U)\supset\dots \mathcal{S}_{M-2}(U)$,
be the $M$-shielding intervals.
We have that
$$\U\subset D_{\theta_0}(\mathcal{S}_0(U))\cup
D_{\lambda \theta_0}(\mathcal{S}_1(U))
\cup\dots\cup D_{\lambda^{M-3}\theta_0}(\mathcal{S}_{M-3}(U))\cup
D_{\alpha_U}(\mathcal{S}_{M-2}(U)),$$
with $\alpha_U=|\mathcal{S}_{M-2}(U)|^{1/3}$.
\end{lem}

\begin{pf}
First we define the number $M$ in the statement of the theorem.
Let $J\owns x$ be a puzzle piece of depth $n$, and
$\delta>0$ be the universal constant so that for all $i\in\mathbb N,$
each $\hat{\mathcal{C}}_x^{i}(Y^0_x)$ is $\delta$-nice.
Let $C$ be the constant given by Lemma~\ref{lem:z2'}.
Let $M$ be so that $C\cdot J\subset \hat{\mathcal{C}}_x^{n-M}(Y^0_x).$

We will prove the lemma by induction of the depth of the puzzle pieces.
Since the top level puzzle pieces are Poincar\'e lens domains,
and puzzle pieces are nested or disjoint,
the initial step of the induction holds automatically for puzzle pieces with depth $M$.
So suppose that the statement holds for puzzle pieces of depth at most
$N\geq M$,
and let us argue that it holds for puzzle pieces of depth $N+1$.

Let $\U$ be any puzzle piece of depth $N+1$.
Let $s$ be so that $\U$ is a component of $F^{-s}(\bm{\mathcal{V}})$, and
let $\{\U_{j}\}_{j=0}^s$ be the chain so that $\U_0=\U$.
Let us divide the chain 
$\{\U_j\}_{j=0}^s$
into pieces based on their depth.
Let $l_0=s$. Let $l_1<s$ be maximal so that
there is a point $x'\in U_{l_1}$ such that
$U_{l_1}\subset\hat{\mathcal{C}}_{x'}(Y^0_{x'})$. Continuing inductively,
if $l_i$ has been defined, let $l_{i+1}<l_i$
be maximal so that there exists a point $x''\in U_{l_{i+1}}$
such that
$U_{l_{i+1}}\subset \hat{\mathcal{C}}_{x''}^{i+1}(Y_{x''}^0)$,
and so on.

By definition, we have that
$\hat{\mathcal{C}}_x^{N+1}(\bm{\mathcal{V}}(x))\subset \U_{0}\subset
\hat{\mathcal{C}}^{N}_x(\bm{\mathcal{V}}(x))\subset
\tilde{\mathcal{C}}^{N}_x(\bm{\mathcal{V}}(x))$
Let us now consider the chain 
$\{\V_j\}_{j=0}^{s'},$ where 
$\V_0=\hat{\mathcal{C}}^{N}_x(\bm{\mathcal{V}}(x))$
and
$\V_{s'}=\tilde{\mathcal{C}}^{N}_x(\bm{\mathcal{V}}(x))$.

\begin{figure}[htp] \hfil
\beginpicture
\dimen0=0.3cm
\setcoordinatesystem units <\dimen0,\dimen0>
\setplotarea x from 15 to 30, y from -10 to 12
\setlinear
\plot -1 0 11 0 /
\circulararc 241 degrees from 10 0  center at  5 3 
\circulararc -241 degrees from 10 0  center at  5 -3 
\circulararc -292 degrees from 8 0  center at 5 -4.5
\circulararc +292 degrees from  8 0 center at  5 4.5 
\circulararc +35 degrees from 5 0 center at 5 30
\circulararc +35 degrees from 5 0 center at 4.5 -30
\circulararc -35 degrees from 5 0 center at 5 30
\circulararc -35 degrees from 5 0 center at 4.5 -30
\put{$D_{\theta_0}(\mathcal{S}_{0}(U))$} at 11 9.5
\put{$D_{\lambda\theta_0}(\mathcal{S}_{1}(U))$} at 14.5 3
\put{$D_{\alpha_U}(\mathcal{S}_{M-2}(U))$} at 23 6
		\endpicture
\caption{The shielding intervals and Poincar\'e disks with $M=4$.
When we pull back, this picture replicates itself at successive depths.
The pullback of
$D_{\theta_0}(\mathcal{S}_{0}(U_{j_1}))$
is contained in
$D_{\lambda\theta_0}(\mathcal{S}_{1}(U_{j_3})),$ 
the pullback of
$D_{\lambda\theta_0}(\mathcal{S}_{1}(U_{j_0})),$
is contained in
$D_{\alpha_U}(\mathcal{S}_{2}(U_{j_3})),$
and 
the pullback of 
$D_{\alpha_U}(\mathcal{S}_{2}(U_{j_0})),$
is contained in 
$D_{\eta\alpha_U}(\mathcal{S}_{2}(U_{j_3}))\cup
D_{\theta_0}(\mathcal{S}_{0}(U_{j_3})),$ with $\eta$ close to 1.
}
\label{fig:shielddisks}
		\end{figure}

\medskip

\noindent\textit{Claim 1.}
There exists $\lambda\in(0,1)$ such that the 
following holds.
Let $0\leq i\leq M-3$, and
let $\tilde{U}_{i}$ be an interval that is well-inside
$\mathcal{S}_i(U_{0})$ and
that contains $V_{s'}$. Then 
there exists an interval $\tilde{U}_{i+2}$ that is well-inside
of $\mathcal{S}_{i+2}(U_{0})$,
such that for any $\theta\in(0,\pi/2)$
$$\comp_{V_0}(F^{-s'}(D_{\theta}(\tilde{U}_i))
\cap\mathbb C_{V_{s'}})\subset D_{\lambda\theta}(\tilde U_{i+2}),$$
where $\mathbb C_{V_{s'}}=\mathbb C\setminus(\mathbb R\setminus V_{s'}).$
\medskip 

\noindent\textit{Proof of Claim 1.}
Choose $x\in U_0$. Let $j_3>0$ be minimal so that
$V_{j_3}\subset \tilde{\mathcal{C}}_x^{N-i+2}(Y^0_x).$
Such a $j_3$ exists since the return to 
$\tilde{\mathcal{C}}_x^{N-i+2}(Y^0_x)$ is non-central.
Let $J^0=\mathcal{S}_i(U_{0})$, 
and 
assuming that $J^i$ has already been constructed
let $J^{i+1}=\mathcal{L}_x(J^i)$.
Recall that there exists $\delta>0$ so that $\mathcal{S}_i(U_{0})$
is $\delta$-nice.
Let $n\geq 0$ be minimal so that 
$R_{J_0}(x)\in J^{n}\setminus J^{n+1}$.
Let $j_1<s$ be maximal so that $V_{j_1}\in J^n\setminus J^{n+1}.$
Note, $j_1>j_3$.
Let $L$ denote the domain of the return mapping to
$J^0$ that contains $R_{J^0}^{n}(V_{j_1})$,
and let $j_1'>j_1$ be minimal so that $V_{j_1'}\subset L.$
It follows from Proposition~\ref{prop:combinatorial depth angle
  control}
and Corollary~\ref{cor:angle control}
 that there exists $\lambda_1'\in(0,1)$ 
and an interval $\tilde{K}_{j_1'}$ that is well-inside $L$,
such that
$$\comp_{V_{j_1'}}(F^{-(s'-j_1')}(D_{\theta}(\tilde{U}_i))\cap\mathbb
C_{V_{s'}})\subset D_{\lambda_1'\theta}(\tilde K_{j_1'}).$$
To see this, observe that we can use 
Proposition~\ref{prop:combinatorial depth angle control}
to control the angle when we 
pull back 
from time $s'$  
to the last visit to $J^0$ before time $j_1'$, and then use 
Corollary~\ref{cor:angle control} to control the angle
when we pull back one more time along the
return mapping of $L$ to $J^0.$

Let $\{L_i\}$ denote the chain with $L_0\subset J^{n}\setminus
J^{n+1},$ the component of  $R_{J_0}^{-n}(L)$ that contains 
$V_{j_1}$ and $L_n=L.$ 
Since the pullback from $L_0$ to $L_n$ is along
a diffeomorphic branch of a central return, 
we have that there exists
an interval 
$\tilde{K}_{j_1}$ well-inside $L^0$
and $\lambda_1\in(0,1)$ so that
$$\comp_{V_{j_1}}F^{-(s'-j_1)} (D_{\theta}(\tilde U_{i})
\cap\mathbb C_{V_{s'}})
\subset D_{\lambda_1\theta}(\tilde{K}_{j_1}).$$

Moreover, $L_0$ is $\delta'$-nice for some $\delta'>0$.
Let $j_2> j_3$ be minimal so that $V_{j_2}\subset
\tilde{\mathcal{C}}_x^{N-i+1}(Y^0_x)\setminus \hat{\mathcal{C}}_x^{N-i+1}(Y^0_x).$ 
Let $j_2'>j_2$ be minimal so that
$V_{j_2'}\subset L_0$. By Proposition~\ref{prop:combinatorial depth
  angle control} there exists $\lambda_2'\in(0,1)$ and an interval
$\tilde{K}_{j_2'}$ that is well-inside $L_0$
such that
$$\comp_{V_{j_2'}}F^{-(s'-j_2')}(D_{\theta}(\tilde U_{i}) \cap\mathbb
C_{V_{s'}})
\subset
D_{\lambda_2'\theta}(\tilde{K}_{j_2'}),$$
and by Proposition~\ref{prop:combinatorial depth angle control}
and Corollary~\ref{cor:angle control}
there exists $\lambda_2$ and an interval $\tilde{K}_{j_2}$
that is well-inside
$\mathcal{L}_{U_{j_2}}(\mathcal{C}_x^{n-i+1}(Y^0_x))$
such that 
$$\comp_{V_{j_2}}F^{-(s'-j_2)} (D_{\theta}(\tilde U_{i}) \cap\mathbb C_{V_{s'}})\subset
D_{\lambda_2\theta}(\tilde{K}_{j_2}).$$
Since $\tilde{K}_{j_2}$ is disjoint from $\hat{\mathcal{C}}_x^{n-i+1}(Y^0_x),$
the pullback of $\tilde K_{j_2}$ is contained in 
$\tilde{\mathcal{C}}_x^{n-i+2}(Y^0_s).$
Arguing as before, by
Proposition~\ref{prop:combinatorial depth angle control} 
and Corollary~\ref{cor:angle control}, 
there exist $\lambda_3\in(0,1)$ and
$\tilde{K}_{j_3}$ well-inside $\tilde{\mathcal{C}}_x^{n-i+2}(Y^0_x)$, so that
$$\comp_{V_{j_3}}F^{-(s'-j_3)} (D_{\theta}(\tilde U_{i}) \cap\mathbb C_{V_{s'}})
\subset
D_{\lambda_3\theta}(\tilde{K}_{j_3}).$$
Finally by Corollary~\ref{cor:angle control},
there exists $\lambda\in(0,1)$
and an interval $\tilde{U}_{i+2}$ that is well-inside $\mathcal{S}_{i+2}(U_{0})$ such
that
$$\comp_{V_{0}}F^{-s'}(D_{\theta}(\tilde U_{i}) \cap\mathbb C_{V_{s'}})
\subset
D_{\lambda\theta}(\tilde{U}_{i+2}).$$
\endpfclaim

\medskip

As above, $\{\V_{j}\}_{j=0}^{s'}$ denotes the chain with
$\V_{s'}=\tilde{\mathcal{C}}_x^{N}(\bm Y^0_x)$
and $\V_0=\hat{\mathcal{C}}_x^{N}(\bm Y^0_x)$.
Now, let us pull back the slit Poincare disk
$D_{\alpha_{\mathcal{S}_{M-2}(U)}}(\tilde U_{M-2})\cap\mathbb C_{V^{s'}}$
back along the chain $\{V_j\}_{j=0}^{s'}$.

\medskip
\noindent\textit{Claim 2.}
There exist $\theta_0\in(0,\pi/2),$
$\alpha\in(|\mathcal{S}_{M-1}(U)|^{1/3},|\mathcal{S}_{M-2}(U)|^{1/3})$
such that
the following holds.
Let $\tilde{U}_{M-2}$ be an interval that is well-inside
$\mathcal{S}_{M-2}(U_{0})$ and
that contains $V_{s'}$. Then
$$\comp_{V_0}F^{-s'}(D_{|\mathcal{S}_{M-2}(U)|^{1/3}}(\tilde U_{M-2})\cap\mathbb C_{V^{s'}})
\subset
D_{\alpha}(\tilde{\mathcal{C}}_x^N(Y^0_x))\cup
D_{\theta_0}(\hat{\mathcal{C}}^{N-M}_x(Y^0_x))\cap \mathbb C_{V_0}.$$

\medskip

\noindent\textit{Proof of Claim 2:}
Recall that
$\tilde{\mathcal{C}}_x^N(Y^0_x)\subset \mathcal{S}_{M-1}(U_{0}).$
Let us start by explaining how we use 
Lemma~\ref{lem:z2'} 
together with Proposition~\ref{prop:combinatorial depth angle control}
to control the pullbacks of the Poincar\'e disk
with small angle along many iterates of the return mapping to
$\mathcal{S}_{M-2}(U_{0})$.

Let $t_1<s'$ be maximal so that
$V_{t_1}\subset \mathcal{S}_{M-2}(U_{0}),$
then since 
there exists $\delta>0$ such that
$(1+2\delta)\tilde{U}_{M-2}\subset \mathcal{S}_{M-2}(U_{0})$, by Lemma~\ref{lem:z2'}, 
Proposition~\ref{prop:Almost Schwarz Inclusion}
and Proposition~\ref{prop:combinatorial depth angle control},
there exists $\eta$ close to 1, $\delta'>0$, and an interval
$(1+2\delta')\tilde{U}_{t_1}'\subset \mathcal{L}_{V_{t_1}}(\mathcal{S}_{M-2}(U_{0})),$
and an angle $\theta_0\in(0,\pi/2)$ such that
$$\comp_{V_{t_1}}F^{-(s'-t_1)}(D_{|\mathcal{S}_{M-2}(U)|^{1/3}}(\tilde{U}_{M-2})
\cap\mathbb C_{V_{s'}})
\subset D_{\eta |\mathcal{S}_{M-2}(U)|^{1/3}}(\tilde{U}_{t_1}')\cup
D_{\theta_0}(\hat{\mathcal{C}}_x^{N-1-M}(Y^0_x)),$$
but now since $\mathcal{S}_{M-2}(U_{0})$ is $\delta$-nice, have that
there exists an interval $\tilde U_{t_1}$ such that $(1+2\delta)
\tilde U_{t_1}\subset \mathcal{S}_{M-2}(U_{0})$
and
$$\comp_{V_{t_1}}F^{-(s'-t_1)}(D_{|\mathcal{S}_{M-2}(U)|^{1/3}}(\tilde{U}_{M-2})
\cap\mathbb C_{V_{s'}})
\subset 
 D_{|\mathcal{S}_{M-2}(U)|^{1/3}}(\tilde{U}_{t_1})\cup
D_{\theta_0}(\hat{\mathcal{C}}_x^{N-1-M}(Y^0_x)).$$
Let us remark that further pullbacks of 
 $D_{\theta_0}(\hat{\mathcal{C}}_x^{N-1-M}(Y^0_x))$
are handled by Claim 1, so we will not treat them in detail in this
part of the proof.

We now argue exactly as in the proof of Claim 1, but we control the
loss of angle differently.
We set $J^0=\mathcal{S}_{M-2}(U_{0}),$ $J^{i+1}=\mathcal{L}_x(J^i)$, and let 
$n\geq 0$ be minimal so that 
$R_{J^0}(x)\in J^{n}\setminus J^{n+1}$.
We define times 
$j_1<s'$ maximal so that $V_{j_1}\in J^n\setminus J^{n+1},$
$j_1'>j_1$ minimal so that $V_{j_1}\subset J^0\setminus J^1.$
We let $L=\mathcal{L}_{V_{j_1'}}(\mathcal{S}_{M-2}(U_{0}))$, and we let
$\{L_i\}_{i=0}^n$
be the chain with $L_n=L$ and $L_0=\comp_{V_{j_1}}R^{-n}_{J^0}(L_n)$.
We let $j_2>0$ be minimal so that
$V_{j_2}\subset\tilde{\mathcal{C}}^{N-1}_x(Y^0_x)\setminus
\hat{\mathcal{C}}^{N-1}_x(Y^0_x)$, and we let $j_2'>j_2$ be minimal so
that $V_{j_2'}\subset L_0$.

We have that there exists
$\eta$ close to 1, and interval $\tilde{U}_{j_1'}$ that is well-inside 
$L_n$ so that  
$$\comp_{V_{j_1'}}F^{-(s'-j_1')}(D_{|\mathcal{S}_{M-2}(U)|^{1/3}}(\tilde{U}_{M-2})
\cap\mathbb C_{V_{s'}})
\subset 
 D_{\eta|\mathcal{S}_{M-2}(U)|^{1/3}}(\tilde{U}_{j_1'})\cup
 D_{\theta_0}(\hat{\mathcal{C}}_x^{N-1-M}(Y^0_x)).$$
Since the pullback from $L_0$ to $L_n$ is by a diffeomorphism, 
there exists $\eta'$ close to one and an interval $\tilde{U}_{j_1}$
that is well-inside $L_0$ so that
$$\comp_{V_{j_1}}R_{J_0}^{-n}(D_{\eta|\mathcal{S}_{M-2}(U)|^{1/3}}(\tilde{U}_{j_1'})
\cap\mathbb C_{V_{s'}})
\subset D_{\eta'|\mathcal{S}_{M-2}(U)|^{1/3}}(\tilde{U}_{j_1}))$$
Now we have that there exists an interval $\tilde U_{j_2'}$
that is well-inside $L_0$ so that
$$\comp_{V_{j_2'}}F^{-(j_1-j_2')}(D_{\eta'|\mathcal{S}_{M-2}(U)|^{1/3}}(\tilde{U}_{j_1})
\cap\mathbb C_{V_{s'}})
\subset D_{\eta'|\mathcal{S}_{M-2}(U)|^{1/3}}(\tilde{U}_{j_2'})\cup
D_{\theta_0}(\hat{\mathcal{C}}^{N-M-1}_x(Y^0_x)),$$
and there exists $\eta''$ close to 1 and an interval $\tilde{U}_{j_2}$
that is well-inside $\mathcal{L}_{V_{j_2}}(\mathcal{C}^{N-1}_x(Y^0_x))$
so that 
$$\comp_{V_{j_2}}F^{-(j_2'-j_2)}(D_{\eta'|\mathcal{S}_{M-2}(U)|^{1/3}}(\tilde{U}_{j_2'})
\cap\mathbb C_{V_{s'}})
\subset D_{\eta''|\mathcal{S}_{M-2}(U)|^{1/3}}(\tilde{U}_{j_2})\cup
D_{\theta_0}(\hat{\mathcal{C}}_x^{N-M-1}(Y^0_x)).$$
Finally, since
$|\mathcal{S}_{M-1}(U_{0})|<(1+2\delta)^{-1}|\mathcal{S}_{M-2}(U_{0})|$
we have that 
there exists an interval $\tilde U$ well-inside
$\hat{\mathcal{C}}_x^{N-1}(Y^0_x)$
such that
$$\comp_{V_{0}}F^{-j_2}(D_{\eta''|\mathcal{S}_{M-2}(U)|^{1/3}}(\tilde{U}_{j_2})
\cap\mathbb C_{V_{s'}})
\subset D_{|\mathcal{S}_{M-1}(U)|^{1/3}}(\tilde{U}_{0})\cup
D_{\theta_0}(\hat{\mathcal{C}}^{N-M}_x(Y^0_x)).$$
\endpfclaim

This concludes the proof.
\end{pf}

\begin{prop}\label{prop:sum of squares, puzzle pieces}
There exists $C>0$ such that for any $N\in\mathbb N$, the following holds. Suppose that
$\U_s$ is a puzzle piece of level 0.
Suppose that $\{\U_j\}_{j=0}^s$ is a chain of puzzle pieces without
central cascades of length greater than $N$ and that $\U_0$ intersects
the real line.
Then 
$$\sum_{j=0}^s (\diam(\U_j))^2
\leq N\cdot C\cdot\mu(\mathcal V)^{1/3},$$
where
$\mu(\mathcal V)=\max_{V\subset\mathcal{V}}\diam(V),$ and the
maximum is taken over connected components $V$ of $\mathcal{V}$.
\end{prop}
\begin{pf}
For each puzzle piece
$\U_j$, let $\mathcal{S}_0(U_{j})\supset \mathcal{S}_1(U_j)
\supset\dots \supset\mathcal{S}_{M-2}(U_j)$ be the shielding intervals for 
$U_j:=\U_j\cap\mathbb R.$

By Lemma~\ref{lem:flower lemma} if for any
$\mathcal{S}_{i_0}(U_{j_0}),$ $0\leq j_0\leq s$ and $0\leq i_0\leq M-2$, there
there are at most $N(M-1)+1$ numbers $j\in\{0,1,2, \dots, s\}$
such that there exists $i, 0\leq i\leq M-2$ with
$\mathcal{S}_i(U_{j})=\mathcal{S}_{i_0}(U_{j_0})$, 
then there exists $C_0>0$ depending on $N$ so that
$$\sum_{j=0}^s\sum_{i=0}^{M-2}\diam(|\mathcal{S}_i(U)_j|)^2<C_0.$$
Since each $$\U_j\subset
\cup_{i=0}^{M-3}D_{\lambda^i\theta_0}(\mathcal{S}_i(U_j))\cup
D_{|\mathcal{S}_{M-2}(U_j)|^{1/3}}(\mathcal{S}_{M-2}(U_j)),$$ we have that there exists a
constants $C,C'>0$ such that
$$\sum_{j=0}^s \diam(\U_j)^2\leq 
C \sum_{j=0}^s\sum_{i=0}^{M-2}\diam(\mathcal{S}_i(U_j))^2
\leq C'.$$ 

So suppose that for some pair $j_0,i_0$ with $0\leq j_0\leq s$
and $0\leq i_0\leq M-2$, we have that
there are at least $N(M-1)+1,$ $j\in\{0,1,2,\dots,s\}$ such that there
exists
$i$, $0\leq i\leq M-2$,
depending on $j$, so that
$\mathcal{S}_i(U_{j})=\mathcal{S}_{i_0}(U_{j_0}).$
Let $k+1$ denote the depth of $U_{j_0}$.
Let $j_{1}, 0\leq j_{1}\leq j_0$ be minimal so that
$\mathcal{S}_{i_0}(U_{j_1})$ is a shielding interval,
with the same subindex, for
$U_{j_{1}}$, and for all $U_{t}$ such that
$U_t\subset\tilde{\mathcal{C}}_c^{k}(U_{j_0})$ we have that 
$\mathcal{S}_{i_0}(U_{j_0})=\mathcal{S}_{i_0}(U_t)$ is a shielding
interval for 
$U_{t}$ and
$\hat k(\mathcal{S}_{i_0}(U_{j_0}),U_t)=\hat k(\mathcal{S}_{i_0}(U_{j_0}),U_{j_0}).$
Choose some $x\in U_{j_1}.$
Let $J^0=\tilde{\mathcal{C}}_x^k(Y_x^0)$,
and let $J^{i+1}=\mathcal{L}_{x}(J^i)$ be the principal nest about
$x$.
Relabeling, if necessary, let $j_0\leq s$ be maximal so that
$\mathcal{S}_{i_0}(U_{j_0})$ is a shielding interval for
$U_{j_0},$ $j_0-j_1=:k_0>N,$
so that $\tilde{\mathcal{C}}_x(J_0)=J^{k_0}$,
and the return from $R_{J_1}\rightarrow J^0$ has a critical point.


Notice that if $j_2<j_1$ is maximal so that
$U_{j_2}\subset \tilde{\mathcal{C}}_x(J_0)=J^{k_0}$,
then $U_{j_2}$ is not contained in the same component of the
landing domain to  $J^{k_0}$
as $U_{j_1}$, so that the depth of $U_{j_2}$
is $k+1$.

By Proposition~\ref{prop:quasidisks},
there exists $\theta\in(0,\pi)$ so that
$\J^2\subset D_{\theta}(J^1)$.
Since the chain, $\{U_{j}\}_{j=0}^s$
does not contain central cascades with length greater than $N$,
there exists $k_1\leq N$ and 
$j_0',$  $j_1<j_0'\leq j_0$, so that
we have that $U_{j_0'}\subset\ J^{k_1}\setminus\ J^{k_1+1}$, and
there exists $\delta'>0$ and $\theta'\in(0,\pi)$ depending on $N$,
so that
$\U_{j_0'}\subset D_{\theta'}(\mathcal{L}_{\U_{j_0'}}(\J^{k_1-1}))$.
We have that
$\comp_{U_{j_1}} F^{-(j_0'-j_1)}(\mathcal{L}_{U_{j_0'}}(\J^{k_1-1}))$,
is a pullback of $\mathcal{L}_{U_{j_0'}}(\J^{k_1-1})$ with bounded
degree
through a monotone branch of a return mapping.
Set $\mu(\J^0)=\max\diam(\U)$ where the maximum is taken
over components $\U$ of $\Dom'(\J^0)$.
Hence by Proposition~\ref{prop:good geometry for central cascades}
we have that there exists a constant $C>0$
such that 
$$
\sum_{j=j_1}^{j_0}\diam(\bm{U}_j)^{2}\leq
\sum_{j=j_0'+1}^{j_0}\diam(\bm{U}_j)^{2}
  +\sum_{j=j_1}^{j_0'}\diam(\bm{U}_j)^{2}
$$
$$
\leq N\cdot\mu(\J^0)+
\max_{j_{1}\leq j\leq j_{0}'}\{\diam(\bm U_{j})\}\cdot
\sum_{j=j_{0}}^{j_{0}'}\diam(\bm U_{j})$$
$$\leq  N\cdot\mu(\J^0) + 
C\cdot \max_{j_{1}\leq j\leq j_{0'}}\diam(\bm U_{j})\leq CN\mu(\J^0).$$
This gives us an estimate on the sum of the 
squares of the diameters whenever we pull back though
a monotone branch of a long central cascade in the principal nest.

Let us now combine the arguments when there are no 
central cascades with the argument for controlling the sum
through a central cascade.
Observe that for $j,$ $j_{0}\leq j\leq s$, we have that
$$\sum_{j=j_{0}}^s \diam(\U_j)^2\leq
N(\sum_{j=j_{0}}^s\sum_{i=0}^{M-3}(\diam
D_{\lambda^{i}\pi/2}(\mathcal{S}_i(U_{j})))^2+(\diam
D_{|\mathcal{S}_{M-2}(U_j)|^{1/3}}(\mathcal{S}_{M-2}(U_j)))^2),$$
and this estimate holds if we take the sum on the left 
over segments of the chain that do not pass though
a cascade of central returns in a principle nest with length greater
than $N$.

Let $\mathcal{S}$ denote the set of all shielding Poincar\'e disks, that
is $\mathcal{S}$ is the set of all Poincar\'e disks, 
$$\bm S_{i,j}=
D_{\lambda^{i}\pi/2}(\mathcal{S}_i(U_{j})),\quad 0\leq j\leq s, 0\leq i\leq M-3,$$
together with
$$D_{|\mathcal{S}_{M-2}(U_j)|^{1/3}}(\mathcal{S}_{M-2}(U_j)),\quad 0\leq j\leq s.$$
where the domains in the union are not included with multiplicity.
Thus we have that
there exists a constant $C'>0$ such that
$$\sum_{j=0}^{s} \diam(\U_j)^2\leq
N\sum_{S\in\mathcal S}  \diam(\bm S)^2 + 
C'\sum_{c\in\crit(f)}\sum_{k}\mu(\tilde{\mathcal{C}}^k_c(\bm Y^0_c)),
$$ 
where the last summand is restricted to 
$k$ such that $\hat{\mathcal{C}}^k_c(\bm Y^0_c))$ is at the start of a
long central cascade,
which is bounded from above. In fact, by
Lemma~\ref{lem:flower lemma},
$M$ is universal,
so we have that there exists a constant $C>0$ so that
$$\sum_{j=0}^{s}\diam(\U_j)^2
\leq N\cdot C\cdot\mu(\mathcal V)^{1/3}.$$
\end{pf}

\subsection{The good nest and complex bounds}
\label{sec:goodnest}
We are going to adapt the construction of the {\em good nest}
\label{def:good nest}
of \cite{AKLS} to multicritical complex box mappings.
Suppose that $F\colon\bm{\mathcal{U}}\rightarrow\bm{\mathcal{V}}$
is a complex box mapping given by 
Theorem~\ref{thm:box mapping persistent},
Theorem~\ref{thm:box mapping persistent infinite branches} or
Theorem~\ref{thm:box mapping reluctant} with $\diam(\bm{\mathcal{V}})$
sufficiently small.
For each critical point $c\in\Crit(F),$ 
let $\V^0_c=\comp_c(\bm{\mathcal{V}})$, and
let $\V^0=\cup_{c'\in\Crit(F)}\V^0_{c'}$.
For each $c\in\Crit(F)$,
let $\W_c^0=\mathcal{L}_c(\V^0)$.
The domain $\W_c^0$ is the
return domain containing $c$ to the union of puzzle pieces, 
$\cup_{c\in\Crit(F)}(\V^0_c).$
Assuming that for each critical point we have
constructed $\V^i_c\supset \W^i_c$, we 
construct $\V^{i+1}_c\supset \W^{i+1}_c:$
Assume that $c\in\Crit(F)$. Let $k_c>0$ be
minimal so that for some $c'\in\Crit(F)$,
$F^{k_c}(c)\in\V^i_{c'}\setminus\W^i_{c'}$,
and let $l_c>k_c$ be minimal so that 
for some $c''\in\Crit(F),$
$F^{l_c}(c)\in \W_{c''}^i$.
Let $\V^{i+1}_c=\comp_c F^{-l_c}(\V^i_{c''})$.
Once we have constructed
$\V^{i+1}_c$ for each $c\in\Crit(F)$, let
$\V^{i+1}=\cup_{c\in\Crit(F)}(\V^{i+1}_c),$
$\W^{i+1}_c=\mathcal{L}_c(\V^{i+1})$,
and $\W^{i+1}=\cup_{c\in\Crit(F)}(\W^{i+1}_c)$ .

The chain from
$\V^{i+1}_c$ onto a component
$\V^{i}_{c''}$ of $\V^{i}$ has order bounded by the 
number of critical points of $F$, and 
the critical value $F^{l_c}(c)$ is contained
in $\W^{i}_{c''}$. To see that it has bounded order
 observe that
the chain from $\V^{i+1}_c$ to $F^{k_c}(\V^{i+1}_c)$ 
can intersect each
critical point at most once and the mapping from 
$F^{k_c}(\V^{i+1}_c)$ to $\V^{i}_{c''}$ is obtained by
iterating diffeomorphic branches of the return mapping to 
$\V^{i}$.
By Proposition~\ref{prop:sum of squares, puzzle pieces} for any $n$
and any $c\in\Crit(F)$,
the mapping from $\V^n_c$ up to the top level has bounded dilatation.

Before we prove complex bounds for this nest, let us state the
following easy lemma.

\begin{lem}
\label{lm:pullbackspace}
For each $K$ and $d$ there exists a constant $\kappa>0$
so that for each $K$-quasi-regular map  $g\colon U\to V$  (which is a branched covering onto $V$)
the following holds.  Let $B\subset V$ a disc and let $A$ be a component of $g^{-1}(B)$. 
Then $$\mod(U\setminus A)\ge \kappa  \mod(V\setminus B)).$$
\end{lem}

\begin{rem} If $g\colon U\to V$ has critical values in $V\setminus B$ then there will be several preimages
of $B$. So in the previous lemma we consider the annulus generated by 
considering $U$ minus just one of the components of $g^{-1}(B)$. 
\end{rem}

\begin{prop}\label{prop:good bounds}
There exists $\delta>0$ so that for all $j>0,$ and $c\in\Crit(f)$,
$$\mod(\V^j_c\setminus \W^j_c)>\delta.$$
\end{prop}
\begin{pf}
We are going to use a similar argument to
the one given in the proof of moduli bounds 
for the enhanced nest, compare
Propostition~\ref{prop:modified delta nice}
and \cite{KvS}.

For each $c\in\Crit(F)$ and each $i$, let 
$p_{i,c}$ be so that $F^{p_{i,c}}(\V^i_c)=\V^{i-1}_{c'}$
for some $c'\in\Crit(F).$
Let $12\#\Crit(F)+5\leq M\leq n$. Let $c_1$ be so that
$F^{p_{n,c}}(\V^n_c)=\V^{n-1}_{c_1},$
and assuming that $c_i$ is chosen, let
$c_{i+1}$ be so that $F^{p_{n-i}}(\V^{n-i}_{c_i})=V^{n-(i+1)}_{c_{i+1}}$.
Let $$P_{n,M,c}=p_{n,c}+p_{n-1,c_1}+\dots+p_{n-(M-1),c_{M-1}}.$$

\medskip

\textbf{Step 1: } There exists a critical point $c_0$ and 
$k_0\leq \#\Crit(F)+2$ such that for some $x\in\mathrm{PC}(F)$,
$F^{P_{n,M,c}}(\W^n_c)\subset \mathcal{L}_x(\V^{n-k_0}_{c_0})$.

\medskip

\textit{Claim 1.}
If $c'$ is so that
$F^{p_{n-1,c}}(\V^{n-1}_c)=\V^{n-2}_{c'}$, then
$r(\W^n_c)\geq p_{n,c}\geq p_{n-1,c}+p_{n-2,c'}.$

\medskip
\noindent\textit{Proof of Claim 1.}
Let us start with the first inequality.
Let $r$ be the return time of $c$ to 
$\V^{n}$. Let $k_c>0$ be minimal so that 
$c\in \V^{n-1}\setminus \W^{n-1}$.
Then 
$\V^{n}_c\subset \mathcal{L}_c(\mathcal{L}_{F^{k_c}(c)}(\V^{n-1})),$
and the chain from  
$ \mathcal{L}_c(\mathcal{L}_{F^{k_c}(c)}(\V^n))$ to 
$\mathcal{L}_{F^{k_c}(c)}(\V^n)$ is disjoint, so 
we must have that $r>k_c$. Since
the return time of $F^{k_c}(c)$ to
$\W^{n-1}$ cannot exceed its
return time to $\V^n$ 
(since $\W^{n-1}\supset\V^n$), we have that
$r\geq p_{n,c}.$

To see the second inequality, observe that as in the first
part since
$\V^{n}_c\subset \mathcal{L}_c(\mathcal{L}_{f^{k_c}(c)}(\V^{n-1}))$
and $\W^{n-2}\supset \V^{n-1},$
we have that $p_{n-1,c}\leq k_c$. But now, since
$F^{k_c}(c)\in\V^{n-1}\setminus \W^{n-1}$,
and $\W^{n-1}\subset \W^{n-2}$, arguing similarly, we have that
$p_{n-1,c}+p_{n-2,c'}\leq p_{n,c}$.
\checkmark

The second inequality of the claim gives us the following recursive bound on the
return times:
$$
\left[
\begin{array}{c}
 p_{n-1,c_1}\\ p_{n-1,c_2}\\
\vdots \\ p_{n-1,c_b}\\ p_{n-2, c_1}\\
p_{n-2,c_2} \\ \vdots \\ p_{n-2,c_b}
\end{array}\right] 
\geq \left[
\begin{array}{cc}
I_b & X \\ I_b & 0
\end{array}
\right]
\left[
\begin{array}{c}
 p_{n-2,c_1}\\ p_{n-2,c_2}\\
\vdots \\ p_{n-2,c_b}\\ p_{n-3, c_1}\\
p_{n-3,c_2} \\ \vdots \\ p_{n-3,c_b}
\end{array}
\right],
$$
where $X$ is a $b\times b$ square matrix with a single 
one in each row. The inequality should be interpreted
entry wise.
Let 
$$A_i=\left[
\begin{array}{cc}
I_b & X_i \\ I_b & 0
\end{array}
\right]$$
and set $\mathcal{A}_n=A_n A_{n-1} \dots A_1,$
then one sees easily that
$$\mathcal{A}_n \left[
\begin{array}{c}
1 \\ \vdots \\ 1
\end{array}
\right]=\left[
\begin{array}{c}
T_{n}\\ \vdots \\ T_n \\ T_{n-1}\\ \vdots \\ T_{n-1}
\end{array}
\right],
$$
where $T_n$ is the $n+1$-st Fibonacci number
$$T_0=1,\; T_1=2,\; T_3=3,\;T_4=5,\dots.$$
Observe that this is the sum of the entries in a
row of $\mathcal{A}_n$.

Let $k\in\mathbb N_{\geq 2}$. By the Pigeonhole Principle
we have that 
for some critical point $c_0\in\Crit(F)$
and some $k_0\in\{k,k+1\}$
the mapping
from $\V^{n-1}_c$ to $\V^{n-2}_{c'}$,
occurs after $T_{k_0-1}/(2\#\Crit(F))$ iterates
of the mapping 
$F^{p_{n-k_0},c_0}\colon\V^{n-k_0}_{c_0}\rightarrow \V^{n-k_0-1}_{c_1}.$

Take $k=6\#\Crit(F)+1,$ so that
$k_0\in\{6\#\Crit(F)+1,6\#\Crit(F)+2\}.$
Then we have that 
the mapping from $\V_c^{n-1}$ to 
$\V_{c'}^{n-2}$ occurs after
six iterates of the mapping from
$\V^{k_0}_{c_0}$ to $\V^{k_0-1}_{c_1}$.
Let $t$ be the entry time of $\V_{c}^{n-1}$ to 
$\V^{n-k_0}_{c_0}$, then we have that 
$F^{t}(\V^{n-k_0}_c)\supset \V^{n-k_0}_{c_0}$. 
Let us now define a sequence of critical points
starting with $c_0$: let $c_{i+1}$ be so that
$F^{p_{n-k_0-i,c_i}}(\V^{n-k_0-i}_{c_i})=\V^{n-k_0-(i+1)}_{c_{i+1}}$.
Set
$$P_{n-k_0,M-k_0,c_0}=p_{n-k_0,c_0}+p_{n-k_0-1,c_1}+p_{n-k_0-2,c_{2}}
+\dots + p_{n-M+1,c_{M-1}}.$$
We have that $F^{P_{n-k_0,M-k_0,c_0}}(\V^{n-k_0}_{c_0})$
is a component of $\V^{n-M}$.

Let $k(\V^{n-k_0}_{c_0})$ denote the minimal return time to
$\V^{n-k_0}_{c_0}$ for $x\in \V^{n-k_0}_{c_0},$ and
$r(\W^{n-k_0}_{c_0})$
be the return time of $\W^{n-k_0}_{c_0}$ to $\V^{n-k_0}_{c_0}.$
Observe that the six iterates of the mapping 
from $\V^{n-k_0}_{c_0}$ to $\V^{n-k_0-1}$ occurs after 
$r(\W^{n-k_0}_{c_0})+4k(\V^{n-k_0}_{c_0})$ iterates of $F$.

\medskip
\textit{Claim 2.}
$P_{n-k_0,M-k_0,c_0}\leq r(\W^{n-k_0}_{c_0})+4k(\V^{n-k_0}_{c_0}),$

\medskip
\noindent\textit{Proof of Claim 2.}
Let $\bm P^0=\V^0_{c_0}$, and consider the modified
principal nest over $c_0:$
$$\bm P^0\Supset \bm Q^0\Supset \bm P^1\Supset \bm Q^1 \dots$$
defined as follows. Let $\bm Q^0=\mathcal{L}_c(\bm P^0)$.
Let $k>0$ be minimal so that
$F^{k}(c_0)\in \bm Q^0\setminus\bm P^0$, and let
$l>k$ be minimal so that $F^l(c_0)\in\bm Q^0$.
Let $\bm P^1=\comp_{c_0}F^{-l}(\bm Q^0)$.
Let $\bm Q^1=\mathcal{L}_{c_0}(\bm P^1)$, and continue
inductively. Let
$t_{j}$ be so that $F^{t_j}(\bm P^{j})=\bm P^{j-1}$,
Observe that
$t_{j}\geq t_{j-1}+t_{j-2}$ since
return of the critical point to $\bm Q^{j-1}$ occurs 
after it passes through $\bm P^{j-1}\setminus\bm Q^{j-1}$, and then
through $\bm P^{j-2}\setminus\bm Q^{j-2}$. 
Assume that $j_0$ is minimal so that 
$\bm P^{j_0}\subset \V^{n-k_0}_{c_0}.$ 
The same reasoning shows that,
$r(\W^{n-k_0}_{c_0})\geq t_{j_0}$.
Then we have that
\begin{equation} 
r(\W^{n-k_0}_{c_0})+4k(\V^{n-k_0}_{c_0})\geq t_{j_0}+4t_{j_0-1}\geq t_{j_0}+ t_{j_0-1}+\dots
+t_1+t_0\geq P_{n-k_0,M-k_0,c_0}\quad\checkmark
\label{eqn:claim2}
\end{equation}

Keeping the same notation as in the proof of the Claim 2
we have that $\bm P^{j_0}\subset
\V^{n-k_0}_{c_0}\subset F^t(\V^{n-k_0}_c)$
and $s=r-P_{n-k_0,M-k_0,c_0}-t\geq 0.$
So $$F^s\circ F^{P_{n-k_0,M-k_0,c_0}}\circ F^t(\W^n_c)$$ is a 
component $\V^n_{c'_n}$ of $\V^n.$
Assuming that $c'_j$ is defined, let
$c'_{j-1}$ be so that 
$F^{p_{j,c_j}}(\V^j_{c_j})=\V^{j-1}_{c_{j-1}}.$
We have that 

$$F^{p_{n-k_0+1},c_{n-k_0+1}}\circ F^{p_{n-k_0+2},c_{n-k_0+2}}
\circ\dots\circ F^{p_{n,c_n}}\circ F^s\circ
 F^{P_{n-k_0,M-k_0,c_0}}\circ F^t(\W^n_c)=\V^{n-k_0}_{c_{n-k_0}}.$$
Therefore, for some $x\in\mathrm{PC}(F)$,
$F^{P_{n,M,c}}(\W^n_c)$ is contained in
$\mathcal{L}_x(\V^{n-k_0}_{c_{n-k_0}}).$

\medskip

\textbf{Step 2:}
Let $i\in\{n-k_0, n-k_0-1,\dots n-M\},$
let $\nu$ be the first entry time of $B=F^{P_{n,M,c}}(\W^n_c)$
to some $\W^i_c.$ 
Let $D_i=\comp_B F^{-\nu}(\V^i_c)$.
The mapping 
$F^{\nu}\colon D_i\rightarrow\V^{i}_c$ has bounded degree,
since after the first entry to $\V^i$ it is obtained by 
iterating the non-critical branches of the return mapping
to $\V^i$.
Let $C_i=\comp_B F^{-\nu}(\W^i_c)$ and set 
$\bm A_i=D_i\setminus C_i.$ Since for all $j,$
$\V^j\supset \W^j\supset \V^{j+1},$ we have that
the $\bm A_i$ are disjoint.

\medskip

\textbf{Step 3:} Let $\mu_{j}$ denote 
$\min(\V^j_c\setminus\W^j_c)$, where the minimum 
is taken over the critical points $c$ of $F$.
Thus we have that 
$$\mod(F^{P_{n,M,c}}(\V^n_c)\setminus F^{P_{n,M,c}}(\W^n_c))\geq
(K_1/K_2)(\mu_{n-M}+\mu_{n-M+1}+\dots\mu_{n-k_0}).$$

\medskip

\textbf{Step 4:} 
It is obvious that the degree of
$F^{P_{n,M,c}}|_{\V^{n}_c}$ is bounded by some constant $D$
that depends on $M$, so we only show that
the degree of $F^{P_{n,M,c}}|_{\W^n_c}$ is bounded by some 
constant $d$ that does not depend on $M$.
From Step 1, we have that there exists $k_0\in\{\#\Crit(F)+1,
6\#\Crit(F)+2\}$ so that
$F^{P_{n,M,c}}(\W^n_c)\subset\mathcal{L}_x(\V^{n-k_0}_{c'}).$
Since $M-k_0\geq 6\#\Crit(F)+2,$ by Claim 2 of Step 1, there exists
 $k_1\in\{\#\Crit(F)+1, 6\#\Crit(F)+2\},$ so that the following holds.
Let us decompose the mapping $F^{P_{n,M,c}}|_{\W^n_c}$
as $F^{P_{n-(k_0+k_1,M-(k_0+k_1),c_1}}\circ
F^{P_{n,k_0+k_1,c}}|_{\W^n_c},$
where $c_1\in\Crit(F),$
is such that $F^{P_{n,k_0+k_1,c}}(\V^n_c)=\V^{n-(k_0+k_1)}_{c_1}$. 
The degree of  $ F^{P_{n,k_0+k_1,c}}|_{\W^n_c}$ is bounded by a constant
which depends only on $\#\Crit(F)$. 
Since $M\geq 12\#\Crit(F)+5>k_0+k_1,$
$F^{P_{n,k_0+k_1,c}}(\W^n_c)$ belongs to some other component
$\mathcal{L}_y(\V^{n-k_0}_c)$. 
Let $r>1$ be such that 
$F^r(\mathcal{L}_y(\V^{n-k_0}_c))=\V^{n-k_0}_c.$
We have that
$r>P_{n-(k_0+k_1),M-(k_0+k_1),c},$ therefore the degree of the mapping
$F^{P_{n-(k_0+k_1),M-(k_0+k_1),c''}}|_{F^{P_{n,k_0+k_1,c}}(\W^n_c)}$ is less than or equal to the degree of
some landing map $R_{\V^{n-k_0}_c}|_{\mathcal{L}_y(\V^{n-k_0}_c)}$,
which is bounded by some constant $d$ which depends only on $b$.

\medskip

\textbf{Step 5:}
Now apply Lemma~\ref{lem:small distortion of thin annuli}
exactly as in the proof of 
Proposition~\ref{prop:modified delta nice}.
Notice that by Lemma~\ref{lm:pullbackspace} there exists $\kappa>0$
such that
$\mod(\V_c^n\setminus \W^n_c) > \kappa\mu_{n-1}$. 
\end{pf}

By Lemma~\ref{lm:pullbackspace}, we have the following corollary:

\begin{cor}\label{cor:Wnice}
There exists $\delta_0>0$ so that each 
$\W^n$ is $\delta_0$-nice. In fact we have that
for each component $\U$ of
$(R_{\V^n}|_{\W^n})^{-1}(\Dom'(\W^n))$ if $\W^n_c$ is the component
of $\W^n$ that contains $\U$ then
$\mod(\W^n_c\setminus\U)\geq\delta_0$.
\end{cor}

Let us now show that the puzzle pieces $\V^n_c$ from the good nest
have $\delta$-bounded geometry. We will require the following
lemmas from \cite{KvS}.

\begin{lem}\label{lem:KvS bg lem 1}
Let $\U$ be a domain that has $\rho$-bounded geometry with
respect to some point $x$ and let $\bm A\subset \U$ be a domain containing
$x$. Then $\U$ has $K\rho$-bounded geometry with respect to 
all $y\in \bm A,$ where $K$ depends only on $\mod(\U\setminus \bm A)$.
\end{lem}

\begin{lem}\label{lem:KvS bg lem 2}
Let $f\colon \U\rightarrow \V$ be a $(1+\eta)$-qc branched covering map
$\bm B\subset \V$ and  $\bm A$ a connected component of $f^{-1}(\bm B)$.
Then if $B$ has $\rho$-bounded geometry with respect to some point
$y$ in $\bm B$ then $\bm A$ has $K\rho$-bounded geometry with respect to 
$x\in \bm A$, where $K$ depends only on $\mod(\V\setminus \bm B)$, the
degree of $f$ and $\eta>0.$
\end{lem}

Let us show that each $\V_c^j$ has bounded geometry at $c\in\Crit(F)$. The proof is the
same as the proof of \cite[Proposition 11.1]{KvS}.
\begin{prop}\label{prop:goodbg}
The domains $\V^j_c$ have bounded geometry at $c$;
that is, there exists $\eta>0$ so that for each $j$,
$\V^j_c \supset D_{\eta\diam(\V^j_c)}(c)$.
\end{prop}

\begin{pf}
Let $\eta_{j,c}>0$ be so that $\V^j_c$ has
$\eta_{j,c}$-bounded geometry at
$c$. 
Let $l_c$ and $c'\in\Crit(f)$ 
be so that $f^{l_c}\colon\V_c^j\rightarrow \V_{c'}^{j-1}$
is onto. 
Since $f^{l_c}(c)\in\W^{j-1}_{c'}$ 
and $\mod(\V^{j-1}_{c'}\setminus\W^{j-1}_{c'})>\delta,$
by Lemma~\ref{lem:KvS bg lem 1} there exists $K_1$
so that $V^{j-1}_{c'}$
has $K_1\eta_{j-1, c'}$-bounded geometry at $f^{l_c}(c)$. 
By Lemma~\ref{lem:KvS bg lem 2}
 there exists $K_2$ so that
$f(\V^{j}_c)$ has $K_1K_2\eta_{j-1,c'}$-bounded geometry at $f(c)$. Thus
we have that $\V^{j}_c$ has $\sqrt{K_1K_2\eta_{j-1,c'}}$
-bounded geometry at
$c$. Thus the $\eta_{j,c}$ are bounded away from 0.
\end{pf}

\section{Dynamics away from critical points: Touching box mappings}
\label{sec:touching box mappings}
	The purpose of this section is to deal with the set of points
	which always stay away from critical points. More precisely,
	given an open set $\mathcal I$ which contains all critical points 
	and the immediate basin of all periodic attractors, 
 	define 
	$$E(\mathcal I)=\{x;  f^n(x)\notin \mathcal{I}\mbox{ for all }n\ge 0\}.\label{notation:escaping}$$
	In this section we will prove the following

	\begin{thm}[QC-rigidity of points staying away from the critical points]
	\label{thm:qc rigidity away from critical points}
		Let $f,\tilde f$ be conjugate maps as in Theorem~\ref{thm:main}.
		Let $I$ be an open set which contains all critical points 
		and the immediate basin of all periodic attractors.
		Then  there exists a quasiconformal mapping $h$ so that
		$h(E(\mathcal I))=E(\tilde {\mathcal I})$ and so that $h|_{E(\mathcal I)}$ is a
		conjugacy between $f$ and $\tilde f$.
	\end{thm}

	The proof of this theorem will be given at the end of this section,
	and  is complicated by the
	fact that we do not want to exclude the presence of parabolic periodic points.
	We will construct suitable ``touching box mappings"
	extending certain certain iterates of the real mapping along diffeomorphic branches.
	These touching box mappings will be different from complex box mappings
	as the boundaries of their domains and the boundaries of their ranges will intersect.
	As mentioned, the reason we need these touching box mappings is because
	we need to deal with the part of the dynamics which does not 
	come close to critical points, and since  we will eventually need to cover
	the entire space by suitable complex domains.

	In Subsection~\ref{subsec:global touching box mappings} we shall give a precise definition of the notion 
	of ``touching box mapping", while in the next subsection we will
	treat the local situation near (not necessarily hyperbolic) repelling periodic points.

	\subsection{Touching box mappings and petals at (possibly parabolic) periodic  points}
		\label{subsec:local touching box mappings}
	To simplify the notation, if $f\in\mathcal{C}$, we will identify $f$ with it's asymptotically holomorphic 
	extension of order 3.
		\begin{prop}[Petals and local QC conjugacies within petals]
		\label{prop:petals}
			Assume that $f\in\mathcal{C}$, 
			that $p$ is a periodic point with $f^{s}(p)=p$ and 
			that there exists a one-sided neighbourhood $[p,q_{0}]$ of $p$ with $Df^{s}(x)\geq 1$. 
			Then there exists $\theta_{0}\in(0,\pi/2)$ and $\delta>0$ 
			so that for each $\theta\in(0,\theta_{0})$ and 
			for any interval $J=[p,q]\subset[p,q_{0}]$ with $|J|<\delta$, 
			such that $f^{s}|_J$ is a diffeomorphism and $J'=f^{s}(J)\supset J$, the following holds.
\begin{enumerate}
\item $\U:=\Comp_{p}f^{-s}(D_{\pi-\theta}(J'))$ 
is contained in $\V:=D_{\pi-\theta}(J')$ 
and they intersect only in the periodic point.
\item $\Comp_{p}f^{-ks}(D_{\pi-\theta}(J'))$ is contained in a set of the form 
$D_{\pi-\theta}(J_{k})$ where $J_{k}\subset J$ and 
$|J_{k}|\rightarrow 0$ as $k\rightarrow \infty$. 
\end{enumerate}
Moreover, assume that $\tilde{f}$ is another mapping 
which is topologically conjugate to $f$ 
with corresponding periodic point $\tilde{p}$ 
and intervals $\tilde{J},\tilde{J}'$, so that either 
$Df^{s}(p)$ and $D\tilde{f}^{s}(p)$ are both greater than one 
or both equal to one, 
and $J$ and $J'$ are so small that 
$\tilde{\U}:=\Comp_{\tilde{p}}\tilde{f}^{-s}(\tilde{\V})$ 
is contained in $\widetilde{\V}:=D_{\pi-\theta}(\tilde{J}')$. 
\begin{itemize}
\item[(3)] Then there exists a $K$-qc homeomorphism 
$H\colon \C\rightarrow\C$ that maps 
$\partial \U\cup\partial \V$ to $\partial\widetilde{\U}\cup\partial\widetilde{\V}$ 
such that $\tilde{f}^{s}\circ H=H\circ f^{s}$ on $\U$. 
Here $K$ depends on $|J|/|J'|,|\tilde{J}|/|\tilde{J}'|$ and $\theta$.  
\end{itemize}
In other words, $H$ maps $\V\setminus \U$ to $\widetilde{\V}\setminus 
\widetilde{\U}$
and is a conjugacy on $\clos(\U)$.
		\end{prop}
		\begin{figure}[htp] \hfil
		\beginpicture
		\dimen0=0.3cm
		\setcoordinatesystem units <\dimen0,\dimen0>
		\setplotarea x from -9 to 30, y from -5 to 4
		\setlinear
		\plot -8 0 8 0 /
		\circulararc 120 degrees from 5 0  center at  0 -3 
		\circulararc -120 degrees from 5 0  center at  0 3 
		\put {$\partial V$} at 5 2 
		\circulararc 110 degrees from 1 0  center at  -2 -2 
		\circulararc -110 degrees from  1 0 center at  -2 2 
		\put {$\partial U$} at 1.5 1 
		\setdots <1 mm>
		\plot -5 0 0 5 /
		\plot -5 0 -2 6 /
		\plot -5 0 5 3 /
		\plot -5 0 0 -5 /
		\plot -5 0 -2 -6 /
		\plot -5 0 5 -3 /
		\setsolid
		\put {$z\mapsto 1/z$} at 10 -2  
		\arrow <3mm> [0.2,0.67] from 7 -3.5 to 13 -3.5
		\put {$\bullet$} at 13 0  
		\plot 13 0 30 0 /
		\plot 17 0 28 4 /
		\plot 17 0 28 -4 /
		\plot 19 0 32.7 5 /
		\plot 19 0 32.7 -5 /
		\setdots <1mm>
		\plot 13 0 32 4 /
		\plot 13 0 32 -4 /
		\plot 13 0 32 2 /
		\plot 13 0 32 -2 /
		\put {$\phi(\partial \V)=l_{\pm}$} at 24 5 
		\put {$\phi(\partial \U)$} at 34 6 
		\endpicture
		\caption{Part of the touching box mapping 
		at a fixed point $p$ of $f^s$, with domain $\U:=\comp_{p} f^{-s}(D_{\pi-\theta}(J'))$ is contained in $\V:=D_{\pi-\theta}(J')$.
		The mapping $z\mapsto 1/z$ sends the curves $\partial \U,\partial \V$  to the wedge-shaped curves
		on the right. \label{fig:touching}}
		\end{figure}

We note that by \cite{MMS}, see also Theorem IV.B in \cite{dMvS},
each $C^3$ mapping of the interval or circle has at
most a finite number of parabolic or attracting periodic points. 

We should note that in the parabolic case, we construct a qc map
which maps the parabolic petal associated to $f(z)=z+z^{n+1}+\dots$
to that of $\tilde f(z)=z+z^{\tilde n}+\dots$ where $n$ and $\tilde n$
are not necessarily the same.  If $n\ne \tilde n$, this qc mapping can definitely 
{\bf not} be extended to a global conjugacy. Since $n$ and $\tilde n$ are not necessarily
the same, we cannot argue exactly as in for example
\cite[Chapter 10]{Milnor}, but have to vary the argument
a little. We find it also useful to 
ensure that $\V$ is a Poincar\'e domain, so the shape of the 'petals'
differs marginally from those usually constructed.
		
That $\U$ is contained in $\V$ does not simply follow from
Proposition~\ref{prop:Almost Schwarz Inclusion}, 
because no loss of angle is allowed
at the points where $\U$ and $\V$ touch. 

\begin{pf}
			It is enough to consider the case that $p$ is a fixed point.
			
			\medskip 
			\noindent
			{\bf Case 1:} Assume that $Df(p)>1$. Recall that 
			$$\frac{\bar{\partial}f}{y^2}\rightarrow 0\ \mathrm{as}\ y\rightarrow0,$$
			and we can express 
			$$f(z)=\lambda z + a z^{2}+O(z^3)+R(x,y),$$
			where $R(x,y)$ is $C^3$ and $D^jR(0,0)=0$ for $j=0,1,2,3$,
			and $R(x,y)=o(\|(x,y)\|^3)$.
			Note that we are abusing notation and denoting both
			the interval mapping and its asymptotically conformal extension by $f$.
			
			Consider the coordinate transformation
			$w:=\phi(z)=1/z$. Since $\partial \V$ consists of two pieces of circles, 
			$\phi(\partial \V)$ consists of two rays of lines $l_\pm$ in the $w$-plane,
			which can be parametrized by 
			$$x=\alpha |y|+x_0\mbox{ where }w=x+iy,$$
			where $\alpha>0$ depends on $\theta$ 
			 (with $\alpha\to \infty$ as $\theta\to 0$)  and where $x_0=1/|J'|$. 
			The mapping $z\mapsto f(z)=\lambda z+ az^2 +O(z^3)$ (where $a\in \R$) corresponds in $\omega$-coordinates
			to $w\mapsto g(w)=(1/\lambda)w-(a/\lambda^2)+O(1/w)$ whereas its preimage
			corresponds to $w\mapsto g^{-1}(w)= \lambda w+(a/\lambda)+O(1/w)$. 
			Since $(\alpha|y|+x_0+i\cdot y)\in l_{\pm}=\partial V$ is mapped by $g^{-1}$ to 
			$(\lambda(\alpha|y|+x_0)+(a/\lambda)+O(1/w))+ i (\lambda y+O(1/w))\in \partial U$,
			we get that $\partial \phi(U)$ is parametrized  by
			$(\alpha|y'|+\lambda x_0+(a/\lambda)+O(1/w)) + i \cdot y')$ (where $y'=\lambda y+O(1/w)$).
			So $\phi(\partial U)$ is also the image of $\phi(\partial V)$ under the map
			$w=(x+i\cdot y)\mapsto (x+(\lambda-1)x_0+a/\lambda +O(1/w)+ i\cdot (y+O(1/w))$.
			That is, $\phi(\partial U)$  is up
			 to  small error a translation of $\phi(\partial V)$
			 to the right (since $x_0<\lambda x_0 + a/\lambda+O(1/w)$ when $x_0$ is large).
			So $\phi(V\setminus U)$ corresponds to the strip drawn on the right in
			Figure~\ref{fig:touching}. The 2nd assertion follows immediately.
			
Assume that $\tilde f$ is a similar map, with corresponding  
$\widetilde{\U},\widetilde{\V}$ and $\tilde g$.
So $\partial \widetilde{\V}$ can also be parametrized by $x=\alpha |y|+\tilde x_0$
			where $w=x+i\cdot y$ and $\tilde x_0=1/|\tilde J|$.  
			Let us now construct a qc mapping $H\colon  \C\to \C$ in the $w$ plane 
			so that
			$\hat H(\phi(\U))=\widetilde{\U}$, $\hat
                        H(\phi(\V))=\phi(\widetilde{\V})$ and 
$\hat H\circ g= \tilde g\circ \hat H$ on $\phi(\U)$.
			To do this, we first define a foliation as follows.
			Take $w=(x+i\cdot y)\in l_+$, i.e. $x=\alpha |y|+x_0$ with $y>0$,
			and consider the line segment $L(w)$ connecting $w$  to $g^{-1}(w)$.
			Let us show that all these line segments are disjoint.
			For this observe that for $w\in l_+$, \,\,
			the vectors $T_w(l_+)$ and $Dg^{-1}T_w(l_+)$ are both
			transversal to and both have the same orientation with respect to the vector
			$g^{-1}(w)-w$.
			(To see this, note that $Dg^{-1}T_w(l_+)=
			Dg^{-1}_w(\alpha+i)=
			(\lambda +O(1/w^2))(\alpha+i)$,
			that $g^{-1}(w)-w= (\lambda-1)w+(a/\lambda)+O(1/w)$ and that $w=x+iy$ is of the
			form $x=\alpha |y|+x_0$.) 
Since the same holds for $w\in l_-$, we obtain a foliation
on $\phi(\V\setminus \U)$ consisting of the line segments $L(w)$, $w\in l_\pm$. 
			For $w\in \phi(\V\setminus \U)$ define $L(w)$ to be the line segments which contains $w$.
			Similarly, we have a foliation on 
			$\phi(\widetilde{\V}\setminus \widetilde{\U})$
                        consisting of the line segments 
$\tilde L(\tilde w)$, $\tilde w\in \tilde l_\pm$
			corresponding to $\tilde f$. Note that $\tilde l_\pm=l_\pm+(\tilde x_0-x_0)$.

			Choose $\hat H\colon  l_\pm \to \tilde l_\pm$ so that it maps
			$(\alpha |t|+x_0+i\cdot t)\in l_\pm$ to 
			$(\alpha | \tilde t|+\tilde x_0+i\cdot \tilde t)\in \tilde l_\pm$ where $\tilde t=t^{\rho}$ and 
			$\rho=\log\tilde\lambda/\log\lambda$. Then define the mapping $\hat H$ mapping the leaves
			of the foliation to the leaves of foliation, interpolating linearly as before.
			Then,  taking $w\in l_{\pm}$ and $\tilde w=\hat H(w)\in \tilde l_\pm$, we have
			$||\tilde g^{-1}(\tilde w)-\tilde w||/||g^{-1}(w)-w||=
			|| (\tilde \lambda-1)\tilde w+(\tilde a/\tilde \lambda)+O(1/\tilde w)||/
			|| (\lambda-1)w+(a/\lambda)+O(1/w)||\approx ||w||^{\rho-1}$. At the same, taking
$w,w'\in l_{\pm}$ nearby we have
$||\hat H(w)-\hat H(w')||/||w-w'||\approx ||w||^{\rho-1}$.
It follows that $\hat H$ is a $K$-quasiconformal map
from  $\phi(\V\setminus \U)$ to $\phi(\widetilde{\V}\setminus \widetilde{\U})$.
Next extend $\hat H$  to the right of $\phi(\partial \V)$ to the right
of 
$\phi(\partial \widetilde{\V})$
by $\hat H(w)=\tilde g^{-n}\circ \hat H \circ g^n(w)$ where $n=n(w)\ge 0$ is minimal
so that $g^n(w)\in \phi(\V\setminus \U)$. 
Next choose $\hat H$ so that for each $0\le x\le x_0$
the point $(\alpha |t|+x+i\cdot t)\in l_\pm$ is mapped to 
$(\alpha | \tilde t|+x(\tilde x_0/x_0)+i\cdot \tilde t)\in \tilde
l_\pm$
where $\tilde t=t^{\rho}$.
Then $\hat H(\alpha |t| + i\cdot t)= (\alpha |t^\rho|+ i\cdot t^\rho)$ for $t\ge 0$.
One can extend $\hat H$ in various ways to the left of the sector
$(\alpha|t| + i\cdot t)$, $t\in \R$, for example defining $\hat H$ for each vector $v$
not in the sector between $(\alpha -i)$ and $(\alpha+i)$, as follows
$\hat H(t\cdot v)= t^\rho\cdot v$.
In this way we get a globally defined qc mapping $\hat H$.
Therefore $H\colon  \C\to \C$ defined by
$H(z)=\phi^{-1}\circ  \hat H \circ \phi(z)$ is 
again a quasiconformal map
such that $H(\U)=\widetilde{\U}$, $H(\V)=\widetilde{\V}$ and 
$\tilde f\circ H=H\circ f$ on $\U$.
							
Observe that the mapping $f|_U$ is monotone,
and $U$ contains no parabolic points.
So that Lemma~\ref{lem:sum of lengths}
implies that 
$$\sum_{n=0}^{\infty}|(f|_U)^{-n}(U)|^2<C|U|^{1/2}.$$
Since we have no loss of angle at periodic points,
and $f\in\mathcal{C}$ this implies that the
dilatations of the
$g^{n}\colon g^{-n}(\phi(\V\setminus \U))
\rightarrow \phi(\V\setminus \U)$ are summable,
and indeed bounded from above by $C|U|^{1/2}.$

\medskip
\noindent 
{\bf Case 2:} Let us now consider the parabolic case, i.e.
$f(z)=z+az^{n+1} + \dots$ with $a>0$.
In this case define the $w=\phi(z)=1/z^n$.
Then $f$ corresponds to  $g(w)= w -(na) + O(1/\sqrt[n]{|w|})$
which has inverse $g^{-1}(w)=w+(na)+O(1/\sqrt[n]{|w|})$. 
The set $\partial \V$ corresponds under $\phi$ to the $n$-th power
of  the half-lines $l_{\pm}$ from Case 1 (note that $\phi(z)=(1/z)^n$).
One can parametrize $l_{\pm}$ by  $\alpha|y|+x_0 + i\cdot y\in \C$
where $y\in \R$ and where $\alpha\to \infty$ when $\theta\to 0$.
The $n$-th power of these two half-lines is equal to 
$(\alpha|y|+x_0+i\cdot y)= y^n(\pm \alpha +x_0/y+ i)^n$,
which for $|y|$ large  is again almost a line.
When $\alpha>0$ is sufficiently large, $\phi(\partial \V)$ is 
contained in $\{w=x+i\cdot y\, ; \, x,y\in \R \mbox{ and }\, x\ge A\}$
with $A$ large.
Moreover, $\phi(\partial \V)$ can be parametrized by 
$\beta(t)+i\cdot t+x_0$ where $\beta$ is a smooth function with
$\beta'(s)\to \pm c$ as $s\to \pm \infty$
with $\alpha\in (0,\infty)$
depends on $\theta$ and $\alpha\to \infty$ when $\theta \to 0$) 
and where $x_0=(1/|J'|)^{n}$.
The set $g^{-1}(\phi(\partial \V))$ is to the right of $\phi(\partial \V)$ since $a>0$.
Hence $\phi(\V\setminus \U)$ corresponds again to a similar strip as before, 
see the right part of Figure~\ref{fig:touching}, and that $\U\subset \V$ with only 
a common point in $p$.
			
Now assume that $\tilde f$ is a similar map, with corresponding  
$\widetilde{\U},\widetilde{\V}$ and $\tilde g$, but not necessarily with the same
order of contact. That is $\tilde f(z)=z+\tilde az^{\tilde n+1} + \dots$
and $\tilde \phi(z)=1/z^{\tilde n}$. 
Let us again define a mapping $\hat H$ similarly as before.
For each $w\in \phi(\partial \V)$
take the line segment $L(w)$ connecting $w$ to $g^{-1}(w)\in \phi( \partial \U)$. 
As before, this defines a foliation in $\phi(\V\setminus \U)$. 
Next define $\hat H\colon  \phi(\partial \V)\to \tilde \phi(\partial \widetilde{\V})$
by mapping $w=\beta(t)+i\cdot t+(1/|J|)^{n}$ to 
$\tilde w=\tilde \beta(t)+i\cdot t+(1/|\tilde J|)^{\tilde n}$.
Next extend $\hat H\colon  \phi(\partial \V)\to \phi(\partial
\widetilde{\V})$
continuously to 
$\hat H\colon  \phi(\V\setminus \U)\to \phi(\widetilde{\V}\setminus \widetilde{\U})$ so that for
each $w\in l_{\pm}$, the line segment $L(w)$ is mapped
to the line segment $\tilde L(w)$ by linear interpolation.
Note that in this case 
$||\tilde g^{-1}(\tilde w)-\tilde w||$ and $||g^{-1}(w)-w||$ 
(with $w\in \phi(\partial \V)$ and $\tilde w\in \tilde  \phi(\partial
\widetilde{\V})$)
are both bounded. 
As defined, we have that $\hat H$ is $K$-qc for some $K\geq 1$.
To see that $\hat{H}$ is $K'$-qc. we need only observe that
backward orbits converge to $0$ at a rate
$$|z_{-n}'|\asymp \Big(\frac{1}{\frac{1}{|z_0|^d}+d n}\Big)^{1/d}$$
and since $f$ is $C^{d+2}$, and the asymptotically conformal extension
is of order
$d+2$ in a neighbourhood of 0 in the complex plane.
So that we can bound that qc constant of $\hat H$
from above by $$K+\sum_{n=0}^\infty \Big(\frac{1}{\frac{1}{|z_0|^d}+d n}\Big)^{(d+1)/d},$$ 
which converges.
It follows that 
$\hat H\colon  \phi(\V\setminus \U)\to \phi(\widetilde{\V}\setminus \widetilde{\U})$
is $K$-quasiconformal.
Next extend $\hat H$  to the right of 
$\phi(\partial \V)$ to the right of $\phi(\partial \widetilde{\V})$
by $\hat H(w)=\tilde g^{-n}\circ \hat H \circ g^n(w)$
where $n=n(w)\ge 0$ is minimal
so that $g^n(w)\in \phi(V\setminus \U)$. 
On the strip between the sector $(c|t|+i\cdot t)$, $t\in \R$ and  $\partial \V$ define
$\hat H$ so it maps to the strip between the sector
$(\alpha|t|+i\cdot t)$, $t\in \R$ $\partial \widetilde{\V}$ by linearly
interpolating horizontal lines.
In particular, $\hat H(\alpha|t|+i\cdot t)=\alpha|t|+i\cdot t$ for all $t\in \R$
and $\hat H$ is $K$-quasiconformal on the sector to the right of these
rays
through $0$.
Since $\phi$ is not univalent anyhow, we will not extend $\hat H$ globally.
			
Let $\beta=\arg(\alpha+i)\in (0,\pi/2)$.
Then $\phi$ maps the sector $\{z\in \C; |\arg(z)|< \beta/n\}$ to the sectorto right of $(\alpha|t|+i\cdot t)$, $t\in \R$. Similarly, 
$\tilde \phi$ maps the sector $\{z\in \C; |\arg(z)|< \beta/\tilde n\}$ to the sector
to right of $(\alpha|t|+i\cdot t)$, $t\in \R$.
So we obtain  that $H(z)=\tilde \phi^{-1}\circ  \hat H \circ \phi(z)$ 
is a well-defined $K$-qc mapping from the sector $\{z\in \C; |\arg(z)|< \beta/n\}$
to the sector $\{z\in \C; |\arg(z)|< \beta/\tilde n\}$.
Note that $H(te^{\pm i\cdot \beta/n})=t^{\tilde n/n}e^{\pm i\cdot \beta/\tilde n}$. 
Finally, one can extend $H$ so that 
$H(te^{\pm i\cdot \gamma})=t^{\tilde n/n}e^{\pm i\cdot \tilde \gamma}$
where $[\beta/n,2\pi-\beta/n] \ni \gamma \mapsto
\tilde \gamma\in [\beta/\tilde n,2\pi-\beta/\tilde n]$
is a linear map.
\end{pf}
		
\subsubsection{Touching box mappings and petals in a complex
                  neighbourhood of (possibly parabolic) periodic
                  points}
\label{subsec:nbhdconjpara}
This subsection is an addendum to Proposition \ref{prop:petals}. We will show how to 
extend conjugacies to neighbourhoods of periodic points in the hyperbolic repelling
case or in the parabolic case when the multiplicities of the parabolic points are the same.
\begin{prop}[Petals and local QC conjugacies within petals] \label{prop:petals2}
Let $f,\tilde f$  be topologically conjugate maps which satisfy the assumptions
of the previous proposition and let $p,\tilde p$ be corresponding
periodic points of, say,  period $s$.
Moreover, assume that either 
\begin{enumerate}
\item $p$ and $\tilde p$ are both hyperbolic repelling periodic points  of $f$, $\tilde f$ 
respectively, or
\item $p$ and $\tilde p$ are both parabolic  and 
$$f^s(x)=p+\lambda (x-p)+a(x-p)^{d+1}+o(|x-p|^{d+1})$$
and 
$$\tilde f^s(x)=\tilde p+\lambda (x-\tilde p)+\tilde a(x-\tilde p)^{d+1}+o(|x-\tilde p|^{d+1})$$
where $\lambda=\tilde \lambda\in \{-1,1\}$, $d=\tilde d$ and $a,\tilde a$ have the same sign.
\end{enumerate}
Then one can obtain a quasiconformal conjugacy $H$ on a complex
neighbourhood of $p$.
\end{prop}

		\begin{figure}[htp] \hfil
		\beginpicture
		\dimen0=0.3cm
		\setcoordinatesystem units <\dimen0,\dimen0> 
		\setplotarea x from -9 to 30, y from -7 to 9
		\setlinear
		\plot -16 0 8 0 /
		\circulararc 120 degrees from 5 0  center at  0 -3 
		\circulararc 110 degrees from 1 0  center at  -2 -2 
		\circulararc -120 degrees from -15 0  center at  -10 -3 
		\circulararc -110 degrees from -11 0  center at  -8 -2

		\put {$\partial \V$} at -5 9 
		\put {$\partial \U$} at -5 4.5 
		\put {$p$} at -5.8 -0.7
		\plot -5 -7 -5 9 /
		\circulararc 125 degrees from  -5 6  center at  -3.5 3 
		\circulararc -125 degrees from -5 6  center at  -6.5 3
 		\put {$l^+_r$} at 0 1.8 
		\put {$l^+_l$} at -10 1.8 
		\circulararc 130 degrees from  -5 8  center at  -3.5 4 
		\circulararc -130 degrees from -5 8  center at  -6.5 4 
		\put {$m^+_r$} at 6.4 1.3 
		\put {$m^+_l$} at -15.7 1.3 
		\setsolid
		
		\put {$z\mapsto 1/z^2$} at 10 -2.3  
		\arrow <3mm> [0.2,0.67] from 7 -1 to 13 -1
		\setcoordinatesystem units <\dimen0,\dimen0> point at -5 0 
		\setplotarea x from -9 to 30, y from -7 to 9

		\put {$\bullet$} at 30 0  
		\plot 13 0 30 0 /
		\plot 28 0 11 6 /
		\plot 28 0 11 -6 /
		\plot 24 0 7 6 /
		\plot 24 0 7 -6 /
		
		\plot 30 1.2    6 5.4 /
		\plot 30 0.6 6 4.8 /
		\plot 30 -1.2    6 -5.4 /
		\plot 30 -0.4 6 -4.8 /
		\put {*} at 20 2.7 
		\put {*} at 14.5 3  
		\put {$\phi(\partial \V)=n^+$} at 14 6 
		\put {$\phi(\partial \V)=n^-$} at 14 -6 
		\put {$\phi(\partial \U)$} at 7 6 
		\put {$\phi(l_r^+)$} at 28 2.5 
		\put {$\phi(l_l^+)$} at 28 -2.5 
		\put {$\phi(m_r^+)$} at 10 3
		\put {$\phi(m_l^+)$} at 10 -3
		\endpicture
\caption{Part of the touching box mapping at a fixed point $p$ of $f^s$, with domain 
$\U:=\comp_{p} f^{-s}(D_{\pi-\theta}(J'))$ is contained in $\V:=D_{\pi-\theta}(J')$. 
In this case $\U$ maps to $\V$. The pre-image  under $z\mapsto 1/z^2$  of the wedge-shaped
region on the right, bounded by half-lines $\phi(m^\pm)$ and $n^\pm$,
have preimages which form part of the boundary of the horizontal and vertical petal. 
The point $0$ in the $w$ plane is marked by a $\bullet$  in the right panel. \label{fig:touching2}}
\end{figure}
		
\begin{pf}
The case when $p$ is hyperbolic and repelling was already considered in the proof of the 
previous proposition. 
Let us next consider the case that $s=1$, $\lambda=1$, $a>0$ and $d=2$. For convenience assume  $p=0$. Then $f$ has a fixed point at $0$ which on $\R$ repels on both sides of $p$. Therefore there exists two repelling  petals along the real line. Let us assume that on these two repelling petals the conjugacy is
already defined.  Since $d=2$, there exists one attracting petal between the repelling  
petals in the vertical direction at $p$. Indeed,  consider $w:=\phi(z)=1/z^2$ and let $g$ be 
the mapping corresponding to $f$ in the $w$-coordinates (taking the branch corresponding
to the positive half-plane in the left-panel of Figure~\ref{fig:touching2}.). Take two  half-lines 
$n^\pm$ through a point in the negative real axis in the $w$-plane as in the figure and consider $g(n^+)$. Then $\phi^{-1}(n^\pm)$ and $\phi^{-1}(g(n^\pm))$ bound two regions $\U$ and  $\V$ (only the components in the upper half-plane are shown in Figure~\ref{fig:touching2}). 
The part of the petal to the right of $p$ above the real axis is bounded by curves $l$ and $m$ 
(and the real axis).  These curves correspond in the $w$-coordinates to 
half-lines $\phi(l^\pm)$ and its image under $g$. 
As in the proof of the previous proposition $g$ is close to a 
a translation. The conjugacy $H$ from the previous proposition is
defined on $\bm W$, the image of $\V$ in $w$-coordinates.
So in the $w$-coordinates $H$ is defined above the line $\phi(m^+)$ and below the line 
$\phi(m^-)$). This mapping  can be extended in the same way to the region in the $w$-plane 
 to the left of $n^\pm$ and above the real line.   Exactly as in the proof of the previous proposition, this gives a quasiconformal conjugacy on the right half of the petal $\V$. 
To construct the conjugacy on the left half of the petal $\V$, we next extend the conjugacy in the same way  on the left petal.
In this way, we obtain in this case a qc conjugacy on a complex neighbourhood of $p$.

If $d>2$ there are more petals, but in the same way we glue the construction of the
conjugacies of the petals along the real line to a conjugacy on a neighbourhood of $p$.

If $d=1$ and $a>0$ then one has an attracting and repelling petal along the real line, and 
there are no other petals. In this case, in order to ensure that 
the petal regions overlap, we introduce a region which 
is no longer lens-shaped: replace the lines 
$x=\alpha |y|+x_0$ with $\alpha>0$ in the $w$-plane, as was done in Figure~\ref{fig:touching},
to lines $x=\alpha |y|+x_0$ with $\alpha<0$ as in Figure~\ref{fig:touching3}. 
The corresponding regions then are no longer lens-shaped, but  
we as before we can glue the construction of the conjugacy in the attracting 
petal on the left with the conjugacy in the repelling petal on the right to obtain 
a quasiconformal conjugacy on a neighbourhood of $p$.
\end{pf}

	\begin{figure}[htp] \hfil
		\beginpicture
		\dimen0=0.3cm
		\setcoordinatesystem units <\dimen0,\dimen0>
		\setplotarea x from -9 to 30, y from -10 to 10
		\setlinear
		\plot -16 0 8 0 /
		\circulararc 240 degrees from 5 0  center at  0 3 
		\circulararc -240 degrees from 5 0  center at  0 -3 
		\put {$\partial \V$} at 5 2 
		\circulararc -250 degrees from 1 0  center at  -2 -2 
		\circulararc +250 degrees from  1 0 center at  -2 2 
		\put {$\partial \U$} at 1.5 1 
		\circulararc -175 degrees from -15 0  center at  -10 -0.3 
		\circulararc -175 degrees from -11 0  center at  -8 -0.3

		\put {$p$} at -5.9 -0.9
		\put {$z\mapsto 1/z$} at 10 -2  
		\arrow <3mm> [0.2,0.67] from 7 -3.5 to 13 -3.5
		\setcoordinatesystem units <\dimen0,\dimen0> point at -9 0 
		\setplotarea x from -9 to 30, y from -7 to 9
		\put {$\bullet$} at 13 0  
		\plot 8 0 25 0 /
		\plot 17 0 5 4 /
		\plot 17 0 5 -4 /
		\plot 19 0 4 5 /
		\plot 19 0 4 -5 /
		\plot 10 0 4 8 / 
			\plot 11 0 5 8 / 
		\put {$l_{\pm}$} at 14 0.8 
		\put {$\phi(\partial \U)$} at 19 1
		\endpicture
		\caption{Part of the touching box mapping 
		at a fixed point $p$ of $f^s$, with domain $\U:=\comp_{p} f^{-s}(D_{\pi-\theta}(J'))$ is contained in $\V:=D_{\pi-\theta}(J')$.
		The mapping $z\mapsto 1/z$ sends the curves 
		$\partial \U,\partial \V$  to the wedge-shaped curves
		on the right. Here $\phi(\partial \V)=l_\pm$ \label{fig:touching3}}.
		\end{figure}

\subsection{Angle control in the presence of parabolic points}

		To prove the existence of the touching box mappings,
		we shall also need the following corollary.

\begin{cor}\label{cor:petals}
Suppose that $f$ is $C^{3}$ and that $p$ is a periodic point of period $s$. 
There exists $\delta>0$ and $\theta_0\in(0,\pi/2),$ 
such that for any $\theta\in(0,\theta_0)$
and all $i$ sufficiently large, we have the following.
Take $J=[p,q_{0}]$ and $J'=f^{s}(J)$, with $|J'|<\delta$.
Assume that $x$ is such that $f^{k}(x)=p$ 
and that $J_{x}$ is a one-sided neighbourhood of $x$ 
so that $f^{k}(J_{x})\rightarrow J$ is a diffeomorphism 
with $Df^{k}(x)\neq 0$.
Then
$$\Comp_{x}f^{-n}(D_{\pi-\theta}(J'))\subset D_{\pi-\theta}(J_{x}),$$
where $n=k+is$.
		\end{cor}

\begin{pf}
When $f$ is analytic, we argue as follows.
By Proposition~\ref{prop:petals},
$\mathrm{Comp}_p f^{-is}(D_{\pi-\theta}(J'))\subset D_{\pi-\theta}(J)$
for all $i> 0$. Moreover, for fixed $\theta$, 
the curvature of the boundary of 
$f^{-is}D_{\pi-\theta}(J'))$ tends to $\infty$ as $i\rightarrow\infty,$
and the diameter of this set shrinks to zero.
At the same time, $f^k$ is a $C^3$ 
diffeomorphism from a real neighbourhood $W_x$ of $J_x$
to a neighbourhood of $J$, the curvature of 
$D_{\pi-\theta}(J)$ is bounded, and
$(f^k|_{W_x})^{-1}(\mathrm{Comp}_p f^{-is}(D_{\pi-\theta}(J')))$ and
$D_{\pi-\theta}(J_x)$ are tangent to each other at $x$.

To make this argument work when $f$ is $C^3$
we need to show that the curvature of $f^{-is}D_{\pi-\theta}(J'))$
tends to $\infty$ as $i\rightarrow\infty$.
Let $J_0=J'$ and set
$J_k= \mathrm{Comp}_p
f^{-ks}(J_0)$.
Let $\Lambda_{J_k}$ be the rescaling from an interval $I$
of unit size to $J_k$, $k=0,\dots, i$ and let
$g=\Lambda^{-1}_{J_0}\circ f^s\circ \Lambda_{J_1}$ near $p.$
By the Fa\`a di Bruno Formula,
we have the following expression for the second derivative
of $g^{k}$.
$$D^2g^{k}(z) = D^2 g\circ g^{k-1}(z) \cdot (Dg^{k-1}(z))^{\otimes 2}
+\sum_{j=1}^{k-1}Dg^{k-j}\circ g^j(z)\cdot D^2 g\circ g^{j-1}(z)\cdot
(Dg^{j-1}(z))^{\otimes 2}.$$
Since the asymptotically holomorphic extension of $f$ is $C^3,$ 
there exists $C_0>0$ be so that $\|f\|_{C^3}<C_0$. In particular, we 
have that $Df^s$ and $D^2f^s$ are uniformly Lipschitz.
Thus we have that if $x_0\in J_k$ and
$y_0\in \mathrm{Comp}_p f^{-ks}(D_{\pi-\theta}(J_0)),$
then 
$$\|Df^s(x_0)-Df^s(y_0)\|\leq C_0|x_0-y_0|\leq C_0|J_k|.$$
So we have that
$$\|Df^s(y_0)\|\leq  C_0|J_k| +\|Df^s(x_0)\|.$$
By the chain rule, we have that
$$\|Df^{ks}(y_0)\|\leq \prod_{i=1}^{k}\big(C_0|J_i| + \|Df^s(f^{si}(x_0))\|\big)$$
$$\leq C_1\prod_{i=1}^{k}\Big(|J_i|+\frac{|J_{i-1}|}{|J_i|}\Big)\leq
C_2 \prod_{i=1}^{k} \frac{|J_{i-1}|}{|J_i|}\Big(1+\frac{|J_i|^2}{|J_{i-1}|}\Big)$$
$$\leq C_3\Big(\frac{J_0}{J_1}\Big)^k\prod_{i=1}^k\big(1+\frac{|J_i|^2}{|J_{i-1}|})
\leq C_4\Big(\frac{|J_0|}{|J_1|}\Big)^k.$$

Now for the rescaled maps, observe that there exists a constant
$C'$ so that
$$\|D (\Lambda_{J_0}^{-1}\circ f^{sk}\circ\Lambda_{J_k})|\|\leq C',$$
and 
$$\|D^2 (\Lambda_{J_{i-1}}^{-1}\circ f^s\circ\Lambda_{J_{i}})(x)\|
\leq C'|J_i|\|D^2f^s|_{D_{\theta}(J_i)}\|.$$
Observe that for a fixed periodic point $p$, 
$\|D^2f^s|_{D_{\theta}(J_1)}\|$ is bounded.
Therefore we have that there exists a constant $C>0$ such that
$$\|D^2g^{k}\|\leq
C\sum_{i=1}^{k}\|D^2(\Lambda_{J_{i-1}}^{-1}\circ f^s\circ\Lambda_{J_{i}})\|
\leq C\sum_{i=1}^{k}|J_i|,$$
which can be made as small as we like by choosing of $J'=J_0$ small.
\end{pf}

%
%
%

\begin{prop}[Pullbacks of Poincar\'e lens domains]  
\label{prop:pullbacklens}
Let $f\in\mathcal{C}$ be a mapping with periodic points.
Let $I$ be a real neighbourhood of the set of critical points of $f$
and let $B_{0}$ be the union of the immediate basins 
of attraction of the periodic attractors of $f$.
There exist constants $\delta>0,$ $C>0$ and $\theta_0\in(0,\pi/2)$,
so that for any $\theta\in(0,\theta_0),$
there exists $\theta'$
with the following properties.
Let $Z'$ and $\mathcal{P'}=M\setminus Z'$ be a set and 
a partition given by 
Lemma \ref{lem:startingpartition},
so that each periodic point in $Z'$
is within distance $\delta$ from a different point in $Z',$
on both sides if the periodic point is repelling on both sides.
Let $J'$ be a component of $\mathcal{P}'$ and
let $J$ be a component of $f^{-n}(J')$ such that
$f^n|_J$ is a diffeomorphism.
Assume that $f^{i}(J)\cap(I\cup B_{0})=\emptyset$ for all $i=0,\dots,n-1$.
Then
$$\Comp_{J}f^{-n}(D_{\pi-\theta}(J'))\subset D_{\pi-\theta'}(J),$$
where $\theta'<C\theta$.
\end{prop}
		
		\begin{pf}
			\noindent\textit{Case 1.}
			Let us first suppose that $J'$
			does not contain a parabolic point in its boundary.
			By Lemma \ref{lem:sum of nice lengths},
			there exists a constant $C$ such that 
			$\sum_{i=0}^n|f^i(J)|^{3/2}<C.$
			Moreover, since $f$ has no wandering domains, there
			exists a constant $C'>0$ such that for all $i$, $|f^{i}(J)|<C'|J'|$.
			So we have that
			$$\sum_{i=0}^n|f^i(J)|^{2}<C'|J'|^{1/2}\sum_{i=0}^n|f^i(J)|^{3/2}<CC'|J'|.$$
			Now, by Lemma \ref{lem:lens angle control}, since
                        $f^n\colon J\rightarrow J'$ is a diffeomorphism,
			$$(f^{n}|_J)^{-1}D_{\pi-\theta}(J')\subset D_{\pi-\theta'}(J),$$
			where $|\theta'-\theta|\rightarrow 0$ as $|J'|\rightarrow 0$.
			
			\medskip
			\noindent\textit{Case 2.}
			Suppose now that $J'$ contains a parabolic periodic point $p$
			in its boundary.
			Let $\hat p$ denote the orbit of $p$.
			For convenience let us assume that $f^s$ is orientation preserving
			near $p$. 
			Let $k_0, 0\leq k_0\leq n$ be minimal so that 
			$^{k_0}(J)$ contains a point of $\hat p$ in its boundary.
			Let $k_1,  k_1\geq k_0$ be minimal so that $f^{k_1}(J)$ contains 
			$p$ in its boundary.
			Then by Proposition \ref{prop:petals}
			$$\mathrm{Comp}_{f^{k_1}(J)}f^{-(n-k_1)}D_{\pi-\theta}(J')\subset D_{\pi-\theta}(f^{k_1}(J)).$$
			By Lemma \ref{lem:lens angle control}, there exists $\theta'$, close to $\theta$ so
			that
			$$\mathrm{Comp}_{f^{k_0}(J)}f^{-(n-k_0)}D_{\pi-\theta}(J')\subset D_{\pi-\theta'}(f^{k_0}(J)).$$
			By the choice of $\mathcal{P}'$, $f^{k_0}(J)$ cannot contain a 
			parabolic periodic point in its boundary. 
			If $k_0=0$ we are done.
			Otherwise, arguing as in Case 1
			finishes the proof.
		\end{pf}

	\subsection{Global touching box mappings}\label{subsec:global touching box mappings}
In order to deal with the part of the dynamics which stays
away from the set of critical points, we shall 
introduce so-called {\em global touching box mappings}.
These are maps for which the domains and ranges touch
at certain points, and where the domain does not
contain critical points. These objects are similar to what was constructed
in Kozlovski's thesis \cite{Koz_ax} and are  defined as follows.
(But in Kozlovski's setting there are big bounds and no parabolic points,
so we have to argue somewhat more carefully here.)

Consider a $C^{3}$ mapping $f\colon M\rightarrow M$,
where $M$ is $[0,1]$ or $S^1$.
Let $B_{0}$ be the union of
immediate basins of attraction of $f$ (possibly one-sided). 
Let $Z'$ be a finite forward invariant set containing 
the boundary points of $B_{0}$,
the boundary points of a nice real neighbourhood
$K$ of the set of critical points of $f$ 
and all parabolic orbits. 
Let $N\geq 1$ be an integer, 
$Z\subset f^{-N}(Z')$ and take $\theta\in(0,\pi/2).$ 
We say that
$$F_{T}\colon \UU_T\rightarrow \VV_T$$
is a \emph{touching box mapping}\label{def:tbm} if the following hold. 
Let $\mathcal{P}$ and $\mathcal{P}'$ 
		be the partitions corresponding to $Z$ and $Z',$ respectively. 
\begin{enumerate}[label=(\alph*)]
\item We have $$\VV_T=\cup_{J'}D_{\pi-\theta}(J')$$
	where the union runs over all components $J'$ of $\mathcal{P}'$.
\item We have $\mathrm{closure}(\mathcal U_T)$ contains
	$$[0,1]\setminus\mathrm{int}(B_{0}\cup K)\quad
\mathrm{or}\quad S^1\setminus\mathrm{int}(B_{0}\cup K).$$
\item For each component $J$ of $\mathcal{P}$ 
	in the complement of $K\cup B_{0}$, 
	there exists an integer $k=k(J)\leq N$ 
	so that $f^{k}|_J$ is a diffeomorphism and
	$J':=f^{k}(J)\in\mathcal{P}'$. 
Correspondingly, there exists a set $\U(J)$ with $\U(J)\cap\mathbb{R}=J$, 
which is mapped by the asymptotically holomorphic extension of 
$f^{k}\colon J\rightarrow J'\in\mathcal{P}'$
diffeomorphically onto $D_{\pi-\theta}(J')$; 
that is, $$f^{k(J)}\colon \U(J)\rightarrow\Comp_{J'}(\VV_T)$$
is a an asymptotically conformal diffeomorphism. 
Define $F_{T}|_{\U(J)}$ to be this extension of $f^{k(J)}$ and let 
$$\UU_T=\bigcup_{J}\U(J)$$
where the union is taken over all these intervals $J$.
\item $\UU_T\subset \VV_T$.
More precisely, each component of $\UU_T$
is either compactly contained in $\VV_T$ or it has a quadratic tangency with $\VV_T$.
Any two components of $\UU_T$ which have a common end point,
do not have a tangency.
\end{enumerate}

Note that by Theorem B of \cite{MMS} there are only
finitely many periodic attractors. 
Note that $\UU_T$ does not contain critical points of $f$, and that
$\UU_T$ 
does not cover the entire interval $[0,1]$ (resp. $S^1$). 
See Figure~\ref{fig:touchingglobal} for a sketch of a
global touching box mapping.

\begin{thm}[Global touching box mappings]\label{thm:touching box map}
Suppose that $f\in\mathcal{C}$ is a mapping with periodic points.
For each $\tau>0,$ 
there exists a nice real neighbourhood $K$ of the set of critical points of $f$,
such that $\diam(K(c))<\tau$ for each $c\in\Crit(f)$,
and finite forward invariant sets $Z$ and $Z'$ as above,
and $\theta\in(0,\pi/2)$ so that these choices give a
touching box mapping as defined above.

Moreover, if $\tilde{f}$ is another mapping of class
$\mathcal{C}$ that is topologically conjugate to $f$
by a conjugacy that maps parabolic points bijectively to parabolic points,
then we can take $Z, Z', \theta$ for $F_{T}$ so that
\begin{itemize}
\item if $\tilde{Z},\tilde{Z}'$ and $\tilde{\theta}=\theta$ are
the corresponding objects for $\tilde{f}$, 
these objects yield a touching box mapping for $\tilde{f}$;
\item there exists a qc mapping $H\colon \C\rightarrow\C$
with $H(\VV_T)=\widetilde{\VV}_T,$
$H(\UU_T)=\widetilde{\UU}_T$ which agrees with the topological conjugacy on $Z$
and which agrees on the petals attached to periodic points in $Z'$
with the conjugacy $H$ defined in Proposition~\ref{prop:petals};
\item inside the components $\V_{i}$ of $\VV_T$, which do not intersect $\UU_T$,
we can choose $H\colon \V_{i}\rightarrow\widetilde{\V}_{i}$,
so that it extends a qs-conjugacy
$h\colon V_{i}\rightarrow\widetilde{V}_{i}$
provided that the following holds:
for each $x\in V_{i}$ so that there exists $f^{n}(x)\in Z'$ one has
$\tilde{f}^{n}\circ h(x)=H\circ f^{n}(x)$.
\end{itemize}
\end{thm}
		
\begin{pf}
			Let $\theta_0\in(0,\pi/2)$ and $\delta>0$ be as in Proposition~\ref{prop:petals}.
			Let $\mathcal P$ be the partition given by  Lemma~\ref{lem:startingpartition}.
			Let $c$ be any critical point of $f$. If $f$ is infinitely renormalizable 
			at $c$ we can choose $\mathcal P(c)$ so that $|\mathcal P(c)|<\tau$;
			moreover, taking $\mathcal P(c)$ smaller if necessary
			we can assume that if $c'$ is any critical point whose orbit 
			accumulates on $c$, then $|\comp_{c'}f^{-n}(\mathcal P(c))|<\tau,$
			where $n$ is minimal so that $f^n(c')\in\mathcal P(c)$.
			If $f$ is finitely renormalizable at $c$, then let $K'=\mathcal P(c)$.
			The interval $K'$ is a periodic interval of period $s\geq 0$
			and the mapping $f^s\colon K'\to K'$ is non-renormalizable. 
			Observe that we do not exclude the possibility that $f$ is
			non-renormalizable at $c$.
			By the last point of Lemma~\ref{lem:startingpartition}
			there exists a nice interval $K(c)$ containing $c$
			such that $|K(c)|<\tau$ for any $c\in\crit(f)$ and
			$\partial K(c)$ consists of preperiodic points.
			Let $Z'$ be a smallest 
			finite forward invariant set 
			that contains the boundary points of all immediate basins of
			periodic attractors of $f$, all parabolic periodic orbits,
			and the points 
			$\cup_{c\in\crit(f)}\partial K(c)$
			and with the property that
			each periodic point $p$ in $Z'$ is within distance $\delta$ 
			from a point $p'$ in $Z'$, with $p'\neq p$
			(on each side if $p$ is repelling on both sides).
			Such a set exists by Lemma \ref{lem:startingpartition}.
			Let $Z_0$ be the collection of periodic points in $Z'$.
			Now choose $\theta\in(0,\theta_0)$ so small that 
			for any interval $T$ and any $n$ for which 
			$f^i(T)\cap(I\cap B_0)=\emptyset$, $i=0,1,\dots,n-1$,
			so that $f^n|_T$ is a diffeomorphism,
			one has
			$$\mathrm{Comp}_T f^{-n}(D_{\pi-\theta}(T'))\subset D_{\pi-\theta_0}(T),$$
			where $T'=f^n(T)$.
			
			Next choose $N_0$ large enough
			so that each point $x\in Z'$ is mapped into 
			$Z_0$ within $N_0$ iterates of $f$, and choose 
			$N>N_0$ so that 
			$N-N_0$ is a multiple of the periods of the periodic orbits of $f$ in $Z_0$ 
			and so large that we can apply Corollary~\ref{cor:petals}
			for each $p\in Z_0$ and each $x\in Z'$;
			in other words,
			when $J_x$ is a component of $\mathcal{P}'$ containing $x$ in its boundary
			and $J_{p}'$ is a component of $\mathcal{P}'$ containing $p=f^N(x)\in Z_0$, then
			$$\mathrm{Comp}_x f^{-N}D_{\pi-\theta}(J_p))\subset D_{\pi-\theta}(J_x).$$
			Now let $$Z=\{x:f^N(x)\in Z'\ \mathrm{and}\ f^i(x)\notin I\cup B_0\ \mathrm{for\ all\ } i=0,1,\dots,N-1\},$$
			and let $\mathcal{P}$ be the corresponding partition.
			For each component $J$ of $\mathcal{P}$, there exists $k(J)\leq N$ such that
			$f^k|_J$ is a diffeomorphism and $J':=f^k(J)\in\mathcal{P}'$.
			If $k<N$, then $f^k(J)$ is equal to a component of $I$ or of $B_0$.
Define $\VV_T=\cup_{J'\in\mathcal{P}'} D_{\pi-\theta}(J')$. 
Let $\UU_T=\cup_{J\in\mathcal{P}} \U(J),$ and
$$F_T(x)=f^{k(J)}(x).$$
			
Now we need to check that $F_T\colon \UU_T\rightarrow \VV_T$
is a touching box mapping.
We must show that $\UU_T\subset \VV_T$ in the required way.
Let $J\in\mathcal{P}$ and let $\hat{J}$
be the component of $\mathcal{P}'$ which contains $J$.
By 
$$\mathrm{Comp}_T f^{-n}(D_{\pi-\theta}(T'))\subset D_{\pi-\theta_0}(T),$$
we have that 
$$\U(J)\subset\mathrm{Comp}_J (f^{-k(J)}D_{\pi-\theta}(J'))
\subset D_{\pi-\theta_0}(J)\subset D_{\pi-\theta}(\hat J),$$
when $J$ and $\hat J$ do not have a common boundary point.
If $J$ has a common boundary point with $J'$ then we have that
$\U(J)\subset \V(J')$.
			
It follows from Proposition~\ref{prop:petals} 
that there is a qc-conjugacy as in the statement of the theorem.
Consider a point $x$ where 
$\partial \UU_T$ and $\partial \VV_T$ are tangent to each other.
Each such point is contained in $Z$ and therefore is mapped by some
iterate 
$f^n$ to a periodic point $f^n(x)=p$.
Hence, provided that $\theta\in(0,\pi/2)$ is chosen sufficiently
small,
we can define $H$ near $x$ by nearly
repeating the argument of Proposition~\ref{prop:petals}. 
			
That $H$ is quasisymmetric at a point $x$ where some
$\V_i$ and $\V_j$ touch can be seen as follows.
Each such point is mapped under $N_0$ iterates to some periodic point
$b\in Z_0$ of say period $s$ with associated petal(s) $\V_{i'}$ and $\V_{j'}$. 
By construction $H$ is a quasiconformal conjugacy on a complex 
neighbourhood of $b$
intersected with the Poincar\'e lens domains $\V_{i'}$ and $\V_{j'}$.
Moreover, if $y$ is near $b$
and $f^n(y)\in Z'$, then 
$\tilde{f}^n(H(y))\in \tilde{Z}'$.
Such points converge no faster than exponentially to $b$
and the rates of approach are the same from both sides of $b$.
Since $f\mapsto f^N(t)$ is a diffeomorphism near $x$,
and the assumption that parabolic periodic points are mapped 
bijectively to parabolic periodic points that
there exists a similar sequence of points approaching from both sides
of
$x$ and the rates of approach are the same.
It follows that $H$ is quasisymmetric at $x$.
		\end{pf}

The proof of Theorem~\ref{thm:qc rigidity away from critical points}
follows immediately from the following corollary to Theorem \ref{thm:touching box map}.
\begin{cor}[Rigidity away from $\crit(f)$]
\label{cor: pullback argument for touching box mappings}
Consider two $\mathcal{C}$ maps $f$ and $\tilde{f}$
as in the Theorem~\ref{thm:main}
that are topologically conjugate, by a topological conjugacy that maps hyperbolic,
respectively parabolic, periodic points of $f$ to hyperbolic,
respectively parabolic, periodic points of $\tilde{f}$.
Assume that associated to these maps we have two touching box mappings
$$F_T\colon \UU_T\rightarrow \VV_T\;\mbox{and}\; 
\tilde{F}_T\colon \widetilde{\UU}_T\rightarrow\widetilde{\VV}_T$$
defined as in the Theorem~\ref{thm:touching box map}.
Then there exists a quasiconformal conjugacy
$H\colon \VV_T\rightarrow\widetilde{\VV}_T$
so that $H(\UU_T)=\widetilde{\UU}_T$, and so that
$\tilde{F}_T\circ H=H\circ F_T$ on $\UU_T$.
Restricted to $\VV_T\setminus \UU_T$,
this $H$ agrees with the $H$ from
Theorem~\ref{thm:touching box map}.
\end{cor}
		
\begin{pf}
Define $H_0$ to be equal to the $K$-qc map
$H\colon \VV_T\rightarrow\widetilde{\VV}_T$
from Theorem~\ref{thm:touching box map}. Since $H(\UU_T)=\widetilde{\UU}_T$,
$\tilde{F}_T\circ H=H\circ F_T$ on $\partial \UU_T$, and 
$F_T\colon \UU_T\rightarrow \VV_T$ and 
$\tilde{F}_T\colon \widetilde{\UU}_T\rightarrow\widetilde{\VV}_T$ are conjugate,
and these maps have no critical points,
we can inductively define
$H_n$ on $\UU_T$ by 
$\tilde{F}_T\circ H_{n}=H_{n-1}\circ F$ on $\UU_T$ and
$H_{n}=H_{n-1}$ on $\VV_T\setminus \UU_T$.

It follows from
Lemma~\ref{lem:sum of nice lengths} as we pull back that the limiting mapping is $K'$-qc:
we can control the sum of diameters along the pullbacks of intervals,
which in turn controls the constant of quasiconformality.
\end{pf}

Note that the Julia set $K(F_T)=\{z\in \UU_T; (f_T)^n(z)\in \UU_T\mbox{ for all }n\ge 0\}$
		of the touching box mapping 
		$F_T \colon  \UU_T\to \VV_T$ is in general {\em not} hyperbolic (if $f$ has parabolic periodic points).

\section{Combinatorial equivalence and external conjugacies}
\label{sec:combinatorial equivalence}

In this section, we will show that two holomorphic complex box mappings that are
non-renormalizable and combinatorially equivalent are 
quasiconformally conjugate, and
in the smooth case, we show that the restrictions of smooth, quasiregular complex
box mappings to their real traces are quasisymmetrically conjugate on
the real line,
provided that there exists an \emph{external}
quasiconformal conjugacy between them.
We will define this notion below and we will also
prove a theorem that tells us when such an external
conjugacy exists (see Theorem \ref{thm:external conjugacy}).
		
We say that two real-symmetric box mappings
$F\colon  \UU\to \VV$ and $\tilde F\colon  \widetilde{\UU}\to\widetilde{\VV}$
are {\em strongly combinatorially equivalent}
if there exist two real-symmetric homeomorphisms
$H_i\colon  \VV\to \widetilde{\VV}$, $i=0,1$	
such that
\begin{itemize}
\item $\tilde F\circ H_1=H_0\circ F$ on $\UU$;
\item $H_1=H_0$ on $\VV\setminus\UU$;
\item $H_1=H_0$ on $\mathcal V=\VV\cap\mathbb R$.
\end{itemize}
	Note that this implies that $H_0$ maps a critical point of $F$ to a critical point
	of $\tilde F$ of the same order (and the post-critical set of $F$ to the post-critical set of $\tilde F$).
	Also note that $H_1$ and $H_0$ are homotopic 
	relative to $\VV\setminus \UU$ and $\VV\cap \R$.
Hence we can find a sequence of real-symmetric
homeomorphisms $H_n\colon  \VV\to \widetilde{\VV}$ with 
$\tilde F\circ H_{n+1}=H_n\circ F$.
In this way, each puzzle piece $\bm P$ of $F$ 
(i.e. component of $F^{-n}(\VV)$) corresponds
to a unique puzzle piece $\widetilde{\bm{P}}$ of $\tilde F$, i.e.,
$H_k(\bm P)=H_n(\bm P)$ is a puzzle piece of level $n$ for
$\tilde F$ for each $k\ge n$.

\begin{figure}[htbp]
\begin{center}
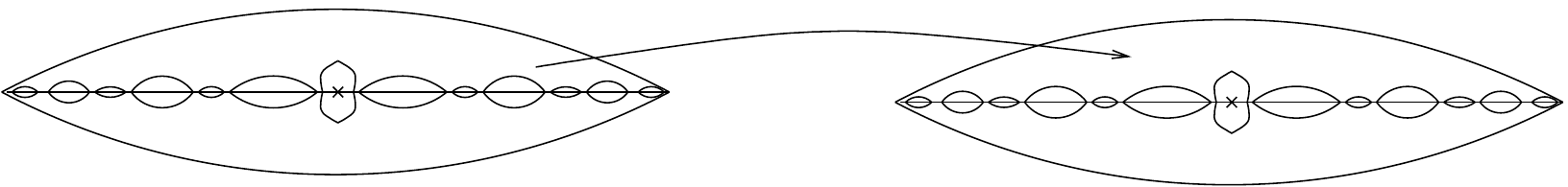
\end{center}
\caption{An external conjugacy. $H_0\colon\bm{\mathcal{V}}\rightarrow
\widetilde{\bm{\mathcal{V}}}$ is a conjugacy on the boundaries of
  the puzzle pieces in $\bm{\mathcal{U}}$.
}
\label{fig:extconj}
\end{figure}

We say that $H_0$ is an \emph{external conjugacy}
\label{def:extconj} if
\begin{enumerate}
\item $H_0(\VV)=\widetilde{\VV}$, 
\item $H_0(\UU)=\widetilde{\UU}$, 
\item $\tilde F\circ H_0=H_0\circ F$ on $\partial \UU.$
\end{enumerate}
See Figure~\ref{fig:extconj}. An external conjugacy provides us with a
$\emph{boundary marking}$ $H$ on $\partial \bm U$
for any puzzle piece $\bm U$ of level $n$ through the formula
$H_0\circ F^n(z)=\tilde F^n\circ H(z)$ for $z\in\partial\bm U,$
when there is a choice of pullback (e.g. when $F^n|_{\bm U}$ is not 
univalent) we choose $H$ so that it maps points on the real line
to points on the real line and preserves the order on the real line.

\subsection{Existence of external conjugacies}
\label{sec:external conjugacy}
\begin{thm}[Conditions for the existence of a quasiconformal external conjuagacy]
\label{thm:external conjugacy}
There exists $\theta\in(0,\pi/2)$ such that the following holds.
Assume that $F\colon  \UU\to \VV$ and 
$\tilde F\colon \widetilde{\UU}\to \widetilde{\VV}$
are complex box mappings which are non-renormalizable and 
strongly combinatorially equivalent that extend the real return maps
to $I$ and $\tilde I$, where $I$ and $\tilde I$ are unions of critical 
puzzle pieces.
Moreover, assume that
\begin{enumerate}
\item $F,\tilde F$ are real-symmetric and 
that all components of $\UU,\widetilde{\UU}$ intersect the real line,
and that the components of $\VV,\widetilde{\VV}$ that intersect $\Omega'$ 
are Poincar\'e discs with angle $\pi-\theta$;
\item \label{itemb}
there exist  touching box mappings $F_T\colon  \UU_T\to \VV_T$  and
$\tilde F_T\colon \widetilde{\UU}_T\to \widetilde{\VV}_T$ with angle
$\pi-\theta'$ satisfying the angle condition
$$\theta\in (0,\theta'/2),$$
where for any critical point $c$ of $F$, $\VV_T(c)\cap\mathbb{R}=I(c),$
and similarly for the corresponding objects associated to $\tilde F$.
\item $F\colon \UU\to \VV$ and $\tilde F\colon  \widetilde{\UU}\to \widetilde{\VV}$
both satisfy the gap and extension properties from 
page~\pageref{def:gapextension}.
\end{enumerate}
Furthermore, assume that there exists a quasiconformal homeomorphism
$h\colon \VV\to \widetilde{\VV}$ so that
	\begin{enumerate}
\setcounter{enumi}{3}
\item $h\colon  \partial (\UU\cap \R)\to \partial (\widetilde{\UU}\cap \R)$
	is a bijection which conjugates $F$ and $\tilde F$ on these sets; 
	\end{enumerate}
Then there exists a quasiconformal external conjugacy $H$ with
$H(\VV)=
\widetilde{\VV}$,
$H(\UU)=\widetilde{\UU}$ and $\tilde F\circ H=H\circ F$ on $\partial \UU$.
	\end{thm}

\begin{pf}
The proof of this theorem is almost the same as the one given in the first part of the
proof of Theorem B' in \cite[page 440--442]{LevStr:invent} except that we shall
use the conjugacy at periodic points associated to the touching box mappings,
see Theorem~\ref{thm:touching box map}.
In order to be complete, we will outline the proof here.

\noindent
Let $\V$ be a component of $\VV$.

\noindent
{\bf Step 1:}  $\bm A^+:=(\V\setminus \UU)\cap \HH^+$
is a quasidisk (i.e. $\partial \bm A^+$ is a quasicircle.)
Here $\HH^{\pm}$ are the upper and lower half-planes. 
For convenience we write $\U^{\pm}=\U\cap \HH^{\pm}$
and $\V^\pm=\V\cap \HH^{\pm}$.
To show that $\bm A^+$ is a quasicircle,
we use the Ahlfors-Beurling Criterion which says a topological circle $\gamma$
is a quasicircle, if there exists a constant
$C$ for any two points $z_1,z_2\in \gamma$
there exists an arc $\gamma_{z_1,z_2}\subset \gamma$ such that
$\diam(\gamma_{z_1,z_2})\le C\dist(z_1,z_2)$.
To verify this condition take $z_1,z_2\in \partial \bm A^+$. If $z_1,z_2$
are both in  $\partial \V^+,$
then the Ahlfors-Beurling Criterion holds because $\partial \V^+$ is piecewise smooth.
Note that by assumption there exists a neighbourhood $\V'$ of $\V$ 
and for each component $\U_i$ of $\UU$ that intersects the orbit of a critical point
 has a neighbourhood $\U_i'$ which is mapped
with bounded degree on a component of $\V'$. 
Observe that by Lemma~\ref{lem:angle control} we may assume that 
any other component of the domain is a Poincar\'e lens domain with angle
bounded uniformly away from $\pi/2$.
Hence the Ahlfors-Beurling Criterion 
also holds on $\partial \U_i^+$ (with a constant which is independent of $i$).
If $z_1,z_2$ are in different components 
$\partial \U_i^+,\partial \U_j^+$ of $\partial \U^+$, then 
the Ahlfors-Beurling Criterion holds uniformly.
		
\noindent
{\bf Step 2:} Let us construct a quasiconformal mapping 
$h\colon  \VV\to \widetilde \VV$
with $h(\UU)=\widetilde{\UU}$.
We start by noting that each $x_i\in \partial (\UU\cap \R)$
is mapped to some $b_i\in \partial \VV$, and that each $b_i\in \partial \VV$
is mapped by $f^{N_0}$  to a periodic points $a_i$ in the finite forward invariant 
set associated to the touching box map.
The angle condition is satisfied, and
there exist complex neighbourhoods $\bm D_i$
around each of these finitely many periodic points $a_i$ and 
neighbourhoods $\bm W_i$ around each $x_i$ so that $f^{N_0}\circ F$
maps $\bm W_i$ univalently onto $\bm D_i$
(and extends univalently onto a neighbourhood of $\bm D_i$)
and so that
$f^{N_0}\circ F(\bm W_i\cap \U_i)\subset \VV_T$. 
Of course around each $a_i\in \partial (\VV\cap \R)$,
there exists also a neighbourhood $\bm W_i$ so that $f^{N_0}$ maps
$\bm W_i$
univalently onto 
$\bm D_i$ and so that $f^{N_0}(\bm W_i\cap \VV)\subset \VV_T$.
For convenience choose $\bm D_i$ so that $H_T(\bm D_i)=\widetilde{\bm{D}}_i$.
Now define $H\colon \partial \VV^+\to \partial \widetilde{\VV}^+$ so that 
$H_T\circ f^{N_0} = \tilde f^{N_0} \circ H$ on $\bm W_i\cap \partial\VV^+$
where $\bm W_i$ are the neighbourhoods of $\partial (\VV^+\cap \R)$ as above.
Extend this mapping on the remainder of $\partial \VV$ by defining
$H\colon \partial \VV\to \partial \widetilde{\VV}$ 
so that it distorts arc-length by a constant factor.
Next define a provisional mapping $H_0\colon  \VV\to \widetilde{\VV}$ which agrees 
with the $H\colon \partial \VV\to \partial \widetilde{\VV}$ we just defined and so that
$H_0(F(c))=\tilde F(\tilde c)$ for the finitely many critical points of $F$.
Choose $H\colon \U_i \to  \widetilde{\U}_i$ so that $\tilde F\circ H=H_0\circ F$,
where the choices of inverses of $\tilde F|_{\widetilde{\U}_i}$ are determined
by the ordering on the real line.  On $\partial \UU\cap \R$, the mapping $H$
is defined by the assumption of the theorem.  Thus $H$
is  defined on $\partial \bm A^+$. 
Let $\phi$ be the restriction of $H$ to $\partial \bm A^+$.
Let us show $\phi$ is a quasisymmetric homeomorphism on this quasicircle.
To do this, we need to show that for any three points
$z_1,z_2,z_3\in \partial \bm A^+$
which are roughly equidistant, $H(z_1),H(z_2),H(z_3)$ are also roughly equidistant.
If $z_i$ are all
in some set $\bm W_j$ as above, then using that $f^{N_0}\circ F|_{\bm W_i}$ is univalent,
by the Koebe Distortion Theorem,
we get that $y_i:=f^{N_0}\circ F(z_i)$, $i=1,2,3$ are also
roughly equidistant. Since $H_T$ is quasiconformal, this implies that 
$H(z_i)=(\tilde f^{N_0}\circ \tilde F)^{-1}\circ  H_T(y_i)$ are also all equidistant. 
If $z_i$ are not all in some set $\bm W_j$ as above, then we argue
exactly as in case (b) in \cite[page 442]{LevStr:invent}.
		
\noindent
{\bf Step 3:} It follows that there exists a quasiconformal map
$H\colon  \bm A^+\to \widetilde{\bm{A}}^+$
(and similarly a quasiconformal mapping  $H\colon \bm A^-\to \widetilde{\bm{A}}^-$.)
Since we already defined a quasiconformal map
$H\colon \bm U_i\to \widetilde{\bm{U}}_i$ we
obtain the required quasiconformal mapping $H\colon  \VV\to \widetilde{\VV}$.
\end{pf}

\subsection{Rigidity of non-renormalizable
analytic complex box mappings}
\label{sec:rigidity analytic}
\begin{prop}[The Spreading Principle (Analytic
  Version)\cite{KSS-rigidity} page 769]
\label{prop:spreading analytic}
Suppose that $F\colon\bm{\mathcal{U}}\rightarrow\bm{\mathcal{V}}$ and
$\tilde F\colon\widetilde{\bm{\mathcal{U}}}\rightarrow\widetilde{\bm{\mathcal{V}}}$
are strongly combinatorially equivalent holomorphic box mappings.
Suppose that $h\colon\mathbb C\rightarrow\mathbb C$ is a $K$-qc external
mapping for $F$ and $\tilde F$.
Let $\bm{\mathcal{O}}\supset\Crit(F)$ be a nice open set consisting of
puzzle pieces.
Let $\phi\colon\bm{\mathcal{O}}\rightarrow\widetilde{\bm{\mathcal{O}}}$ be a
$K$-qc mapping that respects the boundary marking on each
component of $\bm{\mathcal{O}}$.
Then there exists a $K$-qc mapping 
$\Phi\colon\bm{\mathcal{V}}\rightarrow\widetilde{\bm{\mathcal{V}}}$
such that the following hold:
\begin{itemize}
\item $\Phi=\phi$ on $\bm{\mathcal{O}}$;
\item for each $z\in\bm{\mathcal{V}}\setminus \bm{\mathcal{O}},$
$$\tilde F\circ \Phi(z)=\Phi\circ F(z);$$
\item $|\bar \partial \Phi|<k$ on $\bm{\mathcal{V}}\setminus
\Dom(\bm{\mathcal{O}})$,
where $k$ satisfies $K=\frac{1+k}{1-k}$;
\item for each puzzle piece $\bm P$
which is not contained in $\Dom(\bm{\mathcal{O}})$,
$\Phi(\bm P)=\widetilde{\bm{P}}$ and $\Phi\colon\bm P\rightarrow\widetilde{\bm{P}}$
repects the boundary marking.
\end{itemize}
\end{prop}

\begin{thm}[QC rigidity of combinatorially equivalent complex box mappings]
\label{thm: qc rigidity of box mappings}
Suppose that $F\colon \bm{\mathcal{U}}\rightarrow \bm{\mathcal{V}}$
and
$\tilde{F}\colon\widetilde{\bm{\mathcal{U}}}\rightarrow\widetilde{\bm{\mathcal{V}}}$ 
are 
corresponding complex box mappings
given by either Theorem~\ref{thm:box mapping persistent infinite branches}
or Theorem~\ref{thm:box mapping reluctant}.
In particular assume that they satisfy conditions 1, 2 and 3 of
Theorem~\ref{thm:external conjugacy}. 
Assume further that $F$ and $\tilde{F}$ 
are strongly combinatorially equivalent via
homeomorphisms $H_0$ and $H_1$,
where $H_0$ is a quasiconformal external conjugacy.
Then $F|_{\bm{\mathcal{U}}\cap\mathbb R}$ and 
$\tilde{F}|_{\widetilde{\bm{\mathcal{U}}}\cap\mathbb R}$
are quasisymmetrically conjugate;
that is, there exists a qs mapping 
$H\colon \mathcal U\rightarrow\tilde{\mathcal U}$ such that
$H\circ F=\tilde{F}\circ H$ on $\mathcal U$
and
in particular, $H$ maps critical points of $F$ to critical points of $\tilde{F}$ and
$H(F^i(c))=\tilde{F}^{i}(H(c))$. 
\end{thm}
\begin{pf}
This theorem is proved Lemma 6.6 and 6.7 and Section 6.4 of
\cite{KSS-rigidity}.
Lemma 6.6 and 6.7 show that there exist $\delta>0$ and 
arbitrarily small critical
puzzle pieces $\U$ which have $\delta$-bounded geometry 
and so that first the first landing map to $\U$ 
extends to a definite neighbourhood of $\U$,
so that $\U$ is $\delta$-nice. 
Lemma 6.6 deals with the case that F is reluctantly recurrent and
Lemma 6.7 deals with the persistently recurrent case.
The proof of this lemma relies on 
Theorem~\ref{thm:box mapping persistent}.
One then uses the Spreading Principle, as in
In Section 6.4 of \cite{KSS-rigidity}, 
to complete the proof of the theorem. 
\end{pf}

\begin{rem}
In Section~\ref{sec:qc rigidity}, we will show that
the conditions of Theorem~\ref{thm:external conjugacy}
are satisfied for these complex box mappings.
\end{rem}

We will now develop the tools that we will need to
deal with smooth mappings. Quasisymmetric rigidity for 
smooth complex box mappings is proved in
Theorem~\ref{thm:smooth box conj}.

\subsection{Restricting the puzzle to control dilatation}
\begin{figure}
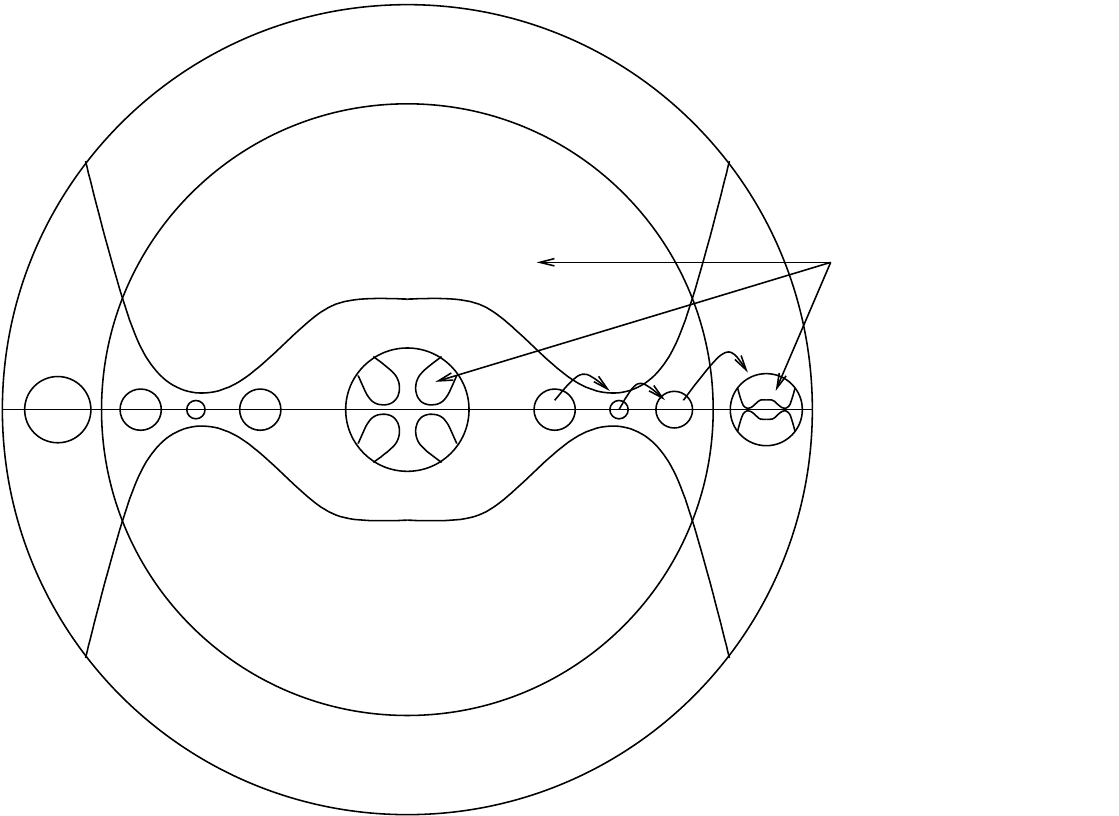
\caption{Exclude puzzle pieces that pass through long central
  cascades,
and their pullbacks from the Yoccoz puzzle.
}
\label{fig:modifiednest}
\end{figure}

Let us give some terminology.
Suppose that $F\colon\bm{\mathcal{U}}\rightarrow\bm{\mathcal{V}}$ is
a complex box mapping.
For $n\in\mathbb N\cup\{0\}$, the \emph{Yoccoz puzzle pieces of level}
$n,$ 
denoted by $\bm{\mathcal{Y}}^n$
is the set of connected components of
$F^{-n}(\bm{\mathcal{V}})$ and $\bm P$ is
a puzzle piece if $\bm P$ a component of $\bm{\mathcal{Y}}^n$ for some
$n$.
The \emph{Yoccoz puzzle} is the set of all puzzle pieces.
In this subsection we will define a sub-puzzle $\bm{\mathcal{M}}$
of $\bm{\mathcal{Y}}$. Roughly, we
construct it from the Yoccoz puzzle by excluding from the Yoccoz
puzzle all puzzle pieces that are
critical puzzle pieces that are deep in long cascades
of central returns or that are complex puzzle pieces that do not
intersect the real line and
are obtained as a ``complex'' pullback of such a puzzle piece 
through a long central cascade,
and we do not refine puzzle pieces that do not
intersect the real line. We will refer to this nest as the
\emph{modified Yoccoz puzzle}. The modified
Yoccoz puzzle does not give a partition of $\bm{\mathcal{U}}$
in the same way as some $\bm{\mathcal{Y}}^n$;
it has ``holes'' off the real line when there is a long cascade of
central returns, see Figure~\ref{fig:modifiednest}.
An important step in the proof 
of rigidity of complex box mappings
is filling in these holes. This is done in Theorem~\ref{thm:central cascades}.
Let us be more precise.

Let $b$ denote the number of critical points of $F$.
Fix $N_0\geq 12b^2-4b$ (notice that $N_0\geq 8$).
We will call central cascades of length longer than $N_0$
\emph{long}. Suppose that $\bm U\owns c,$ $c\in \crit(F)$, is a critical puzzle piece.
Assume that $k\in\mathbb N$ is such that
$\mathcal{C}_c^{k+1}(\bm Y_c^0)\subset \U\subsetneq \mathcal{C}_c^{k}(\bm Y_c^0)$.
Let $\Z^0= \mathcal{C}_c^{k}(\bm Y_c^0)$, and let $$\Z^0\supset \Z^1\supset
\Z^2\supset\dots$$ be the principal nest about $c$.
We will say that $\U$ is 
\emph{deep in a central cascade} if $\U\subset \Z^{N_0}$,
and the return time of $c$ to $\U$ is the same as its return time to
$\Z^{0}.$
It is easy to see that when we exclude puzzle pieces that are 
deep in a sufficiently 
long central cascade, we do not exclude puzzle pieces from 
the enhanced nest,any enhanced nest over a central cascade 
or from the good nest,
from the
modified puzzle:
Provided that $N_0\geq 12b^2-4b$, any puzzle piece $\bm I$ in the enhanced nest
has the property that if $\{\bm U_i\}_{i=0}^s$ is the chain of puzzle
pieces with $\bm U_0=\bm I$ and $\bm U_s=\bm V$,
a component of $\bm{\mathcal{V}}$,
then no critical puzzle piece $\bm U_i$ is deep in a 
central cascade. This follows immediately from the fact that
the mapping $F^{p_i}\colon\II_{i+1}\to\II_i$ 
between two puzzle pieces in the enhanced nest has bounded degree,
see \cite{KSS-rigidity}, Section 8,
and similarly for the puzzle pieces in
the enhanced nest over a long central cascade and for puzzle pieces in
the good nest.


Let $\bm{\mathcal{M}}^0$ be the components of $\bm{\mathcal{V}}.$
These are the puzzle pieces of level $0$ in the modified Yoccoz puzzle.
Assuming that
$\bm{\mathcal{M}}^{n-1}$ is defined, let
$\widehat{\bm{\mathcal{M}}^n}$
denote the collection of all puzzle pieces $(F|_{\bm P'})^{-1}(\bm P)$, for 
$\bm P\in \bm{\mathcal{M}}^{n-1}$ and $\bm P'\in
\bm{\mathcal{M}}^{n-1},$
with $\bm P'\cap\mathbb R\neq\emptyset$.
We have that either
\begin{itemize}
\item $F^{-1}(\bm P)$ intersects the
real line ($\widehat{\bm{\mathcal{M}}^n}$ contains all real puzzle pieces)  or 
\item $F^{-1}(\bm P)$ is not properly contained in 
a connected component of some $\bm{\mathcal{M}}^i$,
$0\leq i\leq n-1$
which does not intersect the real line($\widehat{\bm{\mathcal{M}}^n}$
contains all largest puzzle pieces that do not intersect the real
line).
\end{itemize}
Now we remove certain puzzle pieces which pass through long cascades
of central returns from
$\widehat{\bm{\mathcal{M}}^n}.$
We let $\bm{\mathcal{M}}^n$ denote the set of all puzzle pieces 
$\bm U$ in $\widehat{\bm{\mathcal{M}}^n}$ such that one of the
following holds:
\begin{enumerate}
\item If $\bm U$ is a critical
puzzle piece, $\bm U$ is not deep in central
cascade (we exclude real puzzle pieces that are deep in central cascades).
\item If $\bm U$ is a largest complex puzzle piece, then 
$\bm U\in\bm{\mathcal{M}}^n$ if one of the following holds:
\begin{enumerate}
\item The smallest puzzle piece in the (full) puzzle 
$\bm{\mathcal{Y}}$ that intersects the
real line and that contains $\bm U$ is not deep in a central cascade.
\item The smallest puzzle piece in the (full) puzzle
  $\bm{\mathcal{Y}}$ \label{page:mod2b}
that intersects the
real line and that contains $\bm U$ is deep in a central cascade,
$$\V^0\supset\V^1\supset\V^2\supset\dots\supset\V^N\supset\V^{N+1},$$
and if $1\leq i\leq N$ is such that $\bm U\subset \V^i\setminus
\V^{i+1},$
then $\U$ intersects a component of $(\V^i\setminus
\V^{i+1})\setminus R_{\V^0}^{-i}(V^0\setminus V^1),$ which intersects 
the real line.
\end{enumerate}
(We exclude complex puzzle pieces that are pullbacks of puzzle pieces 
through a long cascade of
central returns, unless these pullbacks are along real branches of the return map.)
\end{enumerate}

Let $\bm{\mathcal{M}}$ denote the collection of all puzzle pieces 
$\bm U$, such
that $\bm U\in\bm{\mathcal{M}}^n$ for some $n$.

\begin{prop}\label{prop:squares}
There exists $C>0$ such that for any $n\geq 0$, 
if $\U_0\in\bm{\mathcal{M}}^n$
then $$\sum_{i=0}^n \Big(\dist\,(F^i(\bm U_0), \mathbb
R)+\diam(F^i(\bm U_0))\Big)^2<
C\cdot\mu(\mathcal V)^{1/3}.$$
\end{prop}

\begin{pf}
Let $\bm U_i=F^i(\bm U_n)$ for $i=0,\dots,n$.
We decompose the chain $\{\U_i\}_{i=0}^n$ into pieces in a  similar 
way as in the proof of 
Lemma~\ref{prop:sum of squares, puzzle pieces}.
If $\bm U_0$ is a puzzle piece that intersects the real line,
then this follows from
Lemma~\ref{prop:sum of squares, puzzle pieces}.
So we can suppose that $\bm U_0$ does not intersect the real line.
Let $m>0$ be minimal so that $\bm U_{m+1}$ is a real puzzle piece.
Then the puzzle pieces $\bm U_0,\dots,\bm U_{m}$ are all largest complex
puzzle pieces, moreover for each of these puzzle pieces either 
condition 2(a) or 2(b) holds, see page~\pageref{page:mod2b}.
If condition 2(a) holds for $\U_m$,
Let $\hat{\bm{U}}_{m}$ be the smallest real puzzle in
$\cup_{i=0}^{n-m}\bm{\mathcal{M}}^i$
that contains $\bm U_{m}$. 
Suppose first that condition 2(a) holds for $\U_m$.
If for all $0\leq i\leq m$ condition 2(a) holds for 
$\U_i$, let $\hat{\bm{ U}}_i=\comp_{\bm U_i}F^{-(m-i)}(\hat{\bm{U}}_m)$.
Observe that since
$\bm U_{m-1}\subset \hat{\bm{U}}_{m-1}$ we must have that
$\hat{\bm{U}}_{m-1}$ intersects the real line,
for if not, then since
$\hat{\bm{U}}_{m-1}$ does not intersect the real line,
$\bm U_{m-1}\notin\bm{\mathcal{M}},$ which is a contradiction.
Moreover, if $\hat{\bm{U}}_{m-1}$ is deep in a central cascade (and
condition 2(b) on page~\pageref{page:mod2b} does not hold),
$\bm U_{m-1}$ would have been excluded from $\bm{\mathcal{M}}^{n-(m-1)},$
and $\bm U_0$ could not be in $\bm{\mathcal{M}}^n$.
Thus we have that $\hat{\bm{U}}_{m-1}$ is a puzzle piece in the modified
Yoccoz puzzle. Arguing inductively, we have that $\hat{\U}_0$ is a
real puzzle piece in the modified Yoccoz puzzle, so
$$\sum_{i=0}^n(\dist\,(F^i(\bm U_0), \mathbb R)+\diam(F^i(\bm U_0)))^2\leq
\sum_{i=0}^{m}\diam(F^i(\hat{\bm{U}}_0))^2+\sum_{i=m+1}^n\diam(F^i(\bm U_{m+1}))^2,$$
which is bounded as required by Lemma~\ref{prop:sum of squares, puzzle pieces}.

So let $k_0\leq m$ be maximal so that condition 2(b) holds for
$U_{k_0}$,
and let
$\hat{\bm{U}}_{k_0}$ be the component of $(\V^i\setminus
\V^{i+1})\setminus R_{\V^0}^{-i}(V^0\setminus V^1)$ that intersects
$\U_{k_0}$. Let $j$ be so that $R_{\V^0}|_{\V^1}=F^j|_{\V^1},$
and let $l_0$ be maximal so that
$\U_{k_0-l_0j}\in\V^{l_0}\setminus\V^{l_0+1}.$
Then there exists a constant $C>0$ such that
$\sum_{k=k_0-l_0j}^{k_0}\diam(\U_k)^2\leq C\diam(\V^0).$
Since the diameters of the puzzle pieces at the tops of long central
cascades decays exponentially, Proposition~\ref{prop:good geometry for
central cascades}, we can control the sum of the squares of the
diameters of the 
subsequences of
$\{\U_i\}_{i=0}^m$ which pass through central cascades.
For the segments of the orbits that are not pullbacks through long
cascades of central returns condition 2(a) holds, and we argue as in the
first case until the next central cascade.
\end{pf}

\subsection{Statement of the QC\textbackslash BG Partition of a
  Central Cascade}
\label{subsec:qcbg partition statement}

Let us describe a theorem that we will use when there
are no long central cascades. If there are no
long central cascades, then this theorem is not required. 
Let $N_0\in\mathbb N$ be the 
lower bound on the length of a long central cascade.
We will assume that $N_0\geq 12b^2-4b$.
Suppose that $c$ is a recurrent critical point of 
$F\colon\bm{\mathcal{U}}\rightarrow\bm{\mathcal{V}}$.
Fix a critical point $c_0\in\omega(c)$ 
chosen so that
$c_0$ has even order if there exists a critical 
point of even order in $\omega(c)$
and let $\II^0$ be the component of $\bm{\mathcal{V}}$
that contains $c_0$.
For each $k\in\mathbb N$,
such that $\hat{\mathcal{C}}^k_{c_0}(\II^0)$ is the first puzzle piece
 in a maximal central cascade of length $N\geq N_0$,
let $\bm{\mathcal{G}}_{c_0}(k)$,
denote the $\delta$-excellent
puzzle piece given by 
Proposition~\ref{prop:good geometry for central cascades}.
Recall that $\bm{\mathcal{G}}_{c_0}(k)$ is a $(1+1/\delta)$-quasidisk.

\begin{defn}Suppose that $I$ is is a real puzzle piece.
A \emph{necklace neighbourhood of} $I$ is
a countable union of pairwise disjoint open sets $\U_i\subset\mathbb C$
such that $I=\overline{\cup U_i}.$
\end{defn}

Recall the definition of a $(K,\delta)$-qc\textbackslash bg mapping
from page~\pageref{page:def qcbg}.
We will say that a mapping $H\colon\bm \Omega\rightarrow\widetilde{\bm{\Omega}}$
is $(K,\delta)$-\emph{qc\textbackslash bg in the upper half-plane} 
if it is $(K,\delta)$-qc\textbackslash bg, and the following holds:
Let $X_1$ and $\tilde X_1$ be the subsets of
$\bm \Omega$ and $\widetilde{\bm{\Omega}}$ given in
 the definition of a $(K,\delta)$-qc\textbackslash bg mapping.
We require additionally that each component
$\bm P$ of $X_1$ that is contained in $\mathbb H^+$
additionally satisfies
$$\mod(\bm{\Omega}\cap\mathbb H^+\setminus \bm P)\geq\delta\mbox{ and, }
\mod(\widetilde{\bm{\Omega}}\cap\mathbb H^+\setminus \widetilde{\bm{P}})\geq\delta.$$
We will refer to this property as $H$ having \emph{moduli bounds in the
upper half-plane.}

Let us assume the following proposition, which we prove in
the Section~\ref{sec:cascades}. It is one of the key ingredients for 
proving rigidity of smooth mappings.
It is designed to be applied with
the QC Criterion in the upper half-plane
to show that there is a quasiconformal
mapping with the same values on the real line.
It can be thought of as
a Carleson Box construction of a quasiconformal mapping,
but where rather than controlling the dilatation on every box, 
we only control the dilatation on a subset, and 
prove that the complementary set satisfies
certain geometric bounds. Recall that
$\mu(\J^0)=\max_{\U\subset\Dom'(\J^0)}\diam(\U),$
where the maximum is taken over connected components $\U$ of 
$\Dom'(\U)$. 

\begin{thm}[QC\textbackslash BG Partition of a Central Cascade]
\label{thm:central cascades}
There exist $K_0,\hat K\geq 1,$ $\delta_0, C>0$ such that for any
$K\geq 1,$ $\delta>0$, 
$c\in\crit(F)$ and $k\in\mathbb N$ with the property that
$\hat{\mathcal{C}}^k_c(\II^0_c)$ is the first puzzle piece
in a maximal central cascade of length $N\geq N_0$
the following holds:
Let $\J^0=\bm{\mathcal{G}}_{c_0}(k)$
and let $\widetilde{\J}^0$ be the corresponding puzzle piece
for $\tilde F$. 
Let $R\colon\J^1\rightarrow \J^0$ be the first return mapping of $\J^1$ to
$\J^0$ and $\tilde R\colon\widetilde{\J}^1\rightarrow \widetilde{\J}^0$ be the
first return mapping of 
$\widetilde{\J}^1$ to
$\widetilde{\J}^0$.
Assume that 
$H_0\colon\J^0\rightarrow\widetilde{\J}^0$ is a $(K,\delta)$-qc\textbackslash bg mapping which
respects the boundary markings on 
$\partial\J^0$ and $\partial\J^1$.
Then there exists a 
$(K',\delta')$-qc\textbackslash bg  
mapping 
$H\colon\J^0\rightarrow\widetilde{\J}^0$ 
with the following properties:
\begin{itemize}
\item $K'=\max\{\hat K,K(1+C\mu(\J^0))(1+C\mu(\widetilde{\J}^0)),$
\item $\delta'=\min\{\delta_0,\delta(1+C\mu(\J^0))^{-1},
\delta(1+C\mu(
\widetilde{\J}^0))^{-1}\},$
\item there exists
a partition of $\J^0$ into 
sets disjoint sets $X_0,$ $X_1,$ $X_2$ and $X_3$,
and similarly for the objects marked with a tilde,
\end{itemize}
such that:
\begin{enumerate}
\item For $i=0,1,2,3,$
$H(X_i)=\tilde X_i.$
\item The set $X_0$ has measure zero
and $H$ satisfies (\ref{eqn:qc bound 0}) at each $x\in X_0$.
\item The set $X_1$ and the corresponding set $\tilde X_1$
are union of pairwise disjoint topological disks 
$Y$, respectively $\tilde Y=H(Y),$
with $\delta'$-bounded geometry
and if $Y$ and $\tilde Y$ are contained in the upper half-plane
$$\mod(\J^0\cap\mathbb H^+\setminus Y)\geq\delta'\mbox{ and, }
\mod(\widetilde{\J}^0\cap\mathbb H^+\setminus \tilde Y)\geq\delta'.$$
\item The set $X_2$ contains a necklace neighbourhood of the real line
  and
$\J^0\setminus \J^3.$ For any
puzzle piece $\bm P\subset X_2$ for $F$ and for each $z\in\partial\bm P$ we have
$H\circ R(z)=\tilde R\circ H(z)$.
\item The set $X_3$ and $\tilde X_3$ are $K_1$-quasidisks that are disjoint from
the set of puzzle pieces that intersect the real line
and $H|_{X_3}$ is $\hat K$-quasiconformal.
\end{enumerate}
We define
$\bm{\mathcal{N}}_{c}(k)=X_2.$
\end{thm}
This will be proved in 
 Section~\ref{sec:cascades}.

\begin{rem}
The definition of a qc\textbackslash bg partition
consists of three sets: a set of measure zero, a set with bounded
geometry and a set where the dilatation is bounded. However,
in the QC\textbackslash BG Partition of a Central Cascade,
the partition consists of four sets. The sets $X_2$, where the mapping
$H$ is defined dynamically, and $X_3,$ which can be decomposed into 
quasidisks, are both included in the set where the dilatation is
bounded.
\end{rem}

\subsection{A Real Spreading Principle}
\label{subsec:restricted spreading}
Let us use Theorem~\ref{thm:central cascades} to prove
a spreading principle.

\begin{figure}
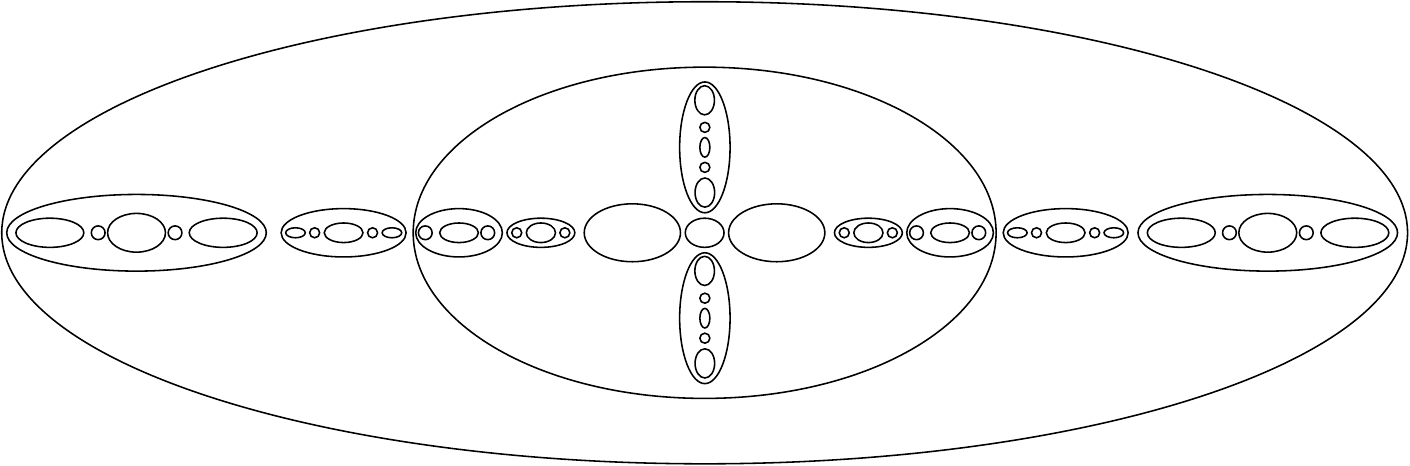
\caption{The puzzle pieces $\V^0$ and $\W^0$ together with some
pullbacks. The pullbacks of components of $\Dom(\V^0)$ which 
do not intersect the real line do not necessarily have moduli bounds
in $\mathbb H^+$. However, such pullbacks of $\Dom'(\W^0)$ do, since
$\mod(\V^0\setminus\W^0)$ is bounded away from zero. The same is true
on deeper levels.}
\label{fig:goodnest1}
\end{figure}

\begin{prop}[A Real Spreading Principle]\label{prop:realspreading}
There exists $\delta_0>0$ and
for any $K\geq 1$ there exists $K_0\geq 1$ such that
the following holds. 
Assume that $F\colon\bm{\mathcal{U}}\rightarrow\bm{\mathcal{V}},
\tilde F\colon\bm{\widetilde{\mathcal{U}}}\rightarrow\widetilde{\bm{\mathcal{V}}},$
are complex box mappings which satisfy the conditions of
Theorem~\ref{thm:external conjugacy},
and let $H_0\colon\mathbb C\rightarrow\mathbb C$ be a $K$-qc external
conjugacy between $F$ and $\tilde F$. Then
there exist arbitrarily small nice
open complex neighbourhoods $\bm{\mathcal{O}}$
of $\Crit(F)$ consisting of puzzle pieces in
$\bm{\mathcal{M}}$ and a mapping $H_{\bm{\mathcal{O}}}\colon\mathbb C\rightarrow\mathbb C$
such that 
\begin{itemize}
\item $H_{\bm{\mathcal{O}}}(\bm{\mathcal{O}})=\widetilde{\bm{\mathcal{O}}};$
\item $H_{\bm{\mathcal{O}}}$ agrees with the
boundary marking on $\partial\bm{\mathcal{O}};$
\item $H_{\bm{\mathcal{O}}}(\Dom'(\bm{\mathcal{O}}))=\Dom'(\widetilde{\bm{\mathcal{O}}});$
\item for each component $\U$ of $\Dom'(\bm{\mathcal{O}}),$
$H_{\bm{\mathcal{O}}}\circ F(z)=\tilde F\circ H_{\bm{\mathcal{O}}}(z)$ for $z\in\partial\U$.
\end{itemize}
Moreover, 
$H_{\bm{\mathcal{O}}}\colon\mathbb{H}^+
\rightarrow \mathbb{H}^+$
is a $(K_0,\delta_0)$-qc\textbackslash bg mapping.
\end{prop}

\begin{pf}
Before we give the details of the proof, let us 
explain our strategy. We start with a 
$K$-qc external conjugacy $H_0:\mathbb C\to\mathbb C$
between 
$F:\UU\to\VV$ and $\tilde F:\widetilde{\UU}\to\widetilde{\VV}$.
We are going to pull $H_0$ back by $F$ and $\tilde F$.
We only pull $H_0$ back along branches of $F$ which intersect
the real line, see Figure~\ref{fig:goodnest1}. 
It is impossible for us to control the
dilatation of high iterates of $F$ along complex branches of $F$,
and
this is one
reason why we cannot repeat the proof of the usual
Spreading Principle.
In Step 1, we pullback to obtain a mapping that is a {\em conjugacy up to
the first landing time to} $\W^0$; that is,
we will show that there exists a qc mapping $\psi_0$ of the plane,
so that for any $z\in\Dom'(\W^0),$ if $s_z\geq 0$ is minimal so 
that $F^{s_z}(z)\in\W^0$, we have that
$\psi_0\circ F^{s_z}(z)=\tilde F^{s_{z}}\circ\psi_0(z).$
In this step we are only pulling back through diffeomorphic
branches of $R_{\V^0} (=F$). The next pullback has a critical point, and
in order to pull $H_0$ back we need to modify it
 so that it maps critical values of
$F$ to critical values of $\tilde F$. This is done in Step 2.
We do not control the dilatation of the new mapping, but
we only modify it in components of $\Dom(\W^0)$ which contain
critical values of $R_{\V^0}$. These puzzle pieces have bounded
geometry. In Step 3, we use Steps 1 and 2 to go down to the 
next level of the good nest. Figure~\ref{fig:cascrep} shows how the
domains are nested.
The mapping that we obtain in this step
is not quite a $(K_0,\delta_0)$-qc\textbackslash bg mapping
on the upper half-plane, since
some of the puzzle pieces where the dilatation is unbounded may
not be well-inside the complex plane. We rectify this by 
refining the mapping in these domains,
see Figure~\ref{fig:wellinside}. 
These puzzle pieces are
components of $F^{-1}(F|_{\W^0})^{-1}(\Dom'(\W^0)).$
Using the fact that the components of 
$(F|_{\W^0})^{-1}(\Dom'(\W^0)),$ are well-inside of $\W^0,$
after we pullback the mapping from $\W^0$ to $\widetilde{\W}^0$
defined in Step 3 to these components of 
$F^{-1}(F|_{\W^0})^{-1}(\Dom'(\W^0)),$ which are mapped
diffeomorphically
onto components of $\W^0$, we obtain a 
$(K_0,\delta_0)$-qc\textbackslash bg mapping.
This is done in Step 4.
In Step 5, we explain how to continue the proof inductively
to obtain a mapping $H_n$ which conjugates the dynamics
of $F$ and $\tilde F$ up to their first entry times to $\W^n$.
We use the QC Criterion to obtain a 
$K_2$-qc mapping from 
$\W^n$ to $\widetilde{\W}^n$ which agrees with the boundary
marking. If there were no long central cascades, we would
have that $H_n$ is a $(K_0,\delta_0)$-qc\textbackslash bg mapping
on the upper half-plane.
Long central cascades are dealt with in Step 6, where we replace the
mapping that we define by pulling back in Steps 1--5,
with the geometrically defined mapping in the
QC\textbackslash BG Partition of a Central Cascade, 
given by Theorem~\ref{thm:central cascades},
see Figure~\ref{fig:cascrep}.
\begin{figure}
\resizebox{\textwidth}{!}
{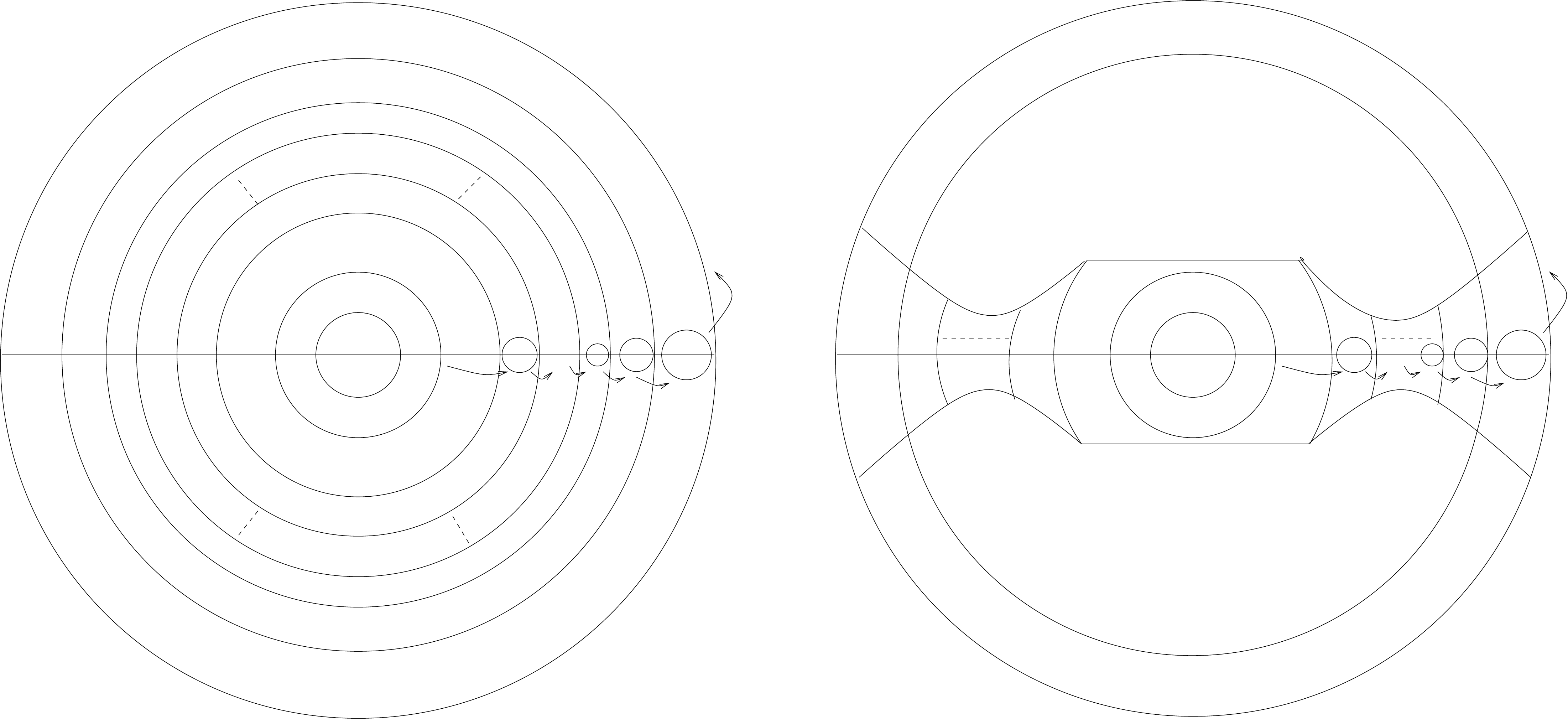}
\caption{Applying Theorem~\ref{thm:central cascades}}
\label{fig:cascrep}
\end{figure}

Let us first give the proof of this result when there
are no long central cascades.

Let $H_0\colon\mathbb C\rightarrow\mathbb C$ be an 
external mapping, so that 
$H_0(\bm{\mathcal{V}})=\widetilde{\bm{\mathcal{V}}},$
$H_0(\bm{\mathcal{U}})=\widetilde{\bm{\mathcal{U}}},$
and for $z\in\partial \bm{\mathcal{U}},$
$H_0\circ F(z)=\tilde F\circ H_0(z).$
For $c\in\Crit(F)$, let 
$\V^0_c=\comp_c\bm{\mathcal{V}},$ and
$\V^0=\cup_{c\in\Crit(F)}\V^0_c$.
Let 
$$\V^0\supset \W^0\supset\V^1\supset\W^1\supset\dots$$
be the good nest of puzzle pieces, see Section~\ref{sec:goodnest}.
Notice that $\W^0=\cup_{c\in\crit(F)}\comp_c(\UU)$.
Recall that by Corollary~\ref{cor:Wnice} and by
Proposition~\ref{prop:goodbg},
we have that there exists $\delta>0$ such that for all $n$
$\W^n$ is $\delta$-nice and each component of 
$\Dom(\W^n)$ has $\delta$-bounded geometry.
We will use these properties to apply the QC Criterion.
Let us show how to pull $H_0$ back through this nest of puzzle pieces.
See Figure~\ref{fig:goodnest1}

\medskip
\textbf{Step 1:  Refine $H_0$ so that it is a conjugacy 
up to the first landing to $\W^0$.} 
For each 
$z$ in $\mathrm{Dom}(\V^0)$ 
(the domain of the first entry mapping to $\V^0$), let $s_z\geq 0$
be minimal so that $F^{s_z}(z)\in\W^0\cup(\V^0\setminus\Dom'(\W^0)),$
where we set $s_z=\infty$ if the orbit of $z$ is contained in  
$\mathrm{Dom}(\V^0)$, but avoids $\W^0$.
The set of points for which $s_z=\infty$ is a hyperbolic
Cantor set of measure zero, and so it is removable.
For each $z$ for which $s_z$ is finite, we define
$\psi_0(z)$ by the formula:
$$H_0\circ F^{s_z}(z)=\tilde F^{s_z}\circ \psi_0(z).$$
Observe that $\psi_0$ is well-defined since the mapping
$F^{s_z}\colon\hat{\mathcal{L}}_z(\W^0)\rightarrow\W^0$ is a
diffeomorphism onto its image, a component of $\W^0.$
For $z\in\mathbb C\setminus\Dom(\V^0),$ we set
$\psi_0(z)=H_0(z).$
Thus, by the complex bounds we have that 
there exists a constant $C>0$ so that
if $H_0$ is $K$-qc, then 
$\psi_0$ is a $K(1+C\mu(\V^0))(1+C\mu(\widetilde{\V}^0))$-qc 
mapping of the plane.

\medskip
\textbf{Step 2. Moving the critical values and pulling back.}
Since $F$ and $\tilde F$ are topologically conjugate,
for each critical point $c$ of $F,$ 
$\psi_0(\hat{\mathcal{L}}_{F(c)}(\W^0))=
\hat{\mathcal{L}}_{\tilde F(\tilde c)}(\widetilde{\W}^0).$
For each component 
$\hat{\mathcal{L}}_{F(c)}(\W^0),$ with $c\in\Crit(F)$,
let $\phi_{c}\colon\mathbb C\rightarrow\mathbb C$ be a quasiconformal
mapping which is the identity on 
$\mathbb C\setminus \hat{\mathcal{L}}_{\tilde F(\tilde c)}
(\widetilde{\W}^0),$
and such that for each critical value $F(c')\in
\hat{\mathcal{L}}_{F(c)}(\W^0),$ we have that
 $\phi_{c}\circ \psi_0(F(c'))=\tilde F(\tilde c')$.
Let $\phi=\phi_{c_1}\circ\cdots\circ\phi_{c_k},$ 
where $\{c_1,\dots,c_k\}\subset\Crit(F)$ is such that
each component of $\Dom'(\W^0)$
 that intersects $F(\Crit(F))$ contains exactly one
of the critical values $F(c_i).$
Then $\phi\circ\psi_0$ may have big quasiconformal distortion. 
However, it still
has quasiconformal distortion bounded by 
$K(1+C\mu(\V^0))(1+C\mu(\widetilde{\V}^0))$
on $\mathbb C\setminus\cup_{i=1}^k\hat{\mathcal{L}}_{F(c_i)}(\W^0).$
We define $\psi_1$ on $\W^0$ by the formula 
$$\phi\circ\psi_0\circ F(z)=\tilde F\circ \psi_1(z).$$
For $z\in\mathbb C\setminus\W^0$, we set
$\psi_1(z)=\psi_0(z)$.
The mapping $\psi_1$ is a 
$K(1+C\mu(\V^0))^{2}(1+C\mu(\widetilde{\V}^0))^{2}$-qc mapping
outside of
$(\Dom'(\W^0)\setminus\W^0)\cup (F|_{\W^0})^{-1}(\Dom'(\W^0))
\supset  \Dom^*(\W^0)$. 

\begin{rem}
Observe that that since $F$ has critical points,
some components
of $F^{-1}(\Dom'(\W^0))$ may
not intersect the real line, and while they have bounded
geometry, they do not necessarily have moduli bounds in the 
upper half-plane, see Figure~\ref{fig:goodnest1}.
We will deal with this in Step 4.
\end{rem}

Let us consider the following multicritical variant of the 
principal nest of puzzle pieces:
Let $\bm T^0=\W^0,$ and for each critical point $c$,
let $\bm T^{i+1}_c=\mathcal{L}_c(\bm T^{i}),$ and set
$\bm T^{i+1}=\cup_{c\in\Crit(F)}\bm T^{i+1}_c.$ 
For each critical point $c$ such that $F(c)\in\bm T^0$,
let 
$N_c\geq 0$ be minimal so that
$F(c)\in\bm T^{N_c}\setminus\bm T^{N_c+1}.$
Let $N=\max_{c\in\Crit(F)}N_c$.

\medskip

\textbf{Step 3: A qc conjugacy up to the first return to 
$\V^1$.}  We will show that
there exist $K_1\geq 1,$ given by Proposition~\ref{prop:squares},
 and a qc mapping 
$\Psi_1\colon\mathbb C\rightarrow\mathbb C$ such that 
$\Psi_1(\V^1)=\widetilde{\V}^1$,
$\Psi_1 (\Dom(\V^1))=\Dom(\widetilde{\V}^1),$
and for each component 
$\U$ of $\Dom(\V^1),$ if $z\in\partial\U$,
$\Psi_1\circ F(z)=\tilde F\circ \Psi_1 (z)$. Moreover, $\Psi_1$ is 
$K_1$-qc outside
$\Dom(\V^1)\cup\bm X^1,$ where 
$\bm X^1\subset\mathbb C\setminus \mathbb R.$

\begin{rem} We will give an explicit description of $\bm X^1$.
For now, let us just point out that
it consists of certain pullbacks of landing domains which do not
intersect the real line. 
\end{rem}

We have that $\psi_1$ conjugates the dynamics of 
$F$ and $\tilde F$ up to the first entry of each point to
$\W^0$. 
Suppose that $c$ is a critical point of $F$ with the property that
$F(c)\in\V^0\setminus\W^0.$ 
Then $\mathcal{L}_{c}(\W^0)\subset\V^1_c\subset\W^0_c.$
Since $\phi$ is the identity outside 
$\cup_{i=1}^k\hat{\mathcal{L}}_{\tilde F(\tilde
  c_i)}(\widetilde\W^0),$
we have that $\psi_1(\V^1_c)=\widetilde{\V}^1_{\tilde c}.$
However, if $F(c)\in\W^0,$ then this need not be
the case. To obtain this, we will 
pullback through the nest 
$$\bm T^0\supset\bm T^1\supset \bm T^2\supset\dots.$$ To
simplify the notation, we will do more 
(that is, pull back further) than we
strictly need to at certain critical points,
but this will not cause problems for us.

Let $\phi_0=\psi_1$. As in Step 1, 
refine $\phi_0$ so that it is a conjugacy up to the first 
landing to $\bm T^1$:
for each $z\in\Dom(\bm T^0),$ let $s_z>0$ be minimal
so that $F^{s_z}(z)\in \bm T^1\cup (\bm T^0\setminus\Dom(\bm T^0))$, and
for
$z\in\Dom(\bm T^0)$ define $\phi_0'(z)$ so that
$$\tilde F^{s_z}\circ\phi_0'(z)=\phi_0\circ F^{s_z}(z).$$
As before $s_z$ is finite except 
at a set of points which is conformally removable.

\begin{figure}
\resizebox{0.5\textwidth}{!}
{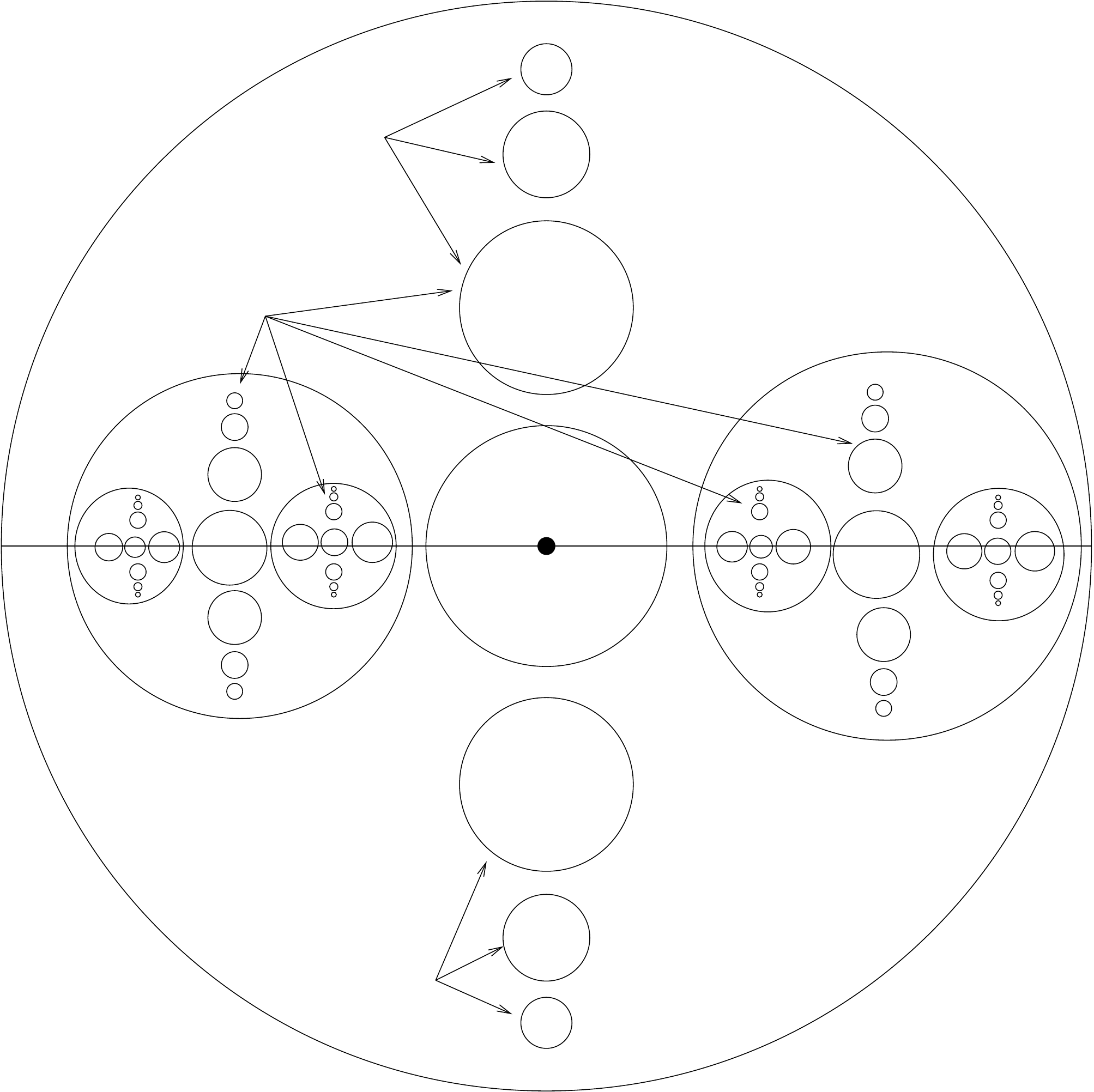}
\caption{Complex pullbacks of landing domains.}
\label{fig:complexdomains}
\end{figure}

Let $$\bm A_0^0\subset \bm T^0\setminus\Dom(\bm T^0)$$
be the union of components of 
$(F|_{\W^0})^{-1}(\Dom'(\W^0))$ that do not intersect
the real line. 
Let $$\bm B_0=\cup_{i=0}^{\infty}(F|_{\Dom'(\bm T^0)\setminus\bm
  T^1})^{-i}(\bm A_0^0).$$ See Figure~\ref{fig:complexdomains} for the
first few pullbacks of $\bm A^0_0.$
For any $i\geq 0,$ the set $(F|_{\Dom'(\bm T^0)\setminus\bm
  T^1})^{-i}(\bm A_0^0)$ is contained in 
$(F|_{\Dom'(\bm T^0)\setminus\bm
  T^1})^{-i}(\bm T^0).$ Each component 
$\bm T'$ of $(F|_{\Dom'(\bm T^0)\setminus\bm
  T^1})^{-i}(\bm T^0)$ intersects the real line, and
the mapping $F^i|_{\bm T'}$ is a diffeomorphism with
 dilatation bounded
independently of $i$ by Proposition~\ref{prop:squares}.

Since $\psi_1$ is $K(1+C\mu(\V^0))^{2}(1+C\mu(\widetilde\V^0))^2$-qc outside of
$(\Dom'(\W^0)\setminus\W^0)\cup (F|_{\W^0})^{-1}\Dom'(\W^0),$
we have that
$\phi_0'$ conjugates the dynamics of 
$F$ and $\tilde F$ up to the first landings 
to $\bm T^1\cup (\bm T^0\setminus\Dom(\bm T^0))$ and it is
$K(1+C\mu(\V^0))^{3}(1+C\mu(\widetilde{\V}^0))^{3}$-qc outside of
$\Dom'(\bm T^1)\cup  \bm B_0$.
Since $F$ and $\tilde F$ are topologically
conjugate, we have that $\phi_0'(F(c))$ 
and $(\tilde F(\tilde c))$ are contained in the same
component of $\Dom'(\widetilde{\bm T}^1)$.
As in Step 2, we modify $\phi_0'$
in each component of 
$\Dom'(\bm T^1)$ that contains a critical value of $F$,
so that $\phi_0'(F(c))=\tilde F(\tilde c)$ for each critical 
point $c$ of $F$. Abusing notation, we will use $\phi_0'$
to denote the modified mapping.
It is important to observe that if
$\U=\comp_{F^{k_c}(c)}(F^{-(l_c-k_c)}(\V^0)),$ where
$c\in\Crit(F),$  then since $\U$ is a pullback of $\V^0$ and
we are only modifying $\phi_0'$ in 
$\Dom'(\bm T^1)$, 
we have that $\phi_0'(\U)=\widetilde{\U}$.

For each $z\in\bm T^1$, let $\phi_1(z)$ be defined by 
$$\phi'_0\circ F=\tilde F\circ\phi_1,$$
and for $z\in\mathbb C\setminus\bm T^1$, set
$\phi_1(z)=\phi_0'(z)$.
Observe that $\phi_1$ is a qc mapping of the plane that is
$K(1+C\mu(\V^0))^{4}(1+C\mu(\widetilde{\V}^0))^{4}$-qc outside of  
$$(\Dom'(\bm T^1)\setminus\bm T^1)\cup \bm B_0\cup 
(F|_{\bm T^1})^{-1}(\Dom'(\bm T^1)\cup \bm B_0)\supset\Dom(\bm T^1),$$
$\phi_1(\bm T^1)=\widetilde{\bm T}^1,$ $\phi_1(\Dom(\bm T^1))
=\Dom(\widetilde{\bm{T}}^1),$
and for any component $\U$ of $\Dom(\bm T^1)$,
if $z\in\partial\U$, $\phi_1\circ F(z)=\tilde F\circ\phi_1(z)$.

Let $\bm A_1^0$ denote the union of the 
components of $(F|_{\bm T^1})^{-1}(\Dom'(\bm T^1))$ that do not intersect the
real line.
Let 
$$\bm A_1=\bm B_0\cup (F|_{\bm T^1})^{-1}(\bm B_0)\cup \bm A_1^0,$$ 
and let 
$$\bm B_1=\cup_{i=0}^{\infty}(F|_{\Dom(\bm T^1)\setminus\bm T^2})^{-i}(\bm A_1).$$
Then after refining $\phi_1$, to obtain $\phi_1',$
which conjugates the dynamics
of $F$ and $\tilde F$ up to the first landing to
$\bm T^2\cup(\bm T^1\setminus\Dom(\bm T^1)),$
we have that $\phi_1'$ is $K(1+C\mu(\V^0))^{5}(1+C\mu(\widetilde{\V}^0))^{5}$-qc
outside of $\Dom'(\bm T^2)\cup \bm B_1$. 

Repeating this argument for each $i\leq N+1$, we obtain
a qc mapping 
\mbox{$\phi_{N+1}'\colon\mathbb C\rightarrow\mathbb C$} 
with the following properties:
\begin{itemize}
\item $\phi_{N+1}'$ is $K_1$-qc outside of 
$\Dom'(\bm T^{N+2})\cup\bm B_{N+1}$, where $K_1$ depends on $N$ and
$\mu(\V^0)$, see 
Proposition~\ref{prop:squares};
\item for $0\leq j\leq N+2$,
$\phi_{N+1}'(\bm T^{j})=\widetilde{\bm T}^{j},$ 
$\phi_{N+1}'(\Dom'(\bm T^{j}))=\Dom'(\widetilde{\bm T}^{j});$
\item for $0\leq j\leq N+1$,
$\phi_{N+1}'(\Dom(\bm T^{j}))=\Dom(\widetilde{\bm{T}}^{j});$
and 
\item if $\U$ is a component of $\Dom(\bm T^j)\cup\Dom'(\bm T^{N+2}),$
for $0\leq j\leq N+1,$ and $z\in\partial\U$, 
then $\phi_{N+1}'\circ F(z)=\tilde F\circ\phi_{N+1}'(z).$ 
\end{itemize}

Since we only ever modify, that is, move, the 
points $\phi_j'(F(c)),$ for $c\in\Crit(F)$ and
$0\leq j\leq N+1,$
inside of $\Dom'(\bm T^j),$ we have that 
for each 
$c\in\Crit(f),$ $\bm T^{N+2}_c\subset \V^1_c,$
$\phi_{N+1}'(\V^1)=\widetilde{\V}^1,$
and if $\U$ is a component of $\Dom'(\V^1)$
and $z\in\partial\U$, 
then $\phi_{N+1}'\circ F(z)=\tilde F\circ\phi_{N+1}'(z).$ 

Let 
$\psi_2(z)=\phi_{N+1}(z).$ 
Let us point out that if every critical
value of $F$ is contained in $\V^0\setminus\W^0$,
then $\psi_2=\psi_1$.
We have that
$\psi_2\colon\mathbb C\rightarrow\mathbb C$ and
$\psi_2(\V^1)=\widetilde{\V}^{1}.$
Let $\bm T=\cup_{c\in\Crit(F)} \bm T^{N_c+1}_c.$
We have that $\psi_2$ conjugates 
$F$ and $\tilde F$
up to the first entry mapping to
$\bm T$.
Now, refine $\psi_2$ in $\Dom(\bm T)$
as in Step 1, so that it
is a conjugacy up to the first landing time to
$\V^1$. 
Let $\bm B^1=\bm B_{N+2}\setminus\Dom'(\V^1)$.
As in Step 2, since $F$ and $\tilde F$ 
are topologically conjugate, corresponding
critical values are contained in corresponding components
of $\Dom(\V^1),$ moving each $\psi_2(F(c))$ so that it agrees with
$\tilde F(\tilde c),$ by a qc mapping which is the identity outside of
$\comp_{\tilde F(\tilde c)}(\Dom(\widetilde{\V}^1))$,
 and pulling back, we obtain a 
a qc mapping $\Psi_1\colon\mathbb C\rightarrow\mathbb C$ which is 
$K_1$-qc outside of 
$\Dom(\V^1)\cup \bm X^1$,
where $\bm X^1=\bm B^1\cup (F|_{\V^1})^{-1}(\bm B^1)$
and $K_1=K_1(N,\mu(\V^0))\geq 1,$ see Proposition~\ref{prop:squares}. 
Moreover, $\Psi_1(\Dom(\V^1))=\Dom(\widetilde{\V}^1)$,
and $\Psi_1$ conjugates $F$ and $\tilde F$ on $\partial\U$ for
each component $\U$ of $\Dom(\V^1)$. 

\medskip

\textbf{Step 4: A qc\textbackslash bg pseudo-conjugacy.}
In this step, we deal with the issue concerning moduli bounds pointed
out in the remark in Step 2.
Let $\U$ be a component of 
$\bm A^0_0$, so that $\V=F(\U)$ is a component of
$\Dom'(\W^0)$. 
We have that $\psi_1|_{\W^0}$ is 
$K(1+C\mu(\V^0))^{2}(1+C\mu(\widetilde{\V}^0))^{2}$-qc on
$\W^0\setminus (F|_{\W^0})^{-1}(\Dom'(\W^0))$
(where $K$ still denotes the dilatation of the external conjugacy).
By Corollary~\ref{cor:Wnice}, we have that 
there exists $\delta>0$, so that for each critical 
point $c$ of $F$ and each component $\U$ of 
$(F|_{\W^0_c})^{-1}(\Dom'(\W^0))$ we have that
$\mod(\W^0_c\setminus\U)\geq \delta.$
Let $s\in\mathbb N$ be so that 
$F^s(\U)$ is a component of $\W^0$. Let
$\psi_3|_{\U}$ be defined by
$$\psi_1 \circ F|_{\U}^s(z)=\tilde F^{s}\circ \psi_3(z).$$
The mapping $\psi_3|_{\U}$ is 
$K(1+C\mu(\V^0))^{3}(1+C\mu(\widetilde{\V}^0))^{3}$-qc on 
$$\U\setminus (F|_{\U})^{-s}(F|_{\W^0_c})^{-1}(\Dom'(\W^0))),$$
and there exists $\varepsilon_0'>0$ so that each component of
$\U'$ of $ (F|_{\U})^{-s}(F|_{\W^0_c})^{-1}(\Dom'(\W^0)))$
has $\varepsilon_0'$-bounded
geometry and $\mod(\U\setminus\U')\geq \varepsilon_0'.$

Thus we have $\psi_3$ defined on $\bm A^0_0$,
where it is a 
$(K(1+C\mu(\V^0))^{3}(1+C\mu(\widetilde\V^0))^{3},\varepsilon_0')$-qc\textbackslash bg
mapping.
For any component $\V$ of $\bm B_0,$ let 
$j\geq 0$ be so that $F^j(\V)$ is a component of 
$\bm A^0_0$, and
define $\psi_3$ on $\V$ by
$\psi_3\circ F^j|_{\V}(z)=\tilde F^j\circ \psi_3(z)$.
This is well defined since $F^j|_{\V}$ is a 
diffeomorphism as $\V\cap\mathbb R=\emptyset$.
Moreover, since we are pulling $\psi_3$ back
along real diffeomorphic branches of $F$, 
we have there exists $\varepsilon_0>0$ so that 
that the resulting mapping is 
$(K(1+C\mu(\V^0))^{4}(1+C\mu(\widetilde{\V}^0))^{4},\varepsilon_0)$-qc\textbackslash bg.

Let $\U$ be a component of $\bm A^0_1$
such that if $\hat \U$ is the component of $F^{-1}(\Dom'(\W_0))$
which contains $\U$, then $\hat \U$ intersects the real line,
see Figure~\ref{fig:wellinside}.
(In the previous step, we defined $\psi_3$ on the components
of $(F|_{\W^0})^{-1}(\Dom'(\W_0))$ which do not intersect the real
line.)
Let $s\in\mathbb N$ be so that 
$F^s(\U)$ is a component of $\bm T^1$.
We have that $\phi_1|_{\bm T^1}$ is 
$K(1+C\mu(\V^0))^{4}(1+C\mu(\widetilde\V^0))^{4}$-qc on 
$ \bm T^1\setminus((F|_{\bm T^1})^{-1}(\Dom'(\bm T^1)\cup \bm B_0)).$
Recall that on each component of $\bm B_0,$ 
$\psi_3$ is $(K(1+C\mu(\V^0))^{4}(1+C\mu(\widetilde\V^0))^{4},\varepsilon_0)$-qc\textbackslash bg.
We define $\psi_3|_{\U}$ by the formula:
$$\phi_1\circ F^s(z)=\tilde F^s\circ\psi_3(z),\mbox{ for } z\in\U.$$
Then $\psi_3|_{\U}$ is $K(1+C\mu(\V^0))^{5}(1+C\mu(\widetilde\V^0))^{5}$-qc on
$\U\setminus (F|_{\U})^{-s}((F|_{\bm T_1})^{-1}(\Dom'(\bm T^1)\cup\bm B_0)).$ 
Moreover, there exists $\varepsilon_1>0$ so that 
each if $\U'$ is a component of the set
$$(F|_{\U})^{-s}((F|_{\bm T^1})^{-1}(\Dom'(\bm T^1))),$$ then $\U'$ has
$\varepsilon_1$-bounded
geometry and $\mod(\U\setminus\U')\geq \varepsilon_1$.
The bound on the geometry follows from 
Lemma~\ref{lem:KvS bg lem 2}, and we obtain the moduli bounds 
from Corollary~\ref{cor:Wnice} and Lemma~\ref{lm:pullbackspace}.
We extend $\psi_3$ to all of $\bm B_1$ by pulling back through 
diffeomorphic branches of $F$ - just as we did to extend it from 
$\bm A^0_0$ to $\bm B_0$.

We repeat this argument for each component of
$(F|_{\bm T^j})^{-1}(\Dom'(\bm T^j))$, $0\leq j\leq N+2$. 
This processes gives us
a $\varepsilon_{N+2}>0,$ depending on $N$, and a 
$(K_1,\varepsilon_{N+2})$-qc\textbackslash bg mapping
$\psi_3\colon\bm X^1\rightarrow\widetilde{\bm{X}}^1,$
where $K_1\geq 1$ is given by Proposition~\ref{prop:squares},
with the following properties.
Suppose that $\U$ is a component of 
$\bm X^1,$ then $\U$ there exists $s\in\mathbb N,$
$c\in\Crit(F)$
and $0\leq j\leq N+2$ so that
$F^s(\U)=\bm T^j_c.$  
Then $\psi_3|_{\U}$ is
$K_1$-qc outside of 
$F^{-s}(F^{-1}(\Dom'(\bm T^j)\cup \bm B_{j-1})),$ and each component
$\U'$ of $F^{-s}(F^{-1}(\Dom'(\bm T^j)\cup \bm B_{j-1}))$ has
$\varepsilon_{N+2}$-bounded geometry and 
$\mod(\U\setminus\U')\geq\varepsilon_{N+2}$.

Define $H_1\colon\mathbb C\rightarrow\mathbb C$ by
$H_1(z)=\Psi_1(z)$ for $z\in\mathbb C\setminus \bm X^1$ and
$H_1(z)=\psi_3(z)$ for $z\in\bm X^1$.
Then $$H_1|_{\mathbb C\setminus\Dom(\V^1)}\colon \mathbb
C\setminus\Dom(\V^1)\rightarrow 
\mathbb C\setminus\Dom(\widetilde{\V}^1)$$ is a 
$(K_1,\varepsilon_{N+2})$-qc\textbackslash bg
mapping.

\medskip

\textbf{Step 5: The inductive step.}
When the the length of long central cascades is bounded, we have that
there exists $\delta_0>0$
such that the following holds:
Repeating Steps 1 through 4 for
each $n\in\mathbb N\cup\{0\}$ we have there is a 
$(K_1,\delta_0)$-qc\textbackslash bg mapping:
$$H_n|_{\mathbb H^+\setminus\Dom'(\W^n)}\colon \mathbb
H^+\setminus\Dom'(\W^n)\rightarrow \mathbb 
H^+\setminus\Dom'(\widetilde{\W}^n),$$
such that $H_n$ agrees with the boundary marking on
each component of $\Dom'(\W^n).$
By the QC Criterion and the complex bounds (see Step 4),
we have that 
there exists $K_2=K_2(K_1,\delta, N)\geq 1$, and $K_2$-qc mapping
$\Psi_2: \W^n\rightarrow
\tilde{\W^n},$
which agrees with the boundary marking on each component of
$\W^n.$ 
Let $\Psi_3$ be the mapping from $\Dom'(\W^n)$ to
$\Dom'(\widetilde{\W}^n)$ obtained by 
pulling this mapping back by the
first entry mapping to $\W^n$. 
Define $H_{\W^n}(z)=H_n(z)$ for $z\in\mathbb C\setminus \Dom(\W^n)$,
and $H_{\W^n}(z)=\Psi_3(z)$ for $z\in\Dom(\W^n)$.
Then we have that there exists $K_0\geq 1$ so that
$H_{\W^n}:\mathbb H^+\rightarrow\mathbb H^+$ is a 
$(K_0,\delta_0)$-qc\textbackslash bg mapping.

Let us just point out that to start the construction of the set off
the real line with bounded geometry at deeper levels we take
$\bm A^0_n$ to be the set $\bm X^{n-1}$ together with the components of 
$(F|_{\W^n})^{-1}(\Dom'(\W^0))$ which do not intersect the real line.
By Proposition~\ref{prop:good bounds}, we have that there exists a
constant $C_1\geq 1,$ so that for any component $\U$ of 
$\cup_{n=1}^\infty\bm X^n,$ if $j>0$ is so that $F^{j}(\U)$
is a component of $\bm A^0_k$ for some $k\geq 0$, then 
$F^k|_{\U}$ is $C_1$-qc, and thus
$\psi_3|_{\U}$ (constructed in Step 4) is a $(K_1,\delta_1)$-qc\textbackslash bg mapping.  

\medskip
\textbf{Step 6: Long central cascades.}
When there is a long central cascade we are unable to control
the dilatation of the
$\phi_n$ constructed in Step 3. To deal with
this problem, we replace the mapping $\phi_n,$ constructed by pulling back, 
with the mapping $H$ given by  Theorem~\ref{thm:central cascades},
which gives us the QC\textbackslash BG Partition of the Central Cascade.
The reason that we cannot do this automatically is that the puzzle 
pieces, $\bm{\mathcal{G}}_c(k)$ are not necessarily contained in
the good nest. We will prove
that  there exist $K_3\geq 1, \delta_3>0,$ 
which do not depend on $c$ or $k$,
such that
for any $\bm{\mathcal{G}}_c(k)$
at the top of a long central cascade, there is a mapping
$\hat H\colon\mathbb C\rightarrow\mathbb C,$
such that $\hat H\colon\mathbb H\rightarrow\mathbb H$ is a
$(K_3,\delta_3)$-qc\textbackslash bg mapping,
$\hat H (\bm{\mathcal{G}}_c(k))=\widetilde{\bm{\mathcal{G}}}_c(k)$,
$\hat H (\mathcal{L}_c(\bm{\mathcal{G}}_c(k)))=
\mathcal{L}_{\tilde c}(\widetilde{\bm{\mathcal{G}}}_{\tilde c}(k)),$
and $\hat H$ respects the boundary marking on 
$\bm{\mathcal{G}}_c(k)$ and $\mathcal{L}_c(\bm{\mathcal{G}}_c(k))$. Then Theorem~\ref{thm:central cascades}
yields a mapping from 
$\bm{\mathcal{G}}_c(k)$ to 
$\widetilde{\bm{\mathcal{G}}}_{\tilde c}(k),$
as in the conclusion of Theorem~\ref{thm:central cascades}.
Once we have this mapping we spread it around dynamically.

Let $n\in\mathbb N$, and let $H_{\W^n}$ be the 
mapping constructed 
in Steps 1--5. We have that
$$H_{\W^n}\colon\mathbb C\rightarrow\mathbb C$$ is such that
\begin{itemize}
\item $H_{\W^n}(\W^n)=\widetilde{\W}^n;$
\item $H_{\W^n}(\Dom(\W^n))=\Dom(\widetilde{\W}^n);$ and
\item for any component $\U$ of $\Dom(\W^n)$,
$H_{\W^n}\circ F(z)=\tilde
  F\circ H_{\W^n}(z),$ for $z\in\partial{\U}$.
\end{itemize}

\medskip
\noindent\textbf{Step 6a.}
Now, let us construct a neighbourhood of the post-critical set as
follows: for each critical point $c_0$ of $F$,
let $\bm{\mathcal{O}}_{c_0}\owns c_0$ be the smallest puzzle piece
from the good nest with the following property:
Suppose that
$\{U_j\}_{j=0}^s$ is the chain with $\U_0=\bm{\mathcal{O}}_c,$
and $\U_s=\V^0_{c}$ for some $c\in\Crit(F)$.
If $\U_j\owns F(c')$, for some $c'\in\Crit(F)$, then
$H_{\W^n}(F(c'))=\tilde F(\tilde c')$.
Let $\bm{\mathcal{O}}=\cup_{c}\bm{\mathcal{O}}_c.$
Note that it follows from equation~\ref{eqn:claim2},
we have that as $n\rightarrow \infty,$
for each $c\in\Crit(F)$, $\diam(\bm{\mathcal{O}}_c)\rightarrow 0$.

Let us regard $H_n\colon\V^0\rightarrow\widetilde{\V}^0$ 
as a conjugacy between $F\colon\W^0\rightarrow \V^0$ and 
$\tilde F\colon\widetilde{\W}^0\rightarrow\widetilde{\V}^0$.

\medskip
\noindent\textbf{Step 6b.}
Set $\hat H_0=H_{\W^n}$.
We pullback $\hat H_0$ as follows:
Let $\bm{\mathcal{Y}}^0=\V^0,$
 $\bm{\mathcal{Y}}^1=\W^0.$
Assuming that $\bm{\mathcal{Y}}^i$ has been defined,
we define $\bm{\mathcal{Y}}^{i+1}$ to be the set of all connected
components
$(F|_{\bm P})^{-1}(\bm Q)$ where $\bm P$ is a puzzle piece in 
$\bm{\mathcal{Y}}^i$ which intersects the real line and 
 $\bm Q$ is a puzzle piece in 
$\bm{\mathcal{Y}}^i$.

For each component $\U$ of $\W^0$,
with the property that if $\U$ contains a critical value $F(c')$,
then $\comp_{c'} F^{-1}(\U)\supset\bm{\mathcal{O}}_c,$
define $\hat H_1$ on $\U$ by
$\tilde F\circ \hat H_1(z)=\hat H_0\circ F(z)$.
Assuming that $\hat H_i$ is defined on $\bm{\mathcal{Y}}^i$, 
on each component $\bm P$ of $\bm{\mathcal{Y}}^{i+1}$,
with the property that if $\bm P$ contains a critical value $F(c')$,
then $\comp_{c'} F^{-1}(\U)\supset\bm{\mathcal{O}}_c,$ we
define $\hat H_{i+1}$ by
$\tilde F\circ \hat H_{i+1}(z)=\hat H_i\circ F(z)$.
By the choice of $\bm{\mathcal{O}}$ 
we can pullback $\hat  H_0=H_n\colon\mathbb C\rightarrow\mathbb C$ to
obtain $\hat H:\mathbb C\rightarrow\mathbb C$ such that
$\hat H(\bm{\mathcal{O}})=\widetilde{\bm{\mathcal{O}}}$.
Now, arguing as in Step 1, we may assume that $\hat H$ is a conjugacy up to the first
landing mapping to $\bm{\mathcal{O}}$.
By Proposition~\ref{prop:squares}, the dilatation of the mapping 
from any component of $\bm{\mathcal{O}}$ onto a component of 
$\V^0$ is bounded. 

\medskip
\noindent\textbf{Step 6c.}
In this step we redefine $\hat H$ inside of central cascades
so that it will be a qc\textbackslash bg
mapping.
By Prop~\ref{prop:good geometry for central cascades},
for any critical point $c$ and any
$\mathcal{C}^k_c(\V^0_c)\supset\bm{\mathcal{O}} $
that is at
the top of a long central cascade ( with length at least
$12b^2-4b,$ where $b$ is the number of critical
points of $F$)
there exists a $\delta$-excellent puzzle piece 
$\bm{\mathcal{G}}_c(k)$ at the ``top'' of the central cascade.

Moreover, unless $\V^n_c$ is at the top of a long central cascade and 
$\bm{\mathcal{G}}_c(k)$ is at the top of the same long central
cascade, we have that $\hat H$ conjugates the dynamics of $F$ and 
$\tilde F$ on the
boundary of a necklace neighbourhood of the real line
contained in
$\bm{\mathcal{N}}_c(k)$, since $\hat H$ conjugates the dynamics of $F$
and $\tilde F$ up to the first landing time to $\bm{\mathcal{O}}$.
By Theorem~\ref{thm:central cascades}, we can redefine $\hat H$ from 
$\bm{\mathcal{G}}_c(k)\setminus \bm{\mathcal{N}}_c(k)
\rightarrow \widetilde{\bm{\mathcal{G}}}_c(k)\setminus \widetilde{\bm{\mathcal{N}}}_c(k)$
obtain a mapping which is $(K',\delta')$-qc\textbackslash bg on $\bm{\mathcal{G}}_c(k)$.

To conclude the proof, we almost repeat the 
argument from Step 6b, except here we do not refine $\hat H$
in $\bm{\mathcal{G}}_c(k)\setminus \bm{\mathcal{N}}_c(k).$
Instead we use the mapping given by 
Theorem~\ref{thm:central cascades}, and then pull this mapping back.
We have that
 $\hat H\colon\V^0\rightarrow\widetilde{\V}^0$ 
is a conjugacy between $F\colon\W^0\rightarrow \V^0$ and 
$\tilde F\colon\widetilde{\W}^0\rightarrow\widetilde{\V}^0$.
Let us abuse notation and redefine the sequence $\hat H_0.$
Let $\hat H_0=\hat H$.
Assuming that $\hat H_i$ is defined on $\bm{\mathcal{Y}}^i$, 
on each component $\bm P$ of $\bm{\mathcal{Y}}^{i+1}$,
with the property that if $\bm P$ contains a critical value $F(c')$,
then $\comp_{c'} F^{-1}(\U)\supset\bm{\mathcal{O}}_c,$ we
define $\hat H_{i+1}$ as follows: 
If neither $\U$ nor the smallest
real puzzle piece that contains $\U$ is deep in a central cascade,
then for $z\in\U$, define $\hat H_{i+1}$ so that
$\tilde F\circ \hat H_{i+1}(z)=\hat H_i\circ F(z)$.
If either $\U$ or the smallest
real puzzle piece that contains $\U$ is deep in a central cascade with
$\bm{\mathcal{G}}_c(k)$ as a top level puzzle piece,
then set  $\hat H_{i+1}(z)=\hat H_0(z)$ 
if  $z\notin \bm{\mathcal{N}}_c(k),$
and if $z\in \bm{\mathcal{N}}_c(k),$
then define $H$ so that
$\tilde F\circ \hat H_{i+1}(z)=\hat H_i\circ F(z)$.
Thus we
obtain a neighbourhood $\bm{\mathcal{O}}$ of $\Crit(F),$
and a $(K_1,\delta_1)$-qc\textbackslash bg mapping
$\hat H_{\bm{\mathcal{O}}}\colon\mathbb C\setminus\Dom(\bm{\mathcal{O}})
\rightarrow \mathbb C\setminus\Dom(\tilde{\bm{\mathcal{O}}}).$

Before concluding the proof, let us point out that we
by the QC Criterion, we have the following
Corollary.

\begin{cor}
For any $N_0$,
there exists $K\geq 1$ such that for any puzzle piece $\bm E$ that is
either from the enhanced
nest, $E=\W^n_c$, from the good nest,
or that is obtained as a pullback of a component of $\bm{\mathcal{V}}$
with order bounded by $N_0$, there exists a $K$-qc mapping from $\bm E$
to $\widetilde{\bm{E}}$
that agrees with the boundary marking. 
\end{cor}

\medskip
\noindent\textbf{Step 6d.}
By the corollary, we have
that there exist $K_3\geq 1$ and arbitrarily small
neighbourhoods $\bm{\mathcal{O}}$
of $\Crit(F)$ such that for each component  
$\bm{\mathcal{O}}_c$ of $\bm{\mathcal{O}}$,
there is a 
$K_4$-qc mapping from 
$\bm{\mathcal{O}}_c\rightarrow \widetilde{\bm{\mathcal{O}}}_c,$
which respects the boundary marking.
Thus redefining $\hat H_{\bm{\mathcal{O}}}$ in $\mathrm{Dom}'(\bm{\mathcal{O}})$,
by pulling back this mapping by the landing map, we 
obtain a $(K_0,\delta_0)$-qc\textbackslash bg mapping
$H_{\bm{\mathcal{O}}}$
of the $\mathbb H$ which agrees with the 
standard boundary marking on 
$\Dom'(\bm{\mathcal{O}})$.
\end{pf}


By the QC Criterion and compactness of the space of $K$-qc mappings we have:
\begin{thm}[Quasisymmetric rigidity for box
  mappings]\label{thm:smooth box conj}
Suppose that $F\colon\bm{\mathcal{U}}\rightarrow\bm{\mathcal{V}}$
$\tilde F\colon\widetilde{\bm{\mathcal{U}}}\rightarrow\widetilde{\bm{\mathcal{V}}}$
are non-renormalizable, strongly combinatorially equivalent box mappings
that are given by either Theorem~\ref{thm:box mapping persistent infinite branches}
or Theorem~\ref{thm:box mapping reluctant}, which in particular
satisfy the conditions of Theorem~\ref{thm:external conjugacy}.
Then restricted to the real line,
the mappings $F\colon\mathcal{U}\rightarrow\mathcal{V}$
and $\tilde F\colon\tilde{\mathcal{U}}\rightarrow\tilde{\mathcal{V}},$ are 
quasisymmetrically conjugate.
\end{thm}
\begin{pf}
We have that there exist arbitrarily small neighbourhoods 
$\bm{\mathcal{O}}$ of $\Crit(F)$ and 
$(K_0,\delta_0)$-qc\textbackslash bg mappings $H:\mathbb H^+
\rightarrow\mathbb H^+$ such that $H$ conjugates $F$
and $\tilde F$ up to the first landing times to 
$\bm{\mathcal{O}}.$ By the QC Criterion, 
there is a $K_1$-qc mapping $\hat H,$ 
which agrees with the restriction of $H$ to 
$\mathbb R$. By the compactness of the space of 
$K_1$-qc mappings, as we take smaller neighbourhoods of 
$\Crit(F),$
there exists a limiting $K_1$-qc mapping which conjugates $F$ and 
$\tilde F$ on their real traces.

\end{pf}

\subsection{Central cascades, the proof of a QC\textbackslash BG
  Partition of a Central Cascade}
\label{sec:cascades}

We will require the following partial ordering
on central cascades:
we say that a central cascade 
$Z^0\supset Z^1\supset\dots\supset Z^{N_1}\supset Z^{N_1+1}\owns c_1$
is \emph{deeper} than a central cascade 
$U^0\supset U^1\supset\dots\supset U^{N_2}\supset U^{N_2+1}\owns c_2$
if $\mathcal{L}_{c_1}(U_0)\supset Z^1.$
Suppose that
$Z^0\supset Z^1\supset\dots\supset Z^{N_1}\supset Z^{N_1+1}\owns c_1$
and
$U^0\supset U^1\supset\dots\supset U^{N_2}\supset U^{N_2+1}\owns c_2$
are two cascades so that neither is deeper than the other.
Then $\mathcal{L}_{c_1}(U^0)\subset Z^1,$ and 
$\mathcal{L}_{c_2}(Z^0)\subset U^1.$
But now, if $x\in Z^0\setminus Z^1$,
then $\mathcal{L}_x(U_0)\subset Z^0\setminus Z^1$ and if
$x\in U^0\setminus U^1$,
then $\mathcal{L}_x(Z_0)\subset U^0\setminus U^1.$
Because of this observation, it will not matter in which order 
we treat central cascades that are unrelated by the partial ordering.
We will prove Theorem~\ref{thm:central cascades}
starting with the shallowest (least deep) central cascades,
and proceeding to deeper central cascades.

We will need the following lemmas:

\begin{lem}\label{lem:central cascade Schwarzian}
If $Z^0$ is sufficiently small, the first return mapping 
$R\colon Z^2\rightarrow Z^1$ has negative Schwarzian derivative.
\end{lem}
\begin{pf}
We will decompose the mapping into a composition of
polynomials and
diffeomorphisms.
Let $\{G_j\}_{j=0}^s$ denote the chain with $G_s=Z^1$ 
and $G_0=Z^2$. Let $0=s_1<s_1<\dots<s_{b-1}<s_b=s$, denote
the indicies such that $G_{s_i}$ contains a critical point.
Label the critical point in $G_{s_i}$ by $c_i$.
We can express $R|_{Z^2}$  as
$$F^{s_b-s_{b-1}}\circ F^{s_{b-1}-s_{b-2}}\circ \dots \circ
F^{s_2-s_1}\circ F^{s_1}|_{Z^2}\colon Z^2\rightarrow Z^1.$$
For each $F^{s_i-s_{i-1}}|_{G_{s_{i-1}}},$ $i=1,2,\dots, b,$
we estimate the Schwarzian derivative as usual.
$$SF^{s_i-s_{i-1}}(x)=SF^{s_{i}-s_{i-1}-1}(F(x))F'(x)^2+SF(x).$$
For $x\in G_{s_{i-1}}$,
we have that 
$$SF(x)\asymp -\frac{1}{(x-c_i)^2}\leq -\frac{1}{|G_{s_{i-1}}|}.$$
For $x\in F(G_{s_{i-1}})$,
$$ SF^{s_i-s_{i-1}-1}(F(x))=\sum_{i=0}^{s_{i}-s_{i-1}-2}
SF(F^{i+1}(x))((F^{i})'(x))$$
$$\leq \frac{C}{|G_{s_{i-1}+1}|^2}\max_{y\in\hat
  Z^0}SF(y)\sum_{s=s_{i-1}+2}^{s_i}|G_j|^2.$$
So we have that
$$SF^{s_i-s_{i-1}}(x)\leq \frac{C}{|G_{s_{i-1}+1}|^2}\cdot\max_{y\in
  \hat
  Z^0}SF(y)\sum_{j=s_{i-1}+1}^{s_i}|G_j|^2\Big(\frac{|G_{s_{i-1}+1}|}{|G_{s_i-1}|}\Big)^2-\frac{1}{|G_{s_{i-1}}|^2}.$$
This is negative when $\sum_{j=s_{i-1}+1}^{s_i}|G_j|^2$ is small,
and this holds when $|Z^0|$ is sufficiently small.
\end{pf}

\begin{lem}
\label{lem:CvST 3.10}
For any $\delta>0$ there exist $\kappa>0$ and $C>0$ with the following
properties.
Let $Z$ be a nice interval (as usual, assumed to be sufficiently small),
having a central cascade $Z:=Z^{0}\supset Z^{1}\supset \dots\supset Z^{N+1}$ with $N\geq 1$ and let $r$ be the return time of $Z^{1}$ to $Z^{0}$. 
\begin{enumerate}
\item
If $|Z^{2}|\geq \delta |Z^{0}|$, then for any critical point $c$ of the mapping $R_{Z}|_{Z^{2}}$ we have
$$|F^{r}(c)-c|\geq\kappa |Z^{0}|\; and \; |DF^{r}(x)|\leq C\ for\ all\ x\in Z^{2}.$$
\item 
If $|Z^{1}|\geq \delta |Z^{0}|$ and we let $\tilde Z=(1+2\delta)Z$ and
$\tilde Z^1=\comp_{Z^1}(F^{-r}(\tilde{Z}))$ the following
holds. Suppose that $F^r$ extends to a mapping from $\tilde Z^1$ to
$\tilde Z,$ with the same set of critical points as $F^r|_{Z^1},$
that can be decomposed into a finite composition of maps with bounded distortion and polynomials. Then for any critical point $c$ of the mapping $R_{Z}|_{Z^{1}}$ we have
$$|F^{r}(c)-c|\geq\kappa |Z^{0}|,\; and \; |DF^{r}(x)|\leq C\ for\ all\ x\in Z^{1}.$$
\end{enumerate}
\end{lem}

Finally we will need the following Lemma of Yoccoz, 
which gives us precise control of the real geometry in saddle-node
cascades. For a proof, we refer the reader to \cite{dFdM 1}.
Suppose that $\Delta_0,\Delta_1,\Delta_2,\dots,\Delta_a$
is a sequence of adjacent intervals in $\mathbb R$.
If $g\colon\cup_{i=0}^{a-1}\Delta_i\rightarrow
\cup_{i=1}^{a}\Delta_i$
is a diffeomorphism with negative Schwarzian derivative such that
$g(\Delta_i)=\Delta_{i+1}$ for $i=0,\dots, a-1$, then we call $g$ an \emph{almost parabolic map
of length} $a$ with fundamental domains $\Delta_i$.

\begin{lem}[Yoccoz Lemma]\label{lem:Yoccoz lemma}
For any $\sigma>0,$ there exists $C> 1$ such that
the following holds.
Let $g\colon J\rightarrow I$ be an almost parabolic mapping of length $a$
with fundamental domains $\Delta_i,$ $i=0,\dots, a$.
If $|\Delta_0|>\sigma|I|$ and $|\Delta_{a-1}|>\sigma|I|$,
then 
$$\frac{1}{C}\frac{|I|}{\min\{i,a-i\}^2}\leq |\Delta_i|\leq C \frac{|I|}{\min\{i,a-i\}^2}.$$
\end{lem}

We will use the following simple lemma to prove that certain domains
are quasidisks.
\begin{lem}\label{lem:qdc}
For any $\alpha\geq 1,\delta>0,n\in\mathbb N$ and $c\in(0,1)$ there exists
$K\geq 1$ such that the following holds.
Suppose that $D\subset\mathbb C$ is a topological disk with marked points
$\{a_1,a_1,\dots,a_{n}\}\in\partial D.$ 
Assume that $D\setminus\{a_1,a_1,\dots,a_{n}\}$
is a union of $n$ piecewise smooth $\alpha$-quaisarcs
$\gamma_i,i=1,2,\dots,n$ with 
\begin{itemize}
\item $\diam(\gamma_i)>\delta\diam(D)$.
\item if $\gamma_i$ and $\gamma_{i'}, i\neq i'$ are not
adjacent in $\partial D$, then $\dist(\gamma_i,\gamma_{i'})\geq
\delta\diam(D)$,
\item If $\gamma_i,\gamma_{i'}$ are adjacent, then for any
$z_1,z_2\in\overline{\gamma_i\cup\gamma_{i'}}$ we have that the 
shortest $\gamma_{[z_1,z_2]}$ path in $\partial D$ connecting $z_1$
and $z_2$
satisfies $c\cdot\diam(\gamma_{[z_1,z_2]})\leq |z_1-z_2|$.
\end{itemize}
Then $D$ is a
$K$-quasidisk.
\end{lem}
\begin{pf}
Let us verify the Ahlfors-Beuring Criterion for any two points
$z_1,z_2\in\partial D$. First suppose that
$z_1,z_2\in\gamma_i.$
Then since $\gamma_i$ is an $\alpha$-quasiarc, we have that the
diameter of the path in $\gamma_i$ that connects $z_1$ and $z_2$
has diameter at most $C\cdot|z_1-z_2|,$ for some $C$ depending only
on $\alpha$. By assumption, the Ahlfors-Beurling Criterion holds
if $z_1$ and $z_2$ are in adjacent components of $\partial
D\setminus\{a_1,\dots,a_n\}$.
So, let us suppose that 
$z_1\in\gamma_i,$ $z_2\in\gamma_{i'},$ and $\gamma_i$ and
$\gamma_{i'}$ are not adjacent. Then we have that 
$\dist(z_1,z_2)>\delta\diam
D>n\delta\max_{i=1,\dots,n}\diam(\gamma_i)\geq n\delta\gamma_{[z_1,z_2]}.$
\end{pf}

\medskip

The following lemma
provides 
the first step in the inductive proof of 
Theorem~\ref{thm:central cascades}.
It is proved using 
the fact that
there are no central cascades lower in the partial ordering on central
cascades, and
the same argument as the one used to 
obtain the Real Spreading Principle.

\begin{lem}\label{lem:qc up to}
There exist $\delta_1>0$ and $K_1\geq 1$ so that
if $\hat{\mathcal{C}}^k_{c_0}(\II^0_{c_0}),$ 
$c_0\in\crit(F)$ and
$k\in\mathbb N\cup\{0\},$ is 
a puzzle piece
at the top of a maximal central cascade of length $N\geq N_0$ which is lowest in the
partial ordering on central cascades, then the following holds.
There exists a mapping $H\colon\mathbb C\rightarrow\mathbb C$ 
such that $H\colon\mathbb H^+\rightarrow\mathbb H^+$  is a 
$(K_1,\delta_1)$-qc\textbackslash bg mapping, 
and a neighbourhood $\bm{\mathcal{O}}$ of $\Crit(F)$
such that $\comp_{c_0} \bm{\mathcal{O}}=\bm{\mathcal{G}}_{c_0}(k),$
$H(\bm{\mathcal{O}})=\widetilde{\bm{\mathcal{O}}},$
$H(\Dom(\bm{\mathcal{O}})=
\Dom(\widetilde{\bm{\mathcal{O}}})$
and $H$ agrees with the boundary markings on $\bm{\mathcal{O}}$ and 
on each component of $\Dom(\bm{\mathcal{O}}),$ where
$\bm{\mathcal{G}}_{c_0}(k)$ is the puzzle piece associated to the long
central cascade by
Proposition~\ref{prop:good geometry for central cascades}.
\end{lem}
\begin{pf}
Since this central cascade is lowest in our ordering, 
we can argue as in the proof of the Real Spreading Principle to see
that there is a $(K_1,\delta_1)$-qc\textbackslash bg mapping of the plane as in the 
statement of the lemma. Specifically,
let $n$ be so that $\W_{c_0}^n$ is at the top of the central cascade,
and let $H_{\W^n}:\mathbb C\rightarrow\mathbb C$ be the mapping 
constructed in Steps 1--5 of Proposition~\ref{prop:realspreading}.
Let $\bm T^0=\W^n$, and let $t>0$ be minimal so that
$\bm T^t\subset\mathcal{L}_{c_0}^7(\bm T^0_{c_0})$. Note that $t$ is
at most $7\#\Crit(F).$
Now, starting with  $H_{\W^n},$ use the argument of Step 3 of
Proposition~\ref{prop:realspreading} to construct 
$\phi_t:\mathbb C\rightarrow\mathbb C,$ which conjugates the 
dynamics of $F$ and $\tilde F$ on  $\Dom(\bm T^t)$.
Starting with $\phi_t,$ use the argument from Step 6b of 
Proposition~\ref{prop:realspreading} to pullback $\phi_t,$
acting as a conjugacy between $F\colon\UU\to\VV$ and 
$\tilde F\colon\widetilde{\UU}\to\widetilde{\VV}$.
Let $r>0$ be so that $F^r(\bm{\mathcal{G}}_{c_0}(k))$ is a component
of $\VV$.
Since by Lemma~\ref{lem:transition times} there are at most 
six visits of $c_0$ to $\bm{\mathcal{G}}_{c_0}(k)$ along the 
orbit $\{f^i(c_0):0< i\leq r$, 
we have that after pulling $\phi_t$ back, we can assume that 
$\phi_t(\bm{\mathcal{G}}_{c_0}(k))=\widetilde{\bm{\mathcal{G}}}_{c_0}(k).$
Now we conclude the proof as in Proposition~\ref{prop:realspreading}.

\end{pf}

Assume that 
$$\Z^0\supset \Z^1\supset\dots\Z^{N+1}$$ is a maximal central cascade
and that
Theorem~\ref{thm:central cascades} 
holds for all maximal central cascades that are lower in the partial ordering
on central cascades than $\Z^0\supset \Z^1\supset\dots\Z^{N+1}$.
Then by the argument used to prove the Real Spreading Principle and
Lemma~\ref{lem:qc up to}, we have that there exists $K\geq 1$,
$\delta>0$ and a mapping 
$H\colon\Z^0\rightarrow\widetilde{\bm{Z}}^0$ with
$H(\bm Z^1)=\widetilde{\bm{Z}}^1$ which is
$(K,\delta)$-qc\textbackslash bg and,
and a $(K,\delta)$-qc\textbackslash bg in the upper half-plane. Let us
now carry out the inductive step of the proof of 
Theorem~\ref{thm:central cascades}.

\begin{figure}[htbp]
\begin{center}
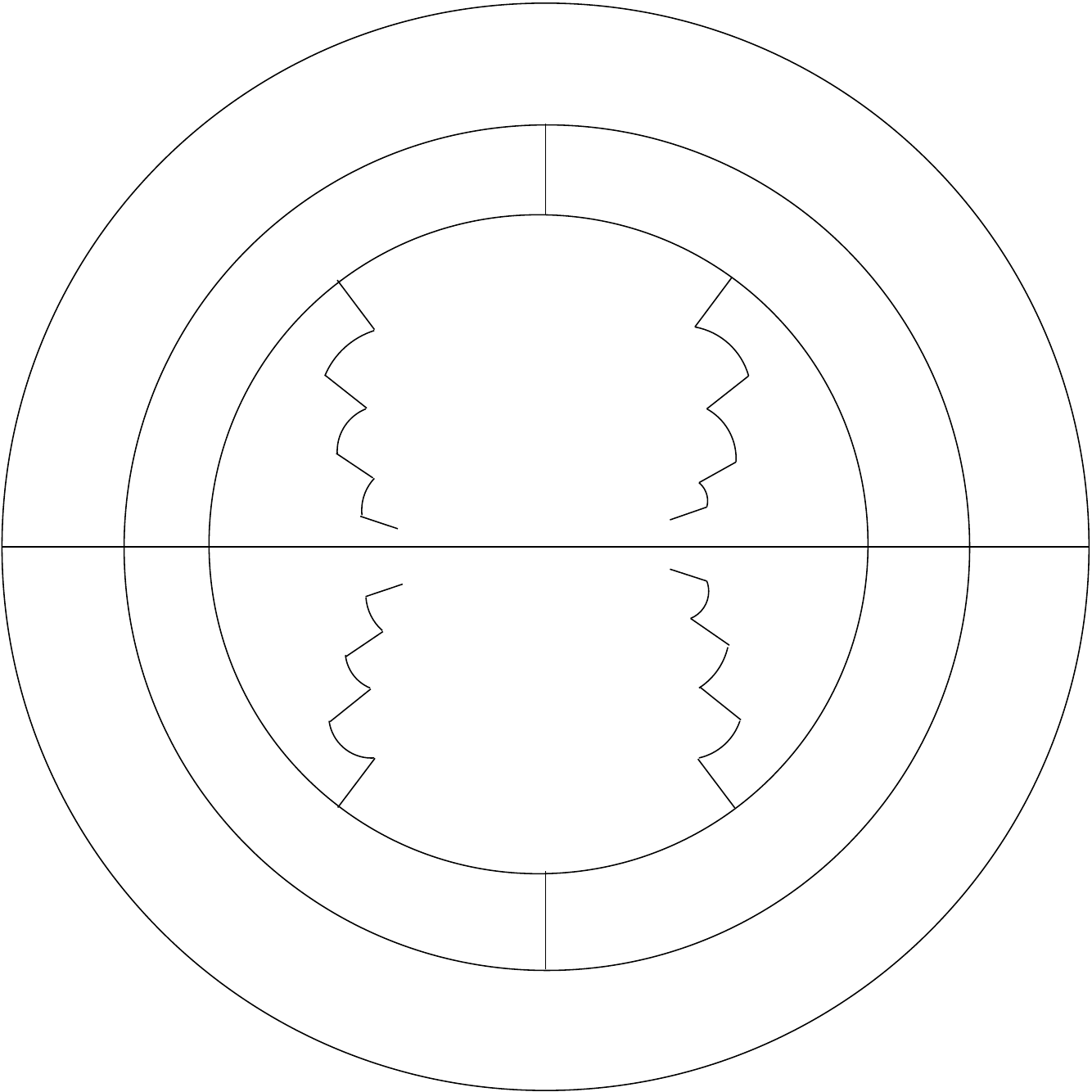
\end{center}
\caption{The first few segments of the paths, $\gamma_i$.}
\label{fig:staircase}
\end{figure}

\medskip
\noindent\textit{The proof of Theorem~\ref{thm:central cascades}}.
Recall that $R\colon\J^1\rightarrow\J^0$ is the first return mapping to
a $\delta$-excellent 
puzzle piece at the ``top'' of a long central cascade.
To begin, let us explain a construction that we will use in all
cases to subdivide $\J^0$ into the sets $X_0,X_1,X_2$ and $X_3$.
We are going to construct
four ``staircases'' which will consist of alternating
pre-images
of the real fundamental domains $J^0\setminus J^1$
under $R^i$ that lie in $\mathbb C\setminus\mathbb R$,
and piecewise smooth arcs in $\partial \J^{i+1},$
see Figure~\ref{fig:staircase}.
These staircases are piecewise smooth curves
and will be part of the set $X_0$.
They
will serve as the boundaries of
artificial, but still dynamically defined,
``puzzle pieces'' that we will use to form the sets
$X_1,$ $X_2$ and $X_3$.

\medskip

Suppose that the degree of $R$ is $d$.
Let $e_0$ be the union of the two fundamental intervals 
comprising
$J^0\setminus J^1$.
Then $R^{-1}(e_0)\subset \J^1\setminus \J^2$
consists of $2d$ components
contained in $\J^1\setminus\J^2$ that connect
$\partial \J^1$ with $\partial\J^2$. Of these components,
$2d-2\geq 2$ lie in $\mathbb C\setminus\mathbb R$. We let $e_1$
denote the union of components $e_1^j$ of $R^{-1}(e_0)$
with the property that there exists a path in
$\J^1\setminus(\overline\J{}^2\cup R^{-1}(e_0))$
that connects $e_1^j$ with the real fundamental domains,
$J^1\setminus J^2$.

Notice that if the degree of $R$ is two, then $e_1$ consists of
exactly two components, and if the degree of $R$ is greater than 
two, $e_1$ consists of four components.
We repeat this construction in each fundamental annulus:
for $i=1,\dots,N$,
we let $e_i$ the union of components of $R^{-i}(e_0)$
with the property that there exists a path in
$\J^i\setminus(\overline\J{}^{i+1}\cup R^{-i}(e_0))$
that connects $e_i^j$ with $J^i\setminus J^{i+1}$.
For all $i>1$, $e_i$ consists of four components.
It is useful to think of these as the components of the
$R^{-i}(J^0\setminus J^1)$ that lie closest to, but not on, the real line.
For $i=1,\dots N$, we let $\alpha_i\subset\partial \J^i$
be the union of the two shortest paths, $\alpha_i^1$,
$\alpha_i^2$,
each crossing the real line that connect
two points of $\partial e_i.$ For $i=1,\dots, N$,
we let $\gamma_i\subset \alpha_i$ be the union of shortest arcs
that connect $\partial e_i$ with $\partial e_{i+1}$.
Then for $i\geq 3$, $\gamma_i$ is a union of four disjoint connected
arcs, $\gamma_i^j,j\in\{1,2,3,4\}$ and
$\cup_{j=1}^4\gamma_{i}^j=\gamma_i.$
We have that
$$e_2\cup\bigcup_{i=3}^N\{\gamma_i,e_{i}\}$$
is a union of four disjoint paths that connect $\partial\J^2$ with
$\partial \J^{N+1}$. 
We will need to adjust this construction below to deal with 
specific cases.

To fix notation, for $i=2,\dots, N$, we let 
$e_i^1, e_i^2$ be the components of $e_i$
that are contained in the upper half-plane, and we let 
$e_i^3,e_i^4$ be the components that are contained in the
lower half-plane. We will also label the components so that 
for all $j\in\{1,2,3,4\}$ we have that 
$e_i^j$ is connected to $e_{i+1}^j$ by the component
$\gamma_{i+1}^j$ of $\gamma_{i+1}.$ 
For $i\geq 2$ and $j\in\{1,2,3,4\}$ we
let $v_{i}^j=e_{i}^j\cap\partial\J^i$ and 
$w_i^j=e_i^j\cap\partial \J^{i+1}.$
Moreover, we will label the components so that
$e_{i}^1$ is connected with $e_{i}^3$ by
$\alpha_i^1$ and $e_{i}^2$ is connected with $e_{i}^4$ by
$\alpha_i^2.$

Let us start by explaining the proof in two cases when
we have bounded geometry and when
there is just one critical point. 
We will then explain how to extend the proof to the
case the $R\colon J^1\rightarrow J^0$ is monotone,
and finally to the general multimodal cases.
We will define the constant separating the
big and bounded geometry cases, when we
treat the cases when $|J^0|$ is much bigger than $|J^1|$.

We divide the proof into the following types of returns.
Suppose that $c$ is a critical point of $F$, $I$ is a nice interval
and
$R\colon I^1\rightarrow I$ is the first return mapping.
We say that the return is \emph{high} if $c\in R(I^1)$
and otherwise we say that the return is \emph{low}.

\begin{figure}[htbp]
\begin{center}
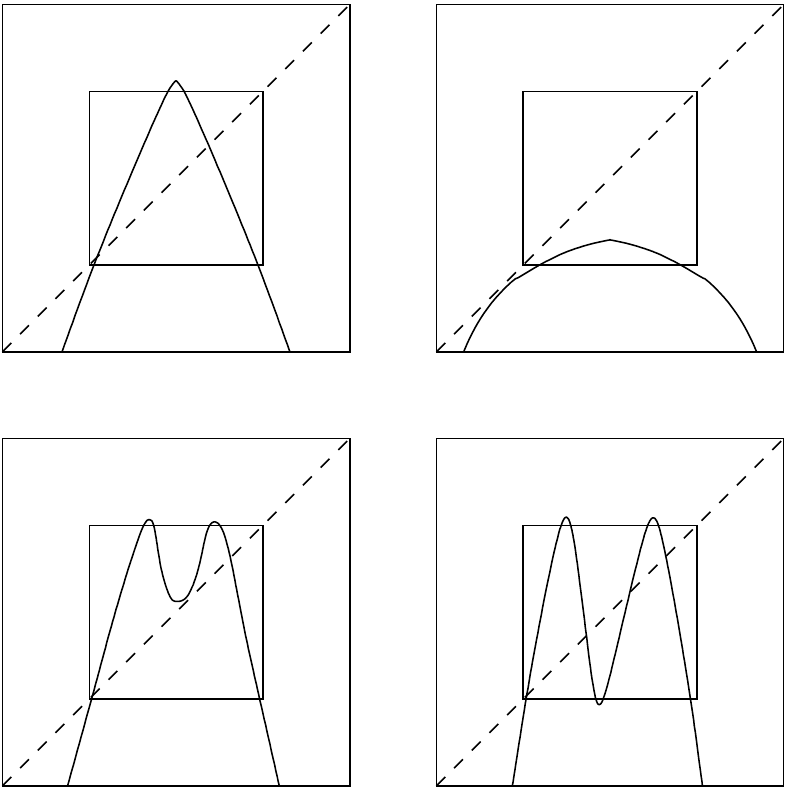
\end{center}
\caption{Types of returns in the setting of long cascades,
clockwise from the top left:
(1) a unimodal high return, (2) a unimodal low return, (3) a high return
followed by a high return, (4) a high return followed by a low return. 
}
\label{fig:returns}
\end{figure}

\medskip
\noindent\textbf{Case 1: Suppose that $R\colon J^1\rightarrow J^0$ is
  unimodal 
mapping with a low return and $|J^0|/|J^1|$ is bounded.}
In this case, $R\colon J^1\rightarrow J^0$ has no fixed points on the real
line,
and $R(J^1)\subset J^0\setminus J^{N+1}$.

\medskip

\noindent\textit {Step 1: Decomposing $\J^0$.}
See Figure~\ref{fig:qcbg} for a sketch of the 
decomposition.
First extend the staircase arcs that we have already constructed
by an extra step.
For $k=2,\dots,\mathrm{degree}(R)$,
let $\bm W_k$ index the components of $R^{-1}(\J^{N+1})$.
Let $e_{N+1}$ be union of the components of 
$R^{-1}(e_N)$ that can be connected to 
$\mathbb R\setminus\{0\}$ by a path in 
$\J^{N+1}\setminus (R^{-1}(e_N)\cup \overline{R^{-1}(\J^{N+1})})$.
Then $e_{N+1}$ consists of four disjoint components
$e_{N+1}^1\cup e_{N+1}^2\cup e_{N+1}^3\cup e_{N+1}^4$
each connecting $\partial\J^{N+1}$ with some $\bm W_k$.
We let $\alpha_{N+1}\subset\partial \J^{N+1}$
be the union of the two shortest paths, $\alpha_{N+1}^1$,
$\alpha_{N+1}^2$,
each crossing the real line that connect
two points of $\partial e_{N+1},$ and we
let $\gamma_{N+1}$ be the union of the four shortest disjoint paths
in $\alpha_{N+1}$ that do not intersect the real line and 
connect $\partial e_{N}$ with $\partial e_{N+1}$. We label 
the components of these sets so that for $j\in\{1,2,3,4\}$
$$\gamma_j=e_3^j\cup\bigcup_{i=4}^{N+1}\{\gamma_i^j\cup e_i^j\}$$
is a path connecting $\partial \J^3$ with $\cup_k \partial\bm W_k$.

\begin{figure}[htbp]
\begin{center}
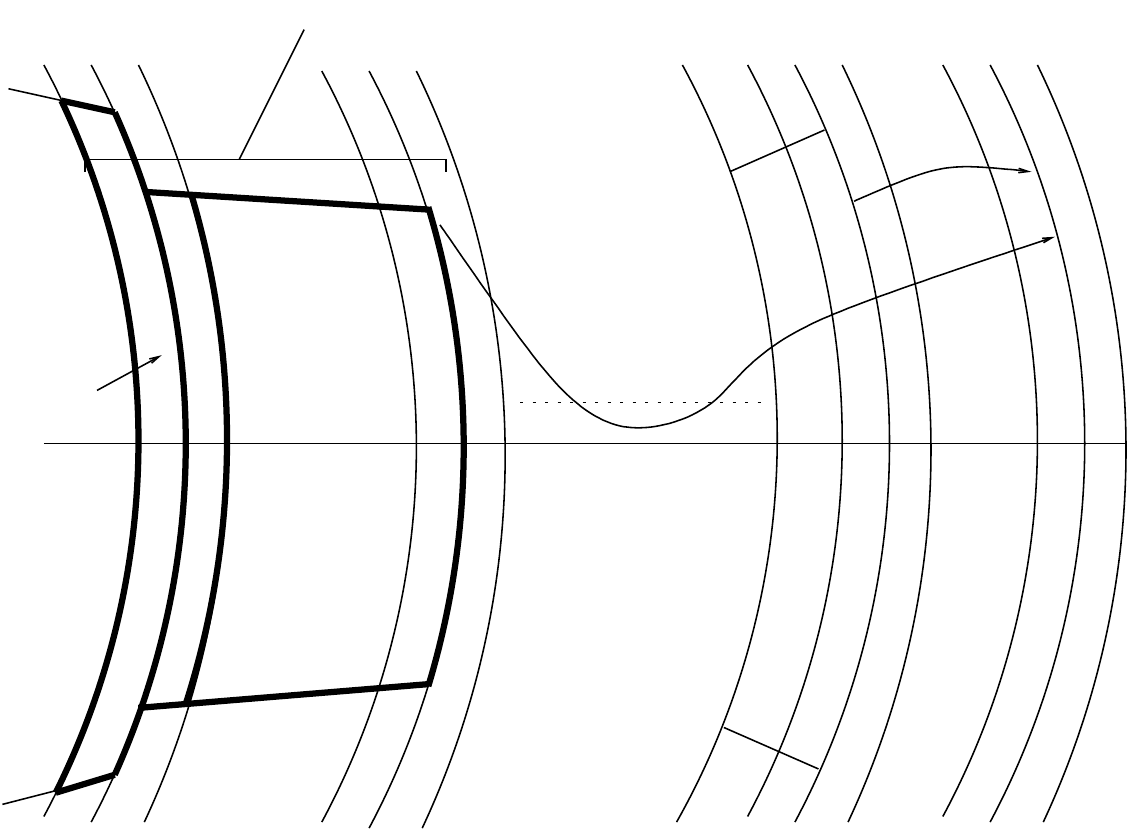
\end{center}
\caption{Choosing a definite neighbourhood of the critical value.
The figure shows one monotone branch of $R$. 
The critical value is contained in $\bm M_1'.$
Depending on its position, we choose the left hand boundary of 
$\bm M_1$, whose preimage under $R$ is off the real line.}
\label{fig:cvshield}
\end{figure}

\medskip
\noindent\textit{A definite neighbourhood of the critical value.}
\label{pageref:dncv}
We will now construct a definite neighbourhood of the critical value.
See Figure~\ref{fig:cvshield} and \ref{fig:saddle critical detail}. 
We need to construct 
a region with bounded geometry and moduli bounds in the upper half
plane that separates a $X_2\owns c$ from $X_3$, in the language of the
qc\textbackslash bg partition. To do this, we
construct a definite neighbourhood of the critical value, where we 
choose the neighbourhood depending on its position.
We will use this idea in other cases as well.

Recall that $\J^0=\II_n^{\Z^0}$ for some $n$, and
let $\bm B_0'= \II_n^{\Z^0}\setminus \II_n^-$,  see Lemma~\ref{lem:nu}.
There are two components of
 $R^{-(N+1)}(\bm B_0')\setminus e_{N+1}$ that intersect the real line.
Let $\bm B_1'$ be the component of $R^{-(N+1)}(\bm B_0')\setminus e_{N+1}$
that intersects the real line
with the property that the region
$\bm M_1'$
bounded by
$\partial\bm B_1'\cup e_{N+1}\cup \gamma_{N+1}\cup e_{N}\cup \alpha_N$
contains the critical value.
Notice that $\bm B_1'$ is mapped onto $\bm B_0'$ 
by a mapping with bounded degree and small dilatation,
so $\bm B_1'$ is a rectangle with bounded modulus,
$m\geq 1$, and there exists $K'\geq 1$ such that
each arc of the boundary of $\bm B_1'$ is 
is a $K'$-quasiarc.
 If $\bm M_1'$ does not contain the critical point,
set $\bm B_1=\bm B_1'$ and $\bm M_1=\bm M_1'.$
Otherwise, if $\bm M_1'$ contains the critical point,
then the critical value is not close to $\partial\bm J^{n+1}$,
and we let  $\bm B_1$ be the rectangle bounded by
$\partial\J^{N+1}\cup e_N\cup R^{-(N+1)}(\partial\II_n^+)$ that is 
closest to the critical value, and we let  
$\bm M_1$ be the region bounded by 
$\partial\bm B_1\cup e_{N}\cup \gamma_{N}\cup e_{N-1}\cup \alpha_{N-1}.$
Let $b_1=\partial \bm B_1\cap\partial\J^{N+1},$
let $b_2$ and $b_4$ be the components of $\partial \bm B_1\cap
\overline e_{N+1},$ and let $b_3=\partial \bm B_1\setminus(b_1\cup
b_2\cup b_4)$.
Let $\phi_{\bm B_1}\colon \bm B_1\rightarrow [0,1]\times[0,1]$ be a $K''$-qc
uniformization of $\bm B_1$
that maps $\bm B_1$ onto  $[0,1]\times[0,1],$
where $K''$-depends only on $m$,
with the property that it maps
$\overline b_1$ to $\{1\}\times[0,1]$ and $\overline b_3$ to
$\{0\}\times[0,1].$
Let us subdivide $\bm B_1$ into $\bm B=\phi_{\bm
  B_1}^{-1}((1/2,1]\times[0,1])$ and  $\bm B'=\phi_{\bm
  B_1}^{-1}[0,1/2]\times[0,1].$
Let $\bm M$ be the region bounded by 
bounded by
$\partial\bm B\cup e_{N+1}\cup \gamma_{N+1}\cup e_{N}\cup \alpha_N.$
Observe that we can carry out this construction for $R$ and $\tilde R$
simultaneously to produce corresponding neighbourhoods of the critical
point for each map.

Let $\Gamma_1$ be the region that intersects the real line and
is disjoint from the critical point that is bounded by
$\partial\bm M\cup\gamma\cup\partial\J^3$,
and let $\Gamma_2=\tau(\Gamma_1)$ where 
$\tau$ is the symmetry about the even critical point.
Let $\Gamma=\Gamma_1\cup\Gamma_2$.
We take $\bm N=\Gamma\cup R^{-1}(\bm M)$.

\medskip

\begin{figure}[htbp]
\begin{center}
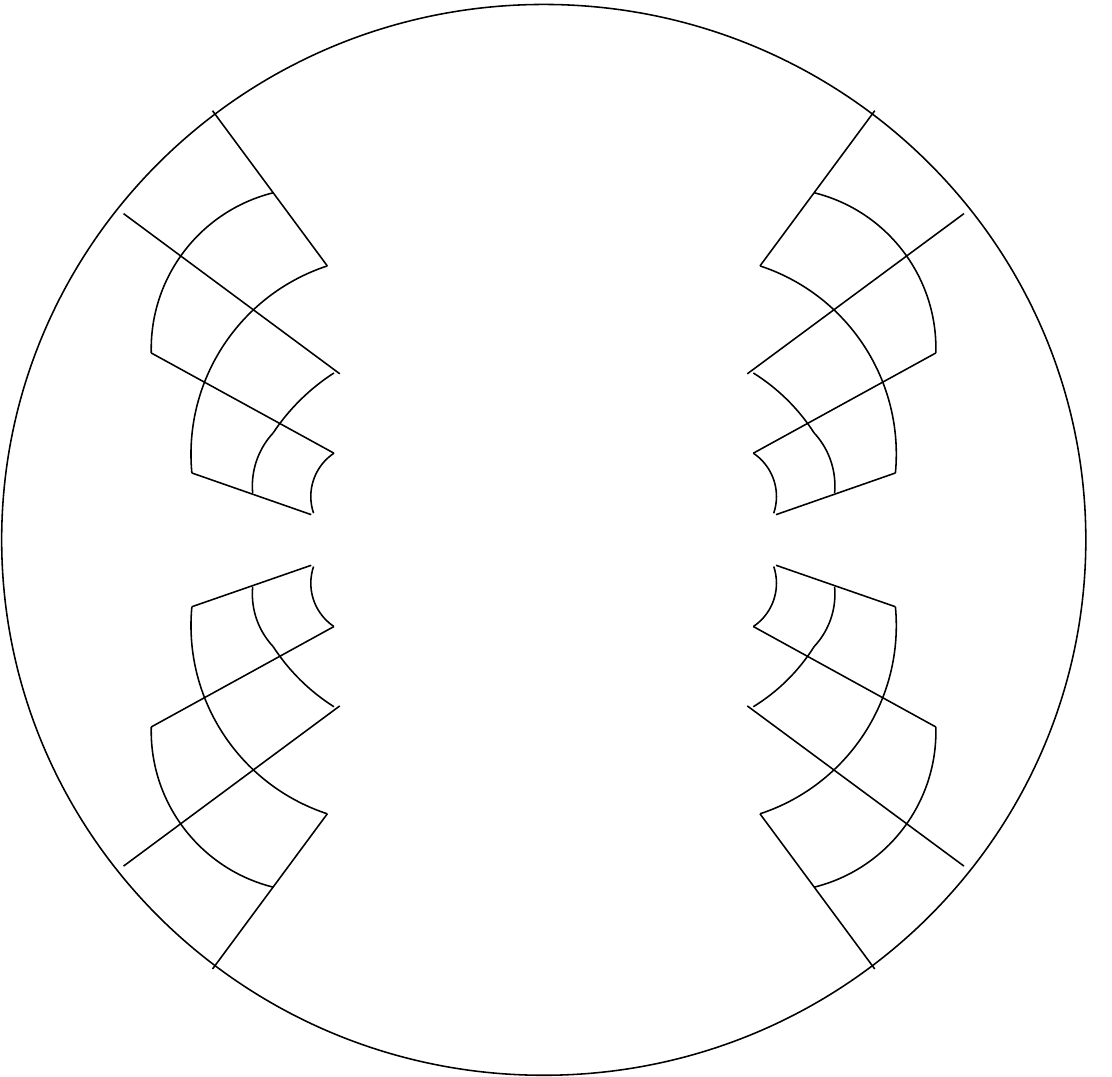
\end{center}
\caption{The first few components of $\bm N'\setminus\bm N\subset
  X_1$ refer to Figure~\ref{fig:qcbg}.}
\label{fig:shield}
\end{figure}

\begin{figure}[htbp]
\begin{center}
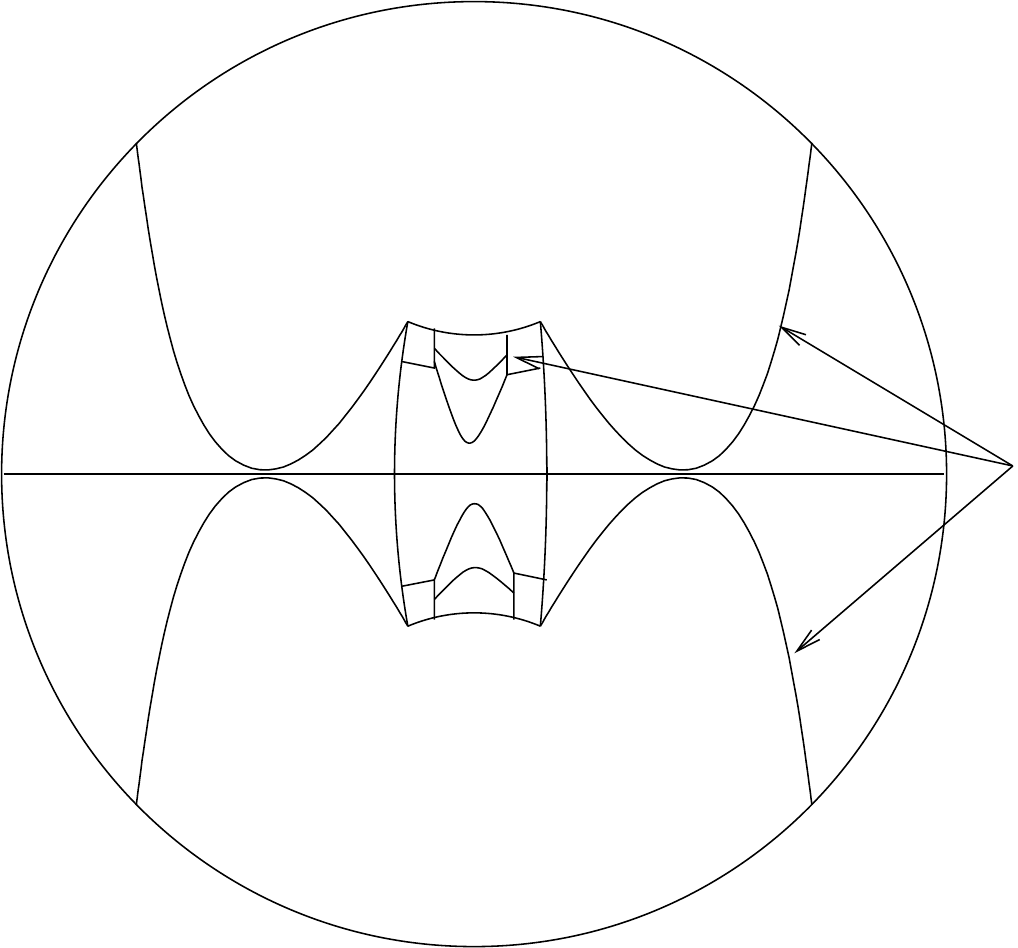
\end{center}
\caption{The components $\bm W_1$ and $\bm W_2$ of $X_3$ 
and detail of the construction around
 the critical point. This picture does not show the detail of $X_1$.
Refer to Figures \ref{fig:qcbg} and \ref{fig:shield} for the detail of $X_1$.
The complex fixed points of $R$ are contained in
$\bm W_1\cup \bm W_2.$
}
\label{fig:saddle critical detail}
\end{figure}

Let us now construct a neighbourhood 
$\bm N'\supset\bm N$, see Figure~\ref{fig:shield}.
As before, we will construct 
staircases which will form 
part of the boundary of $\bm N'$. 
They are constructed in the same way as at the start of the proof, but
the construction is shifted by one level.
For $i=2,\dots,N$, 
let $e_{i}'$ be the components of $R^{-(i-1)}(J^1\setminus J^3)$
that are contained in $\J^i\setminus(\J^{i+1}\cup \mathbb R$),
and that can be connected to $J^i\setminus J^{i+1}$ by a path in 
$\J^i\setminus(\J^{i+1}\cup \overline{R^{-(i-1)}(J^1\setminus J^2)}).$
For $i=N+1$, 
let $e_{N+1}'$ be union of the components of 
$R^{-1}(e_N')$ that can be connected to 
$\mathbb R\setminus\{0\}$ by a path in 
$\J^{N+1}\setminus (R^{-1}(e_N')\cup \overline{R^{-1}(\J^{N+1})})$.
For $i=2,\dots, N$,
$e_i'$ connects $\partial \J^i$ with $\partial \J^{i+1}$,
and $e_{N+1}'$
connects $\partial\J^{N+1}$ with some $\bm W_k$.
Let $v_i'=\partial e_i'\cap\partial \J^i$ and $w_i'=\partial e_i'\cap \partial\J^{i+1}$.
For $i=1,2,\dots, N+1$ let $\gamma_i'\in\partial \J^i$ be the union
of shortest arcs in $\partial\J^i$ connecting $w_{i-1}'$ with $v_i'$.
Let 
$$\gamma'=e_3'\cup\bigcup_{i=3}^{N}\{\gamma_i'\cup  e_{i+1}'\}.$$
$\gamma'$ consists of four arcs connecting $\partial \J^3$
to $R^{-1}(\J^{N+1})$. 
Let us label the four components of 
$\gamma_1', \gamma_2', \gamma_3', \gamma_4'$
of $\gamma'$ so that for $j\in\{1,2,3,4\}$ we have that
$\gamma_j\cap\gamma_j'\neq\emptyset$.

Now we we will close off the region
around the critical value,
see Figure~\ref{fig:saddle critical detail}.
We will give the construction in the case when
the critical value is not contained in
$R^{-N}(\II_n^-)$ - this corresponds to the 
case when it may be close to $\partial\J^{N+1}$.
The modifications to the construction
 when the critical value is 
in $R^{-N}(\II_n^-)$ are the same as in the construction of
$\bm{N}.$
First, we extend the domain 
$\bm B_1$ constructed above:
As before, there are two components of
 $R^{-(N+1)}(\bm B_0)\setminus e_{N+1}'$ that intersect the real line.
Let $\bm B_1'$ be the component of $R^{-(N+1)}(\bm B_0)\setminus e_{N+1}$
that intersects the real line
with the property that the region
$\bm M'$
bounded by
$\partial\bm B_1'\cup e_{N+1}'\cup \gamma_{N+1}'\cup e_{N}'\cup \alpha_N'$
contains the critical value and not the critical point.
Observe that $\bm B_1'\setminus \bm B_1$ consists of three
topological squares with bounded geometry
whose boundaries consist of
arcs in $\alpha_{N+1}',$ the components
of $\partial\J^{N+1}\setminus \overline{e'}_{N+1}$ that intersect the real line.
$\partial \bm B_1'\setminus \partial \bm B_1,$ $e_{N+1}$ and
$e_{N+1}'$, and 
$\phi_{B_1}^{-1}(\{1/2\}\times[0,1])$, and
that 
$\bm M'\setminus\bm M=(\bm B_1'\setminus\bm B_1)\cup\bm B'$.
Let $\Gamma_1'$ be the region that intersects the real line and
disjoint from the critical point that is bounded by
$\partial\bm B_1'\cup\gamma'\cup\partial\J^3$,
and let $\Gamma_2'=\tau(\Gamma_1')$ where 
$\tau$ is the symmetry about the even critical point.
Let $\Gamma'=\Gamma_1'\cup\Gamma_2'$.
Let $\bm N'=\Gamma'\cup R^{-1}(\bm M')$.

The shape of $\Gamma'$ is dictated by the Yoccoz Lemma:
Let $\bm Z_2$ be a topological rectangle bounded by 
$e_2\cup\partial\J^2\cup \partial\J^3\cup(J^2\setminus J^3)$ whose 
intersection with the real line is 
a component of $J^2\setminus J^3$-
a fundamental domain for $R$.
For $i\leq N-2,$ let $\bm Z_{i+2}$ be a component
of $R^{-i}(\bm Z_2)$ that intersects the real fundamental domain
$J^{i+2}\setminus J^{i+3}.$ 
Then $\bm Z_{i+2}$ is mapped with bounded dilatation onto 
$\bm Z_{2}$ by $R^i$. So since the height and width of $\bm Z_2$ are
both comparable to $J^2\setminus J^3$, the 
height and width of $\bm Z_i$ are
both comparable to the length of the real trace of 
$\bm Z_i$, a component of 
$J^i\setminus J^{i+1}.$ By the Yoccoz Lemma
$$|\bm Z_i|\asymp \frac{|J^0|}{\min\{i^2,(N-i)^2\}}.$$
So we see that near the middle of a long saddle-node cascade,
the region $\Gamma'$ is pinched close to the real line.

\medskip
Now we are ready to describe the sets $X_0$, $X_1$, 
$X_2$ and $X_3$.
Set 
$$X_1=(\bm N'\setminus \bm N)\cup (\bm J^3\setminus (\bm J^4\cup \bm
N')).$$ From the construction above, 
we have that $X_1$ is a union of combinatorially defined topological squares
with bounded geometry, 
and moduli bounds as in conclusion (3) from
Theorem~\ref{thm:central cascades}.
To be definite, let us suppose that $\bm S$ is contained in the
upper half-plane.
Then $\bm S$ is a pullback of a component of $X_1$
that is either bounded by $e_i\cup e_i'$ and the shortest arcs in 
$\partial \J^i$ and $\partial \J^{i+1}$ connecting the end points of
$e_i$ and $e_i',$ where $i=3$ or $4$.
Let $\bm S'$ denote this component.
It follows from the complex bounds, Proposition~\ref{prop:good
  geometry for central cascades},
that $\bm S'$ has bounded geometry. We also have that it is contained in 
the set $\hat{\bm{S}}'$ which is contained in the upper
half-plane and is bounded by a component of
$R^{-1}(J^1\setminus J^4)$ lying in $\mathbb C\setminus\mathbb R$,
a component of $J^2\setminus J^5\subset\mathbb R,$ and two arcs
contained in
$\partial\J^2,$ $\partial\J^5,$ respectively. 
We have from the complex bounds 
that $\mod(\hat{\bm{S}}'\setminus\bm S')$ is bounded away from
zero. Moreover, the mapping from $\bm S$ onto $\bm S'$ extends to a 
mapping with small dilatation onto $\hat{\bm{S}}'$. Hence 
there exists $\delta'>0$ depending only on the complex bounds 
such that
$\bm S$ has
$\delta'$-bounded geometry and
$\mod(\J^0\cap\mathbb H^+\setminus \bm S)>\delta.$ 
It follows from this and the construction near the critical value that
$R^{-1}(\bm M'\setminus\bm M)$ has universally
bounded geometry and moduli bounds as well.
Finally, we set 
\begin{itemize}
\item $X_2=\bm N\cup (\J^0\setminus \J^3)$, and 
\item $X_3$ to be the components of 
$\J^0\setminus\overline{X_1\cup X_2}.$
\end{itemize}

There exists a constant
$K_1\geq 1$ such that each component of 
$X_3$ is a $K_1$-quasidisk.
To show this
we use Lemma~\ref{lem:qdc}.
Let $\bm W$ be a component of $X_3$.
The boundary of $\bm W$ consists of four arcs.
$v_1\subset \partial \J^4$, $v_2\subset\Gamma'_1$
$v_3\subset\partial R^{-1}(\bm B_1')$ and
$v_4\subset \Gamma_2'$. 
Since $\bm J^4$ is a $K_0'$-quasidisk,
we have that $v_1$ is a 
$K_1''$-quasiarc.
Similarly, $v_3$
 is a $K_1''$-quasiarc,
since $\partial R^{-1}(\bm B_1')$ is a $K_0'$-quasiarc.
We also have that 
$v_2$ and $v_4$ are
$K_0''$-quasiarcs:
if $\zeta\subset \gamma'_i$ 
is a subarc of $\gamma'_i,$
then $\zeta$ is a union of preimages of
real fundamental domains and the boundaries of puzzle pieces that are
$K_0'$-quasiarcs, the intersections of these paths do not contain any
cusps, so the only interesting case is when $|x-y|$ is small, and
$\zeta$ contains many steps in the staircase.
But now the statement follows since
$\Gamma'\setminus\Gamma$ consists of topological squares
with bounded geometry.
Since $\J^4$ is $\delta'$-nice for some $\delta'>0,$
universal, 
we have that 
$\dist(v_1,v_3),$ 
$\diam(v_2), \diam(v_4)$ are
comparable
to $\diam(\J^0)$, and
Lemma~\ref{lem:CvST 3.10}
and the same analysis as above imply that
$\diam(v_1),\diam(v_3)$,
$\dist(v_2,v_4)$ are comparable to $\diam(\J^0).$
It remains to verify the Ahlfors-Beuring Criterion on 
adjacent components of $\partial \bm W$.
Consider $v_1\cup v_2$. 
Suppose that $z_1,z_2\in v_1\cup v_2.$
Let $\gamma_{[z_1,z_2]}$ be the path connecting $z_1$ and $z_2$
in $v_1\cup v_2$.
We can assume that $z_1\in v_1$ and $z_2\in v_2.$
If $v_2\in\J^5$, then 
$\diam(\gamma_{[z_1,z_2]})$ and $|z_1-z_2|$ are comparable to
$\diam(\J^0),$
so we can assume that $z_2\in e_4,$ but now $e_4$ is a preimage of
of an interval in $J^0\setminus J^1$ that meets $\partial\J^0$ at a
definite angle. So $\diam(\gamma_{[z_1,z_2]})$ and $|z_1-z_2|$ are comparable.
The similar arguments show that the Ahlfors-Beurling Criterion holds
on the
arcs $v_1\cup v_4, v_2\cup v_3$ and $v_3\cup v_4.$

It is important to observe that the boundary of $X_3$ 
consists of pullbacks of the real line and boundaries of puzzle
pieces, so that for each component of $X_3$ there is a corresponding
component of $\tilde{X}_3.$ The set
$X_0$ is the union of the piecewise smooth boundaries of these regions.

\medskip

\noindent\textit{Step 2: Constructing the quasiconformal mapping $H$.}
We define $H$ on $X_2$ dynamically.
For $z\in \Gamma\cup \bm B$, let $0<k_z\leq N+1$ be so that
$R^{k_z}(z)\in \overline \J{}^0\setminus \J^1$,
and define $H(z)$ by
$H_0\circ R^{k_z}(z)=\tilde R^{k_z}\circ H(z)$.
Let $\bm P$ be the component of $\Dom^*(\J^0)$ that 
contains $R^{k_c}(c)$ (the image of the critical point when it escapes
through the cascade). 
For $z\in\bm R^{-1}(\bm M)\setminus R^{-k_c}(\bm P),$
We define $H$ in $R^{-1}(\bm M)$ by
$\tilde R\circ H(z)=H\circ R(z)$
By applying the inductive step,
by the conclusion of Theorem~\ref{thm:central cascades} for 
central cascades that are lower in the partial ordering,
together with the argument used to prove the
Real Spreading Principle, we have that there
exists $K'\geq 1$ and $\delta'>0$ such that 
$H|_{\bm P}$ is a $(K',\delta')$-qc\textbackslash bg
mapping. Thus by the QC Criterion, we have that 
there is a $K''$-qc mapping that
agrees with the
boundary marking on $R^{-k_c}(\bm P).$ 
Observe that this procedure defines $H$ on 
$\partial X_1\cap\partial X_2=\cup_{i=2}^{N+1}\{v_i'\}$.

Let us now define $H$ on $X_3$. 
Let $\bm W$ be a component of $X_3.$
Since $\W$ and $\widetilde{\W}$ are $K_1$-quasidisks,
there exist $K_1$-qc mappings
$\phi_{\bm W }, \phi_{\widetilde{\bm{W}}}\colon\mathbb C\rightarrow\mathbb C$
such that
$\phi_{\bm W} (\mathbb D)=\W$,
and $\phi_{\widetilde{\bm{W}}}(\mathbb D)=\widetilde{\W}$.
We start by partially defining a mapping 
$h_{\bm W}$ from $\partial\bm W$ to $\partial\widetilde{\bm{W}}$.
Observe that $w_i=v_{i+1}'$ and that 
$\gamma\cap\gamma'=\cup_{i=2}^{N+1} \{v_i'\}.$
Observe that we can define $h_{\bm W}$ arbitrarily on
$\partial \J^4\cap \bm W$,
so that it maps $\partial \J^4\cap \bm W$ to 
$\partial \widetilde{\J}^4\cap \widetilde{\bm{W}}$.
We set $h_{\bm W}(v_i')=H(v_i')=H(w_{i-1})$ 
where $H$ is the mapping on $\overline{X}_2$
already defined, and we define 
$h_{\bm W}(w_i')=\tilde w_i'$.

By Theorem A \cite{vSV},
we have that there exists $\delta_0>0,$ universal,
so that $(1+\delta_0)J^1\subset J^0$.
Moreover, since $|J^1|\asymp|J^0|$, we have that
each of the components of $J^0\setminus J^1$ have length
comparable to $|J^0|,$ and
it follows from Lemma~\ref{lem:CvST 3.10},
there exists $\delta'>0$ such
that $|R(c)-c|>\delta'|J^0|.$ 
So we may apply the Yoccoz Lemma. 

Suppose that $x,y,z$ are three points lying in
$\cup_{i=2}^{N+1}\{v_i',w_i'\}.$
Observe that each $e_i^{1'}$
bounds a rectangle $Z_i^{1'}$ in the upper
half-plane bounded by $J^i\setminus J^{i+1}$
and two arcs in $\alpha_i$ and $\alpha_{i+1}$,
respectively.
Moreover, this rectangle is mapped with bounded
dilatation by $R^{i-3}$ onto a rectangle $Z_4^{j'}$ in the upper
half-plane
contained in $J^4\setminus J^{5}$ and
bounded by
$e_4^{j'}$
and two arcs in $\alpha_4$ and $\alpha_{5}$.
Since we have that the distance between the corners of 
$Z_4^{j'}$ are all comparable, we have the same for the corners of
$Z_i^{1'}.$ Hence have that $|w_{i}'-v_i'|$ and
$|v_i'-w_{i-1}'|$ are both comparable to 
$\frac{|J^0|}{\min\{i^2,(N-i)^2\}}.$

For $i\leq N,$
let $\Delta_i$ be the component of $J^i\setminus J^{i+1}$
containing $R^{N-i+1}(c).$
Let $B_1$ be the component of $J^0\setminus J^{N+1}$ that
contains $R(c)$.
By the Yoccoz Lemma, we have that
there exists a constant $C$  for $1\leq i\leq N$,
$$\frac{1}{C}\frac{|B_1|}{\min\{i,N-i\}^2}\leq |\Delta_i|\leq C \frac{|B_1|}{\min\{i,N-i\}^2}.$$
Moreover, we have the same estimates for the objects marked with
tilde, with the same $C$ (since $|J^1|/|J^0|$ is comparable 
to $|\tilde{J}^1|/|\tilde{J}^0|$).
Suppose that $w\in\Delta_{j_w}$, $w\in\{x,y,z\}$
Without loss of generality, we can assume that
$\Delta_{j_x}$ lies to the left of $\Delta_{j_y},$
which lies to the left of $\Delta_{j_z}$, and that
$j_z\leq j_y\leq j_x$.
Then we have that 
$$\frac{1}{|B_1|}\sum_{k=j_y+1}^{j_x-1}|\Delta_{k}|\leq C
\sum_{k=j_y+1}^{j_x-1}\frac{1}{\min\{i,N-i\}^2}
\leq C^2 \frac{1}{|\tilde B_1|}\sum_{k=j_y+1}^{j_x-1}
|\tilde{\Delta}_{k}|.$$
Similarly we have that 
$$\frac{1}{C^2}\frac{1}{|\tilde B_1|}\sum_{k=j_y+1}^{j_x-1}
|\tilde{\Delta}_{k}|\leq\frac{1}{|B_1|}\sum_{k=j_y+1}^{j_x-1}|\Delta_{k}|,$$
and
$$\frac{1}{C^2}\frac{1}{|\tilde B_1|}\sum_{k=j_z+1}^{j_y-1}
|\tilde{\Delta}_{k}|\leq\frac{1}{|B_1|}\sum_{k=j_z+1}^{j_y-1}|\Delta_{k}|\leq 
 C^2 \frac{1}{|\tilde B_1|}\sum_{k=j_z+1}^{j_y-1}
|\tilde{\Delta}_{k}|.$$
Thus we have that 
$$\frac{|\tilde x-\tilde y|}{|\tilde y-\tilde z|}
=\frac{\sum_{k=j_2+1}^{j_1-1}
|\tilde{\Delta}_{k}|}{\sum_{k=j_3+1}^{j_2-1}
|\tilde{\Delta}_{k}|}
\asymp
\frac{
\frac{|\tilde B_1|}{|B_1|}\sum_{k=j_2+1}^{j_1-1}|\Delta_k|}
{\frac{|\tilde B_1|}{|B_1|}\sum_{k=j_3+1}^{j_2-1}|\Delta_k|}=\frac{|x-y|}{|y-z|}.$$
Thus $h_{\bm W}$ is $(C,p)$-qs on $\cup_{i=2}^{N+1}\{v_i',w_i'\},$, 
where $C>0,p\geq 1$ 
depend only on $|J^0|/|J^1|,$ by way of 
Lemma~\ref{lem:CvST 3.10} and
the Yoccoz Lemma.

\medskip

Let us define a mapping $h'\colon \partial\mathbb D\rightarrow\partial\mathbb D$
as follows. For $z\in\cup_{N=2}^{N+1}\phi^{-1}_{\bm W}\{v_i',w_i'\}$
we define $h'(z)=\phi^{-1}_{\widetilde{\bm{W}}}\circ H\circ\phi_{\bm
W}(z)$.
Thus we have that $h'$ is $\kappa'$-quasisymmetric on 
$\cup_{i=2}^{N+1}\phi_{\bm W}^{-1}\{v_i',w_i'\}$,
where $\kappa'\geq 1$ depends only on $K_1$ and $\kappa$.
Let us show that 
$\cup_{i=2}^{N+1}\phi_{\bm
  W}^{-1}\{v_i',w_i'\}$
is $M$-relatively connected.
First, choose $x\in\cup_{i=2}^{N+1}\phi_{\bm
  W}^{-1}\{v_i',w_i'\},$ and $r>0$ so that
$\overline{B(x,r)}\neq S$.
Since the gaps between adjacent points in  
$\cup_{i=2}^{N+1}\phi_{\bm
  W}^{-1}\{v_i',w_i'\}$ are governed by the Yoccoz Lemma,
we have $S$ is $M$-relatively connected
at $x$ with $M\asymp 2$.
Since
$\phi^{-1}(\partial \J^4)\cup\cup_{i=2}^{N+1}\phi_{\bm W}^{-1}\{v_i',w_i'\}$
is $M$-relatively connected,
by Proposition~\ref{prop:Vellis}
we can extend $h'$ to a $\kappa''$-qs mapping
$h'\colon\partial \mathbb D\rightarrow\partial\mathbb D$.
Thus by Ahlfors-Beurling it extends to a $\hat K_2$-qc mapping 
$H'\colon\mathbb D\rightarrow\mathbb D$
and we thus we have that
$h_W\colon W\rightarrow \tilde W$ 
extends to a $\hat K$-qc mapping,
$$H_{\bm W}=\phi_{\widetilde{\bm{W}}}\circ H'\circ \phi_{\bm
  W}^{-1}\colon\bm W\rightarrow\widetilde{\bm{W}}$$
that agrees with $H$ where ever they are both defined.

Finally, we define the mapping on $X_1$. Let $\bm S$ be a component 
of $X_1$, and let $\widetilde{\bm{S}} $ be the corresponding
component of $\tilde X_1$.
Then $\bm S$ and $\widetilde{\bm{S}}$ 
are bounded by piecewise smooth quasiarcs.
The mapping $H$ is already defined on some components of
$\partial \bm S$. We leave $H$ unchanged where it is already
defined, and
extend it to a homeomorphism from
$\partial \bm S$ to
$\partial\widetilde{\bm{S}},$ and
then to a homeomorphism
from the interior of $\bm S$ to 
the interior of $\widetilde{\bm{S}},$
so that it satisfies (\ref{eqn:qc bound 0}) on
$\partial\bm S.$

\medskip

\noindent\textbf{Case 2: Suppose that $R\colon J^1\rightarrow J^0$ is a
  unimodal mapping with a high return and $|J^0|/|J^1|$ is bounded.}
In this case, $R$ has an orientation preserving fixed point
$\beta\in J^{N+1}.$ As usual $\tau(\beta)$ denotes the symmetric point
under the symmetry about the even critical point.
We have that $R(c)\in J^{N}\setminus[\beta,\tau(\beta)]$
and $c\in R((\beta,\tau(\beta)))$.

\begin{figure}[htbp]
\begin{center}
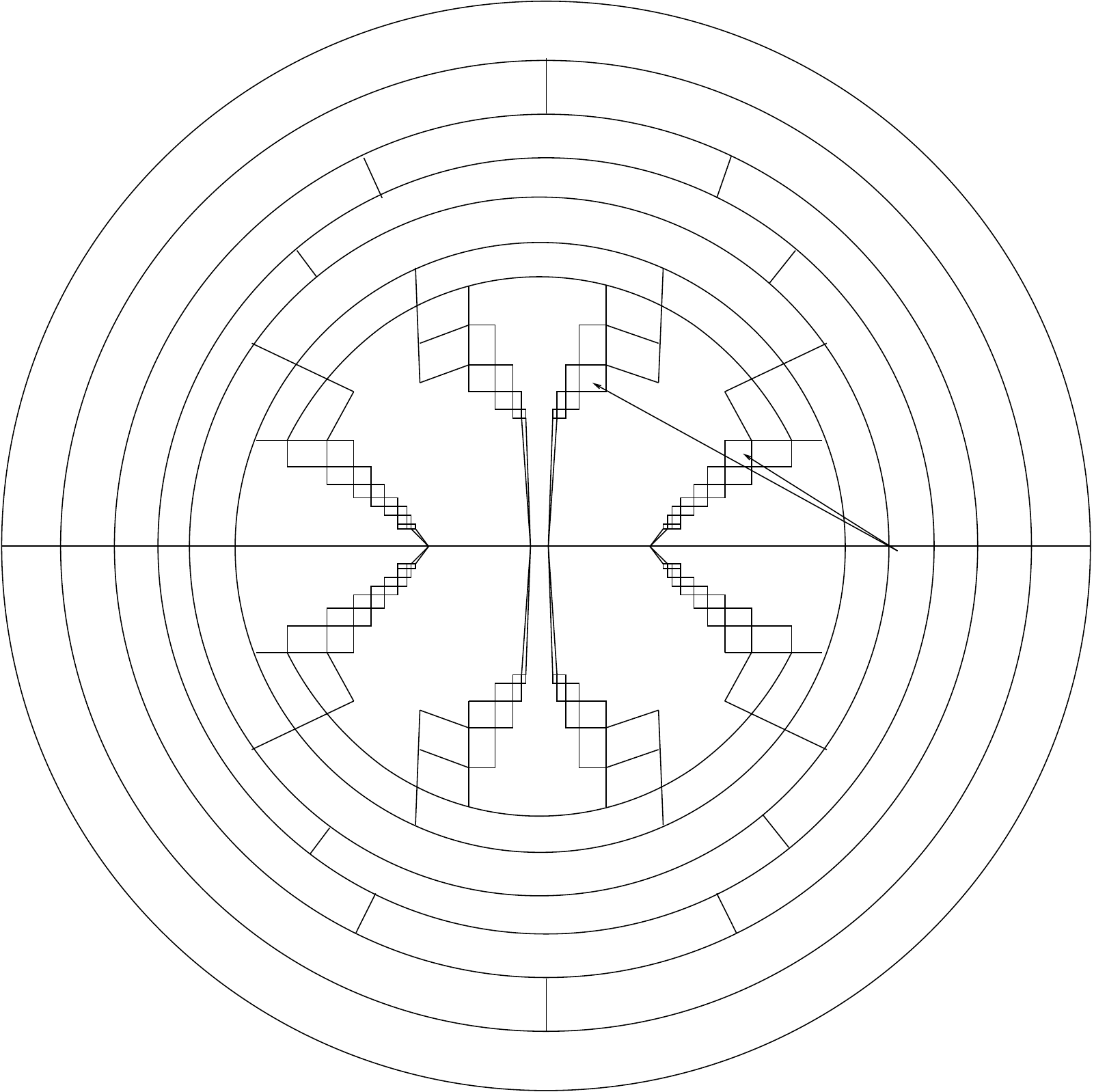
\end{center}
\caption{Domains around the critical point of a high return.}
\label{fig:unimodalhigh}
\end{figure}

\medskip
\noindent\textit{Step 1: Decomposing $\J^0$.}
Let us start by decomposing $\J^0$ into the sets 
$X_0,X_1,X_2,X_3$. See Figure~\ref{fig:unimodalhigh}
for a rough picture of the decomposition.
In this case, the staircases consist of infinitely many 
$e_{i}$, $i=3,4,\dots,$ pullbacks of the real fundamental domains 
$J^0\setminus J^1,$
or $J^1\setminus J^2,$
and
subarcs $\gamma_i$ of the
boundaries of certain puzzle pieces that intersect the real line.

Let us choose the puzzle pieces whose boundaries will contain the $\gamma_i$.
For $0\leq j\leq N+1$, let $\U^j=\J^j$; these are exactly the puzzle
pieces in the principal nest.
For $j> N+1$, 
let $\U^j$ be the components of 
$R^{-j}(\J^0)$ that contain $\beta$ and $\tau(\beta)$.
These are no longer puzzle pieces in the principal nest
(in particular, they do not contain the escaping critical point);
they are the ``outermost''pullbacks of $\J^{N+1}$ by $R$ 
that intersect the real line, 
and they are mapped diffeomorphically onto $\J^{N+1}.$
For $i\in\mathbb N$,
we let $e_i$ the union of components $e_i^j$ of $R^{-i}(e_0)$
with the property that there exists a path in
$\bm U^i\setminus(\overline{\bm{U}}{}^{i+1}\cup R^{-i}(J^0)$
that connects $e_i^j$ with $U^i\setminus (U^{i+1}\cup[\beta,\tau(\beta)])$.
As before, for all $i>1$, $e_i$ consists of four components.
For $i\in\mathbb N$, we let $\alpha_i\subset\partial \U^i\setminus[\beta,\tau(\beta)]$
be the shortest path crossing the real line that connects
two points of $\partial e_i.$
We let $\gamma_i\subset \alpha_i$ be the union of shortest arcs
that connect $\partial e_i$ with $\partial e_{i+1}$.
Then for $i\geq 3$, $\gamma_i$ is a union of four disjoint connected arcs
$\gamma=\cup_{j=1}^4\gamma_{i}^j.$
We have that
$$e_3\cup\bigcup_{i=3}^\infty\{\gamma_i,e_{i}\}$$
union of four disjoint paths that connect $\partial\J^3$ with
$\{\beta,\tau(\beta)\}.$
The arcs
$\partial\J^3\cup\gamma$ bound three regions in 
the complex plane.
Let $\Gamma$ be the union of the two
regions bounded
by $\partial \bm J^3$ and $\gamma$ 
that do not intersect $(\beta,\tau(\beta)),$
see Figure~\ref{fig:unimodalhigh}.
We obtain a neighbourhood of the critical point 
given by $\Comp_cR^{-1}(\Gamma)$.

We construct the neighbourhood 
$\Gamma'\supset\Gamma$ using 
staircases exactly as before.
For $i=2,\dots,\infty$, 
let $e_i'$ be the components of $R^{-(i-1)}(J^1\setminus J^2)$
that are contained in $\mathbb C\setminus\mathbb R$,
and that are not separated from the real line.
For $i\in\{2,3,\dots,\}$,
$e_i'$ connects $\alpha_i'\subset \partial\bm U^i$ with 
$\alpha_{i+1}'\subset\partial\bm U^{i+1}$.
Let $v_i'=e_i'\cap\alpha_i'$ and $w_i'=e_i'\cap \alpha_{i+1}'$.
For $i=1,2,\dots,$ let $\gamma_i'\in\alpha_i'$ be the union
of arcs connecting $w_{i-1}'$ with $v_i'$.
Let 
$$\gamma'=e_3'\cup_{i=4}^{\infty}\{e_i'\cup\gamma_i'\}.$$
$\gamma'$ consists of four arcs connecting $\partial \J^3$
to $\{\beta,\tau(\beta)\}$. 
Let $\Gamma'$ be the union of the
regions bounded
by $\partial \bm J^3$ and $\gamma'$ that intersect the
$J^3\setminus J^4$.
Now, we pull $\Gamma\subset\Gamma'$ back by one iterate of $R$.
See Figure~\ref{fig:unimodalhigh} for the detail of the construction
up to this point.

Let us now partition $\J^0$; this partition will serve as a
qc\textbackslash bg partition. We take 
$$X_2=(\J^0\setminus\J^3)\cup\Gamma\cup\comp_c 
R^{-1}(\Gamma),$$
$$X_3=\J^4\setminus\Big(\Gamma\cup\comp_c 
R^{-1}(\Gamma')\Big),$$
$$X_1=(\bm J^3\setminus(\bm J^4\cup\Gamma'))\cup(\Gamma'\setminus\Gamma) \cup 
(R^{-1}(\Gamma')\setminus R^{-1}(\Gamma)),$$
and the set $X_0$ is the union of their boundaries.
To see that 
there exists $K_1\geq 1$ so that each component of
$X_3$ is a $K_1$-quasidisk, we verify the
conditions of Lemma~\ref{lem:qdc}.
This is essentially the same as the proof in the low return case, so 
we only give a sketch.

\medskip
\noindent\textit{Claim: $X_3$ is a union of quasidisks.}
Let $\bm W$ be a connected component of $X_3$.
We can decompose $\partial\bm W$ into arcs
$\zeta_1,\zeta_2,\zeta_3,\zeta_4$ where 
$\zeta_1\subset \partial\J^4$,$\zeta_2\subset \gamma'$,
$\zeta_3$ consists of an interval in the real-line
and $\zeta_4\subset \comp_c \partial R^{-1}(\Gamma').$
Since $\J^4$ is a $K_0'$-quasidisk,  we have that
$\zeta_1$-is a $K_0''$-quasiarc, and
$\dist(\zeta_1,\zeta_3)$ is bounded from below.
We have that $\dist(\zeta_2,\zeta_4)$  is comparable to $\diam(\J^0),$
since close to the real line, $\zeta_2$ is close to $\beta$ and 
$\zeta_4$ is close to $c_0$ and away from the real line they are
separated by topological squares with bounded geometry and diameters
comparable to $\diam(\J^0)$. It remains for us to show that
$\zeta_2$ is a $K'$-quasiarc, 
\label{page:gamma qa}
and that the
Ahlfors-Beurling Criterion is satisfied on adjacent quasiarcs.
To see that $\zeta_2$ is a  $K'$-quasiarc, we verify the 
Ahlfors-Beurling Criterion. If $z_1,z_2\in\zeta_2$, let
$\gamma_{[z_1,z_2]}$ be the path between them. If $z_1,z_2\in e_i$ for
some $i$, then we have that $\diam(\gamma_{[z_1,z_2]})\leq C|z_1-z_2|$
for some $C$ since $e_i$ is a pullback of a real interval with bounded
dilatation. If $\gamma_{[z_1,z_2]}\subset\alpha_i$ for some $i$, then 
$\diam(\gamma_{[z_1,z_2]})\leq C|z_1-z_2|$ since $\alpha_i$ is a
pullback of a subset of $\partial \J^0$ with bounded dilatation.
In general, we have that $\diam(\gamma_{[z_1,z_2]})\leq C|z_1-z_2|$
since the path between $z_1$ and $z_2$ can be decomposed into 
$e_i$ and pieces of $\alpha_i$ that intersect at approximately right
angles, form the boundaries of squares with bounded geometry,
and the squares bounded by $\alpha_i,\alpha_{i+1}, e_i$ and  $J_i\setminus
J_{i+1}$ also have bounded geometry.
Similarly we can argue that the Ahlfors-Beurling Criterion
holds on $\zeta_3$ and its adjacent arcs.
To see that it holds on
$\zeta_1$ and its adjacent arcs we need only observe that 
$\zeta_2$ and $\zeta_4$ intersect $\zeta_1$ at an almost right
angle, and that $\zeta_i\cap\J^5$ is bounded away from $\zeta_1$ for
$i=2,4$.
\endpfclaim

Also, arguing in the same way as in the 
saddle node case, we have that $X_1$ is a union of
topological squares with bounded geometry and moduli bounds in $\bm
J^0$, and in the upper half-plane. 
Around a component $\bm S'$ of $X_1$ that is contained in 
$\J^4\setminus \J^5$ we can build a topological square
$\hat{\bm{S}}'\supset \bm S'$ such that
$\mod(\hat{\bm{S}}'\setminus \bm S')$ is universally bounded from
below, and any component $\bm S$ of $X_1$ is a pullback of such an
$\bm S'$ and the mapping from $\bm S$ to $\bm S'$ extends to a
diffeomorphism with bounded dilatation onto $\hat{\bm{S}}',$
and the result follows.


\medskip

\noindent\textit{Step 2: Constructing the quasiconformal mapping.}
In $X_2$, the construction is identical to the one used in the 
low return case. Let $\bm P$ denote the puzzle piece that contains 
$R^{k_c}(c),$ then for each $z\in X_2\setminus \comp_c R^{-k_c}(\bm
P)$, there exists $k_z$ such that $R^{k_z}(z)\in \J^0\setminus \J^1$.
We define $H(z)$ by $H_0\circ R^{k_z}(z)=\tilde R^{k_{z}}\circ H(z).$
On  $\comp_c R^{-k_c}(\bm P)$ we define $H$ so that it agrees with the
mapping on the boundary using the QC Criterion.

Now let us define $H$ on $X_3$. 
Let $\bm W$ be a component of
$X_3$, and let $\widetilde{\bm{W}}$
be the corresponding set for $\tilde R$.

\medskip
\noindent\textit{Claim: There exists a $\hat K$-quasiconformal mapping
$H_{\bm W}\colon\bm W\rightarrow\widetilde{\bm{W}}$ such that 
for each $x_i\in\cup_{i=3}^\infty \{v'_i,w'_i\} \cup
R^{-1}(\cup_{i=3}^\infty\{v'_i,w'_i\},$
$H_{\bm W}(x_i)=\tilde x_i$.}
We have that $\bm W$ and $\widetilde{\bm{W}}$ are $K_1$-quasidisks.
Let us start by partially defining a mapping
$h_{\bm W}$ from $\bm W$ to $\widetilde{\bm{W}}$.
For each $x_i\in\{v'_i,w'_i\}, i=2,3,\dots,$, we set $h_{\bm W}(x_i)=\tilde
x_i$.
Pulling back once,
for each point $x\in \comp_c R^{-1}(\Gamma')$ such that
$R(x)\in \cup_{i=N-2}^{\infty}\{v_i',w_i'\}$, we define
$h_{\bm W}$ so that
$h_{\bm W} \circ R(x)=\tilde R\circ h_{\bm W}(x)$.
Let $\{y_i\}_{i=1}^\infty$ denote the set of all
points not in $\partial \bm J^2$ where $h_{\bm W}$
has been defined.
Take any $x,y,z\in \cup\{y_i\}_{i=1}^\infty\subset \partial\bm W$
arranged with $y$ between $x$ and $z$ in
$\partial \bm W.$
Because of the linear behaviour near the $\beta$ 
fixed point, 
there exist $C>0$ and so that
$\frac{|\tilde x-\tilde y|}{|\tilde y-\tilde z|}\leq
C\Big(\frac{|x-y|}{|y-z|}\Big)^\eta.$

Now we extend $h_{\bm W}$ to $\bm W$.
Let $\phi_{\bm W},\phi_{\widetilde{\bm{W}}}\colon\mathbb C\rightarrow\mathbb C$
be $\hat K$-qc mappings of the plane such that
$\phi_{\bm W}(\mathbb D)=\bm W$,
 and $\phi_{\tilde{W}}(\mathbb D)=\widetilde{\bm{W}}$.
Arguing as in the unimodal, low return case, we have that
the mapping 
$\phi^{-1}_{\widetilde{\bm{W}}}\circ h_{\bm W}\circ \phi_{\bm W},$
which maps $\phi_{\bm W}^{-1}(\mathrm{Domain}(h_{\bm W}))$ into 
$\partial\mathbb D$ is $\kappa'$-quasisymmetric, 
and that there exists $M\geq 1$ depending only on $|J^0|/|J^1|$ so that 
$\phi_{\bm W}^{-1}(\mathrm{Domain}(h_{\bm W}))$
is an $M$-relatively connected subset of $\mathbb D$.
So by Proposition~\ref{prop:Vellis},
we have that $\phi^{-1}_{\widetilde{\bm{W}}}\circ h_{\bm W}\circ \phi_{\bm
  W}$
extends to a $\kappa''$-qs mapping from $\partial\mathbb D$ to itself.
Thus by using the 
Ahlfors-Beurling extension
it extends to a $\hat K_1$-qc mapping on the unit disk, and so 
$h_{\bm W}$ extends to a $\hat K$-qc mapping
$H_{\bm W}\colon\bm{W}\rightarrow\widetilde{\bm{W}}$. 
\endpfclaim

\medskip

Now let us extend the mappings that we have
defined to $X_1$. 
Fix a small $\delta'>0$.
Let $\bm S$ be a component of $X_1,$ with corners
$v',w', v$ and $w$.
The mapping $H$ may already be defined at some of the corners.
We extend $H$ to $\overline{\bm{S}}$, 
so that it is a homeomorphism and satisfies
(\ref{eqn:qc bound 0}) on $X_0,$ which contains the boundary of
each component of $\bm S$.

\medskip

Since $H$ is a conjugacy on the boundary of $X_2$,
we can pull this construction back along the two
monotone branches of $R.$ Because $\J^0\setminus \Gamma$ is contained in a
Poincar\'e disk, $D_{\theta}((\beta,\tau(\beta)))$ with $\theta$ bounded away
from zero, and each monotone branch of $R$ is mapped
onto an interval containing $(\beta,\tau(\beta))$ and bounded
by a boundary point of $J^0$ and $\tau(\beta)$, 
we can control the dilatation of $H$ as we pullback it along these
branches. In particular, pulling back along these branches
infinitely many times refines the mapping
$H$ in $X_3$: under $R$, every point in $X_3$ 
either lands in $\overline{\comp_c R^{-1}(\Gamma)}$ or
it escapes through $\J^4\setminus \J^5$,
and the orbit of
almost every point in $[\beta,\tau(\beta)]\setminus \comp_c R^{-1}(\Gamma)$
enters $\comp_c R^{-1}(\Gamma)$. 
Thus we obtain a mapping on the real line that is a
conjugacy up to the first entry time of a point in $N$.
Moreover $H$ is defined dynamically
on $\Gamma$ and on $X_3$ (after refining the original mapping).
Observe that if $H_0$ is quasiconformal on a neighbourhood of the interval
$J^0\setminus J^1$ then so is $H$.

\medskip

\noindent\textbf{Case 3: $R\colon J^1\rightarrow J^0$ is a monotone mapping.}
Here we explain how to deal with the case when $R\colon J^1\rightarrow J^0$ 
is a monotone mapping with possibly several critical points of odd
order. Suppose first that $|J^0|/|J^1|$ is bounded.
To simplify the exposition, we will assume that $R$ is orientation
preserving; the orientation reversing case is similar.
We will focus on the key differences of this case with the 
two unimodal cases already considered.

We construct $\gamma$ and $\gamma'$, by considering
infinitely many pullbacks of the fundamental domains
 $J^0\setminus J^1$ lying in $\mathbb
C\setminus \mathbb R$, which are ``closest to the real line,''
and connecting them by
boundary arcs of puzzle pieces.

Let $\beta\in J^{N+1}$ denote the fixed point of $R$.
By 
Lemma~\ref{lem:CvST 3.10}, and \cite[IV. Theorem B]{dMvS},
if $|J^0|$ is sufficiently small,
there exists a constant $C$ such that $1+1/C<|DR(\beta)|<C$.
Notice that $\gamma_i,\gamma_i',$ 
$i\in\{1,2,3,4\}$ all have a common boundary point at $\beta$.
Let us now construct the regions $\Gamma$ and $\Gamma'$.
There are two cases to consider,
\begin{enumerate} 
\item all critical points of $R$ are contained in a single component of 
$J^0\setminus\{\beta\}$ or 
\item each component of 
$J^0\setminus\{\beta\}$ contains at least one critical point of $R$.
\end{enumerate}
Suppose first that all critical points of $R$ are contained in one
 component of 
 $J^0\setminus\{\beta\}.$
Let $N'$ be maximal so that $J^{N'+1}$ contains a critical point $c'$ of
$R$.
As in the unimodal, low return case, we need to be a little careful
about complex pullbacks of the rectangles bounded by
$\alpha_{i+1}\cup\alpha_i\cup e_i$ or
$\alpha_{i+1}\cup\alpha_i\cup e_i'$
(whose real traces are a fundamental domain in $J^i\setminus
J^{i+1}$), but
observe that $R$ has at most $\crit(F)$ critical points, and
each passes through $J^i\setminus
J^{i+1}$ at most once, so we can construct
a shielding region about each critical point as in the 
unimodal, low return case, see page~\pageref{pageref:dncv}
 and Figure~\ref{fig:cvshield}.
Let us construct a neighbourhood about $c'$ under the
assumption that $R(c')\in R^{-N}(\II_n\setminus \II_n^-)$;
the alternative case is similar.
Let $\Gamma$ be the region bounded by 
$\partial \J^4\cup\gamma\cup\partial\J^{N'+1}$ that contains a
critical value of $R$ and let $\Gamma'$
be the region bounded by 
$\partial \J^3\cup\gamma'\cup\partial R^{-(N+1)}(\II_n^+)$ 
that contains a critical value of $R$.
If each component of 
$J^0\setminus\{\beta\}$ contains at least one critical point of $R$,
then we make the same construction on either side of $\beta$
to obtain a neighbourhoods $\Gamma_1$ 
of the component of $J^0\setminus J^{N'}$ that contains $R(c')$
and a neighbourhood $\Gamma_2$ of 
the component of $J^0\setminus J^{N''}$ that contains $R(c'')$.
We let $\Gamma=\Gamma_1\cup\Gamma_2$, and
define $\Gamma'$ analogously.
We omit the details.

There exists $\delta>0$ such that
$\Gamma'\setminus\Gamma$
consists of a union of topological squares $\bm S$
with
$\delta$-bounded geometry,
so that $\mod(\J^0\setminus\bm S)\geq \delta$ and 
with moduli bounds as in (3) of Theorem~\ref{thm:central cascades}.
To see that these bounds hold, we can argue just as in
the unimodal case.

We let 
\begin{itemize}
\item $X_1=\bm \J^3\cap(\Gamma'\setminus\Gamma)\cup(\J^3\setminus(\J^4\cup\Gamma'))$,
\item $X_2=(\J^0\setminus \J^3)\cup\Gamma$, 
\item$X_3=\J^0\setminus\overline{X_1\cup X_2\cup\mathbb R},$
and 
\item $X_0=\J^0\setminus(X_1\cup X_2\cup X_3).$
\end{itemize}
We define $H$ on $X_2$ by pulling back when we can
and applying the QC Criterion to construct 
$K'$-qc mappings that match the boundary markings
on puzzle pieces that contain critical points. Here are the details.

First we will construct a complex neighbourhood of
$\Crit(F)$.\label{page:mononbhd}
We will use the same construction
in the multicritical saddle-node and high-return cases. 
Let $\{\widetilde{\bm{G}}_j\}_{j=0}^s$ 
be the chain with $\widetilde{\bm{G}}_s=\J^0$,
and $\widetilde{\bm{G}}_0=\J^1$ and 
let $\{\bm{G}_j\}_{j=0}^s$ 
be the chain with $\bm{G}_s=\J^1$,
and $\bm{G}_0=\J^2.$
Let $c$ be a critical point of $F$.
If there exists $i,$ $0< i\leq s$ such that 
$c\in\tilde G_i,$
let $\J_c=\widetilde{\bm{G}}_i=\hat{\mathcal{L}}_c(\J^0).$ 
Otherwise
we define $\J_c$ as follows:
Let $n$ be so that $\W^n_{c_0}$
is at the top of the same long central cascade 
as $\J^0$. If $\W^n_c$ is at the top of a long 
central cascade, then we take $\J_c$ to be
the excellent puzzle piece associated to the long central cascade
by Proposition~\ref{prop:good geometry for central cascades},
otherwise we take $\J_c=\W^n_c$.
Let $T'$ be the set of critical points of $F$
such that $\J_c=W^n_c$.
Let $\hat R$ be the return mapping to $\cup_{c\in\crit(F)}\J_c$.
Let $\bm P_1'$ denote
the components of
$\hat R^{-1}(\Dom(\cup_c\J_c)),$ that intersect
$J^0\setminus J^{1}$.

\medskip

When we pull back a puzzle piece in $\bm P_1'$
through the central cascade,
it could intersect several critical points, we need to 
ensure that the critical values are contained in
puzzle pieces with bounded geometry and moduli bounds.

Let us suppose first that $\bm W$ is a pullback of some 
$\bm J_c$ that is a puzzle piece given by
Proposition~\ref{prop:good geometry for central cascades},
which is at the top of some central cascade.
To be definite, let 
$$\bm W=\comp_{R^{N_1+1}(c_1)}(\Dom(\cup_{c\in\crit(F)}\J_c)).$$
Then $\bm W$ is mapped diffeomorphically onto
$\J_c$.
Consider the chain 
$\{\bm Z_j\}_{j=0}^N$ where $\bm Z_N=\bm W$
and for $j=0,\dots, N-1$, we set $\bm Z_{j}$ to be
the component of $(R_{\J^0}|_{\J^1})^{-1}(\bm Z_{j+1})$ that
intersects the real line. Each $\bm Z_j$ intersects $J^j\setminus J^{j+1}$.
Let $b_0'$ be the number of $\bm Z_{j}, j=0,\dots, N$ that contain a critical
point of $R=R_{\J^0}|_{\bm J^1}$. 
Let us index these critical points by $c'_i$ where $c'_i\in \bm
Z_{j_i}$, 
and $j_1<j_2<\dots<j_{b_0'}$. Let $T=\{c_1',c_2',\dots,c'_{b_0'}\}$.
For each $c_i',$ $j_i\leq N+1$ let $v_i:=R^{j_i}(c_i')$.

Let us now refine $\bm W$:
For any $v_i,$ $i\in\{1,\dots,b_0'\}$,
let $s$ be minimal
so that either
\begin{itemize}
\item $\hat R^{s}(v_1)\in \J_c$ where
$c\in T',$ in which case $\J_c=\W^n_c$ is
a puzzle piece in the good nest
or
\item $\hat R^{s}(v_1)\in \mathcal{L}_c(\J_c),$ for
$c\in\crit(F)\setminus T'$,
in which case
$\J_c$ is given by
Proposition~\ref{prop:good geometry for central cascades}.
\end{itemize}
If $c\in T'$, let $\bm Q_{v_1}=\comp_{v_1}\hat R^{-s}(\J_c)$,
and $\bm Q_{v_1}'=\comp_{v_1}\hat R^{-s}(\comp_{R^s(v_1)}(\Dom(\W^n))),$
and if $c\in\crit(F)\setminus T'$,
let $\bm Q_{v_1}=\comp_{v_1}\hat R^{-s}(\J_c)$,
$\bm Q_{v_1}'=\comp_{v_1}\hat R^{-s}(\mathcal{L}_c(\J_c)).$
By Proposition~\ref{prop:good geometry for central cascades},
Proposition~\ref{prop:good bounds} , Corollary~\ref{cor:Wnice}
and Proposition~\ref{prop:squares},
we have that there exist $K_1\geq 1$ and $\delta_1>0$ 
such that
$\bm Q_{v_1}'$ has 
$\delta_1$-bounded geometry, that
$\mod(\bm Q_{v_1}\setminus \bm Q_{v_1}')\geq\delta_1,$ and
that $H$ is $(K_1,\delta_1)$-qc\textbackslash bg
on each  component of that contains a 
forward image of a critical point.
So by the QC Criterion there exists a $K'=K'(K_1,\delta_1)$-qc
mapping on each $\bm Z_j, j=1,\dots, N$
that agrees with the boundary marking.\label{page:monoconj}
Finally, we can 
extend $H$ to $(\J^0\setminus\J^1)\setminus (\cup\bm Q_v)$
using the same argument given to prove the Real Spreading 
Principle.

\medskip

Once again, we have that $X_3$ is a union of disjoint 
$K_1$-quasidisks. This is checked as in the unimodal, low return case
using Lemma~\ref{lem:qdc}.
We define $H$ on each component $\bm W$of $X_3$
just as we did in the unimodal case:
\begin{itemize}
\item First, we define $h_{\bm W}$
on
$\bigcup_{i=3}^{\infty}\{v_i',w_i'\}$,
so that it matches the boundary marking on this set. 
As in the unimodal case, it is straightforward to
check that $h_{\bm W}$ 
defined on this set is
is $(C,p)$-qs for some $C,p$ depending on $|J^0|/|J^1|$, and the
escape times of the critical points.
Note that when we pass through a long piece of the cascade that contains 
an escaping odd critical point, the cascade behaves according to the Yoccoz
Lemma,as in the unimodal, low return case.
\item There exists a $K_1$-quasiconformal mapping
 of the plane so that
$\phi_{\W}\colon\mathbb D\rightarrow\bm W$. 
\item Just like in the unimodal case,
we have that
$\phi_{\W}^{-1}(\bigcup_{i=3}^{\infty}\{v_i',w_i'\})$
is an $M$-relatively connected subset of $\partial\mathbb D$,
 so by Proposition~\ref{prop:Vellis}, we have that
$\phi_{\widetilde{\W}}^{-1}\circ h_{\bm W}\circ \phi_{\W}$ 
extends to a $\kappa'$-qs mapping on $\partial\mathbb D$.
\item Finally since $\phi_{\widetilde{\W}}^{-1}\circ h_{\bm W}\circ \phi_{\W}$ 
extends to a $\hat K_2$-qc mapping on $\mathbb D$, we have that
$h_{\bm W}$ extends to a $\hat K$-qc mapping
$H_{\bm W}\colon\bm W\rightarrow\bm W$.
\end{itemize}

\medskip

On $X_1$ we define $H$ so that it is a homeomorphims, which
agrees with the mapping already
defined on a subset of 
$\partial X_1$ and so that it 
satisfies (\ref{eqn:qc bound 0}) on $X_0$.
So there exists a $(K',\delta')$-qc\textbackslash bg
mapping $H\colon\mathbb C\rightarrow\mathbb C$ such that 
$H(\J^0)=\widetilde{\J}^0$ as in the statement of 
Theorem~\ref{thm:central cascades}. Now we pull this construction back 
by one iterate of $R$.

The mapping constructed up to this moment is a conjugacy outside of a
neighbourhood
of the $\beta$ fixed point.
Let $\bm W_0=\comp_{\beta}\J^4\setminus R^{-1}(\Gamma)$.
We are going to refine $H$ in
$\bm W_0$, so that we obtain a mapping that 
is a conjugacy in a neighbourhood of $J^1$
(the real trace of the dynamics).
Let $\bm W^1=\comp_{\beta}(R^{-1}(\bm W_0))$.
By definition $H$ conjugates the dynamics of $R$ and
$\tilde R$ on $\partial \bm W^1$, so 
we can pull it back by $R|_{\bm W^1}^i$ for $i=1,2,3,\dots$.
Since $\bm W_0$ is contained in
a Poincar\'e disk with angle bounded away from zero based on its
real trace, the dilatation as we pullback is small, and we obtain a
conjugacy
between $R$ and $\tilde R$ on $\bm W^0$.

\medskip

Suppose now that $|J^0|$ is much bigger than $|J^1|$.
Since $\J^0$ is a $K_0$-quasidisk,
there exists $\theta\in(0,\pi)$ such that
$\J^0\subset D_{\theta}(J^0)$.
So by Lemma~\ref{lem:zd}, for $0\leq i\leq N+1$
we have that there exists $\lambda\in(0,1)$ such that
$\J^i\subset D_{\lambda^i\theta}(J^i)$.
Thus we have that there exists a constant $C$ such that
$\diam(\J^i)\leq C|J^{i}|/\lambda^i\theta$.
So if $|J^{i+1}|$ is much smaller than $|J^i|,$ for $i=0,\dots, i_0$,
with $i_0\in\{0,\dots, N+1\}$,
we have that $|J^{i+1}|$ is much smaller than 
$\lambda|J^i|,$ and so $\frac{|J^{i+1}|}{\lambda^{i+1}\theta}$ is much
smaller than $\frac{|J^{i}|}{\lambda^{i}\theta}$.
Hence there exists a constant $\eta<1$ so that 
$\diam(\J^{i+1})<\eta\diam(\J^{i}).$ So we can assume that
the diameters of the puzzle pieces $\J^i$ decay exponentially, 
for $0\leq i\leq i_0$ and that $|J^{i_0}|$ is comparable to
$|J^{i_0+1}|$ or $i_0=N+1$.

We construct $\bm Q_v\subset\bm Q_v'$ 
just as in the bounded geometry case, and
and refine the conjugacy $H$
exactly as before.
As long as $J^i$ is much bigger than
$J^{i+1}$ we can pull back the conjugacy,
and we deal with puzzle pieces that contain critical points
exactly as in the bounded geometry case.
So we pull back until the moment when the geometry becomes bounded.
Let this moment be $i_0$, so that $|J^{i_0}|$ is comparable to $|J^{i_0+1}|$.
If all critical values of $R$ are contained in $J^{0}\setminus
J^{i_0+1}$ or if both components of $J^{i_0}\setminus J^{i_0+1}$ are
comparable to $J^{i_0}$ then we can conclude the proof using the
argument in the bounded geometry case. It could be the case that one
of the components $F_{i_0}$ of $J^{i_0}\setminus J^{i_0+1}$ is very small. In
this case, we look at the forward images of this fundamental
domain: $F_i=\comp_{R^{i-i_0}(F_{i_0})} J^i\setminus J^{i+1}$.
By the Koebe theorem, we can find $i$ so that $|F_i|\asymp |J^i|$,
and now we can repeat the argument from the bounded geometry case.

\medskip

\noindent\textbf{Case 4: $R\colon J^1\rightarrow J^0$ is a multicritical
  mapping with a low return and $|J^0|/|J^1|$ is bounded.}
This is the case when $R(J^1)$ does not contain the critical
point, and $R$ does not have fix points on the real line,
as in Figure~\ref{fig:returns}(2), but with more critical points.
We will treat the case when the geometry is unbounded 
after we consider the multicritical, high return case.
First, we construct $\bm Q_v\subset\bm Q_v'$
exactly as in the monotone case, see page~\pageref{page:mononbhd}. 
Now, let us show that we can reduce this case to the case when there
is a single turning point. 
For convenience, we will choose the chain $\{G_j\}_{j=0}^s$
so that $G_0=J^2$ contains a critical point $c_0$ of even order, 
and $G_s=J^1$. 
Then $F^s(G_0)\not\owns c_0.$
If $s'<s$ is maximal so that $G_{s'}$ contains a critical point $c'$
of $F$ of even order,
we have that that $F^{s-s'}(c')$ is a global maximum for
$F^s|_{J^1}$ and $F^{s-s'}(G_{s'})\not\owns c_0.$
Thus if we let $s_1>0$ be minimal so that $G_{s_1}$ contains a
critical point $c_1$ of even order,
we have that $F^{s-s_1}(G_{s_1})\not\owns c_0,$
so that the return mapping $F^s\colon G_{s_1}\rightarrow \tilde G_{s_1}$,
has one fewer critical points than $F^{s}|_{G_0}$, and since $F^s|_{G_{s_1}}$
does not have a fixed point, we have that
$F^s(G_{s_1})$ does not contain $c_1$.
Continuing in this fashion we obtain that the return mapping
$F^s\colon G_{s'}\rightarrow \tilde{G}_{s'},$
$\tilde G_{s'}=\comp_{F^{s'}(c_0)}F^{-(s-s')}(J^0),$
 is a unimodal
mapping. Moreover,
$c'$ corresponds to the last critical point of $R$ to
escape the central cascade.

Using the arguments from the previous cases we obtain a
decomposition of
$\bm{J}^0_{s'}$ into sets $X_0,X_1,X_2$ and $X_3$:
We build the paths $\gamma$ and $\gamma'$ just as 
before. Suppose that $N$ is so that
$J^N\setminus J^{N+1}$ contains the even critical point.
Then $\gamma,\gamma'$ are both unions of four paths
that connect $\partial\bm J^3$ with $\partial R^{-1}(J^{N+1})$.
We construct a 
neighbourhood of the critical points just like we did in the
unimodal saddle node and monotone cases, compare Figure~\ref{fig:cvshield},
but just as in the monotone case we have to argue more carefully since
there could be odd critical points in the monotone branches of $R$:
Observe that
there are at most $\#\crit(f)-1$ such odd critical points, 
so $\#\{R^{k_c}(c):c\in\crit(R)\}\cap \bm B_0$ is at most 
$\#\crit(f)-1$. Thus there is a subinterval $J$ of $B_0$ with length
comparable to $|J^0|$ such that if $J'$ is any component of
$R^{-(N+1)}(J)$
that is contained in the real line then
$R^{N+1}\colon J'\rightarrow J$ is a diffeomorphism.
Now as in the unimodal and monotone cases, we can obtain neighbourhoods of
the critical points which are ``shielded'' by topological rectangles
with
bounded geometry and moduli bounds.
With this preparation,
we define $X_i,i\in\{1,2,3,4\},$ exactly as in the unimodal, low return case,
and we can construct $H$ in the
same way as we did in the unimodal
and monotone cases (since there could be critical points in the
monotone branches of $R$).
We will deal with the big geometry case
after we treat cascades with a high return.

\medskip

\noindent\textbf{Case 5: $R\colon J^1\rightarrow J^0$ is a high return.}
First let us consider the case when $|J^0|/|J^1|$ is bounded.
Let $\beta$ denote the outermost orientation
preserving fixed point of $R\colon J^1\rightarrow J^0.$
Let $\lambda_\beta$ be the multiplier of $\beta$.
Provided that $J^0$ is small enough, which we can and will assume,
by 
\cite[Chapter IV. Theorem B]{dMvS} there exists a $\varepsilon_0>0,$
universal, so that
$|\lambda_\beta|>1+\varepsilon_0,$ and by Lemma~\ref{lem:CvST 3.10},
$|\lambda_{\beta}|$ is bounded from above. 

We begin as in the unimodal case by constructing two regions
$\Gamma\subset\Gamma'$ whose boundaries consist of staircases landing
at
$\beta$ and $\tau(\beta)$ respectively and boundary arcs of
$\partial\J^4$.
To construct a conjugacy on the boundaries of puzzle pieces in $\Gamma$ 
that contain critical points, we use the same procedure as in the
monotone case, see page~\pageref{page:mononbhd}. 

Let us now construct a first qc\textbackslash mapping and partition. 
Let 
\begin{itemize}
\item $X_1=\Gamma'\setminus\Gamma\cup(\J^3\setminus(\J^4\cup\Gamma')$,
\item $X_2=\Gamma\cup \J^0\setminus\J^3$
\item $X_3=\J^4\setminus\overline{X_1\cup X_2\cup\mathbb R}$
\item $X_0$ is the union of their boundaries.
\end{itemize}
We construct a qc\textbackslash bg mapping 
$H'$
associated to this partition exactly as in the unimodal case:
Arguing as before, we have that 
there exists a $\hat K$-qc mapping from 
each component
$\bm W$ of $X_3$ to the corresponding
component $\widetilde{\bm{W}}$.
Moreover,
We have that
$(\bm J^3\setminus (\bm J^4\cup\Gamma')$ and
$\Gamma'\setminus\Gamma$, 
are both unions of 
combinatorially defined topological squares with bounded geometry and
moduli bounds as in (3) of Theorem~\ref{thm:central cascades}.
So we can extend the mapping constructed in $X_2$ to a
qc\textbackslash bg mapping $H\colon\J^0\rightarrow\widetilde{\J}^0$.

\begin{figure}[htbp]
\begin{center}\label{fig:multi-high}
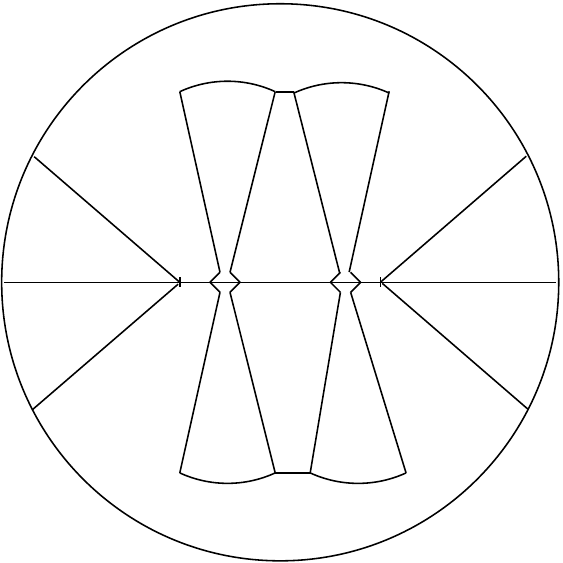
\end{center}
\caption{A sketch of the decomposition for a multicritical high
  return.}
\label{fig:multicriticalhigh}
\end{figure}

Now we pull this picture back by one iterate of $R$.
For each critical point $c$
such that $R(c)\notin[\beta,\tau(\beta)]$, we have that 
$c\in\comp_cR^{-1}(\Gamma)$.
Let $T_0$ denote the critical points of $R$ 
that do not escape $[\beta,\tau(\beta)]$ under $R$. We have that
$T_0\subset J^{N_1+2}$.
If $T_0=\emptyset,$
then we can conclude the argument exactly as in the unimodal case, 
so let us suppose that $T_0\neq\emptyset.$
We will construct neighbourhoods of these critical points.
We let $\bm{\mathcal{W}}'$ be the union of 
components of $\J^3\setminus(\Gamma\cup R^{-1}(\Gamma))$ that
contain the critical points $c\in T_0,$ and
we let $\bm{\mathcal{W}}$ be the union of 
components of $\J^3\setminus (\Gamma\cup R^{-1}(\Gamma'))$ that
contain critical points $c\in T_0$. Observe that
$\bm{\mathcal{W}}\subset \bm{\mathcal{W}}'.$
We let $\bm W_0'$ be a component of $\bm{\mathcal{W}}'$
and  let $\bm W_0$ be the component of $\bm{\mathcal{W}}$
contained in $\bm W_0'$.
As before, it is not hard to use Lemma~\ref{lem:qdc}
to show that $\bm W_0$ and $\bm W_0'$ 
are $K_1$-quasidisks.
We also have that $\bm W_0'\setminus \bm W_0$ is
a union of topological squares with bounded geometry and moduli bounds
as in (3) of Theorem~\ref{thm:central cascades}.
See Figure~\ref{fig:multicriticalhigh}, for a rough idea of
what the decomposition looks like in this case.

For each critical point $c\in T_0$, 
let $k_c\in\mathbb N$ be minimal so that
$R^{k_c}(c)\in [\beta,\tau(\beta)]\setminus \mathcal{W}$.
Order the critical points in $T_0$,
$c_1<c_2<c_3<\dots<c_{|T_0|}$,
so that $k_{c_i}\leq k_{c_{i+1}}$.
Let $i_0<b$ be minimal so that 
$k_{c_{i+1}}-k_{c_i}\geq N_0$.
If no such $i_0$ exists, set $i_0=|T_0|$. 
Let $T_1$ denote the set of critical points
$c_i$ with $i\leq i_0$.
For each critical point $c_i$ with $i\leq i_0$, we have that
$k_{c_i}\leq |T_0|N_0$. Let $N_{2}=k_{c_{i_0}}$.
Let $T_2$ denote the set of critical points that do not escape
$\mathcal W$ under $R^{N_0+N_2}$. This is the same as the set of critical
points $c_i$ with $i>i_0$.

Let $\bm{\mathcal{W}_1'}$ denote the components of $\bm
J^3\setminus(\Gamma\cup R^{-1}(\Gamma)).$
that do not contain a critical point of $R$.
Suppose that $\bm W$ is a component of 
$\bm{\mathcal{W}_1'}$
Then by construction $H$ is a
conjugacy on $\partial\bm W$. So as in the unimodal, high return case,
by pulling $H$ back through infinitely many iterates of $R|_{\bm{\mathcal{W}}_1'}$, 
we construct a conjugacy on the union of such components, up to the
landing time of points in $\bm{\mathcal{W}}\cup\J^0\setminus \J^1.$

When $T_2=\emptyset,$
we can conclude proof by pulling back by a bounded number of steps.
So suppose that $T_2\neq\emptyset$.
For convenience, let us suppose that $c_1\in T_2$.
Then we have that 
$W^0\supset  W^1\supset\dots\supset W^{k_{c_1}-1}\supset W^{k_{c_1}}$
is a long central cascade
and
$R(c_1)\in W^{k_{c_1}-1}\setminus W^{k_{c_1}}.$
There are three cases to consider:
The return mapping to $\bm W^0$ has a low return, a high return or is
monotone. Depending on the case, we one of the arguments we have
already given.
Notice in the low return case, the argument terminates after
one step, since all the critical points escape through the monotone
branches of $R$. In the high return case, if there are more
non-escaping critical points, the argument will not terminate, but
it eventually reduces to the either the unimodal or monotone case.

\medskip

Now let us consider the case when $|J^0|$ is much bigger than $|J^1|$
and $c_0$ is even. Then since $\J^0$ has bounded geometry
and $J^1$ is approximately in the middle of $J^0$,
there exists $C_0$ large so that $|J^0|>C_0|J^1|$
and $\mod(\J^0\setminus J^1)$ is large.

If there exists $\delta_0>0$ such that
for all $n, 0\leq n\leq N$, we have that
$\mod(\J^n\setminus \J^{n+1})\geq\delta_0$,
then, there exists $\delta_0'\in(0,1)$, depending only on $\delta_0$,
so that $\diam(\J^{n+1})<\delta_0'\diam(\J^n)$.
Thus for $n=1,2,\dots, N$, there exists a constant $C>0$ such that
the dilatation of $R|_{\J^{n}\setminus \J^{n+1}}$
is bounded from above by
$$C\max_{0\leq i\leq
  s}\{\diam(\comp_{f^i(c)}f^{-(s-i)}(\J^{n-1}))\},$$
which decays exponentially fast,
and so we can pull the conjugacy back through the entire central
cascade as we did above in the high-return case for the critical
points in $T_1$.

So we can assume that for some $n,$
$\mod(\J^n\setminus\J^{n+1})<\delta'.$
We claim that there exists $C'>1$ such that 
$|J^n|<C'|J^{n+1}|$ and for all $i,0\leq i<n$ we have that
$|J^i|>C'|J^{i+1}|.$ Suppose not, then
$|J^{n+1}|$ is much smaller than $|J^n|$, 
and $\mod(\J^{n}\setminus
  \J^{n+1})\leq D(1+C(\mu(J^{n-1}))\delta_0,$
where $D$ is the degree of $R$.
Since $c_0$ is even, we still have that $J^{n+1}$ is approximately in
the middle of $J^n$.
Thus at least one of $\J^n$ or $\J^{n+1}$ does not have bounded
geometry, even more that
either
the largest Euclidean circle centered at $c$ that
can be inscribed in $\J^n$
has diameter much smaller than $J^n$
or the smallest Euclidean circle centered at $c$ that contains
$\J^{n+1}$ has diameter much larger than $J^{n+1}$.
 But then since for all $i,$ $0\leq i<n$
 $|J^{i}|>C|J^{i+1}|$ the critical value of $F$ is roughly in the
 middle of
 $J^i$, $\J^{i+1}$ cannot be pinched at $c$, which together with the
 moduli bounds for $\mod(\J^i,\setminus \J^{i+1})$ for 
$i=0,\dots, n$,
gives a contradiction.
Thus we can reduce the argument in the big geometry case to
the bounded geometry case. Observe that by Lemma~\ref{lem:qs gluing},
the quasisymmetric mapping that we obtain by
gluing together the mapping $h$ we get by pulling back the
quasiconformal mapping as long as the
geometry is big, and the qs mapping $\hat h$ on $J^n,$
given by argument in the bounded geometry case, is
$\kappa$-qs on $\mathbb R$.
This concludes the proof of 
Theorem~\ref{thm:central cascades}.
\qed

\section{Quasisymmetric rigidity}\label{sec:qc rigidity}
	
Recall that we have a natural partial ordering
	on the set of critical points:   $c_1\ge c_0$ if $c_0\in \omega(c_1)$ or if $c_0=f^n(c_1)$
	for some $n\ge 0$.
	We need to start with critical points which are lowest in
	this ordering because other critical points can accumulate
	on these. We then spread the information to other critical points
	which are higher in this partial ordering. So we start
	with the infinitely renormalizable and attracting cases,
	and proceed to the non-renormalizable cases;
	we first deal with the persistently recurrent and at most finitely renormalizable case,
	and finally treat the reluctantly recurrent case. 
	Recall that the case of critical points that are non-recurrent and do not accumulate on 
	other critical points was 
	treated by Theorem~\ref{thm:qc rigidity away from critical points}
	through the use of touching box mappings.
	We prove that immediate basins of 
	attraction are rigid in
	Proposition \ref{prop:rigidity attracting}.

	To organize things properly, we decompose
	the set of critical points as in Section \ref{subsec:decomp}.
	Let $\Omega$ be a connected component of the graph of critical points.
	Let $\Omega_0\subset\Omega$ be the set of critical points which are in the
	immediate
	basin of a periodic attractor or at which $f$ is infinitely renormalizable.
	Let $\Omega_1\subset\Omega$ be the set of critical points $c$
	for which $f$ is persistently recurrent on $\omega(c)$ and at most finitely renormalizable.
	Let $\Omega_1'\subset \Omega_1$ be
	a smallest subset so that $\cup_{c\in \Omega_1'} \omega(c)= \cup_{c\in \Omega_1} \omega(c)$.
	Let $\Omega_2 = \Omega\setminus(\Omega_0\cup\Omega_1)$ be a set of critical points $c$
	which are reluctantly recurrent or nonrecurrent.
	
	 Let $I_0$\label{page:intervals} consist of the immediate  basins of periodic attractors
	and of  periodic intervals on which the mapping is infinitely renormalizable,
	where if $c'$ is a critical point at which $f$ is infinitely renormalizable, $I_0(c')$ is
	chosen so that if $c$ is any critical point whose orbit does not accumulate on $c'$
	or land on $c'$, then $f^k(c)\notin I_0(c')$ for all $k=0,1,2,\dots$.
	Let $I_1'$ be  nice intervals around critical points $c'\in \Omega_1'$
	again chosen so small that
	if $c$ is any critical point whose orbit does not accumulate on $c'$
	or land on $c'$, then $f^k(c)\notin I_1(c')$ for all $k=0,1,2,\dots.$
	As in Subsection~\ref{subsubsec:persistent infinitely many branches} define
	$$I_1=\bigcup_{c\in \Omega_1\setminus \Omega_1'} \LL_c(I)\cup I'_1$$
	and
	$$I_2=\bigcup_{c\in \Omega_2}I_2(c),$$
	where $I_2(c)=\LL_c(I_0\cup I_1)$ if the latter set is non-empty, and 
	$I_2(c)$ is some arbitrarily small nice interval, disjoint from
	$I_1$ and the other components of $I_2$, otherwise.
	Note that by shrinking $I'_1$ we can choose $I_2$ as small as we like,
	and that $I_0\cup I_{1}\cup I_2$ is by construction a nice set.
	Hence we can assume that iterates of $c\in \Omega_0\cup \Omega_1$
	never enter $I_2$.
	Moreover, we can replace  $I_1$ by an arbitrarily small puzzle piece
	(while keeping $I_2$ as before); the resulting union of intervals will remain nice
	because $c\in \Omega_1$ does not enter $I_2$.
	As usual, if $\tilde f$ is topologically conjugate to $f$, we mark the
	corresponding objects for $\tilde f$ with a tilde.

	\subsection{The immediate basin of attraction of a periodic attractor}
\label{sec:attracting}
Let $\mathcal C$ be the class of maps defined in Section \ref{subsec:class C}.
	\label{subsec:qc rigidity - basin of attraction}
	
			\begin{prop}\label{prop:attracting basin}
			Assume that $f\colon  [0,1]\to [0,1]$ is of class $\mathcal{C}$, $f(\{0,1\})\subset \{0,1\}$
			and that every $x\in (0,1)$ is attracted to a fixed point $p$ of $f$.
			Let $\tilde f\colon  [0,1]\to [0,1]$ be another mapping with the same properties and
			which is topologically conjugate to $f$. 
			Moreover, we require that the order of the maps at the fixed points are the same.
			Then $f$ and $\tilde f$ are quasisymmetrically conjugate.
		\end{prop}

Here we use the convention that we say that the \emph{order} of $f$ at $p$
is equal to zero, if $Df(p)\ne 0$ and otherwise it is the order of the critical point.
			
Although we do not need it in this generality, 
we state this result without requiring that the conjugacy maps
critical points to critical points of the same order (except if the
critical point
is a fixed point). For example, it is allowed to
mapping a critical point of odd order to a regular point,
unless this point is the attracting fixed point of $f$.
Note that if the fixed point is a critical point of $f$ then the conjugacy
can only be quasisymmetric if the orders of the two maps are the same
at these fixed points. Indeed, the rates of convergence near these fixed points
are of the form $\lambda^{d^n}$ and  $\tilde \lambda^{\tilde d^n}$ where
$d,\tilde d$ are the order of these critical points at $p$ resp. $\tilde p$
and this is only compatible with a H\"older conjugacy if $d=\tilde d$.
			
The proof of the Proposition~\ref{prop:attracting basin} also gives the following statement:

		\begin{prop}\label{prop:rigidity attracting}
			Assume that $f,\tilde f$ is of class $\mathcal{C}$ 
			which are topologically conjugate as in Theorem~\ref{thm:main},
			and with one or more periodic attractors.
			Then there exists a quasisymmetric homeomorphism $h$
			such that the following holds.
			\begin{enumerate}
				\item If $p,\tilde p$ are corresponding attracting periodic points,
					then $h\colon  B_0(p)\to B_0(\tilde p)$ and $h$ is a conjugacy between
					$f$ and $\tilde f$ restricted to these sets.
					Here $B_0(p)$ is the immediate basin of $p$ (i.e. the component of the basin of $p$
					which contains $p$).
				\item For each critical point $c$ of $f$ and each $n$ for 
					which  $f^n(c)\in B_0(p)$, one has
					$h(f^n(c))=\tilde f^n(\tilde c)$,
					where $\tilde c$ is the corresponding critical point of $f$.
			\end{enumerate}
		\end{prop}

	\medskip
	\noindent\textit{Proof of Proposition~\ref{prop:attracting basin}}.
	Although the proof of this is more or less standard, we will
			give it to be complete.
			Label corresponding folding points of $f$ and $\tilde f$ by
			$c_1,\dots,c_d$ resp. $\tilde c_1,\dots,\tilde c_d$.
			First assume that $Df(p)\ne 0$. 
			Take an interval $U$ around $p$
			so that by Proposition~\ref{prop:local dynamics}
			$f$ is quasiconformally linearizable 
			on a complex neighbourhood of $p$
			and so that $\tilde f$ has the same property on $h(U)$.
			Take a fundamental  annulus $A\subset U$ 
			which contains an iterate of each critical point.
			Take a corresponding fundamental annulus $\tilde A$
			for $\tilde f$, and take a smooth mapping $H\colon  \clos(A) \to \clos(\tilde A)$
			so that $H\circ f=\tilde f\circ H$ on the outer boundary of $A$
			and so that it maps the $n$-th iterate of $c_j$ to the $n$-th iterate of $\tilde c_j$.
			(We don't care about critical points of odd order.)
			
			Next extend $h\colon  \cup_{n\ge 0} f^n(A)\to \tilde f^n(\tilde A)$ in
			the obvious way by $h=(\tilde f|_{\tilde A})^{n}\circ h \circ (f|_A)^{-n}$ (on $ f^n(A)$).
			This mapping is quasiconformal on a neighbourhood of $p$.
			
			Next extend the qc conjugacy $H$ to a qs conjugacy $H$ defined
			on $(0,1)$ as follows: Take $x\in (\epsilon,1-\epsilon)$,
			take $n$ minimal so that $f^n(x)\in A$, and define $(\tilde f)^{-n}\circ H\circ f^n(x)$.
			Here the choice of the inverse is uniquely dictated by continuity.
			It is easy to see that $H$ is quasisymmetric restricted to an interval 
			of the form $[\delta,1-\delta]$, since then $n$ is uniformly bounded.
			(Restricted to this set, $H$ is possibly only non-smooth points
			in the backward iterates of critical points. There the points are locally
			of the form $t\mapsto t^{\rho}$.)
			But then if $a\in \{0,1\}$ is a fixed point of
                        $f$ then take a 
complex neighbourhood $U_a$
			of $a$ on which $f$ is is linearizable, take a fundamental annulus $A_a$
			and take the corresponding annulus $\tilde A_a$ for $a$.
			Then extend the qs mapping $A_a\cap [0,1]\to \tilde A_a\cap [0,1]$ to a quasiconformal 
			$A_a\to \tilde A_a$, and extend $H$ to a neighbourhood of $a$
			by $(\tilde f)^{n}\circ H\circ (f|_{U_a})^{-n}(x)$.  
			
			Observe that if $a$ is hyperbolic repelling fixed point then 
			backward orbits converge to $a$ at an exponential rate, and if 
			it is parabolic repelling, then
			since $f\in\mathcal{C}$,
			backward orbits converge to $a$ at the rate
			$$|z_{-n}-p'|\asymp \frac{1}{n^{1/d}}$$
			and the 
			asymptotically conformal extension is of order $d+2$ in a neighbourhood 
			of the parabolic periodic point.
			So that 
			in a neighbourhood of $z_{-n}$, $H$ is $K$-quasiconformal with
			$$K<C\sum_{k=1}^{n}\frac{1}{k^{\frac{d+1}{d}}}
			<C\sum_{k=1}^{\infty}\frac{1}{k^{\frac{d+1}{d}}},$$ 
			which converges.
			Thus the resulting limit mapping is quasiconformal.
			If $a$ is not a fixed point,
			then either it is mapped to the other endpoint of $[0,1]$ which is a fixed point
			or both boundary points have period two, and then we argue similarly.
			
		If $Df(p)=0$, then by Proposition~\ref{prop:local dynamics} 
		we can locally conjugate $f$ to $z^d$ by a qc mapping and argue as before.
		\qed

\subsection{The infinitely renormalizable, analytic case}
\label{subsec:qc rigidity - infinitely renormalizable}
Let us first show how to to prove rigidity at critical points at which
$f$ is infinitely renormalizable when $f$ is real-analytic. In this case only
minor modifications of \cite{KSS-rigidity} are needed.
In the infinitely renormalizable case
we have to use an approach which is more `real': we
cannot simply pull back through
the enhanced
nest associated to non-renormalizable complex box mappings. Fortunately, we
can combine a result of \cite{KSS-rigidity} and the touching box mappings.

\begin{thm}\label{thm:infrenorm}
 Let $g\colon  [0,1]\to [0,1]$ and $\tilde g \colon  [0,1]\to [0,1]$ be two real-analytic maps
which are topologically conjugate. Assume that $\omega(c)=\omega(c')$ for each 
two critical points  $c,c'$ of $f$ and that  $f$ is infinitely renormalizable at each
of its critical points.
Then there exists a qs conjugacy between
$g\colon  [0,1]\to [0,1]$ and $\tilde g \colon  [0,1]\to [0,1]$.
\end{thm}

This theorem is a generalization of Theorem 7.1 in \cite{KSS-density},
where it was assumed that $f,\tilde f$ are polynomials
with only have real critical points which
are all of even order,
Here we will show
how to modify this proof and use the touching box mappings to prove
Theorem~\ref{thm:infrenorm}.

First note that by Theorem~\ref{thm:box mapping persistent}, 
if $c$ is a critical point at which $g$ is infinitely renormalizable, then
there exists an interval $I\ni c$ so that $g^s(I)\subset I$
(with the interiors of $I,\dots,g^{s-1}(I)$ disjoint and
with $g^s(\partial I)\subset I$) so that
$g^s\colon  I\to I$ extends to a polynomial-like mapping
$G\colon\U\to \V$  and so that $\mod(\V \setminus \U)\ge \mu$. 
Here $\mu>0$ is \emph{beau} (eventually universal).
We can choose $s$ and $I$,
so that for the corresponding interval $\tilde I$ the map
$\tilde g^s\colon  \tilde I\to \tilde I$
also extends to a polynomial-like mapping
$G\colon \widetilde{\U}\to \widetilde{\V}$.
		
Hence, by Douady and Hubbard's Straightening Lemma,
$g^s\colon \U\to \V$ is 
$K$-qc conjugate to a polynomial
$G\colon  \C\to \C$, where $K$ only depends on $\mu$.
In particular, $g^s\colon  I\to I$ is qs-conjugate to a polynomial
$G\colon  [-1,1]\to [-1,1]$ with $G(\{-1,1\})\subset \{-1,1\}$.  
If $\omega(c)$ contains $b$ critical points, then
this polynomial $G$ is a composition of
$b$ real unicritical polynomials $q_i$,
for which  $q_i(\{-1,1\})\subset \{-1,1\}$, i.e.
$G=q_{b-1}\circ q_{b-2} \circ \dots \circ q_0$.
The mapping $G$ may have critical points off the real line, 
but $G$ is a composition of maps with only real critical values, and in particular
any non-real critical point $\hat c$ is eventually mapped
into the real line and
$\omega(c)=\omega(c')$ for all real critical points $c,c'$ of $G$.

\begin{lem}\label{lem:infreduc}
If $G\colon  [-1,1]\to [-1,1]$ and  $\tilde G\colon  [-1,1]\to [-1,1]$
are qs conjugate,
then $g\colon  [0,1]\to [0,1]$ and $\tilde g \colon  [0,1]\to [0,1]$
are qs conjugate.
\end{lem}		
\begin{pf}
Since $g^s\colon  I\to I$ is qs conjugate to $G\colon  [-1,1]\to [-1,1]$
(and similarly, for $\tilde g^s$ and $\tilde G$), 
if $G\colon  [-1,1]\to [-1,1]$ and $\tilde G\colon  [-1,1]\to [-1,1]$
are qs conjugate, then 
$g^s\colon  I\to I$ and $\tilde g^s\colon  \tilde I\to \tilde I$  are
qs conjugate.
Let $c_0$ be the critical point of $g$ in $I$.
Set $I^s=I$ and 
let $I^i=\comp_{g^i(c_0)}g^{s-i}(I)$ for $s=0,1,\dots s-1$.
Using the conjugacy equation
(and that $g$ and $\tilde g$ have critical points of the same order),
we also get that 
$g^s\colon  I^i\to I^i$ and $\tilde g^s\colon  \tilde{I}^i\to \tilde{I}^i$  
are qs conjugate for each $i=0,1,\dots,s-1$. 
Thus we have obtained a qs conjugacy 
$h\colon  \cup_{i=0,\dots,s-1} I^i \to  \cup_{i=0,\dots,s-1} \tilde{I}^i$
between $g$ and $\tilde g$. 
Let us now show how to extend $h$ to a qs conjugacy defined on $[0,1]$.  
			
To do this, we construct a touching complex box mapping.
Consider a finite, forward invariant set $Z_0'$ containing $\{0,1\}$
and $\partial I^i$, $i=1,\dots,s$.
Let $\delta>0$ be the constant from Proposition \ref{prop:petals}.
Let $N>0$ be so large that each point $z\in Z_0'$ is approximated
by a different point $z'\in g^{-N}(Z_0')\setminus (\cup_{i=0}^s I^i)$
with $|z-z'|<\delta$.
Such an $N$ exists by Lemma \ref{lem:startingpartition}.
			
Let $Z'=g^{-N}(Z_0')$.
Let $\P'=[0,1]\setminus Z'$.
and $\P=[0,1]\setminus g^{-1}(Z')$.
Next, consider a Poincar\'e lens neighbourhood 
(a ``necklace'' neighbourhood)
of $\P'$,
where for each component $J'$ of
$\P'$ the component of the neighbourhood is the 
Poincar\'e lens domain
$D_{\pi-\theta}(J')$ where $\theta$ is chosen close to $0$.
			
Next consider a component $J$ of $\P$ which is not contained
in $\cup g^i(I)$. Since all critical points of $g$ are contained in $\cup g^i(I)$,
$g$ is a diffeomorphism from $J$ to some component $T'$ of $\P'$.
It follows from Proposition~\ref{prop:Almost Schwarz Inclusion}
that the component containing $J$ of $g^{-1}(D_{\pi-\theta}(T'))$  is contained 
in $D_{\pi-\theta'}(J)$ where $\theta'$ is slightly larger than $\theta$. 
In fact, $\theta'/\theta$ can be chosen as close to one as we want
by taking $N$ sufficiently large,
and it can be arranged so that there exists $c<1$ (not depending on
$Z'$)
so that
$|J|/|J'|\le c$ where $J'=J'(J)$ is the component of $\P'$ containing $J$.
It follows that $g^{-s}(D_\theta(T'))$ is compactly contained in $D_{\pi-\theta}(J')$
whenever $J$ is compactly contained in $J'$.
If $J$ has a boundary point in common with $J'$, then
$D_{\pi-\theta'}(J)$ is contained in $D_{\pi-\theta}(J')$,
in view of Proposition~\ref{prop:petals}.
			
Thus we have constructed a touching complex box mapping
$g_T\colon  \UU_T\to \VV_T$.
Here $\UU_T$ is equal to the union of sets of the form
$\comp_{J}g_T^{-1}(D_{\pi-\theta}(T'))$ where
$T'$ is a component of $[0,1]\setminus g^{-1}(Z)$
not contained in $\cup g^i(I)$ and $T'=g(J)$.
Similarly, $\VV_T$ is the union of sets of the form
$D_{\pi-\theta}(T')$ where $T'$
is a component of $[0,1]\setminus Z'$.
Note that $\overline \VV_T$ contains $[0,1]$
and that $\UU_T$ contains $[0,1]\setminus \cup g^i(I)$.
			
Now define a qc mapping $H\colon  \C\to \C$ which is real-symmetric and so that
\begin{itemize}
\item $H$ is an extension of
$h\colon  \bigcup_{i=0,\dots,s-1} g^i(I) \to  \bigcup_{i=0,\dots,s-1} \tilde g^i(\tilde I)$
\item $\VV_T$ is mapped onto $\widetilde{\VV}_T$;
\item $\UU_T$ is mapped onto $\widetilde{\UU}_T$;
\item $\tilde g\circ H=H\circ g$ on $\partial \UU_T$.
\end{itemize}
That this can be achieved follows from Proposition~\ref{prop:petals}
and Theorem~\ref{thm:touching box map}.
Now apply successive pullbacks $H_n$ of $H$. Now take $H_0=H$ and 
define $H_n|_{\UU_T}$ inductively  by $\tilde g_T\circ H_{n+1}=H_n\circ g_T$
and define $H_n$ outside $\UU_T$ to be equal to $H$.
Since $g,\tilde g$ are univalent on each component of $\UU_T,\widetilde{\UU}_T$,
the dilatation of $H_n$ is the same as that of $H$.
Note also that $H_{n+1}$ agrees with $H_n$ outside the
set $\U_n=\{z\in \UU_T; g(z),\dots,g^n(z)\in \UU_T\}$.
It follows that $H_n$ converges to a a qc mapping $H_\infty$
and that $H_\infty$ is a conjugacy between $g_T$ and $\tilde g_T$.
\end{pf}
\medskip

It follows from Lemma~\ref{lem:infreduc}
that it is enough to show that the polynomials $F$ and $\tilde F$
are qc conjugate. Even though not all critical points of $F$ and $\tilde F$
are real, we have as in  Theorem 7.1 of \cite{KSS-rigidity}:

\begin{thm}\label{thm:rigidity of pc inf renorm}
There exists a qs homeomorphism $h\colon  [-1,1]\to [-1,1]$
which maps critical points of $G$ to critical points of $\tilde G$
so that $h(G^i(c))=\tilde G^i(h(c))$.
\end{thm}
\begin{pf}
The proof of this result is based on \cite{KSS-rigidity}.
Since we will give a different proof using the qc\textbackslash bg
partition, we will be relatively brief.
The proof of this Theorem goes essentially as in 
Sections 7.3-7.5  in Section 7 of \cite{KSS-rigidity}.
Indeed, the maps $G,\tilde G$ are  in the same class ${\mathcal T}_b$ as the maps
used in Sections 7.3-7.5  in Section 7 of \cite{KSS-rigidity},
with the only difference being that in our case the maps can have inflection points.
However, from Theorem~\ref{thm:box mapping persistent} 
we have complex bounds 
for maps with inflection points,
and otherwise the modifications needed in the proof are minor
(for example, Fact 7.2 (part 2), in Lemma 7.2 if $F^2\colon  M_1\to M_0$
is monotone, then it does not need to be a diffeomorphism).		
\end{pf}

\begin{cor}\label{cor:plqcconj}
$G$   and $\tilde G$ are qc conjugate.
\end{cor}
		
\begin{pf}
Take a qc mapping $H_0\colon  \C\to \C$ which agrees
with the qs homeomorphism $h$ from
Theorem~\ref{thm:rigidity of pc inf renorm}.
Note that all critical values of $G$  are in the real line,
and have the same $\omega$-limit (and similarly for $\tilde G$).
Hence we can define inductively 
$H_{n+1}$ by $\tilde G\circ H_{n+1}=H_n\circ G$.
Since $h(G^i(c))=\tilde G^i(h(c))$ for all $i\ge 0$, 
this is possible, and $H_n$ has the same dilatation as $H_0$.
It follows that $H_n$ has a convergent limit $H_\infty$, and
that $\tilde G\circ H_\infty=H_\infty\circ G$.
\end{pf}

\medskip
\noindent\textit{Proof of Theorem~\ref{thm:infrenorm}.}
The theorem follows from the Corollary~\ref{cor:plqcconj} and Lemma~\ref{lem:infreduc}.
\qed

\subsection{The infinitely renormalizable, $C^3$ case}
\label{subsec:C3 inf renorm rigidity}
In the smooth case
we cannot make use of the Douady-Hubbard Straightening Theorem.
We will give a proof of quasisymmetric rigidity,
which is new (even in the analytic case), and
which
uses similar ideas as proof of
Theorem~\ref{thm:central cascades}.

Suppose that $c_0$ is a critical point at which $f$ is infinitely renormalizable.
We may assume that $c_0$ is even.
Let $J=J_0\supset J_1\supset J_2\supset J_3\supset\dots$ be successive periodic
intervals about $c_0$
with periods $1<s_1<s_2<s_3<\dots$.
Let $J_n^j=\mathrm{Comp}_{f^j(c_0)}f^{-(s_n-j)}(J_n),$
where $0\leq j<s_n$. 
		We denote the boundary of $J_n$ by $\{\beta,\tau(\beta)\}$ where $\beta$ is a 
		repelling periodic point (for $f$) that is fixed by $R_{J_n}=f^{s_n}|_{J_n}$. 

\begin{prop}[\cite{KSS-rigidity, CvST}]\label{prop:renormbounds}.
		
There exists a positive integer $N=N(f)$ such that if 
$n\geq N$ we have the following.
\begin{enumerate}
\item There exist $0=j_0<j_1<j_2<\dots<j_{b-1}<s_n$
such that each $J_n^{j_i}$ contains a
critical point equivalent to $c_0$
and for any other $i$, $0<i<s_n$, $J_n^i$ is disjoint from $\crit(f)$.
Moreover, $f^{j_{i+1}-j_i}(J_n^{j_i})$ contains the critical point in $J_n^{j_{i+1}}$.
\item Let $\mathcal{J}_n=\cup_{j=0}^{b-1} J_{n}^{j_i}$.
Then the distortion of the 
first landing mapping to $\mathcal{J}_n$  under $f$
restricted to $\cup_{j=0}^{s_n-1}J_n^i$
is bounded from above by a constant depending only on the number of critical
points
of $f$
(in the equivalence class of $c_0$) and their orders.
\item For any $0\leq i\leq s_n$ and $0\leq i'\leq s_{n+1}$,
if $J_{n+1}^{i'}\subset J_n^i$, then
$$(1+\tau)J_{n+1}^{i'}\subset J_n^i,$$
where $\tau>0$ is a constant depending only on the vector of critical points.
\item The derivative of $f^{s_n}\colon J_n\rightarrow J_n$
is bounded from above by a constant that depends
only on the number of critical points of $f$ and
their orders.
\item The multipliers of the periodic points of $f^{s_n}\colon J_n\rightarrow J_n$
are bounded from below by a constant $\rho>1$.
\item There exists a constant $\delta$ such that if $c$ is a critical point 
of $f^{s_n}\colon J_n\rightarrow J_n$, then $|f^{s_n}(c)-c|>\delta|J_n|$.
\item Let $\alpha$ be the fixed point of $f^{s_n}|_{J_n}$ closest to $c_0$, 
$\beta$ be the fixed point of
$f^{s_n}\colon J_n\rightarrow J_n$ in $\partial J_n$ and
$x_0\in f^{-s_n}(\alpha)\cap J_n$ be the point closest to $\beta$.
For $i>0$, let $x_{i+1}$ be the preimage of $x_i$
under $f^{s_n}|_{J_n}$ closest to $\beta$.
There exist constants $C,c>0$ and $0<\lambda<1$ such that
$ c\lambda^n<|x_n-\beta|<  C\lambda^n$ (see \cite[Fact 7.3]{KSS-rigidity}).
\end{enumerate}
\end{prop}

Suppose that $f$ is an infinitely renormalizable mapping as in Theorem~\ref{thm:main}.
		Suppose that $J$ is a periodic interval for $f$ of sufficiently high period
		so that by Theorem~\ref{thm:box mapping persistent}
                there exists a $K$-qr polynomial-like mapping
$F\colon \U\rightarrow \V$ that extends $R_J$ such that
$\mod(\V\setminus \U)$ is bounded away from
		$0$, $\U$ has $\rho$-bounded geometry and $\diam(\U)<C|J|$. 
		Let us first show that we can construct a
		Yoccoz puzzle $\mathcal{Y}$ for $F$, and show that
		the puzzle pieces have good geometric properties.
		The proof of the following is the same as for (holomorphic) polynomial-like maps with
		non-escaping critical points.
		\begin{lem}
		Let $F$ be as above. Let $J_F$ be the set of non-escaping points of $F$. 
		Then $J_F$ is connected.
		\end{lem}
		
Let us show that we can construct external rays for the Julia set of a
$K$-qr polynomial-like mapping $F\colon \U\rightarrow \V$
with a connected Julia set.
We use a method of \cite{Levin-Przytycki}.
Let $K_F$ denote the filled Julia set set of $F$.
Let $B_F\colon \V\rightarrow \bm W_+$ be an isomorphism of 
$\U\setminus K_F$ onto the round annulus $\bm W_+=\{z:1<|z|<R\}$
where $\log R$ is the modulus of $\V\setminus K_F$,
such that $|B_F(z)|\rightarrow 1$ as $z\rightarrow J_F$, so that
$B_F$ is a proper map. 
Let $\bm W'_+=B_F(\U)$.
Let $h_+\colon \bm W'_+\rightarrow \bm W_+$
be the mapping $h_+(x)=B_F\circ F\circ B_F^{-1}(x)$
Let $\sigma\colon  z\mapsto 1/\bar{z}$
be the reflection with respect to the unit circle. 
Let $\bm W_-=\sigma(\bm W_+)$ and $\bm W'_-=\sigma(\bm W'_+)$.
Define $\bm W=\bm W_+\cup S^1\cup \bm W_-$ and
$\bm W'=\bm W'_+\cup S^1\cup \bm W'_{-}$.
By the reflection principle for
qc mappings $h_+$ extends to a $K$-qc mapping
$h\colon \bm W'\rightarrow \bm W$.
We have that $h$ is strongly expanding in a neighbourhood of
$S^1$, and consequently
that $h^{-1}$ is strongly contracting for the Poincar\'e
metric in a neighbourhood of $S^1$.
Let us foliate $\bm V\setminus K_F$ as follows.
Choose $\eta>1$ and let $\phi$ be an
arbitrary smooth conjugacy between $h$ and 
$\psi\colon z\mapsto z^d$ on a neighbourhood
of $\bm W\setminus \bm W'$ to
$\bm \Lambda=\{z:\eta<|z|<\eta^d\}\cup\{z:1/\eta^d<|z|<1/\eta\}$.
Foliate $\bm \Lambda$ by radii,
and extend it to a foliation $\tau'$ of a neighbourhood of $S_1$ by $(\phi\circ h^n)^{1/d^n}$.
Define $\tau^*=\phi^{-1}(\tau')$ and $\tau=B_F^{-1}(\tau^*)$.
Then $\tau$ is a smooth foliation defined on  
$\bm V\setminus K_F$ that is invariant under $F$ with the property that
each leaf of $\tau$ passes through 
$\partial \bm V$ and $\partial \bm U$ and intersects these sets transversally.
Note the $\tau$ is smooth since all the critical points of $F$ are contained in $J_F$.
We will call any leaf of $\tau$ an external ray of $J_F$.
Since $h$ is strongly expanding,
each leaf of $\tau^+$ intersects $S^1$ in a unique point
		$e^{2\pi i t}$.
Thus we extend $\tau^*$ to a foliation of $\bm W$.
From the construction of $\tau^*$ via $\tau'$ we have that 
for each $t\in[0,1)$ only one leaf of $\tau^*$ passes through $e^{2\pi i t}$.
For each external ray $\ell$ we let $\ell^*$ denote the leaf of $\tau^*$ that
contains $B_F(\ell)$.
		\begin{lem}
		Keeping the same notation as above. Suppose that $p$ is a real repelling periodic point of $F\colon \U\rightarrow \V$.
		Then at least one ray lands at $p$.
		\end{lem}
		 \begin{pf}
		 Suppose that the period of $p$ is $s$, then by considering the mapping $F^s\colon p \mapsto p$, we can suppose that
		 $p$ is a fixed point of $F$, so that is what we will do. Let $\lambda=F'(p)$.
		 Let $N$ denote a complex neighbourhood of $p$
		 on which $F$ is quasiconformally conjugate to $x\mapsto \lambda x$. Let $g=F^{-1}|_N$.
		 Then there exists a constant $C$ such that for $x\in N$
		 $$\frac{1}{C\lambda^n}<|(g^{n})'(x)|<\frac{C}{\lambda^n}.$$
		 For any ray $R$ that is fixed by $F$, let $R_t$ denote the segment on the ray $R$
		 lying between $t$ and $F(t)$, and including $t$. Then $R_t$ is a fundamental domain
		 for the dynamics on $R$.
		 It follows from the facts that $h$ is strongly expanding near $S^1$ and
		that $B_F$ is proper, that the Euclidean length of $R_t$ tends to 0 as $t\rightarrow 0.$
		 Let $S'\subset N$ be a complex neighbourhood of $p$ and let $S=g^{-1}(S')$.
		 Let $$\varepsilon=\inf_{z\in\partial S, z'\in\partial S'}|z-z'|.$$
		 There exists $t_{\varepsilon}$, such that for all $t\in(0,t_{\varepsilon})$,
		 the Euclidean length of $R_t$ is less than $\varepsilon$.
		 Take $t$ so small that $R_t$ lies in $S'$.
		 Then for all $n$, $R_{g^n(t)}\subset S'$ and
		 $$R_{g^n(t)}\subset B\Big(p,\frac{C\diam(S')}{\lambda^n}\Big),$$
		 which implies that that the ray $R$ lands.
		 \end{pf}

Let $\alpha$ be the orientation reversing fixed point in $J_n$ 
that is furthest from $c\in J_n$.
From the above and by slightly modifying the 
arguments of Lemmas 5.1 and 5.2
of \cite{KSS-rigidity}, we have two rays,
symmetric about the real line, landing at $\alpha$,
and pulling them back by one iterate of $F$, we obtain two rays landing at
$\tau(\alpha)$.	
Let $M_0=(\alpha,\tau(\alpha))$.
Inductively define $M_{n+1}$ to be the pullback of $M_n$
that contains $\alpha$ in its boundary.
Assuming that the first renormalization of $F$ is 
not of intersection type, we have $M_1\neq M_0$.
If the mapping $F^2\colon M_1\rightarrow M_0$ is monotone,
then it has a unique fixed point, which we denote by $\rho$.
Just as above, we have two rays landing at each of $\rho$ and $\tau(\rho)$.
If the mapping $F^2\colon M_1\rightarrow M_0$ is not monotone, then
the mapping $F^2\colon M_2\rightarrow M_1$ is. In this case we can 
apply the argument on pages 796-797 of \cite{KSS-rigidity}
to reduce to the previous case. 
So let us just consider the case where 
$F^2\colon M_1\rightarrow M_0$ is monotone.
		
Let $\bm Y'$ be the puzzle piece bounded by the rays through
$\alpha$ and $\tau(\alpha)$
and $\partial \bm V$ containing $c_0$, and let
$\bm Y$ be the puzzle piece bounded by the rays through
$\rho$ and $\tau(\rho)$
and the curve $\partial U$ containing $c_0$.
Let $\bm{\mathcal{Y}}_0$ be the union of puzzle pieces 
bounded by 
$\partial \V$, $\partial \U$ and the rays landing at
$\alpha,\tau(\alpha), \rho$ and $\tau(\rho)$
Just as for polynomials, the pullbacks of these
puzzle pieces are either nested or disjoint.	
				
\begin{lem}\label{lem:top level geometry}
There exists $\tau>0$ and $\theta\in(0,\pi/2)$ such that 
$\bm Y$ has $\theta$-bounded shape and $\mod(\bm Y'\setminus \bm Y)\geq \tau$.
\end{lem}
\begin{pf}
The relative length in $I$ of each component of
$I\setminus \{\alpha,\tau(\alpha),\rho,\tau(\rho)\}$
is uniformly bounded away from zero.
The multipliers of $\alpha$ and $\rho$
are uniformly bounded from above and away from one,
$\bm Y\Subset \bm Y'$,
and there exists a constant $C$ so that
$F$ is $(1+C\max_{0\leq j<s_n}(|J_n^j|)$-quasiregular
 on $\U\cap\mathbb{R}$, so
$\mod(\bm Y'\setminus \bm Y)$ is bounded away from 0.
Let $T\owns \rho$ be the maximal open interval such that 
$F^4|_T$ is monotone.
Let $L$ be the component of 
$T\setminus \{\rho\}$
that is contained in $(\rho,c_0)$
Let $R$ be the other one.
The relative lengths of $L$ and $R$ are bounded away from zero.
For some constant $\theta>0,$
$\partial \bm Y\cap(\mathbb C\setminus F^{-4}(\V))$
is contained in $D_{\theta}(L)\cap D_{\theta}(R)$,
and is disjoint from $D_{\pi-\theta}((\rho,-\rho)),$
so, the external rays landing at $\rho$ are contained in 
$D_{\theta}(L)\cap D_{\theta}(R)$. Similarly around the preimage of $\rho$.
This is enough to show that $Y$ has $\theta$-bounded shape.
\end{pf}

Let $f,\tilde{f}$ be maps as in Theorem~\ref{thm:main}. 
Suppose that $f$ is infinitely renomalizable at $c_0$, and let
$F\colon \U\rightarrow \V$ and 
$\tilde F\colon \widetilde{\U}\rightarrow\widetilde{\V}$ 
		be the $K$-qr polynomial-like mappings that extend the return maps to 
		corresponding periodic intervals $I\owns c_0$, $\tilde{I}\owns \tilde{c}_0$ 
		of sufficiently high period given by
		Theorem~\ref{thm:box mapping persistent}.
		Let $\bm{\mathcal{Y}}_0$ and $\widetilde{\bm{\mathcal{Y}}}_0$ be the puzzles constructed as above
		for $F$ and $\tilde{F}$, and
		let $H_0=B_{\tilde{F}}^{-1}\circ B_F.$ Then
$H_0\colon \partial \bm{\mathcal{Y}}_0\rightarrow \partial \widetilde{\bm{\mathcal{Y}}}_0$
		 conjugates $F$ and $\tilde F$ on $\partial\bm{\mathcal{Y}}_0$, and 
		since $F,\tilde{F}$ are
		asymptotically conformal on the real slices of their domains
		there exists $K'>1$ so that $H_0$ is a $K'$-qc in a
                neighbourhood of
$\bm{\mathcal{Y}}_0$.
		Mimicking the case of polynomials, we refer to $H_0$
		as a \emph{boundary marking}.

\begin{lem}
[Puzzle geometry control for $K$-qr polynomial-like mappings,
c.f. Lemma 7.4 \cite{KSS-rigidity}]
\label{lem:puzzle geometry}
Let $f,\tilde{f}$ be maps as in Theorem~\ref{thm:main}.
Suppose that $f$ is infinitely renomalizable at $c_0\in\Crit(f)$,
and let
$F\colon \U\rightarrow \V$ 
be a $K$-qr polynomial-like mapping that extends the return mapping to 
a periodic interval $I\owns c_0$ of sufficiently high period $s$ given by
Theorem~\ref{thm:box mapping persistent}.
Suppose that $F$ is not immediately renormalizable (equivalently, the
next renormalization of $F$ is not of intersection type).
		Let $J_0$ be the maximal properly periodic interval for $F|_I$ that contains $c_0$.
Then there exist two puzzle pieces $\bm P'\supset \bm P$ for $F$ that contain $J_0$
		and are contained in the Yoccoz puzzle $\bm{\mathcal{Y}}$ for $F$ such that 
\begin{enumerate}
\item the first return time of $c_0$ to $\bm P'$ 
is equal to the first renormalization period $s$;
\item $\bm P'\setminus \bm P$ is disjoint from $\omega(c_0)$
\item $\mod(\bm P'\setminus \bm P)\geq \eta$;
\item $\bm P$ is an $\eta$-excellent puzzle piece, see page~\pageref{page:defgood},
\end{enumerate}
where $\eta>0$ is a constant independent of $F$.
Moreover, if we replace the mapping $F$ and the puzzle pieces
$\bm P'$, $\bm P$ by the corresponding
objects for $\tilde f$, the statements remain true.
\end{lem}
\begin{pf}
We will explain how to modify the proof of Lemma 7.4 of \cite{KSS-rigidity}.		
First, by the Complex Bounds of \cite{CvST}, 
and Lemma~\ref{lem:sum of lengths},
the existence of domains $\V$ and $\V'$
as in the start of the proof of \cite{KSS-rigidity}, Lemma 7.4 imply
the result.
Since we are using $\V$ to denote the range of a 
qr polynomial-like mapping, we will deviate from the notation
in \cite{KSS-rigidity} and denote $\V$ and $\V'$ by $\bm Q$ and $\bm Q'$ respectively. 
Moreover, Lemma~\ref{lem:top level geometry}
gives us puzzle pieces $\bm Y,\bm Y'$ as in Fact 7.4 of \cite{KSS-density}.
Proceeding as in \cite{KSS-rigidity}, we introduce an extended mapping as follows.
Let $\{c_0,\dots,c_{b-1}\}$ be the set of critical points in $\omega(c_0)$
ordered so that $\U\owns c_0, f^{s_1}(\U)\owns c_1,
f^{s_2}(\U)\owns c_2, \dots, f^{s_{b-1}}(\U)\owns c_{b-1}$
with $0=s_0<s_1<s_2<\dots s_{b-1}$.
For each $c_i$, let $\U_i=f^{s_i}(\U)$.
For each $i$,$0<i<b-1$, if $z\in \U_i$
define $\bm{F}(z,i)=(f^{s_{i+1}-s_i}(z),i+1)$.
We say that an interval $K\times\{i\}\subset [0,1]\times\{i\}$
or $S^1\times\{i\}$ is an $\bm{F}$-pullback of
$M_0\times 0$ of \emph{depth} $k$ if
it is a component of $\bm{F}^{-k}(M_0\times\{0\})\cap(\mathbb{R}\times\{i\})$
so that $i+k=0\;\mod(b)$.
Let $m$ be the maximal positive integer so that 
$M_0\times\{0\}$ has a unicritical
$\bm{F}$-pullback, $K\times\{i\}$. 
For each $(x,i)\in\mathrm{PC}(\bm{F})\cap (K\times\{i\})$,
we have $\bm{F}^m(x,i)\in(M_0\setminus(M_1\cup(-M_1)))\times\{0\}$.
This follows exactly as on page 795 \cite{KSS-rigidity}.
Let $\bm W=\comp_{c_i}(\bm{F}^{-m}(\bm Y\times\{0\}))$
and $\bm W'=\comp_{c_i}(F^{-m}(\bm Y'\times\{0\}))$.
Then $\bm W'\setminus \bm W$ is disjoint from the post-critical set of $\bm{F}$,
and since there exists a constant $C>0$ such that the mapping
$$\bm{F}^m\colon (\bm W',\bm W)\rightarrow (\bm Y'\times\{0\},\bm Y\times\{0\})$$
is a degree $d_{c_i}$ unicritical,
$(1+C\max_{0\leq i\leq m-1}(\diam(\bm G_i))^{1/2})$-qr branched
covering,
where $\{\bm G_i\}_{i=0}^m$ is the chain with $\bm G_m=\bm
Y'\times\{0\},$
and $\bm G_0=\bm W'$, we have
$$\mod(\bm W'\setminus \bm W)\geq \frac{\mod(\bm Y'\setminus \bm Y)}
{d (1+C\max_{0\leq i\leq m-1}(\diam(\bm G_i))^{1/2})},$$
which is bounded from below.
To see that $\bm W$ has bounded geometry follows
from Lemmas \ref{lem:sum of lengths}, \ref{lem:zd} and \ref{lem:zd lower bound}.
Now we define $\bm Q$, respectively $\bm Q '$, to be the 
topological disk that is the component of the domain of the
first landing mapping to $\bm W$, respectively $\bm W'$, under $F$ that contains $c_0$.
		\end{pf}

The following gives qc\textbackslash bg partition for a
qr polynomial-like mapping. It is analogous to 
Theorem~\ref{thm:central cascades}.
for a qr polynomial-like renormalization of $f$ at an even
critical point $c_0$. For an image showing part of the construction
see Figure~\ref{fig:ren2}.

\begin{figure}
\input{renormsmall.pdf_t}
\caption{The qc\textbackslash bg partition of a polynomial-like
  mapping.}
\label{fig:ren2}
\end{figure}

\begin{lem}\label{lem:pl decomp}
There exist $\delta_0>0,$ $K_1\geq 1$ and $\hat K\geq 1$,
such that for any $K,K_0\geq 1,\delta>0$ the following holds.
Suppose that $c_0$ is an even critical point at which 
$f$ is infinitely renormalizable. 
Suppose that
$F\colon\U\rightarrow\V$ is a $(1+\mu(\V))$-qr polynomial-like mapping that extends the 
first return mapping to a periodic interval of sufficiently high
period,
where $\partial\V$ is a $K_0$-quasidisk.
Assume that there exists
a $(K,\delta)$-qc\textbackslash bg mapping $H_0\colon\mathbb C\rightarrow\mathbb C$
such that $H_0(\V)=\widetilde{\V}$ and $H_0(\U)=\widetilde{\U}.$
Assume that the set $X_1^{H_0}$ (the set with bounded geometry)
has moduli bounds in $\V$ and $\V\cap\mathbb H$.
Then there exists a $(K',\delta')$-qc\textbackslash bg mapping
$H\colon\mathbb C\rightarrow\mathbb C,$ 
with $K'=\max\{\hat K,(1+C\mu(\V)),(1+C\mu(\widetilde{\V}))K\}$ and
$\delta'=\min\{\delta_0,(1+C\mu(\V))^{-1}\delta,
(1+C\mu(\widetilde{\V}))^{-1}\delta\},$
where $C>0$ is a universal constant,
a partition of $\V$ into combinatorially defined
disjoint sets
$X_0,X_1,X_2$ and $X_3,$ and a partition of 
$\widetilde{\V}$ into combinatorially defined disjoint sets
$\tilde X_0,\tilde X_1,\tilde X_2$ and $\tilde X_3$
such that the following holds:
\begin{enumerate}
\item $H(X_i)=\tilde X_i,$ for $i=0,1,2,3.$
\item $X_1$ is a union of combinatorially defined topological squares
  $\bm S$ such that each $\bm S$ has $\delta'$-bounded geometry and satisfies the following
  moduli bounds:
$$\mod(\U\setminus\bm S)>\delta',$$ and if $\bm S\subset\mathbb H^+,$
then 
$$\mod(\U\cap\mathbb H\setminus\bm S)>\delta'.$$


\item Each component of $X_3$ is a $K_1$-quasidisk and $H$ is $\hat
  K$-qc on $X_3.$
\item $H$ is $K'$-qc on $X_2$ moreover, $H$ is defined dynamically on
  $X_2$ and $X_2$ includes $\bm V\setminus F^{-2}(\U),$ and a
  combinatorally defined neighbourhood of the rays landing at
  $\beta,\alpha,\rho, \tau(\beta),\tau(\alpha),$ and $\tau(\rho)$.
\item The measure of $X_0$ is zero and $H$ satisfies (\ref{eqn:qc
    bound 0}) at each $x\in X_0$.
\end{enumerate}
\end{lem}
\begin{rem}
We will prove in Lemma~\ref{lem:infrenqd} that there exists
$K_0\geq 1$ such that
that $\partial\bm V$ is always a $K_0$-quasidisk.
\end{rem}
\begin{pf}
To simplify the notation, let $\II=\II^0=\V$
and $\II^1=\U$, so that
$F\colon \II^1\rightarrow \II^0$,
and as usual,
we define the puzzle pieces in the principle nest about
$c_0$ by $\II^i=\comp_{c_0}F^{-i}(\II^0)$.
In this case, the puzzle pieces shrink to the filled Julia set of $F$.

We will only consider the bounded geometry case. We can reduce the big
geometry case to the bounded geometry case exactly as we did at the
end of the proof of Theorem~\ref{thm:central cascades}.

Let us start by partitioning $\II$ into sets $X_i, i\in\{0,1,2,3\}.$
Our construction is similar to the one used to deal with 
long central cascades when there is a high return; however, we need to
partition $\II$ into new ``puzzle pieces''. 

Let $\beta$ denote the outermost orientation preserving fixed point of
$F$ and let $\tau(\beta)$ denote its image under the even symmetry
about $c_0$. As in the long saddle node case,
let $\Gamma_\beta\subset \Gamma_{\beta}'$ denote the regions intersecting
$\II^0\setminus[\beta,\tau(\beta)]$, which are bounded by
$\partial\II^4$ and the staircases $\gamma,\gamma'$, respectively,
that connect $\partial\II^4$ with $\beta$ and $\tau(\beta)$.

We will now construct regions about the rays landing at 
$\alpha,\tau(\alpha),\rho,\tau(\rho)$.
Order the components of
$(\II^i\setminus\II^{i+1})\setminus F^{-i}(I^0\setminus I^0)$
counterclockwise
by $S^j_i,$ $j=\{0,\dots, 2d^{i}\}$ where $S_i^0$ 
 lies in the upper half-plane
is bounded by the component of
$I_i\setminus I_{i+1}$ that is contain in the image of $I^1$.
For each $i\in\mathbb N\cup\{0\}$ let
$S_{i}^{j_\alpha}$ be the component of
$\cup_{j=0}^{2d^i}S_i^j$ that is contained in the upper
half-plane and
intersects one of the two rays landing at $\alpha$, then
$S_{i}^{2d^i+1-j_{\alpha}}$ is the component of 
$(\II^i\setminus\II^{i+1})\setminus F^{-i}(I^0\setminus I^0)$
that lies in the lower half-plane and 
intersects a ray landing at $\alpha$.
Let $S_{i,\alpha}=S_{i}^{j_\alpha}\cup S_{i}^{2d^i+1-j_{\alpha}}.$
We construct the ``shielding'' region about these topological squares
by setting
$$S'_{i,\alpha}=\Big(\bigcup_{k=-2}^2 S_{i}^{j_\alpha+k}\Big)\bigcup \Big(\bigcup_{k=-2}^s S_{i}^{2d^i+1-j_{\alpha}+k}\Big).$$
We let 
$$\Gamma_\alpha=\bigcup_{i=0}^\infty S_{i,\alpha}\quad\mbox{ and } \quad\Gamma_\alpha'=\bigcup_{i=0}^\infty S_{i,\alpha}'.$$ 
Carrying out identical constructions, we construct regions 
$\Gamma_a\subset\Gamma_a'$ about the rays landing at $a$ for
$a\in\{\tau(\alpha),\rho,\tau(\rho)\}$. We set
$$\Gamma'=\bigcup_{a\in\{\alpha,\beta,\rho,\tau(\alpha),\tau(\beta),\tau(\rho)\}}\Gamma'_a\quad\mbox{
  and }\quad \Gamma=\bigcup_{a\in\{\alpha,\beta,\rho,\tau(\alpha),\tau(\beta),\tau(\rho)\}}\Gamma_a.$$

We partition $\V$ as follows:

\begin{itemize}
\item $X_1=(\Gamma'\setminus\Gamma)\cup(\II^3\setminus(\II^4\cup\Gamma')$


\item $X_2=\Gamma\cup\II^0\setminus\II^3$

\item $X_3=\II^{4}\setminus\overline{\Gamma}$

\item $X_0=\II^{0}\setminus\overline{X_1\cup X_2\cup X_3}$.
\end{itemize}

We saw in the proof of Lemma~\ref{lem:top level geometry} that
there exists a real neighbourhood of $W$
of $\alpha$ such that $F|_W$ is a diffeomorphism and
$|W|\asymp\diam(\U).$ Thus there exists $\theta\in(0,\pi)$ such that
$D_{\theta}(W)$ contains all the components of 
$\II^2\setminus(\II^3\cup F^{-3}(V\setminus U)$ except those whose
boundaries
intersect the real line. Since the lengths of
$W^i=\comp_{\alpha}F^{-i}(W)$
shrink exponentially quickly,
we have that the dilation of $F^i$ is bounded, independently of
$i$, on
$\Gamma_{\alpha}\cap(\II^i\setminus\II^{i+1})$.
The proof is similar in the other components of 
$\Gamma\cap(\II^i\setminus\II^{i+1})$.

Suppose that $\bm W$ is a component of
$(\II^4\setminus\Gamma')\cap\mathbb H$.
We will use Lemma~\ref{lem:qdc} to show that $\bm W$ is a
$\hat K$-quasidisk for some universal $\hat K$. 
To be concrete we will carry out the proof in a specific case. The
proofs for the other components are identical.
We assume $\bm W$ is a topological rectangle 
with $\partial\bm
W=\overline{\eta_1\cup\eta_2\cup\eta_3\cup\eta_4}$ where
$\eta_1\subset \partial\II^4,$ $\eta_2\subset \partial\Gamma_{\alpha}',$
$\eta_3\subset \mathbb R$ and $\eta_4\subset \partial\Gamma_{\tau(\beta)}'$.

Since $\partial \V$ is a $K_0$-quasicircle,
we have that $\eta_1$ is a $k$-quasiarc.
We verify that both $\eta_2$ and $\eta_4$ are $k$-quasiarcs,
exactly as we showed that $\gamma$ is a quasiarc in the central
cascades, high return argument, see page~\pageref{page:gamma qa}. We have that 
$\eta_3$ is a 1-quasiarc since it is an interval in the
real line. Theorem~\ref{thm:box mapping persistent}
implies that there exists $\delta_1>0$ such that
$\dist(\eta_1,\eta_3)\geq \delta_1\diam(\II^4),$
and that for each $i\in\{1,2,3,4\}$, we have 
$\diam(\eta_i)>\delta_1\diam(\II^4)$.
A compactness argument shows that there exists $\delta_2>0$ such that
$\dist(\eta_2,\eta_4)\geq \delta_2\diam(\II^4).$

So it remains for us to show that the Ahlfors-Beurling Criterion
holds on adjacent boundary arcs.
The arguments
are similar 
near the intersections $\bar\eta_1\cap\bar\eta_2$ and
$\bar\eta_1\cap\bar\eta_4,$ and near
$\bar\eta_3\cap\bar\eta_2$ and
$\bar\eta_3\cap\bar\eta_4.$ So let us only consider
the cases when $z_1\in\eta_1$, $z_2\in\eta_2$ and when
$z_1\in\eta_3$ and $z_2\in\eta_2.$
First suppose that $z_1\in\eta_1$, $z_2\in\eta_2$
and $z_2$ is not contained in $F^{-4}(I^0\setminus I^2),$
then, there exists a constant $c>0$ such that
$|z_1-z_2|>c\cdot\diam(\II^4),$ and so $|z_1-z_2|\geq
c\cdot\diam(\gamma_{[z_1,z_2]}),$
so we can assume that $z_2\in F^{-4}(I^0\setminus I^2).$
Then we have that $\gamma_{[z_1,z_2]}$ consists of an arc in
$F^{-4}(I^0\setminus I^2)$ and an arc in $\alpha_1$. 
Since $\II^0$ is a real-symmetric $K_0$-quasidisk,
we have that for any $x\in I^0\setminus I^1$ and $z\in\partial\II^0,$
there exists $c'>0$ such that 
$|x-z|\geq c'\cdot\diam(\eta')$ where $\eta'$ is the shortest
path in $\partial\II^0\cup\mathbb R$ connecting $x$ and $z$.
Since there exists a definite neighbourhood $\bm U$ of 
$\gamma_{[z_1,z_2]}$ such that $F^4$ has small dilatation on $\U,$
$F^4(z_1)\in\II^0$ and  $F^4(z_2)\in(I^0\setminus I^2)$,
there exists $c>0$ such that 
$|z_1-z_2|\geq
c\cdot\diam(\gamma_{[z_1,z_2]}).$
In case $z_1\in\eta_3$ and $z_2\in\eta_2$
the estimate follows exactly as in the proof that
$\bm Y$ has bounded shape.

We construct $H$ using the same method as in the proof of 
Theorem~\ref{thm:central cascades}. 
We define it on
$X_2$ dynamically, by pulling $H_0$ back. 
We have shown that each component $\bm W$ of $X_3$ is 
a $K_1$-quasidisk. Using the same argument as in 
the high return case, we define $H\colon\bm W\rightarrow\widetilde{\bm{W}}$
so that it is $\hat K$-qc mapping which maps each point of
$\partial\bm W\cap\Gamma'\cap \partial \II^i$ to the corresponding 
point for $\tilde F,$ for $i=4,5,6,\dots.$ 
Finally, we extend $H$ so that it is a homeomorphism on $\V$,
maps each topological square
in $\Gamma'\setminus\Gamma$ to the corresponding topological
square for $\tilde F$ and satisfies (\ref{eqn:qc bound 0}) on $X_0,$
which contains the boundary of each such topological square.
\end{pf}

The following proposition combines 
Lemmas \ref{lem:puzzle geometry} and \ref{lem:pl decomp}
in the case when $F\colon\U\rightarrow\V$ is immediately
renormalizable. 

\begin{lem}\label{lem:renintersection}
There exist  $K_1,\hat K\geq 1,$ and $\delta_0,\delta_1>0,$
such that for any $K_0\geq 1,$ $K\geq 1$ and $\delta>0$ the following holds.
Suppose that $F\colon\U\rightarrow\V$ is a qr polynomial-like mapping
that extends the return mapping to a periodic interval $I\owns c_0$
of sufficiently high period.
Assume that $H_0\colon\mathbb C\rightarrow\mathbb C$ is a 
$(K,\delta)$-qc\textbackslash bg mapping such that
$H_0(\V)=\widetilde{\V}$, $H_0(\U)=\widetilde{\U}$, and
for $z\in\partial\U,$ $H_0\circ F(z)=\tilde F\circ H_0(z)$.
Let $X_1^{H_0}$ be the set with bounded geometry for $H_0$.
Suppose that $F$ is immediately
renormalizable, and
let $J=(\alpha,\tau(\alpha))=(\beta',\tau(\beta'))$
be the periodic interval of period two
under $F$. Then there exists a qr polynomial-like renormalization
$F'\colon\U'\rightarrow \V'$ such that $\V'$ is combinatorially defined and
$\mod(\V'\setminus\U')>\delta_1$.
Moreover, 
there exist a $(K',\delta')$-qr\textbackslash bg mapping
$H\colon\V\rightarrow\widetilde{\V},$
with $K'=\max\{\hat K,(1+\mu(\V))(1+\mu(\widetilde{\V}))K\}$ and
$\delta'=\min\{\delta_0,(1+C\mu(\V))^{-1}\delta, (1+C\mu(\widetilde{\V}))^{-1}\delta\},$
where $C>0$ is a universal constant,
a partition of $\V$ into sets
$\V=X_0\cup X_1\cup X_2\cup X_3,$ 
and similarly for the objects marked with a tilde,
such that 
\begin{enumerate}
\item $H(X_i)=\tilde X_i,$ for $i=0,1,2,3$.
\item $H(\V')=\widetilde{\V}',$ $H(\U')=\widetilde{\U}'$ and if
  $x\in\partial\U'$ then $H(F^i(z))=\tilde F^i(H(z)), i=1,2.$ 
\item $X_1$ is a union of combinatorially defined topological squares
  $\bm S$ such that each $\bm S$ has $\delta'$-bounded geometry and
satisfies the following moduli bounds:
$$\mod(\U\setminus\bm S)>\delta',$$ and if $\bm S\subset\mathbb H^+,$
then 
$$\mod(\U\cap\mathbb H\setminus\bm S)>\delta'.$$
\item Each component of $X_3$ is a $K_1$-quasidisk and $H$ is $\hat
  K$-qc on $X_3.$
\item $H$ is $K'$-qc on $X_2$ moreover, $H$ is defined dynamically on
  $X_2$ and $X_2$ includes $\bm V^0\setminus F^{-1}(\U),$ and a
  combinatorally defined neighbourhood of the rays landing at
  $\beta,\alpha,\zeta, \tau(\beta),\tau(\alpha),$ and $\zeta'$,
where $\zeta,\zeta'\in U$ are the real preimages of 
$\tau(\alpha).$ 
\item The measure of $X_0$ is zero and $H$ satisfies (\ref{eqn:qc
    bound 0}) on $X_0$.
\end{enumerate}
\end{lem}

\begin{pf}
The first part of the proof of this lemma is almost the same as the
proof of Lemma~\ref{lem:pl decomp}; however, we are going to construct
larger shielding regions about the rays. Instead of using two
``squares'' on either side of the square containing the ray, we will
use three. This will let us use the outer boundary of this ``shielding
region'' as part of the boundary of  $\V'$.

Let $\II^0=\V$, and as usual,
let $\II^{i+1}=\comp_{c_0}(\II^i)$.
We will first construct regions about the rays landing at 
$\alpha,\tau(\alpha),\zeta$ and $\zeta'$.
As before, order the components of
$(\II^i\setminus\II^{i+1})\setminus F^{-i}(I^0\setminus I^0)$
counterclockwise
by $S^j_i,$ $j=\{0,\dots, 2d^{i}\}$. 
For each $i\in\mathbb N\cup\{0\}$ let
$S_{i}^{j_\alpha}$ be the component of
$\cup_{j=0}^{2d^i}S_i^j$ that is contained in the upper
half-plane and
intersects a ray landing at $\alpha$.
Let $S_{i,\alpha}=S_{i}^{j_\alpha}\cup S_{i}^{2d^i+1-j_{\alpha}}.$
We set
$$S'_{i,\alpha}=\cup_{k=-3}^3 S_{i}^{j_\alpha+k}\bigcup \cup_{k=-2}^s S_{i}^{2d^i+1-j_{\alpha}+k}.$$
We let 
$$\Gamma_\alpha=\cup_{i=0}^\infty S_{i,\alpha}\mbox{ and } \Gamma_\alpha'=\cup_{i=0}^\infty S_{i,\alpha}'.$$ 
Carrying out identical constructions, we construct regions 
$\Gamma_a\subset\Gamma_a'$ about the rays landing at $a$ for
$a\in\{\tau(\alpha),\zeta,\zeta'\}$. We set
$$\Gamma'=\bigcup_{a\in\{\alpha,\beta,\tau(\alpha),\tau(\beta),\zeta,\zeta'\}}\Gamma'_a\quad\mbox{
  and}\quad \Gamma=\bigcup_{a\in\{\alpha,\beta,\tau(\alpha),\tau(\beta),\zeta,\zeta'\}}\Gamma_a.$$

Let us now construct little topological disks about $\alpha$ and its preimages.
Let $\alpha'$ be the orientation reversing fixed
point of $F^2$ closest to $c_0.$ Let $x_0=\alpha',$ and let 
$x_{i+1}$ be the point in $F^{-1}(x_i)$ that is contained in the
monotone branch of $F$ which contains $\alpha$.
Then the $x_i$ converge to $\alpha$ from either side at a 
definite exponential rate. We can choose $i_0$ so that
$|\alpha-x_{i_0}|$ and $|\alpha-x_{i_0+1}|$ are both comparable 
to $|U|.$ Now, choose $D_{\alpha}$ so that
 $\partial D_{\alpha}\cap\partial\Gamma'_{\alpha}$
are combinatorially defined points lying the $\partial\II^{i_0'}\cap
F^{-i_0'}(\overline V\setminus \overline U).$ 
By doing this, we can carry out this construction simultaneously
for $F$ and $\tilde F$.

Let $D_{\alpha}'=\comp_{\alpha}F^{-1}(D_{\alpha})$
By \cite{dMvS}, Chapter IV, Theorem B 
$|DF(\alpha)|>1+\eta,$
so $\mod(D_{\alpha}\setminus
 D_{\tau(\alpha)})$ is bounded away from zero.
Let $D_{\tau(\alpha)}=\comp_{\tau(\alpha)}F^{-1}(D_{\alpha})$ 
and let $D_{\zeta}=\comp_{\zeta}F^{-1}(D_{\tau(\alpha)}),$ and
 $D_{\zeta'}=\comp_{\zeta'}F^{-1}(D_{\tau(\alpha)}).$
Set $$\Delta'=\cup_{a\in\{\alpha,\tau(\alpha),\zeta,\zeta'\}} D_{a}.$$

Let us now define an intermediate qc\textbackslash bg
mapping $H'\colon\V\rightarrow\widetilde{\V}$ that agrees with $H_0$
on $\partial\U$. This is identical to the construction of $H$ in
Lemma~\ref{lem:pl decomp}, so we will be brief. 
We define $H'$ in the annuli, 
$F^{-i}(\V\setminus\U)$ for $i=1,2,3$, by pulling back $H_0$.
The components of $\II^4\setminus\overline{\Gamma'\cup \Delta'}$
are $K_1$-quasidisks, which can be seen by using 
Lemma~\ref{lem:qdc},
so as in the proof of 
Theorem~\ref{thm:central cascades} in the high return case, 
for any component $\bm W$ of 
$\II^4\setminus\overline{\Gamma'\cup \Delta'}$
we can define a $\hat K$-qc mapping from $\bm W$
to $\widetilde{\bm W}$ that is a conjugacy on $\overline{\Gamma'}
\cap\partial\II_i\cap\partial\bm W$ for $i=5,6,\dots$.
We define $H'$ on $\Gamma$ and in $\Delta$ by pulling back
the mappings already constructed. We extend it to each topological
square $\overline{\bm{S}}$ of 
$\Gamma'\setminus\Gamma$ so that it is a homomorphism
and satisfies (\ref{eqn:qc bound 0}) at each point $x\in\partial\bm S$,
for each $\bm S$.

We obtain a polynomial-like mapping with complex bounds that extends
$F^2|_{\alpha,\tau(\alpha)}$ as follows:
Let
\begin{itemize}
\item $v_{1}$ be the connected component of  
$\partial D\setminus\Gamma'_\alpha$ that does not
intersect $(\alpha,\tau(\alpha))$,
\item  $v_2$ be the two arcs of 
$\partial\Gamma'_{\alpha}\setminus\partial\U$ that
intersect $\bar v_\alpha$,
\item  $v_{3}$ be the connected component of $\partial
  D_{\tau(\alpha)} \setminus\Gamma'_{\tau(\alpha)}$ that does not
intersect $(\alpha,\tau(\alpha))$,
\item $v_4$ be two arcs of $\partial\Gamma'_{\tau(\alpha)}\setminus\partial\U$ that
intersect $\bar v_{\tau(\alpha)}$, and 
\item$v_5$ be the two arcs of $\partial \U$ lying in the upper and lower
half-planes, respectively, that connect  $v_{2}$ and
$v_{4}.$
\end{itemize}

Let $\V'$ be the region bounded by $\overline{v_1\cup v_2\cup v_3\cup
  v_4\cup v_5},$
and let $\U'=\comp_{c_0}F^{-2}(\V')$.
By construction we have that $\U'\Subset\V'$.
Since $\diam(J_F)\asymp\diam(\V')$ and 
$\dist(\partial\V', J_{F^2})$ is bounded away from zero, we have that
$\mod(\V'\setminus\U')$ is also bounded away from 0. 
We have already
defined $H'$ on $\V$ so that $H'(\V')=\widetilde{\V}'$.
Now we are going to define $H$. Outside of $\U'$, we
define $H$ by pulling back $H'$. 
To define $H$ in $\U'$, we use the same ideas as before:
In the sets $\Delta'\cap\U'$ and $\Gamma\cap\U',$ $H$ 
is already defined and the complementary regions in $\U$ 
consist of a $K_1$-quasidisk,
and a region consisting of topological
squares with bounded geometry and with moduli bounds 
as stated in conclusion (2) of the lemma.
\end{pf}

\begin{lem}\label{lem:infrenqd}
For any $K_0'\geq 1$ there exists
$K_0\geq 1$ such that the following holds.
Suppose that $f$ is infinitely renormalizable at an even critical
point $c_0.$
Suppose that $F_0\colon\U_0\rightarrow\V_0$ is a qr polynomial-like mapping
that extends the first return mapping to a periodic interval $J_0\owns
c_0$ of sufficiently high period $s$. 
Let $s_0=s$ and for $i\in\mathbb N$, let $s_i>s_{i-1}$ be minimal 
so that there exists a periodic interval $J_i\owns c_0$ of period $s_i$.
Assume that $\V_0$ is a
$K_0'$-quasidisk.
Then there exists a qr polynomial-like mapping
$F_i\colon\U_i\rightarrow\V_i$ that extends $f^{s_i}\colon J_i\rightarrow J_i$
such that $\V_i$ is a $K_0$-quasidisk.
\end{lem}
\begin{pf}
We will use
Lemma~\ref{lem:qdc}.
Let us consider the
case
when there are several immediately renormalizable levels in a row.
Let $F_i\colon\U_i\rightarrow\V_i$ denote this sequence of renormalization
levels $i=1,2,3,\dots, k.$
Suppose that $z_1,z_2\in\partial\V_k$, and let
$\gamma_{[z_1,z_2]}$ denote the shortest path connecting $z_1$ and
$z_2$ in $\partial\V_k$.
The boundary of $\V_k$ consists of arcs in $\II^6_{k-1},$ arcs in 
$\partial\Delta_k$ and arcs in $\partial\Gamma'_k$.
If $\gamma_{[z_1,z_2]}\subset \partial \Delta_k\cup\partial\Gamma'_k,$
then arguing exactly as for the staircases landing at the $\beta$
fixed point in the high return case, we have 
that there exists a constant $C$ such that
$\diam(\gamma_{[z_1,z_2]})\leq C|z_1-z_2|.$
Thus we can assume that
$\gamma_{[z_1,z_2]}\cap\partial\II^6_{k-1}\neq\emptyset.$
If $\gamma_{[z_1,z_2]}\cap\II^7_{k-1}\neq\emptyset$, then by the moduli
bounds, we have that there exists $C$, depending only on
the moduli bounds, such that
$|z_1-z_2|\geq \frac{1}{C}\diam(\V_{k})\geq
\frac{1}{C}\diam(\gamma_{[z_1,z_2]})$.
So we can assume that 
$\gamma_{[z_1,z_2]}\cap\partial\II^{6}_{k-1}\neq\emptyset$
and $\gamma_{[z_1,z_2]}\cap\II^7_{k-1}=\emptyset.$
Now we have that $F_{k-1}^6(\gamma_{[z_1,z_2]})\subset\partial\V_{k-1}\cup
V_{k-1}\setminus U_{k-1}.$
If $F_{k-1}^6(\gamma_{[z_1,z_2]})\cap V_{k-1}\setminus U_{k-1}\neq\emptyset,$
then, since $\partial \Delta_{k-1}\cap\partial \V_{k-1}$ consists of 
$K_2$-quasiarcs (for a universal $K_2$),
and
the diameters 
and the distances between
the arcs in 
$\partial \Delta_{k-1}\cap\partial\V_{k-1},$ and
$\diam(\partial \Gamma'_{k-1}\cap\partial\V_{k-1})$
are all comparable to $\diam(\V_{k-1}),$ 
we have that there exists $C'$ such that
$C'|F_{k-1}^6(z_1)-F_{k-1}^6(z_2)|\geq \diam(F_{k-1}^6(\gamma_{[z_1,z_2]})),$
and so there exists $C$ such that $C|z_1-z_2|\geq \gamma_{[z_1,z_2]},$
since $F_{k-1}^6$ has small dilatation in a definite neighbourhood of 
$\gamma_{[z_1,z_2]}.$
Thus we can assume that
$F_{k-1}^6(\gamma_{[z_1,z_2]})\subset\partial\V_{k-1}.$
Now we can repeat the argument that we just gave to to show that
$F_{k-1}^6(\gamma_{[z_1,z_2]})\cap\partial\II^{6}_{k-2}\neq\emptyset$
and $F_{k-1}^6(\gamma_{[z_1,z_2]})\cap\II^7_{k-2}=\emptyset,$
and again if 
$F_{k-2}^6(F_{k-1}^6(\gamma_{[z_1,z_2]}))\cap V_{k-2}\setminus
U_{k-2}\neq\emptyset$, we have that there exists
$C'$ such that
$C'|F_{k-2}^6\circ F_{k-1}^6(z_1)-F_{k-2}^6\circ F_{k-1}^6(z_2)|\geq 
\diam(F_{k-2}^6(F_{k-1}^6(\gamma_{[z_1,z_2]}))),$ and again we are
done since
$F_{k-2}^6\circ F_{k-1}^6$ has bounded dilatation in a
definite neighbourhood of $\gamma_{[z_1,z_2]}.$
So we can assume that 
$F_{k-2}^6(F_{k-1}^6(\gamma_{[z_1,z_2]}))\subset\partial\V_{k-2}.$
We can repeat this argument infinitely many times since,
in the circumstances where we repeat the argument, 
by the moduli bounds, the dilatation of 
$F_{i}^6\circ F_{i+1}^6\circ\dots\circ F_{k-1}^6$  is bounded in a
neighbourhood of $\gamma_{[z_1,z_2]}.$ 

The proof is the almost the same in the general case,
when there are both immediate and non-immediate renormalizations
of $f$ about a critical point $c_0$. 
However, if $F_{i-1}\colon\U_{i-1}\rightarrow\V_{i-1}$ is not immediately
renormalizable,
then $F_i\colon\U_i\rightarrow\V_i$ is the first return mapping
from $\U_i=\mathcal{L}_{c_0}(\bm P)$ to  $\V_i=\bm P,$
where $\bm P$ is the puzzle piece given by 
Lemma~\ref{lem:puzzle geometry}. The domain $\bm P=\II_n$ is the
first terminating interval in the enhanced nest about $c_0,$
starting with $\II_0=\bm Q$, see the proof of 
Lemma~\ref{lem:puzzle geometry}.
So $F_{i-1}^{p_{n-1}+p_{n-2}+\dots+p_0}\colon\bm P\rightarrow\bm Q$
has bounded dilatation, and we can argue as in
Proposition~\ref{prop:quasidisks} in between renormalization levels, 
and we can argue as in the immediately renormalizable case
to pass from one renormalization to the next.
\end{pf}

\medskip
\noindent\textit{Concluding the infinitely renormalizable case.}
Let $c_0$ be an even critical point at which $f$ is infinitely
renormalizable.
To conclude the infinitely renormalizable case, 
let us describe how to use these lemmas to build a 
qc\textbackslash bg mapping
that gives a partial conjugacy between $f$ and $\tilde f$;
that is, a conjugacy between $f$ and $\tilde f$ defined on a
subset of $M$ up to the landing
time of a point to a neighbourhood of $\crit(f)\cap\omega(c_0)$.
Let $J_0\owns c_0$ be a sufficiently small periodic interval that
contains $c_0$. Let $s_0$ be the period of $J_0$. Inductively define
$s_{i+1}>s_i,$ minimal, so that there exists a periodic interval
$J_{i+1}\owns c_0$ with period $s_{i+1}$.
Let
$F_i\colon\bm U_i\rightarrow\bm V_i$
be the qr polynomial-like mapping that extends the return mapping to
$J_i,$ where we choose the qr polynomial-like mappings
$F_i\colon\U_i\rightarrow\V_i$ 
as follows:
The mapping $F_0\colon\U_0\rightarrow \V_0$ is given by
Theorem~\ref{thm:box mapping persistent}.
Inductively, we define $F_{i+1}\colon\U_{i+1}\rightarrow\V_{i+1},$ 
if $F_i$ is not immediately renormalizable 
$F_{i+1}\colon\U_{i+1}\rightarrow\V_{i+1}$ is given by the first return mapping
to the first terminating interval in the enhanced nest about $c_0$
starting with $\II_0=\bm Q,$ where $\bm Q$ is the puzzle piece 
from the proof of Lemma~\ref{lem:puzzle geometry},
and if $F_i$ is immediately renormalizable,
$F_{i+1}\colon\U_{i+1}\rightarrow\V_{i+1}$ is the qr polynomial-like mapping
given by Lemma~\ref{lem:renintersection}.
Let $\{J_i^j\}_{j=0}^{s_i}$ be the chain of intervals with
$J_i^{s_i}=J_{i}$ and $J_{i}^j=\comp_{f^j(c_0)}f^{-(s_i-j)}(J_i)$.
Let $\{\U_i^j\}$ be the chain of domains associated to the first return mapping
$F_i\colon\U_i\rightarrow\V_i$, so that $\U_i^j\supset J_i^j$.
Since there are only finitely many components of $\{\U_0^j\}$
it is easy to build by hand a qc mapping $H_0$
of the plane so that 
for $0\leq j<s_0-1,$
$H_0(\U_0^j)=\widetilde{\U}_0^j$,
$z\in\partial\U_0^j$ we have
$H_0\circ f(z)=\tilde f\circ H_0(z)$.
Now,
Lemma~\ref{lem:pl decomp} allows us to extend $H_0$
as a qc\textbackslash bg mapping to a mapping 
$H'_0\colon\V_0\rightarrow \widetilde{\V}_0$ so that it
preserves some additional dynamical information inside of $\V_0$,
see the conclusion of Lemma~\ref{lem:pl decomp}.
If the next renormalization is not of intersection type, then 
by Theorem~\ref{thm:central cascades} and the (proof of)
the Real Spreading Principle, we obtain a
qc\textbackslash bg mapping $H_1$ of the plane 
such that $H_0(\U_1^j)=\widetilde{\U}_1^j$, if
$z\in\partial\U_0^j$ we have
$H_0\circ f(z)=\tilde f\circ H_0(z),$
and for $z\notin\V_1$ we have
$H_1\circ f(z)=\tilde f\circ H_1(z)$.
Thus we can continue on to the next renormalization level.

If the next renormalization level is of intersection type, then we have
from
Lemma~\ref{lem:renintersection}, a qc\textbackslash bg
mapping $H'_0$ of the plane such that 
$H_0'(\U_1^j)=\widetilde{\U}_1^j$,
and for $0\leq j<s_1-1, z\in\partial\U_1^j$ 
$H_0'\circ f(z)=\tilde f\circ H_0'(z).$
We need to extend $H_0'$ to the union of intervals
$\cup_{j=0}^{s_0-1}V_0^j,$ so that it is a conjugacy outside of
a neighbourhood of the critical points.
In this case, all the critical points of
$F_0\colon\U_0\rightarrow\V_0$ are contained in 
$(-\alpha,\tau(\alpha))\cup(\alpha,\zeta),$
where we choose $\zeta$ to be the preimage of $\tau(\alpha)$ that is
closest to
to $\alpha$.
First, we extend $H_0'$ to the component $\bm X$ 
of $\V_0\setminus\Gamma'$
that has $\beta_0$ in its boundary.
Since this branch contains no critical points, we define $H_0'$
dynamically. If  $z\in \bm X$ and $k>0$ is minimal so that 
$F_0^k(z)\in\comp_{c_0}(\V_0\setminus\Gamma')\cup\Gamma'_{\tau(\alpha)},$
then we define $H_0'(z)$ by $\tilde F_0^k \circ H_{0}'(z)=H_0'\circ
F_0^k(z).$ We define $H_0'$ on the component
of $\V_0\setminus\Gamma'$ that has $\tau(\beta)$ in its boundary
by pulling the extension of $H_0'$ just constructed back by one
iterate of $F$.
Now, we spread this around to the 
each $\V_0^j$ by pulling back by the 
landing mapping $f^{s_0-j}\colon\V_0^j\rightarrow \V_0$,
except that we extend $H'_0$ to each component of the partition that
contains a critical point using the fact that it is a quasidisk, that
is bounded by topological squares with bounded geometry and moduli
bounds as in the conclusion of Lemma~\ref{lem:renintersection}, and
the same argument to build a $\hat K$-qc mapping that we have used
many times.

Observe that the qc constants and geometric bounds are controlled as 
we pullback from one level to the next since $\diam(\V_i)$
decay exponentially.

\subsection{The persistently recurrent,
finitely renormalizable case}
\label{subsec:qc rigidity-persistent} 
In the remaining sections,
we will consider the smooth and analytic cases together.
In the analytic case, we obtain a qc mapping that conjugates
the dynamics of $f$ and $\tilde f$ in a necklace neighbourhood of the
real line, but in the smooth case we only obtain a conjugacy on the
real-line. It is important to observe that this conjugacy is obtained 
as the restriction to the real-line of a qc\textbackslash bg mapping
in the upper half-plane, which is quasisymmetric by the QC Criterion.

In this section we treat critical points
$c\in\mathrm{Crit}(f)$ such that $\omega(c)$ is persistently recurrent 
and $f$ is not infinitely renormalizable at $c$.
Fix such a critical point $c$ and let $\tilde c$
be the corresponding critical point for $\tilde f$.
In this case, Theorem~\ref{thm:box mapping persistent}
implies that $c$ and $\tilde c$ possess arbitrarily
small combinatorially defined nice real neighbourhoods
$I$ and $\tilde I$, such that the first return maps to
$I\cup \cup_{c'\in\omega(c)}\mathcal{L}_{c'}(I)$ and
$\tilde I\cup \cup_{\tilde c'\in\omega(\tilde c)}\mathcal{L}_{\tilde c'}(\tilde I)$ 
restricted to components that intersect
$\omega(c)$ and $\omega(\tilde c)$, respectively,
extend to complex box mappings 
$$F\colon \bm{\mathcal{U}}\rightarrow \bm{\mathcal{V}}\;
\mathrm{and}\;\tilde{F}\colon \widetilde{\bm{\mathcal{U}}}\rightarrow
\widetilde{\bm{\mathcal{V}}}$$
whose ranges are $\delta$-nice and have $\delta$-bounded geometry.
Moreover, these box mappings are strongly combinatorially equivalent.
		
If $F$ and $\tilde F$ are analytic,
by Theorem~\ref{thm: qc rigidity of box mappings}, these mappings are
qc-conjugate provided that there exists a quasiconformal external conjugacy
between the maps. Since both $\bm{\mathcal{U}}$ and 
$\widetilde{\bm{\mathcal{U}}}$
consist of just finitely many components,
whose boundaries are piecewise smooth Jordan curves 
with a finite number of corners,
such a conjugacy is simple to construct.
On the other hand, if $F$ and $\tilde F$ are $C^3$,
then we still have that there is an external conjugacy between $F$ and
$\tilde F$;however, we do not have that $F$ and $\tilde F$ are
qc-conjugate.
Instead by Theorem~\ref{thm:central cascades} and the argument
used to prove the Real Spreading Principle, we have that
there exist $K\geq 1$ and $\delta>0$ and a
$(K,\delta)$-qc\textbackslash bg mapping in the upper half-plane
$H\colon\bm{\mathcal{V}}\rightarrow\widetilde{\bm{\mathcal{V}}}$
that conjugates $F$ and $\tilde F$ on their real traces.

In the analytic case, we immediately get that the box mappings
$F\colon \bm{\mathcal{U}}\rightarrow \bm{\mathcal{V}}$ and
$\tilde{F}\colon \widetilde{\bm{\mathcal{U}}}\rightarrow\widetilde{\bm{\mathcal{V}}}$
are qc-conjugate by a
quasiconformal homeomorphism 
$h\colon\bm{\mathcal{V}}\rightarrow\widetilde{\bm{\mathcal{V}}}$
that is a conjugacy on the set
$$K(I_1)=\{z:(R_{I_1})^n(z)\in \Dom_{\Omega_1}(I_1)\ \mathrm{for\ all}\ n\geq 0\}\supset \omega(c),$$
where $\Dom_{\Omega_1}(I_1)$ is as defined on Page
\pageref{notation:D(I)}.
In the smooth case, we have that $h$ is qc\textbackslash bg
in the upper half-plane, so again, by the QC Criterion, we have that
$h|_{\mathbb R}$ is a quasisymmetric mapping that
is a conjugacy on
$K(I_1).$
		
Recall that $I_1=\cup_{c\in\Omega_1\setminus\Omega_1'}\mathcal{L}_c(I)\cup I'_1$,
$\Dom_{\Omega_1}(I_1)$
is the collection of components of $\Dom(I_1)$ that intersect
$\cup_{c\in\Omega_1}\omega(c)$ and $\Dom_{\mathrm{diff}}(I_1)$ is the union of
components, $J$, of $\mathrm{Dom}(I_1)$ satisfying
\begin{enumerate}
\item $J\subset I'_1$;
\item $J\cap (\cup_{c\in\Omega_1}\omega(c))=\emptyset$;
\item if $s$ is the first return time of $J$ to $I_1$,
then $f^i(J)\cap\mathcal{I}=\emptyset$ for all $i=0,1,\dots, s-1$.
\end{enumerate}
Consider the first entry map
$$R_{\hat{I}_*}\colon \Dom_{\Omega_1}(I_1)\cup
\Dom_{\mathrm{diff}}(I_1)\cup (\Dom(I_0)\cap I)
\rightarrow I_1\cup I_0.$$
By definition if $J$ is a component of the domain and $s=s(J)$ is its return time,
then $f^i(J)\cap\mathcal{I}=\emptyset$ for $i=1,\dots,s-1$.
By Theorem~\ref{thm:box mapping persistent infinite branches},
there exists an integer $m$ such that the $m$-th iterate of this map
has an extension to a complex box mapping 
$$F_{it}\colon\bm{\mathcal{U}}_{it}\rightarrow\bm{\mathcal{V}}_{it}$$
with $\bm{\mathcal{V}}_{it}\cap\mathbb{R}=I_1\cup I$ 
such that $\bm{\mathcal{V}}_{it}$ is $\delta$-nice and any component of either
$\bm{\mathcal{U}}_{it}$ or $\bm{\mathcal{V}}_{it}$
has $\delta$-bounded geometry,
and so that the mapping $F_{it}\colon\Dom_{\Omega_1}(\bm{\mathcal{V}}_{it})
\rightarrow\bm{\mathcal{V}}_{it}$;
that is, the mapping $F_{it}$
restricted to components of the domain that
intersect $\omega(c)$, for some $c\in\Omega_1,$ satisfies the gap
and extension conditions.
From this it is easy to construct a 
box mapping $F_{it}'\colon\bm{\mathcal{U}}_{it}'\rightarrow\bm{\mathcal{V}}_{it}'$
where $\bm{\mathcal{V}}_{it}'=\cup_{c\in\Omega_1'}\bm{\mathcal{V}}_{it}(c)$,
and each component of $\bm{\mathcal{U}}_{it}'$ 
is a component of $\bm{\mathcal{U}}_{it}$
(if necessary, simply post-compose $F_{it}|_{\U}$ by a first landing onto
some 
component $\V\owns c$ of $\bm{\mathcal{V}}_{it}$ with
$ c\in\Omega_1'$).
As usual, we mark the corresponding object for $\tilde f$ with a tilde.
Moreover, by Theorem~\ref{thm:qc rigidity away from critical points}, 
the set 
$$E(\mathcal{I})=\{z\in [0,1]:f^{n}(z)\notin\mathcal{I},\ \mathrm{for}\ n=0,1,\dots\}$$
is quasiconformally rigid. 
Since $\partial \Dom_{\Omega_1}(I_1)\cup 
\partial\Dom_{\mathrm{diff}}(I_1)\subset E(\mathcal{I}),$
it follows that we may apply Theorem~\ref{thm:external conjugacy}
to obtain
a quasiconformal external conjugacy between $F_{it}'$ and $\tilde{F}_{it}'.$
In the analytic case, we have that Theorem~\ref{thm: qc rigidity of box mappings} implies that
$F_{it}'\colon \bm{\mathcal{U}}_{it}'\rightarrow
                \bm{\mathcal{V}}_{it}'$ 
and the corresponding box mapping
		$\tilde{F}_{it}'\colon
                \widetilde{\bm{\mathcal{U}}}_{it}'\rightarrow
\widetilde{\bm{\mathcal{V}}}_{it}'$
are quasiconformally conjugate. Arguing as before in the smooth case
we have that there exists a qc\textbackslash bg mapping
in the upper half-plane that conjugates $F_{it}'$ and $\tilde F_{it}'$
on their real traces.
It follows that the set
		$$K_m(I_0\cup I_1\cup \mathcal{I}^{\mathrm{Comp}})=
		\{z:(R_{\hat{I}_*})^{mn}(z)\in I_0\cup I_1\cup E(\mathcal{I})\ \mathrm{for\ all\ }n\geq 0\}$$
is quasiconformally rigid.
By pulling back at most $m$ times along the branches of 
		$$R_{\hat{I}_*}\colon \Dom_{\Omega_1}(I_1)\cup\Dom_{\mathrm{diff}}(I_1)\rightarrow I_1$$
one can obtain a
quasisymmetric homeomorphism that is a conjugacy on the set
		\begin{equation}\label{eqn:rigidityhatI}
		K(I_0\cup I_1\cup \mathcal{I}^{\mathrm{comp}})=
		\{z:(R_{\hat{I}_*})^n(z)\in I_0\cup I_1\cup E(\mathcal{I})\ \mathrm{for\ all\ }
		n\geq 0\}.
		\end{equation}

\subsection{The reluctantly recurrent case}
\label{subsec:qc rigidity-reluctant} 
Finally, we treat the critical points;
$c\in\Omega_2$ such that
$f$ is reluctantly recurrent at $c.$
Let $\Omega_r\subset\Omega_2$
be a non-trivial block of critical points such that
each recurrent critical point $c\in\Omega_r$
is reluctantly recurrent, and with the property that
if $\Omega$ is the component of the graph on the critical points
that contains $\Omega_r$, then 
$f$ is not infinitely renormalizable at any $c'\in\Omega$.
Since we have rigidity away from the critical set, we may as well assume
that if $c\in\Omega_r$ is non-recurrent, then the orbit of $c$
accumulates on a critical point in $\Omega_r$.
Let $\Omega_r'\subset\Omega_r$ denote a smallest set of critical points
such that each critical point in $\Omega_r'$ is recurrent, and
any critical point $c\in\Omega_r\setminus \Omega_r'$
accumulates on some critical point $c'\in\Omega_r'$.
		
Let $I_r$ denote a small real neighbourhood of
$\Omega_r$ with the properties given by
Proposition~\ref{prop:real bounds reluctant},
and let $I_r'=\cup_{c\in\Omega_r'}I_r(c).$
		
Let
$$R_{I_r}\colon \Dom_{\Omega_r}(I_r)\cup
\Dom_{I_0\cup I_1'}(I_r)\rightarrow I_r\cup I_0\cup I_1'$$
denote the first entry map, where
$\Dom_{\Omega}(I_r)$ denotes the collection of domains of the return mapping to $I_r$
that intersect $\omega(c)$ for some critical point $c\in \Omega_r$
and
$\Dom_{I_0\cup I_1'}(I_r)$ consists of the intervals in $I_r$ that
are eventually mapped into $I_0\cup I_1'$.
		
By Theorem~\ref{thm:box mapping reluctant}
and our choice of neighbourhood of
$\crit(f)$, we can assume that
this mapping has has an extension to a box mapping
$$\hat F\colon \UU\rightarrow \VV,$$
where for components
$I_{0,i}$ of $I_0$ the set
$\mathrm{Comp}_{I_{0,i}}\VV$ is a Poincar\'e disk $D_{\pi-\theta_0}(I_{0,i})$,
for components
$I_{1,i}'$ of $I_1'$ the set
$\mathrm{Comp}_{I_{1,i}'}\VV$ is a Poincar\'e disk $D_{\pi-\theta_1}(I_{1,i}')$
and for components $I'_{r,i}$ of $I_r'$
the set
$\mathrm{Comp}_{I_{r,i}'}\VV$ is a Poincar\'e disk $D_{\pi-\theta}(I_{r,i}')$,
		where $0<\theta_0<\theta_1<\theta$,
and restricted to components that intersect 
the orbit of a critical point $c\in\Omega_r$
this box mapping has the gap and extension properties.
The angles $\theta_0,\theta_1$ are chosen as in the proof of
Theorem~\ref{thm:box mapping persistent infinite branches}.
		
\medskip		
\noindent		
\textit{Claim:} $\partial \UU\cap\mathbb R$ is qs-rigid.

\medskip
\noindent		
\textit{Proof.} 
Let
$$E(I_r)=\{x\in[0,1]: f^n(x)\notin I_r\ \mathrm{for\ all\ }n\geq 0\},\ E'(I_r)
			=\{x\in[0,1]:f^{n}(x)\notin I_r\ \mathrm{for\ all\ }n\geq 1\}$$
	and
$$CV(f)=\{f(c):c\in\Omega_r\ \mathrm{and}\ f(c)\notin E(I_r)\}.$$
Notice that $\partial \UU\cap\mathbb{R}\subset E'(I_r)$.
Consider a touching box mapping as in
Subsection \ref{subsec:global touching box mappings}
$$F_T\colon \UU_T\rightarrow \VV_T,$$
so that the components of $\VV_T\cap\mathbb{R}$ that contain critical points 
in $\Omega_0\cup\Omega_1\cup\Omega_r$
are the components of 
$I_0, I_1$ and $I_r$.
Recall that $\UU_T\cap I_i=\emptyset, i=0,1,r$,
	and let 
$$\tilde{F}\colon \widetilde{\UU}_T\rightarrow\widetilde{\VV}_T$$
be the corresponding touching box mapping for $\tilde{f}$.
By Theorem~\ref{thm:touching box map} 
and the qs-rigidity of the set
$$
K(I_0\cup I_1\cup \mathcal{I}^{\mathrm{comp}})=
\{z:(R_{\hat{I}_*})^n(z)\in I_0\cup I_1\cup E(\mathcal{I})\ \mathrm{for\ all\ }
n\geq 0\},$$ 
as above,
there exists a qc-homeomorphism $H_0$ that maps
$\VV_T$ to $\widetilde{\VV}_T$ with 
$H_0(\UU_T)=\widetilde{\UU}_T$ that agrees with the given topological conjugacy
between $f$ and $\tilde f$. 
As in Corollary~\ref{cor: pullback argument for touching box mappings}
we get a qc-homeomorphism $H$ with
			$H(\UU_T)=\widetilde{\UU}_T$
with $\tilde{F}_T\circ H= H\circ F_T$
on $\UU_T$ and which is a conjugacy on the closed set $E(I_r)$.
We modify $H$ on the finite set, $CV(f)$, so that it sends $CV(f)$ to $CV(\tilde{f}).$
Pulling back $H$ once, we obtain a qc-conjugacy on the set $E'(I_r),$
and since $\partial \UU\cap\mathbb{R}\subset E'(I_r)$, we have 
			completed the proof of the claim. \endpfclaim

It follows from Theorem~\ref{thm:external conjugacy},
that there exists a quasiconformal external conjugacy
$H$ between $F$ and $\tilde{F}$ so that
$H(\cup_{c\in\Omega_r'}\VV(c))=\cup_{c\in\tilde{\Omega}_r'}\widetilde{\VV}(c),$
$H(\UU)=\widetilde{\UU}$ and
$H$ is a conjugacy on $\partial \UU$.
Thus, in the analytic case, Theorem~\ref{thm: qc rigidity of box mappings}
implies that
$F\colon \UU\rightarrow \cup_{c\in\Omega_r'}\VV$
			and
			$\tilde F\colon \widetilde{\UU}\rightarrow 
\cup_{c\in\tilde \Omega_r'}\widetilde{\VV}(c)$
are qc conjugate. In the smooth case we have that the mappings
are as conjugate on their real traces by Theorem~\ref{thm:smooth
  box conj}.
Thus
we obtain a quasisymmetric homeomorphism that is a conjugacy on 
$$K(I_r)=\{z\in\I_r: \hat R_{I_r}(z)\in I_r\mbox{ for all } n\geq 0\} .$$
We now repeat this argument for each block of critical points $\Omega_r.$

\subsection{Conclusion}
Since almost every point in the interval eventually lands in $I_0\cup I_1\cup I_2$,
we extend the conjugacies near the critical points and in immediate basins of attraction
almost everywhere in the interval
by considering the first entry mappings to $I_0\cup I_1\cup I_2$,
and we do the same for the objects marked with a tilde.
This yields a globally defined conjugacy on $M$, since
for any component $J$ of the first landing map
			to $I_0\cup I_1\cup I_2$, 
			the first landing mapping extends to a 
			complex domain $\U_J$, and 
			near the boundary of $\U_J$ 
			conjugacy is defined as the pullback
			of the conjugacy between the touching box mappings
			$F_T\colon \UU_T\rightarrow \VV_T$ and
			$\tilde{F}_T\colon\widetilde {\UU}_T\rightarrow \widetilde{\VV}_T$.
			
					
			\begin{figure}[htp] \hfil
			\beginpicture
			\dimen0=0.3cm
			\setcoordinatesystem units <\dimen0,\dimen0> 
			\setplotarea x from -9 to 30, y from -5 to 4
			\setlinear
			\plot -12 0 26 0 /
			\circulararc 120 degrees from 0 0  center at  -5 -3 
			\circulararc -120 degrees from 0 0  center at  -5 3  
			\circulararc 50 degrees from -2 0  center at  -3 -2 
			\circulararc -50 degrees from  -2 0 center at  -3 2 
			\circulararc -90 degrees from -9 0  center at  -7.5 -1.5 
			\circulararc 90 degrees from  -9 0 center at  -7.5 1.5 
			\circulararc -90 degrees from -5.5 0  center at  -5 -0.5 
			\circulararc 90 degrees from  -5.5 0 center at  -5 0.5 
			\put {$I_2$} at -7 4
			
			

			\circulararc 130 degrees from 15 0  center at  8.5 -3 
			\circulararc -130 degrees from 15 0  center at  8.5 3 
			\circulararc -90 degrees from 4 0  center at  5.5 -1.5 
			\circulararc 90 degrees from  4 0 center at 5.5 1.5
			 \arrow <3mm> [0.2,0.67] from 5.5 1 to -3 3
			\circulararc 125 degrees from 12 0  center at  10 -1 
			\circulararc -125 degrees from  12 0 center at 10 1 
			\arrow <3mm> [0.2,0.67] from 10 1.2 to 10 4
			\arrow <3mm> [0.2,0.67] from -2.3 -0.5 to 4 -0.2
			
			\circulararc 50 degrees from 14.5 0  center at  13.5 -2 
			\circulararc -50 degrees from  14.5 0 center at  13.5 2 
			 \arrow <3mm> [0.2,0.67] from 13.5 0.5 to 12 2.5
			\put {$I_1=I$} at 11 5
			
			\circulararc -90 degrees from 19 0  center at  20.5 -1.5 
			\circulararc 90 degrees from  19 0 center at 20.5 1.5
			\put {$I_0$} at 20.5 2
			
			\endpicture
			\caption{The box mapping associated to the non-minimal case.
			Some domains of the first return mapping to $I_2$ can spend a lot of time in 
			$I_1$. We also need to consider domains which eventually enter into $I_0$. \label{fig:global}}
			\end{figure}

\section{Table of notation and terminology}
${}$\\

\medskip
\begin{tabular}{cc}
\begin{tabular}{|l|l|}
\hline
$\mathrm{Forward}(c)$, $\mathrm{Back}(c)$ & page~\pageref{page:forward}\\
\hline
$\mathrm{Dom}(P),$ $\mathrm{Dom}_{\Omega}(P)$ & page \pageref{notation:Dom}\\
\hline
$\mathrm{Dom}^*(P)$, $\mathrm{Dom}'(P)$ & page~\pageref{page:return}\\
\hline
$\mathcal{L}_{x}(I)$ & page \pageref{notation:entry domain} \\
\hline
$\hat{\mathcal{L}}_{x}(V)$ & page \pageref{landing domain}\\
\hline
well-inside & page \pageref{well-inside}\\
\hline
$\mathcal{C}(M)=\mathcal C$ & page \pageref{def:class of maps}\\
\hline
\end{tabular} &
\begin{tabular}{|l|l|}
\hline
$\Omega,\Omega_0,\Omega_1,\Omega_2$ & page \pageref{notation:Omega}\\
\hline
$I_0$, $I_1'$, $I_1$, $I_2$ & pages \pageref{notation:critical neighbourhoods}\\
\hline
$\Omega_r, $ $\Omega_{r,e}$, $\Omega_{r,o}$ & page
                                              \pageref{notation:Omega_reluct},
                                              \pageref{notation:reluct
                                              decomp}\\
\hline
$D_{\theta}(I)$ & page \pageref{page:poincare}\\
\hline
$\mathcal{I}$ & page \pageref{notation:all intervals}\\
\hline
$\mathrm{Dom}_{\mathrm{diff}}(I_1)$ & page \pageref{notation:dom_diff}\\
\hline
$E(\mathcal{I})$ & page \pageref{notation:escaping}\\
\hline
\end{tabular}
\end{tabular}

\input{bib.tex}
\end{document}

%% file: prepicte.tex
\catcode`@=11 \catcode`!=11

\expandafter\ifx\csname fiverm\endcsname\relax
  \let\fiverm\fivrm
\fi
  
\let\!latexendpicture=\endpicture 
\let\!latexframe=\frame
\let\!latexlinethickness=\linethickness
\let\!latexmultiput=\multiput
\let\!latexput=\put
 
\def\@picture(#1,#2)(#3,#4){%
  \@picht #2\unitlength
  \setbox\@picbox\hbox to #1\unitlength\bgroup 
  \let\endpicture=\!latexendpicture
  \let\frame=\!latexframe
  \let\linethickness=\!latexlinethickness
  \let\multiput=\!latexmultiput
  \let\put=\!latexput
  \hskip -#3\unitlength \lower #4\unitlength \hbox\bgroup}

\catcode`@=12 \catcode`!=12

%% file: postpict.tex
\catcode`@=11 \catcode`!=11
  
\let\!pictexendpicture=\endpicture 
\let\!pictexframe=\frame
\let\!pictexlinethickness=\linethickness
\let\!pictexmultiput=\multiput
\let\!pictexput=\put

\def\beginpicture{%
  \setbox\!picbox=\hbox\bgroup%
  \let\endpicture=\!pictexendpicture
  \let\frame=\!pictexframe
  \let\linethickness=\!pictexlinethickness
  \let\multiput=\!pictexmultiput
  \let\put=\!pictexput
  \let\input=\@@input   
  \!xleft=\maxdimen  
  \!xright=-\maxdimen
  \!ybot=\maxdimen
  \!ytop=-\maxdimen}

\let\frame=\!latexframe

\let\pictexframe=\!pictexframe

\let\linethickness=\!latexlinethickness
\let\pictexlinethickness=\!pictexlinethickness

\let\\=\@normalcr
\catcode`@=12 \catcode`!=12

%% file: macros.tex
\setcounter{secnumdepth}{4}

\newcommand{\proclaim}[2]{\medbreak {\bf #1}{\sl #2} \medbreak}

\newcommand{\ntop}[2]{\genfrac{}{}{0pt}{1}{#1}{#2}}

\let\newpf\proof \let\proof\relax \let\endproof\relax
\newenvironment{pf}{\newpf[\proofname]}{\qed\endtrivlist}

\def\PC{{\mathrm{PC}}}
\def\R{{\mathbb{R}}}
\def\C{{\mathbb{C}}}
\def\N{{\mathbb{N}}}
\def\H{{\mathbb{H}}}
\def\P{{\mathcal{P}}}
\def\Crit{{\mathrm{Crit}}}
\def\Per{{\mathrm{Per}}}
\def\Y{{\mathcal{Y}}}
\def\Comp{{\mathrm{Comp}}}
\def\dbar{{\bar{\partial}}}
\def\Im{{\mathrm{Im}}}
\def\diam{{\mathrm{diam}}}
\def\Dom{{\mathrm{Dom}}}
\def\dist{{\mathrm{dist}}}
\def\gap{{\mathrm{Gap}}}
\def\Space{{\mathrm{Space}}}
\def\spa{{\mathrm{Space}}}
\def\cen{{\mathrm{Cen}}}
\def\Land{{\hat{\mathcal{L}}}}

\def\comp{{\mathrm{Comp}}}
\def\LL{{\mathcal{L}}}
\def\clos{{\mathrm{Cl}}}
\def\crit{{\mathrm{Crit}}}
\def\mod{{\mathrm{mod}}}
\newcommand{\I}{\mathcal{I}}
\newcommand{\II}{\bm{I}}
\newcommand{\dombf}{\operatorname{\bf D}}
\newcommand{\DD}{\operatorname{\Dom}}
\newcommand{\U}{\bm{U}}
\newcommand{\J}{\bm{J}}
\newcommand{\W}{\bm{W}}
\newcommand{\V}{\bm{V}}
\newcommand{\PP}{\bm{P}}
\newcommand{\Z}{\bm{Z}}
\newcommand{\G}{\bm{G}}
\newcommand{\A}{\mathcal{A}}
\newcommand{\HH}{\mathbb{H}}
\newcommand{\frm}{R}
\newcommand{\compl}{\operatorname{compl}}
\newcommand{\interior}{\operatorname{int}}
\newcommand{\UU}{\bm{\mathcal{U}}}
\newcommand{\VV}{\bm{\mathcal{V}}}

\newtheorem{thm}{Theorem}[section]
\newtheorem{thmx}{Theorem}
\renewcommand{\thethmx}{\Alph{thmx}}
\newtheorem{cor}[thm]{Corollary}
\newtheorem{lem}[thm]{Lemma}
\newtheorem{lemma}[thm]{Lemma}
\newtheorem{claim}[thm]{Claim}
\newtheorem{prop}[thm]{Proposition}
\newtheorem{conj}[thm]{Conjecture}
\newtheorem{schw}{Schwarz Lemma} 
\newtheorem{sectl}{Sector Lemma}
\theoremstyle{remark}
\newtheorem{rem}{Remark}[section]
\newtheorem{notation}{Notation}
\newtheorem{example}{Example}[section]

\numberwithin{equation}{section}
\newcommand{\thmref}[1]{Theorem~\ref{#1}}
\newcommand{\propref}[1]{Proposition~\ref{#1}}
\newcommand{\secref}[1]{\S\ref{#1}}
\newcommand{\lemref}[1]{Lemma~\ref{#1}}
\newcommand{\corref}[1]{Corollary~\ref{#1}} 
\newcommand{\figref}[1]{Fig.~\ref{#1}}

\theoremstyle{definition}
\newtheorem{defn}{Definition}[section]
\newtheorem{definition}{Definition}[section]
\def\proof{\bn {\bf Proof.} }

%% file: lensmapping.pdf_t
\begin{picture}(0,0)%
\includegraphics{lensmapping.pdf}%
\end{picture}%
\setlength{\unitlength}{3947sp}%
\begingroup\makeatletter\ifx\SetFigFont\undefined%
\gdef\SetFigFont#1#2#3#4#5{%
  \reset@font\fontsize{#1}{#2pt}%
  \fontfamily{#3}\fontseries{#4}\fontshape{#5}%
  \selectfont}%
\fi\endgroup%
\begin{picture}(4899,1468)(1114,-4092)
\put(4641,-2724){\makebox(0,0)[lb]{\smash{{\SetFigFont{6}{7.2}{\familydefault}{\mddefault}{\updefault}{\color[rgb]{0,0,0}$\bm{\mathcal{V}}$}%
}}}}
\put(3923,-3555){\makebox(0,0)[lb]{\smash{{\SetFigFont{6}{7.2}{\familydefault}{\mddefault}{\updefault}{\color[rgb]{0,0,0}$\bm{\mathcal{U}}$}%
}}}}
\put(4263,-4046){\makebox(0,0)[lb]{\smash{{\SetFigFont{6}{7.2}{\familydefault}{\mddefault}{\updefault}{\color[rgb]{0,0,0}$c\in\crit(f)$}%
}}}}
\end{picture}%

%% file: qcbgsmall.pdf_t
\begin{picture}(0,0)%
\includegraphics{qcbgsmall.pdf}%
\end{picture}%
\setlength{\unitlength}{3947sp}%
\begingroup\makeatletter\ifx\SetFigFont\undefined%
\gdef\SetFigFont#1#2#3#4#5{%
  \reset@font\fontsize{#1}{#2pt}%
  \fontfamily{#3}\fontseries{#4}\fontshape{#5}%
  \selectfont}%
\fi\endgroup%
\begin{picture}(3610,3446)(289,-1708)
\put(2350,-1477){\makebox(0,0)[lb]{\smash{{\SetFigFont{5}{6.0}{\familydefault}{\mddefault}{\updefault}{\color[rgb]{0,0,0}$X_3$}%
}}}}
\put(1912,787){\makebox(0,0)[lb]{\smash{{\SetFigFont{5}{6.0}{\familydefault}{\mddefault}{\updefault}{\color[rgb]{0,0,0}$X_3$}%
}}}}
\put(3350,1141){\makebox(0,0)[lb]{\smash{{\SetFigFont{5}{6.0}{\familydefault}{\mddefault}{\updefault}{\color[rgb]{0,0,0}$\bm{J}^3$}%
}}}}
\put(2842, 60){\makebox(0,0)[lb]{\smash{{\SetFigFont{5}{6.0}{\familydefault}{\mddefault}{\updefault}{\color[rgb]{0,0,0}$X_2$}%
}}}}
\put(1821, 60){\makebox(0,0)[lb]{\smash{{\SetFigFont{5}{6.0}{\familydefault}{\mddefault}{\updefault}{\color[rgb]{0,0,0}$X_2$}%
}}}}
\put(3884, 37){\makebox(0,0)[lb]{\smash{{\SetFigFont{5}{6.0}{\familydefault}{\mddefault}{\updefault}{\color[rgb]{0,0,0}$X_1$}%
}}}}
\end{picture}%

%% file: renormsmall.pdf_t
\begin{picture}(0,0)%
\includegraphics{renormsmall.pdf}%
\end{picture}%
\setlength{\unitlength}{3947sp}%
\begingroup\makeatletter\ifx\SetFigFont\undefined%
\gdef\SetFigFont#1#2#3#4#5{%
  \reset@font\fontsize{#1}{#2pt}%
  \fontfamily{#3}\fontseries{#4}\fontshape{#5}%
  \selectfont}%
\fi\endgroup%
\begin{picture}(3814,3632)(1143,-2740)
\put(2059,-895){\makebox(0,0)[lb]{\smash{{\SetFigFont{5}{6.0}{\familydefault}{\mddefault}{\updefault}{\color[rgb]{0,0,0}$\beta$}%
}}}}
\put(2625,-867){\makebox(0,0)[lb]{\smash{{\SetFigFont{5}{6.0}{\familydefault}{\mddefault}{\updefault}{\color[rgb]{0,0,0}$\tau(\alpha)$}%
}}}}
\put(3416,-895){\makebox(0,0)[lb]{\smash{{\SetFigFont{5}{6.0}{\familydefault}{\mddefault}{\updefault}{\color[rgb]{0,0,0}$\alpha$}%
}}}}
\put(3810,-867){\makebox(0,0)[lb]{\smash{{\SetFigFont{5}{6.0}{\familydefault}{\mddefault}{\updefault}{\color[rgb]{0,0,0}$\tau(\beta)$}%
}}}}
\put(4178,461){\makebox(0,0)[lb]{\smash{{\SetFigFont{5}{6.0}{\familydefault}{\mddefault}{\updefault}{\color[rgb]{0,0,0}$\bm I^2$}%
}}}}
\put(2851,-499){\makebox(0,0)[lb]{\smash{{\SetFigFont{5}{6.0}{\familydefault}{\mddefault}{\updefault}{\color[rgb]{0,0,0}$X_3$}%
}}}}
\put(2172,-585){\makebox(0,0)[lb]{\smash{{\SetFigFont{5}{6.0}{\familydefault}{\mddefault}{\updefault}{\color[rgb]{0,0,0}$X_3$}%
}}}}
\put(2200,-1320){\makebox(0,0)[lb]{\smash{{\SetFigFont{5}{6.0}{\familydefault}{\mddefault}{\updefault}{\color[rgb]{0,0,0}$X_3$}%
}}}}
\put(1634,-867){\makebox(0,0)[lb]{\smash{{\SetFigFont{5}{6.0}{\familydefault}{\mddefault}{\updefault}{\color[rgb]{0,0,0}$X_2$}%
}}}}
\put(2851,-1489){\makebox(0,0)[lb]{\smash{{\SetFigFont{5}{6.0}{\familydefault}{\mddefault}{\updefault}{\color[rgb]{0,0,0}$X_3$}%
}}}}
\put(3613,-1291){\makebox(0,0)[lb]{\smash{{\SetFigFont{5}{6.0}{\familydefault}{\mddefault}{\updefault}{\color[rgb]{0,0,0}$X_3$}%
}}}}
\put(4094,-1150){\makebox(0,0)[lb]{\smash{{\SetFigFont{5}{6.0}{\familydefault}{\mddefault}{\updefault}{\color[rgb]{0,0,0}$X_2$}%
}}}}
\put(2822,-103){\makebox(0,0)[lb]{\smash{{\SetFigFont{5}{6.0}{\familydefault}{\mddefault}{\updefault}{\color[rgb]{0,0,0}$X_1$}%
}}}}
\put(2822,-245){\makebox(0,0)[lb]{\smash{{\SetFigFont{5}{6.0}{\familydefault}{\mddefault}{\updefault}{\color[rgb]{0,0,0}$X_2$}%
}}}}
\put(4942,-952){\makebox(0,0)[lb]{\smash{{\SetFigFont{5}{6.0}{\familydefault}{\mddefault}{\updefault}{\color[rgb]{0,0,0}$X_1$}%
}}}}
\end{picture}%

%% file: goodnest1.pdf_t
\begin{picture}(0,0)%
\includegraphics{goodnest1.pdf}%
\end{picture}%
\setlength{\unitlength}{3947sp}%
\begingroup\makeatletter\ifx\SetFigFont\undefined%
\gdef\SetFigFont#1#2#3#4#5{%
  \reset@font\fontsize{#1}{#2pt}%
  \fontfamily{#3}\fontseries{#4}\fontshape{#5}%
  \selectfont}%
\fi\endgroup%
\begin{picture}(9852,7895)(1264,-8173)
\put(6376,-5611){\makebox(0,0)[lb]{\smash{{\SetFigFont{12}{14.4}{\familydefault}{\mddefault}{\updefault}{\color[rgb]{0,0,0}$\U'$}%
}}}}
\put(8851,-736){\makebox(0,0)[lb]{\smash{{\SetFigFont{12}{14.4}{\familydefault}{\mddefault}{\updefault}{\color[rgb]{0,0,0}$\V^n$}%
}}}}
\put(7051,-1111){\makebox(0,0)[lb]{\smash{{\SetFigFont{12}{14.4}{\familydefault}{\mddefault}{\updefault}{\color[rgb]{0,0,0}$\W^n$}%
}}}}
\put(8626,-1636){\makebox(0,0)[lb]{\smash{{\SetFigFont{12}{14.4}{\familydefault}{\mddefault}{\updefault}{\color[rgb]{0,0,0}$\mathcal{L}_{R_{\V^n}(c_0)}(\V^n)$}%
}}}}
\put(10276,-4261){\makebox(0,0)[lb]{\smash{{\SetFigFont{12}{14.4}{\familydefault}{\mddefault}{\updefault}{\color[rgb]{0,0,0}$\mathcal{L}_{R_{\V^n}(c_0)}(\W^n)$}%
}}}}
\put(11101,-2011){\makebox(0,0)[lb]{\smash{{\SetFigFont{12}{14.4}{\familydefault}{\mddefault}{\updefault}{\color[rgb]{0,0,0}$\U=\mathcal{L}_x(\W^n)$}%
}}}}
\put(7726,-1711){\makebox(0,0)[lb]{\smash{{\SetFigFont{12}{14.4}{\familydefault}{\mddefault}{\updefault}{\color[rgb]{0,0,0}$R_{\V^n}$}%
}}}}
\put(6226,-1486){\makebox(0,0)[lb]{\smash{{\SetFigFont{12}{14.4}{\familydefault}{\mddefault}{\updefault}{\color[rgb]{0,0,0}$\U'$}%
}}}}
\put(7876,-5611){\makebox(0,0)[lb]{\smash{{\SetFigFont{12}{14.4}{\familydefault}{\mddefault}{\updefault}{\color[rgb]{0,0,0}$R_{\W^n}$}%
}}}}
\end{picture}%

%% file: rigidityscheme1.pdf_t
\begin{picture}(0,0)%
\includegraphics{rigidityscheme1.pdf}%
\end{picture}%
\setlength{\unitlength}{3947sp}%
\begingroup\makeatletter\ifx\SetFigFont\undefined%
\gdef\SetFigFont#1#2#3#4#5{%
  \reset@font\fontsize{#1}{#2pt}%
  \fontfamily{#3}\fontseries{#4}\fontshape{#5}%
  \selectfont}%
\fi\endgroup%
\begin{picture}(15915,8874)(-11,-8023)
\put(8926,-6811){\makebox(0,0)[lb]{\smash{{\SetFigFont{12}{14.4}{\familydefault}{\mddefault}{\updefault}{\color[rgb]{0,0,0}Propositions~\ref{prop:real bounds persistent} and \ref{prop:real bounds reluctant}}%
}}}}
\put(1238,-1264){\makebox(0,0)[lb]{\smash{{\SetFigFont{12}{14.4}{\familydefault}{\mddefault}{\updefault}{\color[rgb]{0,0,0}Propostition \ref{prop:rigidity attracting}}%
}}}}
\put(1238,-1030){\makebox(0,0)[lb]{\smash{{\SetFigFont{12}{14.4}{\familydefault}{\mddefault}{\updefault}{\color[rgb]{0,0,0}hyperbolic attractors,}%
}}}}
\put(1238,-800){\makebox(0,0)[lb]{\smash{{\SetFigFont{12}{14.4}{\familydefault}{\mddefault}{\updefault}{\color[rgb]{0,0,0}Rigidity in basins of }%
}}}}
\put( 24,-2268){\makebox(0,0)[lb]{\smash{{\SetFigFont{12}{14.4}{\familydefault}{\mddefault}{\updefault}{\color[rgb]{0,0,0}Rigidy at infinitely renormalizable}%
}}}}
\put( 24,-2503){\makebox(0,0)[lb]{\smash{{\SetFigFont{12}{14.4}{\familydefault}{\mddefault}{\updefault}{\color[rgb]{0,0,0}critical points, Sections~\ref{subsec:qc rigidity - infinitely renormalizable} and \ref{subsec:C3 inf renorm rigidity}}%
}}}}
\put(8926,-2236){\makebox(0,0)[lb]{\smash{{\SetFigFont{12}{14.4}{\familydefault}{\mddefault}{\updefault}{\color[rgb]{0,0,0}Complex bounds, Theorem \ref{thm:box mapping persistent}}%
}}}}
\put(8926,-2461){\makebox(0,0)[lb]{\smash{{\SetFigFont{12}{14.4}{\familydefault}{\mddefault}{\updefault}{\color[rgb]{0,0,0}Real bounds, Proposition~\ref{prop:renormbounds}}%
}}}}
\put(4576,-961){\makebox(0,0)[lb]{\smash{{\SetFigFont{12}{14.4}{\familydefault}{\mddefault}{\updefault}{\color[rgb]{0,0,0}Local behaviour at periodic points, }%
}}}}
\put(4576,-1186){\makebox(0,0)[lb]{\smash{{\SetFigFont{12}{14.4}{\familydefault}{\mddefault}{\updefault}{\color[rgb]{0,0,0}Proposition \ref{prop:local dynamics}}%
}}}}
\put(4576,-2461){\makebox(0,0)[lb]{\smash{{\SetFigFont{12}{14.4}{\familydefault}{\mddefault}{\updefault}{\color[rgb]{0,0,0}asymptotically holomorphic,}%
}}}}
\put(4576,-2686){\makebox(0,0)[lb]{\smash{{\SetFigFont{12}{14.4}{\familydefault}{\mddefault}{\updefault}{\color[rgb]{0,0,0}polynomial-like mappings}%
}}}}
\put(4576,-2236){\makebox(0,0)[lb]{\smash{{\SetFigFont{12}{14.4}{\familydefault}{\mddefault}{\updefault}{\color[rgb]{0,0,0}A qc\textbackslash bg partition  for}%
}}}}
\put(4576,-2911){\makebox(0,0)[lb]{\smash{{\SetFigFont{12}{14.4}{\familydefault}{\mddefault}{\updefault}{\color[rgb]{0,0,0}Section \ref{subsec:C3 inf renorm rigidity}}%
}}}}
\put(8934,-5129){\makebox(0,0)[lb]{\smash{{\SetFigFont{12}{14.4}{\familydefault}{\mddefault}{\updefault}{\color[rgb]{0,0,0}Complex bounds,}%
}}}}
\put(8926,-5386){\makebox(0,0)[lb]{\smash{{\SetFigFont{12}{14.4}{\familydefault}{\mddefault}{\updefault}{\color[rgb]{0,0,0}Theorems \ref{thm:box mapping persistent}, \ref{thm:box mapping persistent infinite branches} and \ref{thm:box mapping reluctant}}%
}}}}
\put(8926,-5611){\makebox(0,0)[lb]{\smash{{\SetFigFont{12}{14.4}{\familydefault}{\mddefault}{\updefault}{\color[rgb]{0,0,0}and Proposition~\ref{prop:good bounds}}%
}}}}
\put(8926,-3886){\makebox(0,0)[lb]{\smash{{\SetFigFont{12}{14.4}{\familydefault}{\mddefault}{\updefault}{\color[rgb]{0,0,0}Complex bounds, Proposition~\ref{prop:modified delta nice}}%
}}}}
\put(8926,-4141){\makebox(0,0)[lb]{\smash{{\SetFigFont{12}{14.4}{\familydefault}{\mddefault}{\updefault}{\color[rgb]{0,0,0}and quasidisks, Proposition~\ref{prop:quasidisks}}%
}}}}
\put(4576,-3886){\makebox(0,0)[lb]{\smash{{\SetFigFont{12}{14.4}{\familydefault}{\mddefault}{\updefault}{\color[rgb]{0,0,0}A qc\textbackslash bg partition of a central}%
}}}}
\put(4576,-4141){\makebox(0,0)[lb]{\smash{{\SetFigFont{12}{14.4}{\familydefault}{\mddefault}{\updefault}{\color[rgb]{0,0,0}cascade, Theorem~\ref{thm:central cascades}}%
}}}}
\put(4576,-6436){\makebox(0,0)[lb]{\smash{{\SetFigFont{12}{14.4}{\familydefault}{\mddefault}{\updefault}{\color[rgb]{0,0,0}Existence of  QC external conjugacies,}%
}}}}
\put(4576,-6661){\makebox(0,0)[lb]{\smash{{\SetFigFont{12}{14.4}{\familydefault}{\mddefault}{\updefault}{\color[rgb]{0,0,0}Theorem \ref{thm:external conjugacy}}%
}}}}
\put(4576,-7711){\makebox(0,0)[lb]{\smash{{\SetFigFont{12}{14.4}{\familydefault}{\mddefault}{\updefault}{\color[rgb]{0,0,0}Touching box mappings,}%
}}}}
\put(4576,-7936){\makebox(0,0)[lb]{\smash{{\SetFigFont{12}{14.4}{\familydefault}{\mddefault}{\updefault}{\color[rgb]{0,0,0}Theorem \ref{thm:touching box map}}%
}}}}
\put(15751,389){\makebox(0,0)[lb]{\smash{{\SetFigFont{12}{14.4}{\familydefault}{\mddefault}{\updefault}{\color[rgb]{0,0,0}T}%
}}}}
\put( 76,347){\makebox(0,0)[lb]{\smash{{\SetFigFont{12}{14.4}{\familydefault}{\mddefault}{\updefault}{\color[rgb]{0,0,0}basins of hyperbolic attractors, Theorem \ref{thm:qc rigidity away from critical points}}%
}}}}
\put( 76,578){\makebox(0,0)[lb]{\smash{{\SetFigFont{12}{14.4}{\familydefault}{\mddefault}{\updefault}{\color[rgb]{0,0,0}Rigidty away from $\mathrm{Crit}(f)$ and }%
}}}}
\put(4576,614){\makebox(0,0)[lb]{\smash{{\SetFigFont{12}{14.4}{\familydefault}{\mddefault}{\updefault}{\color[rgb]{0,0,0}Touching box mappings,}%
}}}}
\put(4576,389){\makebox(0,0)[lb]{\smash{{\SetFigFont{12}{14.4}{\familydefault}{\mddefault}{\updefault}{\color[rgb]{0,0,0}Theorem \ref{thm:touching box map}}%
}}}}
\put( 84,-4066){\makebox(0,0)[lb]{\smash{{\SetFigFont{12}{14.4}{\familydefault}{\mddefault}{\updefault}{\color[rgb]{0,0,0}Rigidity at finitely renormalizable,}%
}}}}
\put( 84,-4290){\makebox(0,0)[lb]{\smash{{\SetFigFont{12}{14.4}{\familydefault}{\mddefault}{\updefault}{\color[rgb]{0,0,0}persistently recurrent critical points,}%
}}}}
\put( 76,-4561){\makebox(0,0)[lb]{\smash{{\SetFigFont{12}{14.4}{\familydefault}{\mddefault}{\updefault}{\color[rgb]{0,0,0}Section \ref{subsec:qc rigidity-persistent}}%
}}}}
\put( 76,-5911){\makebox(0,0)[lb]{\smash{{\SetFigFont{12}{14.4}{\familydefault}{\mddefault}{\updefault}{\color[rgb]{0,0,0}Rigidity at reluctantly recurrent}%
}}}}
\put( 76,-6136){\makebox(0,0)[lb]{\smash{{\SetFigFont{12}{14.4}{\familydefault}{\mddefault}{\updefault}{\color[rgb]{0,0,0}critical points, Section \ref{subsec:qc rigidity-reluctant}}%
}}}}
\put(8926,-6586){\makebox(0,0)[lb]{\smash{{\SetFigFont{12}{14.4}{\familydefault}{\mddefault}{\updefault}{\color[rgb]{0,0,0}Real bounds,}%
}}}}
\put(4575,-5150){\makebox(0,0)[lb]{\smash{{\SetFigFont{12}{14.4}{\familydefault}{\mddefault}{\updefault}{\color[rgb]{0,0,0}Rigidity for complex box mappings,}%
}}}}
\put(4575,-5440){\makebox(0,0)[lb]{\smash{{\SetFigFont{12}{14.4}{\familydefault}{\mddefault}{\updefault}{\color[rgb]{0,0,0}Theorems \ref{thm: qc rigidity of box mappings} and \ref{thm:smooth box conj}}%
}}}}
\end{picture}%

%% file: circlemap.tex
\begin{figure}[htb] \hfil \beginpicture
\dimen0=4 cm
\setcoordinatesystem units <\dimen0,\dimen0>
\setplotarea x from 0 to 1, y from 0 to 1
\axis bottom invisible ticks
 withvalues 0 1 / at 0  1 / /
\axis left   invisible ticks
  withvalues 0 1 / at 0  1 / /
\grid 1 1
\setlinear \plot
0 0 1 1 /  \plot
0.000  0.300
0.002  0.305
0.004  0.310
0.006  0.314
0.008  0.319
0.010  0.324
0.012  0.329
0.014  0.334
0.016  0.338
0.018  0.343
0.020  0.348
0.022  0.353
0.024  0.357
0.026  0.362
0.028  0.367
0.030  0.372
0.032  0.376
0.034  0.381
0.036  0.386
0.038  0.391
0.040  0.395
0.042  0.400
0.044  0.405
0.046  0.410
0.048  0.414
0.050  0.419
0.052  0.424
0.054  0.428
0.056  0.433
0.058  0.437
0.060  0.442
0.062  0.447
0.064  0.451
0.066  0.456
0.068  0.460
0.070  0.465
0.072  0.469
0.074  0.474
0.076  0.478
0.078  0.483
0.080  0.487
0.082  0.492
0.084  0.496
0.086  0.501
0.088  0.505
0.090  0.509
0.092  0.514
0.094  0.518
0.096  0.522
0.098  0.527
0.100  0.531
0.102  0.535
0.104  0.539
0.106  0.544
0.108  0.548
0.110  0.552
0.112  0.556
0.114  0.560
0.116  0.564
0.118  0.568
0.120  0.573
0.122  0.577
0.124  0.581
0.126  0.585
0.128  0.588
0.130  0.592
0.132  0.596
0.134  0.600
0.136  0.604
0.138  0.608
0.140  0.612
0.142  0.615
0.144  0.619
0.146  0.623
0.148  0.627
0.150  0.630
0.152  0.634
0.154  0.637
0.156  0.641
0.158  0.645
0.160  0.648
0.162  0.652
0.164  0.655
0.166  0.658
0.168  0.662
0.170  0.665
0.172  0.669
0.174  0.672
0.176  0.675
0.178  0.678
0.180  0.682
0.182  0.685
0.184  0.688
0.186  0.691
0.188  0.694
0.190  0.697
0.192  0.700
0.194  0.703
0.196  0.706
0.198  0.709
0.200  0.712
0.202  0.715
0.204  0.718
0.206  0.720
0.208  0.723
0.210  0.726
0.212  0.728
0.214  0.731
0.216  0.734
0.218  0.736
0.220  0.739
0.222  0.741
0.224  0.744
0.226  0.746
0.228  0.749
0.230  0.751
0.232  0.753
0.234  0.756
0.236  0.758
0.238  0.760
0.240  0.762
0.242  0.765
0.244  0.767
0.246  0.769
0.248  0.771
0.250  0.773
0.252  0.775
0.254  0.777
0.256  0.779
0.258  0.781
0.260  0.782
0.262  0.784
0.264  0.786
0.266  0.788
0.268  0.789
0.270  0.791
0.272  0.793
0.274  0.794
0.276  0.796
0.278  0.797
0.280  0.799
0.282  0.800
0.284  0.802
0.286  0.803
0.288  0.804
0.290  0.806
0.292  0.807
0.294  0.808
0.296  0.810
0.298  0.811
0.300  0.812
0.302  0.813
0.304  0.814
0.306  0.815
0.308  0.816
0.310  0.817
0.312  0.818
0.314  0.819
0.316  0.820
0.318  0.821
0.320  0.822
0.322  0.822
0.324  0.823
0.326  0.824
0.328  0.825
0.330  0.825
0.332  0.826
0.334  0.826
0.336  0.827
0.338  0.828
0.340  0.828
0.342  0.829
0.344  0.829
0.346  0.829
0.348  0.830
0.350  0.830
0.352  0.831
0.354  0.831
0.356  0.831
0.358  0.831
0.360  0.832
0.362  0.832
0.364  0.832
0.366  0.832
0.368  0.832
0.370  0.832
0.372  0.832
0.374  0.833
0.376  0.833
0.378  0.833
0.380  0.833
0.382  0.832
0.384  0.832
0.386  0.832
0.388  0.832
0.390  0.832
0.392  0.832
0.394  0.832
0.396  0.831
0.398  0.831
0.400  0.831
0.402  0.831
0.404  0.830
0.406  0.830
0.408  0.830
0.410  0.829
0.412  0.829
0.414  0.829
0.416  0.828
0.418  0.828
0.420  0.827
0.422  0.827
0.424  0.826
0.426  0.826
0.428  0.825
0.430  0.825
0.432  0.824
0.434  0.824
0.436  0.823
0.438  0.823
0.440  0.822
0.442  0.821
0.444  0.821
0.446  0.820
0.448  0.820
0.450  0.819
0.452  0.818
0.454  0.818
0.456  0.817
0.458  0.816
0.460  0.815
0.462  0.815
0.464  0.814
0.466  0.813
0.468  0.812
0.470  0.812
0.472  0.811
0.474  0.810
0.476  0.809
0.478  0.809
0.480  0.808
0.482  0.807
0.484  0.806
0.486  0.806
0.488  0.805
0.490  0.804
0.492  0.803
0.494  0.802
0.496  0.802
0.498  0.801
0.500  0.800
0.502  0.799
0.504  0.798
0.506  0.798
0.508  0.797
0.510  0.796
0.512  0.795
0.514  0.794
0.516  0.794
0.518  0.793
0.520  0.792
0.522  0.791
0.524  0.791
0.526  0.790
0.528  0.789
0.530  0.788
0.532  0.788
0.534  0.787
0.536  0.786
0.538  0.785
0.540  0.785
0.542  0.784
0.544  0.783
0.546  0.782
0.548  0.782
0.550  0.781
0.552  0.780
0.554  0.780
0.556  0.779
0.558  0.779
0.560  0.778
0.562  0.777
0.564  0.777
0.566  0.776
0.568  0.776
0.570  0.775
0.572  0.775
0.574  0.774
0.576  0.774
0.578  0.773
0.580  0.773
0.582  0.772
0.584  0.772
0.586  0.771
0.588  0.771
0.590  0.771
0.592  0.770
0.594  0.770
0.596  0.770
0.598  0.769
0.600  0.769
0.602  0.769
0.604  0.769
0.606  0.768
0.608  0.768
0.610  0.768
0.612  0.768
0.614  0.768
0.616  0.768
0.618  0.768
0.620  0.767
0.622  0.767
0.624  0.767
0.626  0.767
0.628  0.768
0.630  0.768
0.632  0.768
0.634  0.768
0.636  0.768
0.638  0.768
0.640  0.768
0.642  0.769
0.644  0.769
0.646  0.769
0.648  0.769
0.650  0.770
0.652  0.770
0.654  0.771
0.656  0.771
0.658  0.771
0.660  0.772
0.662  0.772
0.664  0.773
0.666  0.774
0.668  0.774
0.670  0.775
0.672  0.775
0.674  0.776
0.676  0.777
0.678  0.778
0.680  0.778
0.682  0.779
0.684  0.780
0.686  0.781
0.688  0.782
0.690  0.783
0.692  0.784
0.694  0.785
0.696  0.786
0.698  0.787
0.700  0.788
0.702  0.789
0.704  0.790
0.706  0.792
0.708  0.793
0.710  0.794
0.712  0.796
0.714  0.797
0.716  0.798
0.718  0.800
0.720  0.801
0.722  0.803
0.724  0.804
0.726  0.806
0.728  0.807
0.730  0.809
0.732  0.811
0.734  0.812
0.736  0.814
0.738  0.816
0.740  0.818
0.742  0.819
0.744  0.821
0.746  0.823
0.748  0.825
0.750  0.827
0.752  0.829
0.754  0.831
0.756  0.833
0.758  0.835
0.760  0.838
0.762  0.840
0.764  0.842
0.766  0.844
0.768  0.847
0.770  0.849
0.772  0.851
0.774  0.854
0.776  0.856
0.778  0.859
0.780  0.861
0.782  0.864
0.784  0.866
0.786  0.869
0.788  0.872
0.790  0.874
0.792  0.877
0.794  0.880
0.796  0.882
0.798  0.885
0.800  0.888
0.802  0.891
0.804  0.894
0.806  0.897
0.808  0.900
0.810  0.903
0.812  0.906
0.814  0.909
0.816  0.912
0.818  0.915
0.820  0.918
0.822  0.922
0.824  0.925
0.826  0.928
0.828  0.931
0.830  0.935
0.832  0.938
0.834  0.942
0.836  0.945
0.838  0.948
0.840  0.952
0.842  0.955
0.844  0.959
0.846  0.963
0.848  0.966
0.850  0.970
0.852  0.973
0.854  0.977
0.856  0.981
0.858  0.985
0.860  0.988
0.862  0.992
0.864  0.996
0.866  1.000 /
\plot
0.867  0.000
0.868  0.004
0.870  0.008
0.872  0.012
0.874  0.015
0.876  0.019
0.878  0.023
0.880  0.027
0.882  0.032
0.884  0.036
0.886  0.040
0.888  0.044
0.890  0.048
0.892  0.052
0.894  0.056
0.896  0.061
0.898  0.065
0.900  0.069
0.902  0.073
0.904  0.078
0.906  0.082
0.908  0.086
0.910  0.091
0.912  0.095
0.914  0.099
0.916  0.104
0.918  0.108
0.920  0.113
0.922  0.117
0.924  0.122
0.926  0.126
0.928  0.131
0.930  0.135
0.932  0.140
0.934  0.144
0.936  0.149
0.938  0.153
0.940  0.158
0.942  0.163
0.944  0.167
0.946  0.172
0.948  0.176
0.950  0.181
0.952  0.186
0.954  0.190
0.956  0.195
0.958  0.200
0.960  0.205
0.962  0.209
0.964  0.214
0.966  0.219
0.968  0.224
0.970  0.228
0.972  0.233
0.974  0.238
0.976  0.243
0.978  0.247
0.980  0.252
0.982  0.257
0.984  0.262
0.986  0.266
0.988  0.271
0.990  0.276
0.992  0.281
0.994  0.286
0.996  0.290
0.998  0.295
1.000  0.300
/
\setdots <2pt>
\plot 0.282 0.8 0.718 0.8 /
\setdashes \plot 0.866  1.000 0.866  0.000 /
\put  {$f_t$} at 0.38 0.76
\put {$f$} at 0.15 0.531
\endpicture
\caption[ ]{A map $f\colon S^1\to S^1$ of degree one which is not monotone, and the
corresponding monotone maps $f_t$ with plateaus.\label{fig:circlemap}}
\end{figure}

%% file: cen1.pdf_t
\begin{picture}(0,0)%
\includegraphics{cen1.pdf}%
\end{picture}%
\setlength{\unitlength}{3947sp}%
\begingroup\makeatletter\ifx\SetFigFont\undefined%
\gdef\SetFigFont#1#2#3#4#5{%
  \reset@font\fontsize{#1}{#2pt}%
  \fontfamily{#3}\fontseries{#4}\fontshape{#5}%
  \selectfont}%
\fi\endgroup%
\begin{picture}(7902,1452)(286,-280)
\put(7726,539){\makebox(0,0)[lb]{\smash{{\SetFigFont{12}{14.4}{\familydefault}{\mddefault}{\updefault}{\color[rgb]{0,0,0}$I(c)\owns c$}%
}}}}
\put(5251,539){\makebox(0,0)[lb]{\smash{{\SetFigFont{12}{14.4}{\familydefault}{\mddefault}{\updefault}{\color[rgb]{0,0,0}$\hat{J}(c)\owns f(c)$}%
}}}}
\put(301,539){\makebox(0,0)[lb]{\smash{{\SetFigFont{12}{14.4}{\familydefault}{\mddefault}{\updefault}{\color[rgb]{0,0,0}$I(c_0)\owns c_0$}%
}}}}
\put(3076,-136){\makebox(0,0)[lb]{\smash{{\SetFigFont{12}{14.4}{\familydefault}{\mddefault}{\updefault}{\color[rgb]{0,0,0}$f(I(c))$}%
}}}}
\put(2176,989){\makebox(0,0)[lb]{\smash{{\SetFigFont{12}{14.4}{\familydefault}{\mddefault}{\updefault}{\color[rgb]{0,0,0}$f^k$ - the first landing map}%
}}}}
\put(6901,-211){\makebox(0,0)[lb]{\smash{{\SetFigFont{12}{14.4}{\familydefault}{\mddefault}{\updefault}{\color[rgb]{0,0,0}$f$}%
}}}}
\end{picture}%

%% file: cen2.pdf_t
\begin{picture}(0,0)%
\includegraphics{cen2.pdf}%
\end{picture}%
\setlength{\unitlength}{3947sp}%
\begingroup\makeatletter\ifx\SetFigFont\undefined%
\gdef\SetFigFont#1#2#3#4#5{%
  \reset@font\fontsize{#1}{#2pt}%
  \fontfamily{#3}\fontseries{#4}\fontshape{#5}%
  \selectfont}%
\fi\endgroup%
\begin{picture}(8277,1357)(961,-2005)
\put(2926,-1036){\makebox(0,0)[lb]{\smash{{\SetFigFont{12}{14.4}{\familydefault}{\mddefault}{\updefault}{\color[rgb]{0,0,0}$f^k$ - the first landing map}%
}}}}
\put(976,-1111){\makebox(0,0)[lb]{\smash{{\SetFigFont{12}{14.4}{\familydefault}{\mddefault}{\updefault}{\color[rgb]{0,0,0}$I(c_0)\owns c_0$}%
}}}}
\put(5251,-1861){\makebox(0,0)[lb]{\smash{{\SetFigFont{12}{14.4}{\familydefault}{\mddefault}{\updefault}{\color[rgb]{0,0,0}$f(I(c))$}%
}}}}
\put(6226,-1186){\makebox(0,0)[lb]{\smash{{\SetFigFont{12}{14.4}{\familydefault}{\mddefault}{\updefault}{\color[rgb]{0,0,0}$\hat{J}(c)\owns f(c)$}%
}}}}
\put(8851,-1186){\makebox(0,0)[lb]{\smash{{\SetFigFont{12}{14.4}{\familydefault}{\mddefault}{\updefault}{\color[rgb]{0,0,0}$I(c)\owns c$}%
}}}}
\put(8176,-1936){\makebox(0,0)[lb]{\smash{{\SetFigFont{12}{14.4}{\familydefault}{\mddefault}{\updefault}{\color[rgb]{0,0,0}$f$}%
}}}}
\end{picture}%

%% file: extmap.pdf_t
\begin{picture}(0,0)%
\includegraphics{extmap.pdf}%
\end{picture}%
\setlength{\unitlength}{3947sp}%
\begingroup\makeatletter\ifx\SetFigFont\undefined%
\gdef\SetFigFont#1#2#3#4#5{%
  \reset@font\fontsize{#1}{#2pt}%
  \fontfamily{#3}\fontseries{#4}\fontshape{#5}%
  \selectfont}%
\fi\endgroup%
\begin{picture}(7820,925)(518,-1018)
\put(4426,-211){\makebox(0,0)[lb]{\smash{{\SetFigFont{7}{7.2}{\familydefault}{\mddefault}{\updefault}{\color[rgb]{0,0,0}$H_0$}%
}}}}
\put(2440,-725){\makebox(0,0)[lb]{\smash{{\SetFigFont{7}{6.0}{\familydefault}{\mddefault}{\updefault}{\color[rgb]{0,0,0}$\bm{\mathcal{U}}$}%
}}}}
\put(2931,-157){\makebox(0,0)[lb]{\smash{{\SetFigFont{7}{6.0}{\familydefault}{\mddefault}{\updefault}{\color[rgb]{0,0,0}$\bm{\mathcal{V}}$}%
}}}}
\put(7395,-208){\makebox(0,0)[lb]{\smash{{\SetFigFont{7}{6.0}{\familydefault}{\mddefault}{\updefault}{\color[rgb]{0,0,0}$\tilde{\bm{\mathcal{V}}}$}%
}}}}
\put(6904,-788){\makebox(0,0)[lb]{\smash{{\SetFigFont{7}{6.0}{\familydefault}{\mddefault}{\updefault}{\color[rgb]{0,0,0}$\tilde{\bm{\mathcal{U}}}$}%
}}}}
\end{picture}%

%% file: modifiednest.pdf_t
\begin{picture}(0,0)%
\includegraphics{modifiednest.pdf}%
\end{picture}%
\setlength{\unitlength}{3947sp}%
\begingroup\makeatletter\ifx\SetFigFont\undefined%
\gdef\SetFigFont#1#2#3#4#5{%
  \reset@font\fontsize{#1}{#2pt}%
  \fontfamily{#3}\fontseries{#4}\fontshape{#5}%
  \selectfont}%
\fi\endgroup%
\begin{picture}(5359,3902)(1458,-3325)
\put(5489,-711){\makebox(0,0)[lb]{\smash{{\SetFigFont{7}{8.4}{\familydefault}{\mddefault}{\updefault}{\color[rgb]{0,0,0}Excluded puzzle pieces}%
}}}}
\put(5489,-861){\makebox(0,0)[lb]{\smash{{\SetFigFont{7}{8.4}{\familydefault}{\mddefault}{\updefault}{\color[rgb]{0,0,0}deep in a central cascade}%
}}}}
\put(5489,-1011){\makebox(0,0)[lb]{\smash{{\SetFigFont{7}{8.4}{\familydefault}{\mddefault}{\updefault}{\color[rgb]{0,0,0}and their pullbacks}%
}}}}
\end{picture}%

%% file: goodnest2.pdf_t
\begin{picture}(0,0)%
\includegraphics{goodnest2.pdf}%
\end{picture}%
\setlength{\unitlength}{3947sp}%
\begingroup\makeatletter\ifx\SetFigFont\undefined%
\gdef\SetFigFont#1#2#3#4#5{%
  \reset@font\fontsize{#1}{#2pt}%
  \fontfamily{#3}\fontseries{#4}\fontshape{#5}%
  \selectfont}%
\fi\endgroup%
\begin{picture}(6764,2232)(518,-1311)
\put(4862,464){\makebox(0,0)[lb]{\smash{{\SetFigFont{7}{8.4}{\familydefault}{\mddefault}{\updefault}{\color[rgb]{0,0,0}$\W^0$}%
}}}}
\put(6131,709){\makebox(0,0)[lb]{\smash{{\SetFigFont{7}{8.4}{\familydefault}{\mddefault}{\updefault}{\color[rgb]{0,0,0}$\V^0$}%
}}}}
\end{picture}%

%% file: cascrep.pdf_t
\begin{picture}(0,0)%
\includegraphics{cascrep.pdf}%
\end{picture}%
\setlength{\unitlength}{3947sp}%
\begingroup\makeatletter\ifx\SetFigFont\undefined%
\gdef\SetFigFont#1#2#3#4#5{%
  \reset@font\fontsize{#1}{#2pt}%
  \fontfamily{#3}\fontseries{#4}\fontshape{#5}%
  \selectfont}%
\fi\endgroup%
\begin{picture}(21129,9681)(-28,-8863)
\put(6076,-2461){\makebox(0,0)[lb]{\smash{{\SetFigFont{17}{20.4}{\familydefault}{\mddefault}{\updefault}{\color[rgb]{0,0,0}$\bm T^{N+1}$}%
}}}}
\put(19801,-886){\makebox(0,0)[lb]{\smash{{\SetFigFont{17}{20.4}{\familydefault}{\mddefault}{\updefault}{\color[rgb]{0,0,0}$\V^0$}%
}}}}
\put(16951,-3136){\makebox(0,0)[lb]{\smash{{\SetFigFont{17}{20.4}{\familydefault}{\mddefault}{\updefault}{\color[rgb]{0,0,0}$\V^{1}$}%
}}}}
\put(5476,-3061){\makebox(0,0)[lb]{\smash{{\SetFigFont{17}{20.4}{\familydefault}{\mddefault}{\updefault}{\color[rgb]{0,0,0}$\V^1$}%
}}}}
\put(6376,-2011){\makebox(0,0)[lb]{\smash{{\SetFigFont{17}{20.4}{\familydefault}{\mddefault}{\updefault}{\color[rgb]{0,0,0}$\bm T^N$}%
}}}}
\put(6751,-1636){\makebox(0,0)[lb]{\smash{{\SetFigFont{17}{20.4}{\familydefault}{\mddefault}{\updefault}{\color[rgb]{0,0,0}$\bm T^2$}%
}}}}
\put(7051,-1261){\makebox(0,0)[lb]{\smash{{\SetFigFont{17}{20.4}{\familydefault}{\mddefault}{\updefault}{\color[rgb]{0,0,0}$\bm T^1$}%
}}}}
\put(7426,-811){\makebox(0,0)[lb]{\smash{{\SetFigFont{17}{20.4}{\familydefault}{\mddefault}{\updefault}{\color[rgb]{0,0,0}$\W^0=\bm T^0$}%
}}}}
\put(7876,-286){\makebox(0,0)[lb]{\smash{{\SetFigFont{17}{20.4}{\familydefault}{\mddefault}{\updefault}{\color[rgb]{0,0,0}$\V^0$}%
}}}}
\put(5176,-3436){\makebox(0,0)[lb]{\smash{{\SetFigFont{17}{20.4}{\familydefault}{\mddefault}{\updefault}{\color[rgb]{0,0,0}$\bm T^{N+2}$}%
}}}}
\end{picture}%

%% file: complexdomains.pdf_t
\begin{picture}(0,0)%
\includegraphics{complexdomains.pdf}%
\end{picture}%
\setlength{\unitlength}{3947sp}%
\begingroup\makeatletter\ifx\SetFigFont\undefined%
\gdef\SetFigFont#1#2#3#4#5{%
  \reset@font\fontsize{#1}{#2pt}%
  \fontfamily{#3}\fontseries{#4}\fontshape{#5}%
  \selectfont}%
\fi\endgroup%
\begin{picture}(9624,9614)(-11,-8768)
\put(3376,-7861){\makebox(0,0)[lb]{\smash{{\SetFigFont{17}{20.4}{\familydefault}{\mddefault}{\updefault}{\color[rgb]{0,0,0}$\bm A^0_0$}%
}}}}
\put(3001,-436){\makebox(0,0)[lb]{\smash{{\SetFigFont{17}{20.4}{\familydefault}{\mddefault}{\updefault}{\color[rgb]{0,0,0}$\bm A^0_0$}%
}}}}
\put(1876,-2011){\makebox(0,0)[lb]{\smash{{\SetFigFont{17}{20.4}{\familydefault}{\mddefault}{\updefault}{\color[rgb]{0,0,0}$\bm B_0$}%
}}}}
\put(3826,-3136){\makebox(0,0)[lb]{\smash{{\SetFigFont{17}{20.4}{\familydefault}{\mddefault}{\updefault}{\color[rgb]{0,0,0}$\bm T^1$}%
}}}}
\put(4801,-4261){\makebox(0,0)[lb]{\smash{{\SetFigFont{17}{20.4}{\familydefault}{\mddefault}{\updefault}{\color[rgb]{0,0,0}$c$}%
}}}}
\put(7951,-286){\makebox(0,0)[lb]{\smash{{\SetFigFont{17}{20.4}{\familydefault}{\mddefault}{\updefault}{\color[rgb]{0,0,0}$\bm T^0$}%
}}}}
\end{picture}%

%% file: init.pdf_t
\begin{picture}(0,0)%
\includegraphics{init.pdf}%
\end{picture}%
\setlength{\unitlength}{3947sp}%
\begingroup\makeatletter\ifx\SetFigFont\undefined%
\gdef\SetFigFont#1#2#3#4#5{%
  \reset@font\fontsize{#1}{#2pt}%
  \fontfamily{#3}\fontseries{#4}\fontshape{#5}%
  \selectfont}%
\fi\endgroup%
\begin{picture}(6398,6388)(589,-4943)
\put(4759,-1542){\makebox(0,0)[lb]{\smash{{\SetFigFont{5}{6.0}{\familydefault}{\mddefault}{\updefault}{\color[rgb]{0,0,0}$\gamma_5^1$}%
}}}}
\put(4762,-22){\makebox(0,0)[lb]{\smash{{\SetFigFont{7}{8.4}{\familydefault}{\mddefault}{\updefault}{\color[rgb]{0,0,0}$\bm J^2$}%
}}}}
\put(5337,244){\makebox(0,0)[lb]{\smash{{\SetFigFont{7}{8.4}{\familydefault}{\mddefault}{\updefault}{\color[rgb]{0,0,0}$\bm J^1$}%
}}}}
\put(6090,509){\makebox(0,0)[lb]{\smash{{\SetFigFont{7}{8.4}{\familydefault}{\mddefault}{\updefault}{\color[rgb]{0,0,0}$\bm J^0$}%
}}}}
\put(6488,-1660){\makebox(0,0)[lb]{\smash{{\SetFigFont{7}{8.4}{\familydefault}{\mddefault}{\updefault}{\color[rgb]{0,0,0}$e_0$}%
}}}}
\put(3567,421){\makebox(0,0)[lb]{\smash{{\SetFigFont{7}{8.4}{\familydefault}{\mddefault}{\updefault}{\color[rgb]{0,0,0}$e_1$}%
}}}}
\put(4388,-923){\makebox(0,0)[lb]{\smash{{\SetFigFont{5}{6.0}{\familydefault}{\mddefault}{\updefault}{\color[rgb]{0,0,0}$\gamma_1$}%
}}}}
\put(4787,-1387){\makebox(0,0)[lb]{\smash{{\SetFigFont{5}{6.0}{\familydefault}{\mddefault}{\updefault}{\color[rgb]{0,0,0}$e_4^3$}%
}}}}
\put(2857,-914){\makebox(0,0)[lb]{\smash{{\SetFigFont{5}{6.0}{\familydefault}{\mddefault}{\updefault}{\color[rgb]{0,0,0}$\gamma_2$}%
}}}}
\put(2882,-2529){\makebox(0,0)[lb]{\smash{{\SetFigFont{5}{6.0}{\familydefault}{\mddefault}{\updefault}{\color[rgb]{0,0,0}$\gamma_4$}%
}}}}
\put(4472,-2481){\makebox(0,0)[lb]{\smash{{\SetFigFont{5}{6.0}{\familydefault}{\mddefault}{\updefault}{\color[rgb]{0,0,0}$\gamma_3$}%
}}}}
\put(4808,-361){\makebox(0,0)[lb]{\smash{{\SetFigFont{5}{6.0}{\familydefault}{\mddefault}{\updefault}{\color[rgb]{0,0,0}$e_2^1$}%
}}}}
\put(4905,-577){\makebox(0,0)[lb]{\smash{{\SetFigFont{5}{6.0}{\familydefault}{\mddefault}{\updefault}{\color[rgb]{0,0,0}$\gamma_3^1$}%
}}}}
\put(4905,-890){\makebox(0,0)[lb]{\smash{{\SetFigFont{5}{6.0}{\familydefault}{\mddefault}{\updefault}{\color[rgb]{0,0,0}$e_3^1$}%
}}}}
\put(4905,-1109){\makebox(0,0)[lb]{\smash{{\SetFigFont{5}{6.0}{\familydefault}{\mddefault}{\updefault}{\color[rgb]{0,0,0}$\gamma_4^1$}%
}}}}
\end{picture}%

%% file: returns.pdf_t
\begin{picture}(0,0)%
\includegraphics{returns.pdf}%
\end{picture}%
\setlength{\unitlength}{3947sp}%
\begingroup\makeatletter\ifx\SetFigFont\undefined%
\gdef\SetFigFont#1#2#3#4#5{%
  \reset@font\fontsize{#1}{#2pt}%
  \fontfamily{#3}\fontseries{#4}\fontshape{#5}%
  \selectfont}%
\fi\endgroup%
\begin{picture}(3774,3774)(1189,-4123)
\put(1276,-2611){\makebox(0,0)[lb]{\smash{{\SetFigFont{6}{7.2}{\familydefault}{\mddefault}{\updefault}{\color[rgb]{0,0,0}$(4)$}%
}}}}
\put(1276,-511){\makebox(0,0)[lb]{\smash{{\SetFigFont{6}{7.2}{\familydefault}{\mddefault}{\updefault}{\color[rgb]{0,0,0}$(1)$}%
}}}}
\put(3376,-511){\makebox(0,0)[lb]{\smash{{\SetFigFont{6}{7.2}{\familydefault}{\mddefault}{\updefault}{\color[rgb]{0,0,0}$(2)$}%
}}}}
\put(3376,-2611){\makebox(0,0)[lb]{\smash{{\SetFigFont{6}{7.2}{\familydefault}{\mddefault}{\updefault}{\color[rgb]{0,0,0}$(3)$}%
}}}}
\end{picture}%

%% file: cvshieldsmall.pdf_t
\begin{picture}(0,0)%
\includegraphics{cvshieldsmall.pdf}%
\end{picture}%
\setlength{\unitlength}{3947sp}%
\begingroup\makeatletter\ifx\SetFigFont\undefined%
\gdef\SetFigFont#1#2#3#4#5{%
  \reset@font\fontsize{#1}{#2pt}%
  \fontfamily{#3}\fontseries{#4}\fontshape{#5}%
  \selectfont}%
\fi\endgroup%
\begin{picture}(5417,3976)(664,-2319)
\put(1051,1364){\makebox(0,0)[lb]{\smash{{\SetFigFont{5}{6.0}{\familydefault}{\mddefault}{\updefault}{\color[rgb]{0,0,0}$\J^{N+1}$}%
}}}}
\put(1925,613){\makebox(0,0)[lb]{\smash{{\SetFigFont{5}{6.0}{\familydefault}{\mddefault}{\updefault}{\color[rgb]{0,0,0}$e_N'$}%
}}}}
\put(5076,869){\makebox(0,0)[lb]{\smash{{\SetFigFont{5}{6.0}{\familydefault}{\mddefault}{\updefault}{\color[rgb]{0,0,0}$R^2$}%
}}}}
\put(3628,-352){\makebox(0,0)[lb]{\smash{{\SetFigFont{5}{6.0}{\familydefault}{\mddefault}{\updefault}{\color[rgb]{0,0,0}$R^N$}%
}}}}
\put(4253,812){\makebox(0,0)[lb]{\smash{{\SetFigFont{5}{6.0}{\familydefault}{\mddefault}{\updefault}{\color[rgb]{0,0,0}$e_2'$}%
}}}}
\put(959,-267){\makebox(0,0)[lb]{\smash{{\SetFigFont{5}{6.0}{\familydefault}{\mddefault}{\updefault}{\color[rgb]{0,0,0}$\bm B_1'$}%
}}}}
\put(2067,1579){\makebox(0,0)[lb]{\smash{{\SetFigFont{5}{6.0}{\familydefault}{\mddefault}{\updefault}{\color[rgb]{0,0,0}$\bm M_1'$}%
}}}}
\put(5388,1408){\makebox(0,0)[lb]{\smash{{\SetFigFont{5}{6.0}{\familydefault}{\mddefault}{\updefault}{\color[rgb]{0,0,0}$\J^0$}%
}}}}
\put(5388,1494){\makebox(0,0)[lb]{\smash{{\SetFigFont{5}{6.0}{\familydefault}{\mddefault}{\updefault}{\color[rgb]{0,0,0}$\II_n$=}%
}}}}
\put(5644,1408){\makebox(0,0)[lb]{\smash{{\SetFigFont{5}{6.0}{\familydefault}{\mddefault}{\updefault}{\color[rgb]{0,0,0}$\II_n^+$}%
}}}}
\put(5190,1380){\makebox(0,0)[lb]{\smash{{\SetFigFont{5}{6.0}{\familydefault}{\mddefault}{\updefault}{\color[rgb]{0,0,0}$\bm I^-_n$}%
}}}}
\put(976,1214){\makebox(0,0)[lb]{\smash{{\SetFigFont{5}{6.0}{\familydefault}{\mddefault}{\updefault}{\color[rgb]{0,0,0}$e_{N+1}'$}%
}}}}
\end{picture}%

%% file: shield2.pdf_t
\begin{picture}(0,0)%
\includegraphics{shield2.pdf}%
\end{picture}%
\setlength{\unitlength}{3947sp}%
\begingroup\makeatletter\ifx\SetFigFont\undefined%
\gdef\SetFigFont#1#2#3#4#5{%
  \reset@font\fontsize{#1}{#2pt}%
  \fontfamily{#3}\fontseries{#4}\fontshape{#5}%
  \selectfont}%
\fi\endgroup%
\begin{picture}(5320,5152)(21,-4315)
\put(3751,-2836){\makebox(0,0)[lb]{\smash{{\SetFigFont{7}{8.4}{\familydefault}{\mddefault}{\updefault}{\color[rgb]{0,0,0}$\gamma_3'$}%
}}}}
\put(751,-1111){\makebox(0,0)[lb]{\smash{{\SetFigFont{7}{8.4}{\familydefault}{\mddefault}{\updefault}{\color[rgb]{0,0,0}$\gamma_2$}%
}}}}
\put(1501,-661){\makebox(0,0)[lb]{\smash{{\SetFigFont{7}{8.4}{\familydefault}{\mddefault}{\updefault}{\color[rgb]{0,0,0}$\gamma_2'$}%
}}}}
\put(676,-2461){\makebox(0,0)[lb]{\smash{{\SetFigFont{7}{8.4}{\familydefault}{\mddefault}{\updefault}{\color[rgb]{0,0,0}$\gamma_4$}%
}}}}
\put(1426,-2836){\makebox(0,0)[lb]{\smash{{\SetFigFont{7}{8.4}{\familydefault}{\mddefault}{\updefault}{\color[rgb]{0,0,0}$\gamma_4'$}%
}}}}
\put(3676,-661){\makebox(0,0)[lb]{\smash{{\SetFigFont{7}{8.4}{\familydefault}{\mddefault}{\updefault}{\color[rgb]{0,0,0}$\gamma_1'$}%
}}}}
\put(4351,-1111){\makebox(0,0)[lb]{\smash{{\SetFigFont{7}{8.4}{\familydefault}{\mddefault}{\updefault}{\color[rgb]{0,0,0}$\gamma_1$}%
}}}}
\put(5326,-1486){\makebox(0,0)[lb]{\smash{{\SetFigFont{7}{8.4}{\familydefault}{\mddefault}{\updefault}{\color[rgb]{0,0,0}$\bm J^3$}%
}}}}
\put(4351,-2461){\makebox(0,0)[lb]{\smash{{\SetFigFont{7}{8.4}{\familydefault}{\mddefault}{\updefault}{\color[rgb]{0,0,0}$\gamma_3$}%
}}}}
\end{picture}%

%% file: saddle1.pdf_t
\begin{picture}(0,0)%
\includegraphics{saddle1.pdf}%
\end{picture}%
\setlength{\unitlength}{3947sp}%
\begingroup\makeatletter\ifx\SetFigFont\undefined%
\gdef\SetFigFont#1#2#3#4#5{%
  \reset@font\fontsize{#1}{#2pt}%
  \fontfamily{#3}\fontseries{#4}\fontshape{#5}%
  \selectfont}%
\fi\endgroup%
\begin{picture}(4877,4550)(14,-3725)
\put(4876,-1411){\makebox(0,0)[lb]{\smash{{\SetFigFont{6}{7.2}{\familydefault}{\mddefault}{\updefault}{\color[rgb]{0,0,0}$X_1$}%
}}}}
\put(2672,-1322){\makebox(0,0)[lb]{\smash{{\SetFigFont{5}{6.0}{\familydefault}{\mddefault}{\updefault}{\color[rgb]{0,0,0}$\partial \bm J^{N+1}$}%
}}}}
\put(2194,-1163){\makebox(0,0)[lb]{\smash{{\SetFigFont{5}{6.0}{\familydefault}{\mddefault}{\updefault}{\color[rgb]{0,0,0}$R^{-1}(\bm B)$}%
}}}}
\put(2194,-686){\makebox(0,0)[lb]{\smash{{\SetFigFont{5}{6.0}{\familydefault}{\mddefault}{\updefault}{\color[rgb]{0,0,0}$\partial R^{-1}(\bm B_1')$}%
}}}}
\put(1655,-82){\makebox(0,0)[lb]{\smash{{\SetFigFont{5}{6.0}{\familydefault}{\mddefault}{\updefault}{\color[rgb]{0,0,0}$\bm W_1$}%
}}}}
\put(1686,-3038){\makebox(0,0)[lb]{\smash{{\SetFigFont{5}{6.0}{\familydefault}{\mddefault}{\updefault}{\color[rgb]{0,0,0}$\bm W_2$}%
}}}}
\put(383,-1672){\makebox(0,0)[lb]{\smash{{\SetFigFont{5}{6.0}{\familydefault}{\mddefault}{\updefault}{\color[rgb]{0,0,0}$\Gamma_2'$}%
}}}}
\put(3911,-1703){\makebox(0,0)[lb]{\smash{{\SetFigFont{5}{6.0}{\familydefault}{\mddefault}{\updefault}{\color[rgb]{0,0,0}$\Gamma_1'$}%
}}}}
\put(4388,-464){\makebox(0,0)[lb]{\smash{{\SetFigFont{5}{6.0}{\familydefault}{\mddefault}{\updefault}{\color[rgb]{0,0,0}$\bm J^3$}%
}}}}
\end{picture}%

%% file: highdetailsmall2.pdf_t
\begin{picture}(0,0)%
\includegraphics{highdetailsmall2.pdf}%
\end{picture}%
\setlength{\unitlength}{3947sp}%
\begingroup\makeatletter\ifx\SetFigFont\undefined%
\gdef\SetFigFont#1#2#3#4#5{%
  \reset@font\fontsize{#1}{#2pt}%
  \fontfamily{#3}\fontseries{#4}\fontshape{#5}%
  \selectfont}%
\fi\endgroup%
\begin{picture}(7673,7666)(-11,-5620)
\put(3601,1364){\makebox(0,0)[lb]{\smash{{\SetFigFont{6}{7.2}{\familydefault}{\mddefault}{\updefault}{\color[rgb]{0,0,0}$e_1$}%
}}}}
\put(6324,-1851){\makebox(0,0)[lb]{\smash{{\SetFigFont{5}{6.0}{\familydefault}{\mddefault}{\updefault}{\color[rgb]{0,0,0}$X_1$}%
}}}}
\put(5611,-115){\makebox(0,0)[lb]{\smash{{\SetFigFont{5}{6.0}{\familydefault}{\mddefault}{\updefault}{\color[rgb]{0,0,0}$\bm J^4$}%
}}}}
\put(4925,-1730){\makebox(0,0)[lb]{\smash{{\SetFigFont{5}{6.0}{\familydefault}{\mddefault}{\updefault}{\color[rgb]{0,0,0}$\Gamma\subset X_2$}%
}}}}
\put(4214,-1387){\makebox(0,0)[lb]{\smash{{\SetFigFont{5}{6.0}{\familydefault}{\mddefault}{\updefault}{\color[rgb]{0,0,0}$X_3$}%
}}}}
\put(3310,-115){\makebox(0,0)[lb]{\smash{{\SetFigFont{5}{6.0}{\familydefault}{\mddefault}{\updefault}{\color[rgb]{0,0,0}$\comp_c R^{-1}(\Gamma)$}%
}}}}
\put(3499,-269){\makebox(0,0)[lb]{\smash{{\SetFigFont{5}{6.0}{\familydefault}{\mddefault}{\updefault}{\color[rgb]{0,0,0}$\subset X_2$}%
}}}}
\put(3094,-1356){\makebox(0,0)[lb]{\smash{{\SetFigFont{5}{6.0}{\familydefault}{\mddefault}{\updefault}{\color[rgb]{0,0,0}$X_3$}%
}}}}
\put(3094,-2413){\makebox(0,0)[lb]{\smash{{\SetFigFont{5}{6.0}{\familydefault}{\mddefault}{\updefault}{\color[rgb]{0,0,0}$X_3$}%
}}}}
\put(4214,-2413){\makebox(0,0)[lb]{\smash{{\SetFigFont{5}{6.0}{\familydefault}{\mddefault}{\updefault}{\color[rgb]{0,0,0}$X_3$}%
}}}}
\put(1851,-1730){\makebox(0,0)[lb]{\smash{{\SetFigFont{5}{6.0}{\familydefault}{\mddefault}{\updefault}{\color[rgb]{0,0,0}$\Gamma\subset X_2$}%
}}}}
\put(4881,-410){\makebox(0,0)[lb]{\smash{{\SetFigFont{5}{6.0}{\familydefault}{\mddefault}{\updefault}{\color[rgb]{0,0,0}$\zeta_1$}%
}}}}
\put(4676,-1378){\makebox(0,0)[lb]{\smash{{\SetFigFont{7}{8.4}{\familydefault}{\mddefault}{\updefault}{\color[rgb]{0,0,0}$\zeta_2$}%
}}}}
\put(4123,-1721){\makebox(0,0)[lb]{\smash{{\SetFigFont{7}{8.4}{\familydefault}{\mddefault}{\updefault}{\color[rgb]{0,0,0}$\zeta_3$}%
}}}}
\put(4054,-963){\makebox(0,0)[lb]{\smash{{\SetFigFont{7}{8.4}{\familydefault}{\mddefault}{\updefault}{\color[rgb]{0,0,0}$\zeta_4$}%
}}}}
\put(7076,271){\makebox(0,0)[lb]{\smash{{\SetFigFont{6}{7.2}{\familydefault}{\mddefault}{\updefault}{\color[rgb]{0,0,0}$\bm J^0$}%
}}}}
\put(6660,160){\makebox(0,0)[lb]{\smash{{\SetFigFont{6}{7.2}{\familydefault}{\mddefault}{\updefault}{\color[rgb]{0,0,0}$\J^1$}%
}}}}
\put(6268, 57){\makebox(0,0)[lb]{\smash{{\SetFigFont{6}{7.2}{\familydefault}{\mddefault}{\updefault}{\color[rgb]{0,0,0}$\J^2$}%
}}}}
\put(5946,-35){\makebox(0,0)[lb]{\smash{{\SetFigFont{6}{7.2}{\familydefault}{\mddefault}{\updefault}{\color[rgb]{0,0,0}$\J^3$}%
}}}}
\put(5301,242){\makebox(0,0)[lb]{\smash{{\SetFigFont{6}{7.2}{\familydefault}{\mddefault}{\updefault}{\color[rgb]{0,0,0}$e_3$}%
}}}}
\put(4886,841){\makebox(0,0)[lb]{\smash{{\SetFigFont{6}{7.2}{\familydefault}{\mddefault}{\updefault}{\color[rgb]{0,0,0}$e_2$}%
}}}}
\put(7355,-1924){\makebox(0,0)[lb]{\smash{{\SetFigFont{6}{7.2}{\familydefault}{\mddefault}{\updefault}{\color[rgb]{0,0,0}$e_0$}%
}}}}
\end{picture}%

%% file: multihighsmall.pdf_t
\begin{picture}(0,0)%
\includegraphics{multihighsmall.pdf}%
\end{picture}%
\setlength{\unitlength}{3947sp}%
\begingroup\makeatletter\ifx\SetFigFont\undefined%
\gdef\SetFigFont#1#2#3#4#5{%
  \reset@font\fontsize{#1}{#2pt}%
  \fontfamily{#3}\fontseries{#4}\fontshape{#5}%
  \selectfont}%
\fi\endgroup%
\begin{picture}(2690,2690)(606,-2456)
\put(2916,-942){\makebox(0,0)[lb]{\smash{{\SetFigFont{5}{6.0}{\familydefault}{\mddefault}{\updefault}{\color[rgb]{0,0,0}$\Gamma$}%
}}}}
\put(1879,-918){\makebox(0,0)[lb]{\smash{{\SetFigFont{5}{6.0}{\familydefault}{\mddefault}{\updefault}{\color[rgb]{0,0,0}$\bm W^0$}%
}}}}
\put(2916,-123){\makebox(0,0)[lb]{\smash{{\SetFigFont{5}{6.0}{\familydefault}{\mddefault}{\updefault}{\color[rgb]{0,0,0}$\bm J^3$}%
}}}}
\put(2167,-1931){\makebox(0,0)[lb]{\smash{{\SetFigFont{5}{6.0}{\familydefault}{\mddefault}{\updefault}{\color[rgb]{0,0,0}$\mathrm{Comp}_c R^{-1}(\Gamma)$}%
}}}}
\end{picture}%